\documentclass[11pt]{article}
\usepackage{amssymb}
\usepackage{amsmath}
\usepackage{datetime}
\usepackage{enumerate}
\usepackage{mathrsfs}
\usepackage{graphics}
\usepackage{graphicx}
\IfFileExists{stmaryrd.sty}{\usepackage{stmaryrd}}{}
\usepackage{xcolor}
\usepackage{subfigure}
\usepackage[T1]{fontenc}
\usepackage{latexsym,amssymb,amsmath,amsfonts,amsthm,txfonts}
\IfFileExists{thmtools.sty}{\usepackage{thmtools}}{}

\usepackage[colorlinks=true, linkcolor=blue, citecolor=blue]{hyperref}

\usepackage{amsfonts}
\usepackage{mathtools}
\usepackage{esint}
\usepackage{enumitem}
\IfFileExists{dsfont.sty}{\usepackage{dsfont}}{\let\mathds\mathbf}

\usepackage[T1]{fontenc}
\usepackage[toc,page]{appendix}
\usepackage{verbatim}

\usepackage{setspace}
\IfFileExists{fouridx.sty}{\usepackage{fouridx}}{}

\newcommand{\dom}{\mathscr{O}}

\newcommand{\ab}{a}
\newcommand{\bb}{b}
\newcommand{\xb}{x}
\newcommand{\ub}{u}
\newcommand{\vb}{\mathrm{v}}
\newcommand{\zb}{z}

\newcommand{\fb}{f}
\newcommand{\jb}{j}
\newcommand{\kb}{k}
\newcommand{\wb}{w}

\newcommand{\eb}{e}

\newcommand{\tor}{\mathbb{T}^3}
\newcommand{\rA}{\mathrm{A}}
\newcommand{\loc}{\mathrm{loc}}

\usepackage{todonotes}

\usepackage{hyperref}

\newcommand\exxp[1]{e^{#1}}

\newcommand\dela[1]{}

\newcommand{\N}{{\mathbb N}}

\newcommand{\be}{\begin{eqnarray}}
\newcommand{\ben}{\begin{eqnarray*}}
\newcommand{\en}{\end{eqnarray}}
\newcommand{\enn}{\end{eqnarray*}}

\newcommand{\eps}{\varepsilon}

\newcommand{\Om}{\Omega}

\newtheorem{theorem}{Theorem}[section]

\newtheorem{lemma}[theorem]{Lemma}
\newtheorem{prop}[theorem]{Proposition}
\newtheorem{cor}[theorem]{Corollary}
\newtheorem{remark}[theorem]{Remark}

\newtheorem{definition}[theorem]{Definition}

\newtheorem{assumption}[theorem]{Assumption}

\numberwithin{equation}{section}

\begin{document}
\renewcommand{\theequation}{\arabic{section}.\arabic{equation}}
\begin{titlepage}
\title{Global well-posedness and backward uniqueness of stochastic 3D Burgers equation with randomized $\mathbb{L}^{2}$ initial data}
\author{ Zdzis{\l}aw Brze{\'z}niak$^{1},$\ \ Zhao Dong$^{2},$\ \  Shijia Zhang$^{3}$ and  Guoli Zhou$^{4}$\thanks{Corresponding author}\\
{\small $^1$ Department of Mathematics, University of York, United Kingdom}\\
{\small $^{2,3}$  RCSDS, Academy of Mathematics and Systems Science, Chinese Academy of Sciences, Beijing 100190, China}\\
{\small $^{4}$ School of Statistics and Mathematics, Chongqing University, Chongqing 400044, China}\\
({\small zdzislaw.brzezniak@york.ac.uk, dzhao@amt.ac.cn,
shijiazhang@amss.ac.cn, zhouguoli736@126.com })}
\date{}
\end{titlepage}
\maketitle
\noindent\textbf{Abstract}:\
Global well-posedness for 3D deterministic Burgers equation with $\mathbb{L}^{2}(\tor;\mathbb{R}^{3})$-valued initial data is unknown.
Here, by a suitable randomization for the initial data, we prove almost-sure global existence and pathwise uniqueness for the stochastic 3D Burgers equation with the randomization of each fixed deterministic initial data in $\mathbb{L}^{2}(\tor;\mathbb{R}^{3})$. Furthermore, we establish the backward uniqueness for the stochastic 3D Burgers equation.



\renewcommand{\contentsname}{}
\tableofcontents
\newpage
\section{Introduction}
In this paper, we are concerned with 3D viscous stochastic Burgers equation with random initial data on the three
dimensional torus $\tor:=\mathbb{R}^{3}/2\pi \mathbb{Z}^{3}$.  To be more precise,  let $(\Omega,\mathscr{F},\mathbb{F}, \mathbb{P})$ be a filtered probability space satisfying the usual assumptions, where the filtration is $\mathbb{F}=\{\mathscr{F}_t\}_{t\in [0,\infty)}$.
We use $\mathbb{E}$ to denote the expectation with respect to $\mathbb{P}$.
We assume that $B=(B(t): t \geq 0)$ is an $\mathbb{R}$-valued Wiener process with respect to $\mathbb F$.
Let $b\in\mathbb R$ be a fixed constant.
We will denote by $\nu>0$ the viscosity which is a fixed constant throughout the whole paper. If $E$ is a metric space, then by
$\mathscr{B}(E)$ we will denote the  Borel $\sigma$-field on $E$.
If $p \in [1,\infty]$, then by $\mathbb{L}^{p}(\tor)$ we denote the Banach space of all (equivalence classes) of $p$-integrable
Lebesgue measurable $\mathbb{R}^{3}$-valued functions endowed with the classical norm.
We assume that   $u_0 \in \mathbb{L}^{2}(\tor)$ is a  deterministic function and
\[\fb: [0,\infty)\times \tor\rightarrow \mathbb{R}^{3}
\]
is  Lebesgue measurable satisfying some integrability conditions.
We consider the following Cauchy problem with the initial data
 \begin{align}\label{eqn-Burgers-intro}
&d\ub(t)=(\nu\Delta \ub-(\ub\cdot \nabla) \ub)(t)dt+\fb(t)dt+\ub(t)\circ b dB(t),\ \ \mathrm{on}\ [0,\infty)\times \tor,\\
&\ub(0, \omega)=\overrightarrow{\ub_{0}},\quad  \omega\in \Omega ,\nonumber
\end{align}
where the initial data $\overrightarrow\ub_{0}$ is  a "randomization" of  $u_0 $ given by \eqref{eqn-randomization of u_0} using a sequence $r_{\kb}^j$ of independent, $\mathscr F_0$-measurable random variables satisfying \autoref{A1}.
A solution of equation \eqref{eqn-Burgers-intro} is  an $\mathbb{L}^{2}(\tor)$-valued progressively measurable
stochastic process $\ub:= (u_1, u_2,u_3): [0, \infty)  \times \Omega \ni (t,\omega) \mapsto \ub (t,\omega) \in  \mathbb{L}^{2}(\tor)$ which satisfies equation \eqref{eqn-Burgers-intro} in the weak form, see Definition \ref{def-solution}.
\par
The Burgers equation describes the interaction of nonlinear transport and viscous dissipation. It was introduced by Bateman \cite{Beteman_1915} and studied by Burgers \cite{Burgers_1939} as a model for turbulent fluid motion. In three dimensions, it retains convective acceleration and diffusion without the pressure term or the incompressibility constraint. In our model, the scalar multiplicative noise acts uniformly in space on the velocity amplitude, while the randomization describes uncertainty in the Fourier coefficients of a prescribed finite-energy initial velocity.

\par
One dimensional stochastic Burgers equation has been fairly well studied. By an adaptation of the celebrated Hopf-Cole transformation, Bertini, Cancrini, and Jona-Lasinio \cite{Bertini+Cancrini+Jona-Lasinio_1994} solved the one dimensional modified Burgers equation with additive space-time white noise.
Later, Da Prato, Debussche and Temam \cite{Prato+Debussche+Temam_1994} established well-posedness for the one dimensional stochastic Burgers equation on bounded domains by semigroup methods and proved the existence of an invariant measure.
Furthermore, Gy\"ongy and Nualart \cite{Gyongy+Nualart_1999} extended the results of Da Prato
and Debussche to the Burgers equation defined on the whole line.

 \par
For multidimensional Burgers equations, the maximum principle of Kiselev and Ladyzhenskaya \cite{Kiselev+Ladyzhenskaya_1957} plays an important role in the global theory. Pooley and Robinson \cite{Pooley+Robinson_2016} proved global well-posedness for 3D Burgers equations in $C([0,T];\mathbb{H}^{\frac{1}{2}}(\tor))$.
Utilising the maximum principle, Brze{\'z}niak, Goldys and Neklyudov \cite{Brz+Goldys+Neklyudov_2014} established global existence and uniqueness of mild solutions to multidimensional Burgers equations with additive noise for $L^p$ initial data, $p>d$, where $d$ is the spatial dimension. By the Doss--Sussman transformation, Dong, Guo, Wu and Zhou \cite{Dong+Guo+Wu+Zhou_2023} obtained global well-posedness for 3D Burgers equations with linear multiplicative noise in $C([0,T];\mathbb{H}^{\frac{1}{2}}(\tor))$. For It\^{o} linear multiplicative noise, Dong, Wu and Zhou \cite{Dong+Wu+Zhou_2024} studied regularity and ergodicity.

 \par

\par
If the initial data and the force are of potential type, the so-called gradient solutions are studied. There are a number of works about the potential Burgers equation, see \cite{Bakhtin+Li_2019, Boritchev_2016,  E+Khanin+Mazel+Sinai_2000, Gomes+Iturriaga+Khanin+Padilla_2005, Iturriaga+Khanin_2003, Khanin+Zhang_2003} and other references.

\par
In this article, we establish global well-posedness for the stochastic 3D Burgers equation with the randomization of each fixed deterministic initial datum in $\mathbb L^2(\tor)$. More precisely, under the assumptions of \autoref{thm-main result}, there exists a unique global weak solution $\ub$ satisfying
\[
\ub\in C([0,\infty);\mathbb{L}^2(\mathbb{T}^3))
\cap L^4_{\loc}([0,\infty);\mathbb{L}^6(\mathbb{T}^3))
\cap L^2_{\loc}([0,\infty);\mathbb{H}^{1,3}(\mathbb{T}^3)).
\]

The main results and features of this paper are summarized as follows.
(i) We establish probabilistic estimates for randomized $\mathbb L^2(\tor)$ initial data under uniformly bounded second moments of independent, centered random multipliers, without imposing any higher-order moment conditions. These estimates include the spatial $\mathbb L^6$ estimate and the small-time mixed-norm estimates for the randomized heat evolution and its gradient, and remain valid when the multipliers have infinite sixth moments.
(ii) We give a direct construction of almost-sure global weak solutions in $Z_\infty$ for the stochastic 3D Burgers equation with randomized $\mathbb L^2(\tor)$ initial data. We combine the small-time estimates with maximal regularity, parabolic regularization, the maximum principle and a continuation argument to prove global existence and pathwise uniqueness. Moreover, we establish continuous dependence in probability in $C([0,T];\mathbb L^2(\tor))$ with respect to the deterministic initial data, using the same random multipliers, Brownian motion and external force.
(iii) We establish backward uniqueness for these global weak solutions without imposing any potential structure on the velocity field or external force. More precisely, two solutions driven by the same Brownian motion and external force cannot coincide at a positive time unless they agree at all earlier times.

The second-moment assumption concerns integrability in probability. Under \autoref{A1}, the randomized initial datum belongs to $L^2(\Omega;\mathbb L^6(\tor))$, see \eqref{eq-randomization-L6}, but need not possess positive Sobolev regularity. Thus the $\mathbb H^{\frac12}$ theories in \cite{Pooley+Robinson_2016, Dong+Guo+Wu+Zhou_2023} do not cover the full class of initial data considered here. The $L^p$ theory in \cite{Brz+Goldys+Neklyudov_2014} is also relevant; its global theorem on the torus treats additive noise under stronger space-time regularity assumptions on the force, see \cite[Theorem 4.1]{Brz+Goldys+Neklyudov_2014}.

\begin{remark}[Relation with the $L^p$ theory]\label{rem-L6-Lp-theory}
By \eqref{eq-randomization-L6}, we infer that $\overrightarrow{\ub_0}\in\mathbb L^6(\tor)$, $\mathbb P$-a.s..
Since $6>3$, after the Doss--Sussman transformation, the local $L^p$ argument in \cite{Brz+Goldys+Neklyudov_2014} can be adapted to \eqref{ee-o} with $p=6$.
Under \autoref{A2}, positive-time regularization and the maximum principle then provide another way to obtain global existence and uniqueness.

However, this argument does not directly give the full solution class $Z_\infty$.
In particular, the time-weighted gradient regularity in \cite[Theorem 3.4]{Brz+Goldys+Neklyudov_2014} does not imply
\[
\nabla\ub\in L^2([0,T];\mathbb L^3(\tor)),\qquad T>0.
\]

Our method provides this regularity  directly.
Indeed, by \eqref{eq-623-1}, we deduce
\[
\mathbb E\int_0^T\Vert\nabla e^{-tA}\overrightarrow{\ub_0}\Vert_{L^3}^2\,dt\leq cK_T^2,\qquad T\in(0,1],
\]
where $K_T\to0$ as $T\downarrow0$ for every deterministic $\ub_0\in\mathbb L^2(\tor)$.
This estimate requires only \autoref{A1}, without any positive Sobolev regularity of $\ub_0$.
Together with the estimates for $\wb$ and the nonlinear remainder $\vb$, it yields the required mixed-norm regularity up to the initial time and allows us to construct the solution directly in $Z_\infty$.

Moreover, using the same random multipliers, Brownian motion and external force, we prove continuous dependence in probability in $C([0,T];\mathbb L^2(\tor))$ with respect to the deterministic initial data $\ub_0$.
We also prove backward uniqueness for the resulting solutions, including equality at the initial time by continuity in $\mathbb L^2(\tor)$.
\end{remark}

For random initial Cauchy problems, see \cite{Burq+Tzvetkov_2008, Burq+Tzvetkov_2008_2, Burq+Tzvetkov_2014, Deng+Cui_2011, Nahmod+Pavlovic+Staffilani_2013, Nahmod+Rey-Bellet+Sheffield+Staffilani_2011, Mohammed+Zhang_2013}.

\par
Here, we give a brief description of the route to solving (\ref{eqn-Burgers-intro}).
For this outline, let $\nu=1$. Let \[\alpha(t)=\exp(bB(t)), t\geq 0,\] and set $$ \wb=\wb(t)=\wb(t,\xb) = \int_{0}^{t} e^{(t-s)\Delta}\alpha^{-1}(s)\fb(s)\,ds, (t,\xb) \in [0,\infty)\times \tor. $$ Let $\ub=\alpha\tilde{\vb}+\alpha\wb. $
Then the stochastic 3D Burgers equation with random initial value (\ref{eqn-Burgers-intro}) can be transformed into the following random case:
\begin{align}\label{ee-o}
&\partial_{t} \tilde{\vb}(t)-\Delta \tilde{\vb}(t)+\alpha(t)\left((\tilde{\vb}+\wb\right)\cdot \nabla )(\tilde{\vb}+\wb)(t)=0,\ \mathrm{on}\ (0,\infty)\times \tor,\\
&\tilde{\vb}_0=\tilde{\vb}(0),\nonumber
\end{align}
where
\[
\tilde{\vb}_0=\overrightarrow{\ub_{0}} \in \mathbb{L}^{2}(\tor)
\]
and
$\overrightarrow{\ub_{0}} $ is the randomization of $\ub_{0}\in \mathbb{L}^{2}(\tor)$ given in \eqref{eqn-randomization of u_0} and satisfies \autoref{A1}.  To obtain the well-posedness for (\ref{ee-o}), we decompose it into a nonlinear partial differential equation with zero initial data (\ref{eqn-Burgers-intro.4}) and a linear part $\mathbf{z}$.
\begin{align}\label{eqn-Burgers-intro.4}
&\partial_{t} \vb(t)-\Delta \vb(t)+\alpha(t)((\vb+\mathbf{z})\cdot \nabla )(\vb+\mathbf{z})(t,\xb)=0,\ \mathrm{on}\ (0,\infty)\times \tor,\\
&\vb(0)=0,\nonumber
\end{align}
where function $\mathbf{z}:[0,\infty)\times \tor\mapsto \mathbb{R}^3$ is defined by
\begin{align}\label{eq-z}
\mathbf{z}(t, \overrightarrow{\ub_{0}}):=e^{t\Delta}\overrightarrow{\ub_{0}}+\wb(t),
\end{align}
and satisfies
\begin{align}
&d\mathbf{z}(t)-\Delta\mathbf{z}(t)dt=\alpha^{-1}(t)\fb(t)dt, \ \ \mathrm{on}\ (0,\infty)\times \tor,\\
&\mathbf{z}(0)=\overrightarrow{\ub_{0}}\in \mathbb{L}^{2}(\tor).\nonumber
\end{align}
For general vector fields, the convection term has no $\mathbb L^2$ energy cancellation. By \autoref{prop-ouz}, $\mathbf z\in L^4(0,T;\mathbb L^6(\tor))$ and $\nabla\mathbf z\in L^2(0,T;\mathbb L^3(\tor))$, and these norms tend to zero as $T\to0+$, almost surely. In particular, $(\mathbf z\cdot\nabla)\mathbf z\in L^{\frac43}(0,T;\mathbb L^2(\tor))$. These bounds also control the mixed terms and allow us to apply maximal regularity and a fixed-point argument to the nonlinear remainder $\vb$, see \autoref{thm-local l2} and \autoref{thm-local l2-weak}. Observing that $$ \tilde{\vb}(t)=  \vb(t)+ e^{t\Delta}\overrightarrow{\ub_{0}},  t\in [0, \tau],$$ the local well-posedness of (\ref{ee-o}) is obtained in  \autoref{thm-local weak solution}. By \autoref{cor-ae-H2}, for almost every $\omega\in\Omega$, there exists $\tau^\ast\in(0,\tau(\omega))$ such that
\[
\tilde{\vb}(\tau^\ast,\omega)\in\mathbb H^2(\tor).
\]
The next question is to extend the local solution $(\tilde{\vb},\tau)$ globally, using the strong solution starting at $\tau^\ast$. To justify the maximum-principle estimates at this regularity, we first study a globally well-posed regularized system and then pass to the limit. Finally, positive-time strong regularity allows us to apply the classical Dirichlet-quotient argument for backward uniqueness; continuity in $\mathbb L^2(\tor)$ extends the conclusion to the initial time.

\par
The rest of the paper is organised as follows. \autoref{section-2} contains the notation, assumptions and main result. Local existence and uniqueness are proved in \autoref{section-3}. Global existence, continuous dependence and backward uniqueness are established in \autoref{section-4}. Auxiliary proofs are collected in the appendices.

Throughout the paper, we use $c > 0$ for a generic constant with possibly different values at each appearance, and
denote by $c(a)>0$ a constant which depends on parameter $a$.

Without loss of generality, we assume that $\nu=1$ in the following.
Indeed, for general $\nu>0$, set $s=\nu t$, $v(s,x)=\nu^{-1}u(s/\nu,x)$ and $\widetilde B(s)=\sqrt{\nu}B(s/\nu)$. Then $\widetilde B$ is a standard Brownian motion and the transformed equation has viscosity one, with the corresponding rescaled coefficients and force. Thus it is sufficient to consider $\nu=1$.

\section{Notations, assumptions and main result}\label{section-2}

\subsection{Sobolev spaces on the torus}

By $\dom$ we will denote the three dimensional torus $\tor$.

For $1 \leq p \leq \infty$, we denote by $\mathbb{L}^p(\dom)$ the usual Lebesgue space of $\mathbb{R}^3$-valued measurable functions $\ub : \dom \to \mathbb{R}^3$ with norm
\[
\Vert \ub\Vert_{L^p} := \left( \int_{\dom} |\ub(x)|^p \, dx \right)^{1/p}, \quad 1 \leq p < \infty,
\]
and
\begin{equation}\label{eqn-L^infty}
\Vert\ub\Vert_{L^\infty(\dom)}:= \mathrm{ess} \sup\limits_{x\in \dom}|\ub(x)|.
\end{equation}
For $p=2$, the space $\mathbb{L}^2(\dom)$ is a real Hilbert space with the inner product
\begin{equation}\label{eqn-L^2-real inner product}
\langle \ub, \vb \rangle := \int_{\dom} \ub(x) \cdot \vb(x) \, dx.
\end{equation}

We use the following real Fourier basis; the complex formulation and its relation to the real basis are recalled in Appendix~\ref{sec-Fourier}.
\begin{theorem}\label{thm-ONB}
Let $\mathbb{N}=\{0,1,2,\ldots\}$ and $\mathbb Z_0^3=\mathbb Z^3\setminus\{(0,0,0)\}$. Define
\begin{align}
    K_0 &:= \left\{ k \in \mathbb{Z}_0^3 : k_1 > 0 \;\text{or}\; (k_1 = 0, k_2 > 0) \;\text{or}\; (k_1 = k_2 = 0, k_3 > 0) \right\},\\
    K&:=K_0 \cup \{(0,0,0)\}.
\end{align}

If $\kb \in K_0$, then we define  $\mathbb{R}^3$-valued functions $\eb_\kb^j$, $j=1,\cdots,6$,  on $\dom$ by
\begin{equation}\label{eq-ONB-real}
\begin{aligned}
     \eb_\kb^1(x)&:=(4\pi^3)^{-\frac12}\left(\sin(\kb\cdot x),0,0\right), \;\;x \in \dom,
     \\
     \eb_\kb^2(x)&:=(4\pi^3)^{-\frac12}\left(0,\sin(\kb\cdot x),0\right), \;\; x \in \dom,
     \\
     \eb_\kb^3(x)&:=(4\pi^3)^{-\frac12}\left(0,0,\sin(\kb\cdot x)\right), \;\; x \in \dom,\\
     \eb_\kb^4(x)&:=(4\pi^3)^{-\frac12}\left(\cos(\kb\cdot x),0,0\right), \;\; x \in \dom,
     \\
     \eb_\kb^5(x)&:=(4\pi^3)^{-\frac12}\left(0,\cos(\kb\cdot x),0\right), \;\; x \in \dom,
     \\
     \eb_\kb^6(x)&:=(4\pi^3)^{-\frac12}\left(0,0,\cos(\kb\cdot x)\right), \;\; x \in \dom,
\end{aligned}
\end{equation}
If $\kb=0=(0,0,0)\in  K$, then we define  $\mathbb{R}^3$-valued functions $\eb_0^j$, $j=1,\cdots,3$,  on $\dom$ by
\begin{equation}\label{eq-ONB-real-2}
\begin{aligned}
          \eb_0^1(x)&:=(2\pi)^{-\frac32}\left(1,0,0\right), \;\;x \in \dom,
     \\
     \eb_0^2(x)&:=(2\pi)^{-\frac32}\left(0,1,0\right), \;\; x \in \dom,
     \\
     \eb_0^3(x)&:=(2\pi)^{-\frac32}\left(0,0,1\right), \;\; x \in \dom,
\end{aligned}
\end{equation}

We put
\begin{equation}\label{eqn-mathcal K}
\mathcal{K} := (K_0 \times \{1,2,3,4,5,6\})\cup(\{0\}\times \{1,2,3\})
\end{equation}

Then, the sequence $(\eb_\kb^j:(\kb,j) \in \mathcal{K})$
is an orthonormal basis of the real Hilbert space  $\mathbb{L}^2(\dom; \mathbb{R}^3)$.
In particular, if $\ub \in \mathbb{L}^2(\dom; \mathbb{R}^3)$, then there exists a unique real sequence
$\left(\ub_k^j\right)_{ (k,j) \in \mathcal{K}}$ such that
\begin{align}\label{eqn-u-real}
\ub(x) =  \sum_{(k,j) \in \mathcal{K}} u_k^j \eb_k^j(x),
\end{align}
with convergence in $\mathbb{L}^2(\dom; \mathbb{R}^3)$, see \autoref{thm-Fourier-real}.
Moreover, this sequence is determined by
\begin{equation}\label{eqn-u_k-real}
\ub_\kb^j := \langle \ub , \eb_k^j \rangle.
\end{equation}
Moreover, the Parseval identity holds:
\begin{align}\label{eqn-Parseval-real}
\Vert \ub\Vert_{L^2}^2 =  \sum_{(k,j) \in \mathcal{K}} |\ub_k^j|^2.
\end{align}

Moreover, there exists a constant $C>0$, independent of $(k,j)\in\mathcal K$, such that
\[
\Vert e_k^j\Vert_{L^6}\leq C,
\qquad
\Vert\nabla e_k^j\Vert_{L^3}\leq C|k|.
\]
\end{theorem}

For $s \geq 0$, we define $\mathbb{H}^{s}(\dom; \mathbb{R}^{3})$, the real Sobolev space  of order $s$    by the  following
\begin{equation}\label{eqn-H^2d-real}
\mathbb{H}^{s}(\dom; \mathbb{R}^{3}) := \left\{ \ub \in \mathbb{L}^2(\dom; \mathbb{R}^3) : \sum_{(k,j) \in \mathcal{K}}(1+ |\kb|^{2s})|\ub_{\kb}^j|^{2}    < \infty\right\},
    \end{equation}
with norm
\begin{equation}\label{eqn-norm H^s}
\Vert \ub\Vert_{H^s}^2:=  \sum_{(k,j) \in \mathcal{K}}(1+ |\kb|^{2s})|\ub_{\kb}^j|^{2}  .
\end{equation}
It is known that $\mathbb{H}^{s}(\dom; \mathbb{R}^{3})$ is a separable real Hilbert  space. For $s=0$, we set $\mathbb{H}^{0}(\dom;\mathbb{R}^{3})=\mathbb{L}^{2}(\dom;\mathbb{R}^{3})$ with the $L^2$-norm.

We define the
negative Sobolev space as the dual of the positive Sobolev space, i.e. for $s>0$
\begin{align}\label{eq-negative Sobolev space}
 \mathbb{H}^{-s}(\dom):=
(\mathbb{H}^{s}(\dom))^{\prime},
\end{align}
with norm
\begin{align}\label{eq-negative norm}
\Vert \vb\Vert_{H^{-s}}:=
\sup_{\Vert u\Vert_{H^{s}}=1}|\langle \vb,u\rangle|,\;\;\vb\in \mathbb{H}^{-s}.
\end{align}
Equivalently, writing $v_k^j:=\langle\vb,\eb_k^j\rangle$, we have
\[
\Vert\vb\Vert_{H^{-s}}^2=\sum_{(k,j)\in\mathcal K}
\frac{|v_k^j|^2}{1+|k|^{2s}},\qquad s>0.
\]

\begin{definition}\label{def-Laplace}
Let us define operator $\rA$  in the Hilbert space $H=\mathbb{L}^2(\dom,\mathbb{R}^3)$ by
\begin{equation}\label{eqn-Laplace}
 \begin{aligned}
D(\rA):=\mathbb{H}^2(\dom;\mathbb{R}^3)\\
\rA u := (-\Delta u_1, -\Delta u_2,-\Delta u_3) \end{aligned}
\end{equation}
\end{definition}
The following result is well known.
\begin{prop}\label{prop-eigenvectors of A}
The operator $\rA$ is a non-negative and self-adjoint  in the Hilbert $H$. Moreover,
for every  $(k,j) \in \mathcal{K}$, the function $e_k^j$  is an eigenvector of $\rA$ with
\begin{equation}\label{eqn-Laplace-1}
\rA e_k^j =\vert k \vert^2 e_k^j.
\end{equation}
\end{prop}
Define $\Lambda$ to be the square root of $\rA$, i.e.
\begin{equation}
    \label{eqn-Laplace-square root}
\Lambda:= \rA^{\frac12}
\end{equation}

It follows from the definition \eqref{eqn-Laplace-square root} and interpolation theory  that
\begin{equation}\label{eqn-Laplace-2}
    \begin{aligned}
D(\Lambda):=\mathbb{H}^1(\dom;\mathbb{R}^3)\\
\Lambda (\Lambda u)=\rA u, \;\; u \in D(\rA).
\end{aligned}
\end{equation}

Since $\Lambda$ is also a  non-negative self-adjoint operator in the Hilbert $H$, the fractional powers
$\Lambda^{s}$ are well defined. They can be characterised as follows.

For $s> 0$, assume that
$\ub \in \mathbb{H}^{s}(\dom)$, we define
\begin{equation} \label{def-Lambda^s}
	(\Lambda^{s}\ub)(\xb):=\sum_{(\kb,j)\in \mathcal{K}}|\kb|^{s}\ub_{\kb}^je_k^j, \;\; \xb \in \dom.
\end{equation}
For $s=0$, we define $\Lambda^0=I$.
One can show that
\[
\Lambda^{s}\ub \in \mathbb{L}^2(\dom).
\]
Indeed, by the Parseval Theorem,
\begin{align}
\Vert \Lambda^{s}\ub\Vert_{L^2}^{2}
	&= \sum_{(\kb,j)\in \mathcal{K}}
	 |\kb|^{2s}|\ub_{\kb}^j|^2\leq  \sum_{(\kb,j)\in \mathcal{K}} (1+ \vert \kb \vert^{2s}) \vert \ub_{\kb}^j \vert^2 <\infty,\;\;s>0.
\end{align}

For $s>0$, denote by $\Vert \cdot\Vert_{s}$ the seminorm on $\mathbb{H}^{s}(\dom)$
\begin{equation}
 \Vert u\Vert_{s}:=\Vert \Lambda^{s} u \Vert_{L^2}, \;\; u \in \mathbb{H}^{s}(\dom).
\end{equation}
For $s=0$, we define
\begin{equation}
\label{eqn-norm H^s-homogenous}
 \Vert \ub\Vert_{0}^2= \sum_{(\kb,j)\in \mathcal{K},\; k\neq0} |\ub_{\kb}^j|^{2} .
\end{equation}
Note that
\begin{equation}
\label{eqn-1011}
\Vert \ub-\hat{\ub}_0\Vert_{L^2}^2= \Vert \ub\Vert_{0}^2
\end{equation}
and
\begin{equation}
\label{eqn-1012}
 \Vert \ub\Vert_{0}^2 \leq  \Vert \ub\Vert_{s}^2,\;\;s\geq 0.
\end{equation}
We deduce that for each $s\geq 0$,
\begin{equation}
\label{eqn-1013}
\Vert \ub-\hat{\ub}_0\Vert_{L^2}^2 \leq  \Vert \ub\Vert_{s}^2.
\end{equation}
It is known that $\ub \in \mathbb{H}^s(\dom)$ if and only if $\ub \in \mathbb{L}^2(\dom)$ and $\Lambda ^s \ub \in \mathbb{L}^2(\dom)$.
Indeed, it follows from \eqref{eqn-H^2d-real} that  for $s>0$
\begin{align}\label{eq-homogeneous norm}
\Vert \cdot\Vert_{H^s}^2=\Vert \cdot\Vert_{L^2}^2+\Vert \cdot\Vert_{s}^2.
\end{align}

Note that for $0\leq s_{1}\leq s_{2}$,  we have $\Vert \ub\Vert_{{s_{1}}}\leq \Vert \ub\Vert_{s_{2}}$.

Note that for $s,s_1,s_2\geq 0$ satisfying $s=s_1+s_2$, we have
\begin{align}\label{eq-Lambda}
 A^{\frac{s}{2}}=A^{\frac{s_1}{2}}A^{\frac{s_2}{2}}.
\end{align}

Note that for $s,s_1,s_2\geq 0$ satisfy $s=s_1+s_2$, we have
\begin{align}\label{eq-Lambda-1}
 \Lambda^{s}=\Lambda^{s_1} \Lambda^{s_2}.
\end{align}

For heat-semigroup estimates with integer-order spatial derivatives, see \cite[Proposition 3.3]{Brz+Goldys+Neklyudov_2014}.
\begin{lemma}\label{lem3.3-1}
Let $p \in (1, \infty)$, $q \in [p, \infty)$ and $r\in (p,\infty]$ such that
\begin{equation}\label{eqn-pqr}
    \frac{1}{r}=\frac{1}{p}-\frac{1}{q}.
\end{equation}
Then for every $\gamma\geq 0$, there exists a constant $C>0$ independent of $h$ such that for  $0<t\leq 1$
\begin{align}\label{eqn-semigroup-1}
\left\Vert \Lambda^\gamma e^{-tA} h\right\Vert_{L^q} \leq C t^{-\frac{\gamma}{2}-\frac{3}{2 r}}\Vert h\Vert_{L^p}, \quad h \in \mathbb{L}^p(\dom).
\end{align}

\end{lemma}
Appendix~\ref{app-heat-semigroup-proofs} proves the estimate for every real $\gamma\geq0$.

\begin{remark}\label{rem-lem3.3}
All these and similar inequalities which are stated for $t \in (0,1]$ in \autoref{lem3.3-1} are also valid, with different constants, for $t \in (0,T]$, whenever $T>0$ is fixed.
\end{remark}

\begin{lemma}\label{lem-heat semigroup on homogenous spaces}
Let $s_1\in[0,\infty)$ and $s_2\in [s_1,\infty)$.
Suppose that $\gamma$ is defined by
\[
\gamma=s_2 - s_1
\]
then for $0< r\leq1$
\begin{trivlist}
\item[(i)]  there exists $C>0$ such that for every $g \in {\mathbb{H}}^{s_1}(\dom)$
\begin{equation}\label{eqn-3.101}
\begin{aligned}
\left\Vert  e^{-rA} g \right\Vert_{H^{s_2}}& \leq (C r^{-\frac{\gamma}{2}}+1 )\Vert  g \Vert_{H^{s_1}},
\end{aligned}
\end{equation}
\item[(ii)]  there exists $C>0$ such that for every $g \in {\mathbb{H}}^{-s_2}(\dom)$
\begin{equation}\label{eqn-3.102}
\begin{aligned}
\left\Vert  e^{-rA} g \right\Vert_{H^{-s_1}}& \leq C (r^{-\frac{\gamma}{2}}+1) \Vert  g \Vert_{H^{-s_2}}.
\end{aligned}
\end{equation}
\end{trivlist}
\end{lemma}
The proof is given in Appendix~\ref{app-heat-semigroup-proofs}.

For the real valued sequence $\ub=\{\ub_\kb^j\}_{(\kb,j)\in \mathcal{K}}$, $\ub_{\kb}^j\in \mathbb{R}$, $(\kb,j)\in \mathcal{K}$, we define
\[
\Vert \ub\Vert_{l^{2}}:=\Big{(}\sum_{(\kb,j)\in \mathcal{K}}|\ub_{\kb}^j |^{2}\Big{)}^{\frac{1}{2}}.
\]
Similarly,

for real valued sequence $a=\{a_{\kb}\}_{\kb\in K}$, $a_{\kb}\in \mathbb{R}, \kb\in K, $
\[
\Vert a\Vert_{l^{2}}:=\Big{(}\sum_{\kb\in K} |a_{\kb} |^{2}\Big{)}^{\frac{1}{2}}.
\]

\subsection{Spaces of time dependent functions}
Similarly, let $T > 0$ be fixed. For $1 \leq p \leq \infty$, we denote by $L^p([0,T])$ the standard Lebesgue space of measurable, real-valued functions on the time interval $[0,T]$, equipped with the norm
\[
\Vert f\Vert_{L^p([0,T])} :=
\begin{cases}
\left( \int_0^T |f(t)|^p \, dt \right)^{1/p}, & \text{if } 1 \leq p < \infty, \\
\mathrm{ess} \sup\limits_{t\in [0,T]} |f(t)|, & \text{if } p = \infty.
\end{cases}
\]

We introduce the following mixed norms, see e.g. \cite{Benedek+Panzone_1961}.

For $1 \leq p,q < \infty$, we also define mixed function spaces  $L^p([0,T]; \mathbb{L}^q(\dom))$ consists of measurable functions $f: [0,T] \times \dom \to \mathbb{R}$ such that
\[
\Vert f\Vert_{L^p([0,T]; \mathbb{L}^q(\dom))} := \left( \int_0^T \left( \int_{\dom} |f(t,x)|^q \, dx \right)^{p/q} \,dt \right)^{1/p}< \infty.
\]
For $p=\infty$, we define spaces $L^\infty([0,T]; \mathbb{L}^q(\dom))$ as follows:
\[
    L^\infty([0,T]; \mathbb{L}^q(\dom)) := \left\{ f : [0,T] \times \dom \to \mathbb{R} \,\middle|\,
    \mathrm{ess} \sup\limits_{t\in [0,T]} \left( \int_{\dom} |f(t,x)|^q dx \right)^{1/q} < \infty \right\}, \quad 1 \leq q < \infty.
    \]
In addition, we introduce the following commonly used spaces.
For $s\geq 0$ and $1\leq p<\infty$, we denote by
$L^p([0,T];\mathbb H^s(\dom))$ the space consisting of all
measurable functions $f:[0,T]\to \mathbb H^s(\dom)$ such that
\[
\Vert f\Vert_{L^p([0,T];\mathbb H^s(\dom))}
:=
\left(\int_0^T \Vert f(t)\Vert_{H^s(\dom)}^p\,dt\right)^{1/p}
<\infty.
\]
For $s \geq 0$, the space $L^\infty([0,T]; \mathbb{H}^s(\dom))$ consists of functions $f: [0,T] \times \dom \to \mathbb{R}$ such that $f(t,\cdot)\in \mathbb{H}^s(\dom)$ and
\[
\Vert f\Vert_{L^\infty([0,T]; \mathbb{H}^s(\dom))} := \mathrm{ess} \sup\limits_{t\in [0,T]}\Vert f(t)\Vert_{H^s(\dom)} < \infty.
\]
The space $L^\infty([0,T] \times \dom)$ consists of all measurable functions $f: [0,T] \times \dom \to \mathbb{R}$ such that
\[
\Vert f\Vert_{L^\infty([0,T] \times \dom)} := \mathrm{ess} \sup\limits_{(t,x) \in [0,T] \times \dom} |f(t,x)| < \infty.
\]

\subsection{Assumptions}

\begin{assumption}\label{A1}

We assume that $\{r_{\kb}^j\}_{(\kb,j)\in\mathcal K}$ is a sequence of independent, zero-mean, real, $\mathscr F_0$-measurable random variables on $(\Omega,\mathscr F,\mathbb P)$, that the family is independent of $B$, and that

\begin{align}\label{eq-random initial condition}
M_2:=\sup_{(\kb,j)\in\mathcal K}\mathbb E|r_{\kb}^j|^2<\infty.
\end{align}
\end{assumption}
\begin{assumption} \label{A2}
There exists $\delta>\frac{3}{2}$ such that the external force
$\fb:[0,\infty)\times \dom\rightarrow \mathbb{R}^{3}$ in (\ref{eqn-Burgers-intro}) satisfies
\begin{align}\label{eqn-f-assumption}
&\fb \in L_{\loc}^{2}([0,\infty); \mathbb{H}^{\delta}(\dom)).
\end{align}
\end{assumption}

\begin{remark}\label{remark-random variables}
In fact, many examples satisfy \autoref{A1}. For instance, let $\{r_{\kb}^j\}_{(\kb,j)\in\mathcal K}$ be a sequence of independent standard Gaussian or symmetric Bernoulli random variables, which is $\mathscr F_0$-measurable and independent of $B$.
\end{remark}

\begin{definition}\label{def-randomization-u}
Assume that $\ub\in \mathbb{L}^{2}(\dom)$  and
\[
u=\sum_{(\kb,j)\in \mathcal{K}}\ub_{\kb}^je_{\kb}^j.
\]
 We define the  randomization of   $u$
 to be the following   $\mathbb{L}^{2}(\dom)$-valued random variable
\begin{align}\label{eqn-randomization-u}
\overrightarrow\ub&:=\sum_{(\kb,j)\in \mathcal{K}}r_{\kb}^j\ub_{\kb}^je_{\kb}^j, \;\; \omega \in \Omega .
\end{align}
\end{definition}

Indeed, by \autoref{A1} and the Parseval identity, for every finite subset $\mathcal K_N\subset\mathcal K$,
\[
\mathbb E\Big\Vert\sum_{(\kb,j)\in\mathcal K_N}r_{\kb}^j\ub_{\kb}^je_{\kb}^j\Big\Vert_{L^2}^2
=\sum_{(\kb,j)\in\mathcal K_N}\mathbb E|r_{\kb}^j|^2|\ub_{\kb}^j|^2
\leq M_2\Vert\ub\Vert_{L^2}^2.
\]
Hence the series in \eqref{eqn-randomization-u} converges in $L^2(\Omega;\mathbb L^2(\dom))$. Let $\{\mathcal K_N\}_{N\in\mathbb N}$ be an increasing sequence of finite subsets of $\mathcal K$ such that $\mathcal K_N\uparrow\mathcal K$. For $M>N$, by the same calculation and the Parseval identity, we infer that
\begin{align*}
\mathbb E\Big\Vert\sum_{(\kb,j)\in\mathcal K_M\setminus\mathcal K_N}r_\kb^j\ub_\kb^je_\kb^j\Big\Vert_{L^2}^2
\leq M_2\sum_{(\kb,j)\in\mathcal K\setminus\mathcal K_N}|\ub_\kb^j|^2\to0,\qquad N\to\infty.
\end{align*}
Hence, the series in \eqref{eqn-randomization-u} converges in $L^2(\Omega;\mathbb L^2(\dom))$ and
\begin{align}\label{eq-randomization-L2}
\Vert\overrightarrow\ub\Vert_{L^2(\Omega;\mathbb L^2)}\leq M_2^{\frac12}\Vert\ub\Vert_{L^2}.
\end{align}

\begin{remark}\label{rem- randomization}
 The following has been proved in  \cite[Lemma B.1]{Burq+Tzvetkov_2008}.
 \begin{trivlist}
\item[(i)] If $\ub\in \mathbb{L}^{2}(\dom)$, then
  $\mathbb{P}$-almost surely, $\overrightarrow\ub \in \mathbb{L}^{2}(\dom)$.
  \item[(ii)] Suppose that there exist $c,p>0$ such that $\mathbb{P}(|r_k^j|>c)\ge p$ for every $(k,j)\in\mathcal K$.
If $u\in L^2(\mathcal O)$ and $u\notin H^\delta(\mathcal O)$
for every $\delta>0$, then
$\overrightarrow{u}\notin H^\delta(\mathcal O)$
for every $\delta>0$, $\mathbb{P}$-a.s..
 \end{trivlist}

\end{remark}

\par

\subsection{Main result}\label{subsec-main-result}

Throughout the paper, $B$ is a Wiener process with respect to $\mathbb F$. By \autoref{A1},
\[
\mathscr G:=\sigma(r_{\kb}^j:(\kb,j)\in\mathcal K)\subset\mathscr F_0,
\]
and $\mathscr G$ is independent of $B$. In particular, $\overrightarrow{\ub_0}$ is $\mathscr F_0$-measurable. Such a filtration can be obtained by the usual augmentation of $\bigl(\mathscr G\vee\sigma(B(s):0\leq s\leq t)\bigr)_{t\geq0}$.

We state the main result after specifying the solution class. Set
\begin{equation}\label{eqn-Z_infty}
 Z_\infty \coloneqq    C([0,\infty); \mathbb{L}^{2}(\dom)) \cap L^{4}_{\loc}([0,\infty); \mathbb{L}^{6}(\dom))\cap L^{2}_{\loc}([0,\infty); \mathbb{H}^{1,3}(\dom)).
\end{equation}

\begin{definition}[\textbf{Global weak solutions to (\ref{eqn-Burgers-intro})}]\label{def-solution}
Let $\ub_0\in\mathbb{L}^2(\dom)$ be a deterministic function and let $\overrightarrow{\ub_0}$ be its randomization defined by \eqref{eqn-randomization-u}. A stochastic process $\ub:[0,\infty)\times\Omega\to\mathbb{L}^2(\dom)$ is called a global weak solution to \eqref{eqn-Burgers-intro} if it is an $\mathbb{L}^2(\dom)$-valued continuous $\mathbb{F}$-progressively measurable process satisfying the following conditions:

(i) There exists $\Omega_0\subset\Omega$ with $\mathbb{P}(\Omega_0)=1$ such that, for every $\omega\in\Omega_0$,
\[
\ub(\cdot,\omega)\in Z_\infty.
\]

(ii) For every $\eta\in\mathbb{H}^1(\dom)$, the real-valued process $\langle\ub(t),\eta\rangle$ is a continuous semimartingale, and the following equality holds for every $t\geq0$, $\mathbb{P}$-a.s.:
\begin{align}\label{eq-global-weak-solution}
\langle\ub(t),\eta\rangle={}&\langle\overrightarrow{\ub_0},\eta\rangle-\nu\int_0^t\langle\Lambda\ub(s),\Lambda\eta\rangle\,ds-\int_0^t\langle(\ub(s)\cdot\nabla)\ub(s),\eta\rangle\,ds\nonumber\\
&+\int_0^t\langle\fb(s),\eta\rangle\,ds+b\int_0^t\langle\ub(s),\eta\rangle\circ dB(s).
\end{align}
\end{definition}

The following theorem gives global existence, uniqueness, continuous dependence and backward uniqueness.
\begin{theorem}\label{thm-main result}
Let $\ub_0 \in \mathbb{L}^{2}(\dom)$ be a deterministic function and let
$\overrightarrow{\ub_{0}}$ be its randomization given by \eqref{eqn-randomization-u}, using a sequence $\{r_{\kb}^j\}_{(\kb,j)\in\mathcal K}$ satisfying \autoref{A1}. Suppose that $f$ satisfies \autoref{A2}, i.e., for some $\delta>\frac{3}{2}$,
\[\fb \in L^{2}_{\loc}([0,\infty); \mathbb{H}^{\delta}(\dom)).
\]
Then there exists a unique global weak solution $\ub$ to (\ref{eqn-Burgers-intro}) in the sense of \autoref{def-solution}.

 Moreover, let $u_{0,n}\to u_0$ in $L^2(\mathcal O)$. Randomize all
initial data with the same multipliers in \eqref{eqn-randomization-u}, and use the same
$B$ and $f$. Then, for every $T>0$, the corresponding solutions
satisfy
\[
u_n\longrightarrow u
\quad\text{in probability in }C([0,T];L^2(\mathcal O)).
\]

Moreover, the backward uniqueness for equation (\ref{eqn-Burgers-intro}) holds as well. More precisely, let $\ub_0^1,\ub_0^2\in\mathbb L^2(\dom)$ be two deterministic functions, and let $\overrightarrow{\ub_0^1},\overrightarrow{\ub_0^2} \in L^2(\Omega;\mathbb L^2(\dom))$ be their randomizations defined by \eqref{eqn-randomization-u}, using the sequences $\{r_{\kb}^{j,1}\}_{(\kb,j)\in\mathcal K}$ and $\{r_{\kb}^{j,2}\}_{(\kb,j)\in\mathcal K}$, respectively, satisfying \autoref{A1}. Assume that both random initial values are $\mathscr F_0$-measurable.

Suppose that $u^1$ and $u^2$ are the corresponding global weak solutions to (\ref{eqn-Burgers-intro}), driven by the same $B$, $b$ and $f$.
Then, for every $T>0$, there exists $\Omega_T\in\mathscr F$ with $\mathbb P(\Omega_T)=1$ such that,
for every $\omega\in\Omega_T$,
\[
u^1(T,\omega)=u^2(T,\omega)
\quad\Longrightarrow\quad
u^1(t,\omega)=u^2(t,\omega),
\qquad t\in[0,T].
\]
\end{theorem}

The proof is given in Section~\ref{section-4}.

\section{Local existence of the solutions to (\ref{eqn-Burgers-intro.4})}\label{section-3}

 \subsection{Some lemmata}

If $X$ and $Y$ are some normed vector spaces both embedded into a common vector space, then
$X\cap Y$ is a normed vector space endowed with norm
\begin{equation}\label{eqb-X cap Y}
    \Vert u \Vert_{X\cap Y}:=\Vert u \Vert_{X}+\Vert u \Vert_{Y}, \;\; u \in X\cap Y.
\end{equation}

For $q\in[1,\infty]$ and $s\in\mathbb R$, define $L^q([0,T];\mathbb H^s(\dom))$ as the following space of equivalence classes of strongly measurable maps, where maps equal almost everywhere are identified:
\[
L^{q}\left([0,T] ; {\mathbb{H}}^{s}(\dom)\right)
:= \left\{ f:[0,T]\to\mathbb H^s(\dom): \Vert f\Vert_{L^q([0,T];\mathbb H^s(\dom))}<\infty \right\},
\]
where
\[
\Vert f\Vert_{L^q([0,T]; {\mathbb{H}}^{s}(\dom))} :=
\begin{cases}
\left( \displaystyle\int_0^T \Vert f(t)\Vert_{H^{s}}^q\, dt \right)^{1/q}, \;\; q < \infty, \\
\operatorname*{ess\,sup}_{t \in (0,T)} \Vert f(t)\Vert_{H^{s}}, \;\; q = \infty,
\end{cases}
\]
and $\Vert \cdot\Vert_{H^{s}}$ is defined in \eqref{eqn-norm H^s} and \eqref{eq-negative norm}.

 \begin{lemma}\label{lem-moments of rvs}
Assume that $\{r_{\kb}^j\}_{(\kb,j)\in\mathcal K}$ is a sequence of independent, zero-mean, real random variables satisfying
\begin{align}\label{eq-moments condition}
M_2:=\sup_{(\kb,j)\in\mathcal K}\mathbb E|r_{\kb}^j|^2<\infty.
\end{align}
Then, for every $\{c_{\kb}^j\}_{(\kb,j)\in\mathcal K}\in\ell^2(\mathcal K;\mathbb C)$, the series $\sum_{(\kb,j)\in\mathcal K}c_{\kb}^jr_{\kb}^j$ converges in $L^2(\Omega)$ and satisfies
\begin{align}\label{eqn-boundedness in L^p}
\Big\Vert\sum_{(\kb,j)\in\mathcal K}c_{\kb}^jr_{\kb}^j\Big\Vert_{L^2(\Omega)}
\leq M_2^{\frac12}\Vert c\Vert_{\ell^2}.
\end{align}
Moreover, let $T>0$ and let $E=\mathbb L^6(\dom)$, $E=L^4([0,T];\mathbb L^6(\dom))$ or $E=L^2([0,T];\mathbb L^3(\dom))$. Then, for every sequence $\{F_{\kb}^j\}_{(\kb,j)\in\mathcal K}\subset E$ satisfying
\[
\sum_{(\kb,j)\in\mathcal K}\Vert F_{\kb}^j\Vert_E^2<\infty,
\]
the series $\sum_{(\kb,j)\in\mathcal K}r_{\kb}^jF_{\kb}^j$ converges in $L^2(\Omega;E)$ and there exists a constant $c>0$ such that
\begin{align}\label{eqn-vector-valued-second-moment}
\Big\Vert\sum_{(\kb,j)\in\mathcal K}r_{\kb}^jF_{\kb}^j\Big\Vert_{L^2(\Omega;E)}
\leq cM_2^{\frac12}\Big(\sum_{(\kb,j)\in\mathcal K}\Vert F_{\kb}^j\Vert_E^2\Big)^{\frac12}.
\end{align}
The same conclusion holds for the corresponding finite-dimensional vector-valued and matrix-valued function spaces.
\end{lemma}

The scalar estimate follows from independence and zero means. For the vector-valued estimate, symmetrization, the Khintchine inequality and the mixed-norm Minkowski inequality give, for every finite $\mathcal K_N\subset\mathcal K$,
\[
\mathbb E\Big\Vert\sum_{(k,j)\in\mathcal K_N}r_k^jF_k^j\Big\Vert_E^2
\leq cM_2\sum_{(k,j)\in\mathcal K_N}\Vert F_k^j\Vert_E^2.
\]
Applying this to the tails gives convergence in $L^2(\Omega;E)$. The full proof is given in Appendix~\ref{app-second-moment-estimate}.

In particular, by the uniform bound $\Vert e_\kb^j\Vert_{L^6}\leq c$ in \autoref{thm-ONB}, applying \eqref{eqn-vector-valued-second-moment} with $E=\mathbb L^6(\dom)$ and $F_\kb^j=\ub_{0,\kb}^je_\kb^j$ yields
\begin{align}\label{eq-randomization-L6}
\mathbb E\Vert\overrightarrow{\ub_0}\Vert_{L^6}^2
\leq cM_2\sum_{(\kb,j)\in\mathcal K}|\ub_{0,\kb}^j|^2\Vert e_\kb^j\Vert_{L^6}^2
\leq cM_2\Vert\ub_0\Vert_{L^2}^2.
\end{align}
The randomized series therefore converges in $L^2(\Omega;\mathbb L^6(\dom))$. Since $\mathbb L^6(\dom)\hookrightarrow\mathbb L^2(\dom)$, this limit agrees with the randomization defined in \autoref{def-randomization-u}.

According to \cite[Theorem in Section 3.5.4 and Theorem 1 in Section 3.6.1]{Schmeisser+Triebel_1987}, for $q \in (1,\infty)$ and $s\geq 0$
\begin{equation}\label{eqn-rit}
\left(\mathbb{H}^{s}(\dom),\mathbb{H}^{s+2}(\dom)\right)_{1-\frac1q,q}= B^{s+2-\frac{2}{q}}_{2,q}.
    \end{equation}

The following Theorem can be found in \cite[Theorem 1.1]{GallaratiVeraar2017}.
\begin{theorem}\label{thm-maximal regularity}
Assume that $T>0$. Assume that $q \in (1,\infty)$ and $s\geq 0$.
Let us define an  operator $\Xi$  by
\begin{equation}\label{eqn-Xi}
    [\Xi f](t):=\int_0^t e^{-(t-s)A}f(s)\,ds, \;\; t \in [0,T]
\end{equation}
where $f \in L^{q}\left([0,T] ; {\mathbb{H}}^{s}(\dom)\right)$. Then
$\Xi f$ belongs to $L^{q}\left([0,T] ; {\mathbb{H}}^{s+2}(\dom)\right)$,
$\frac{d}{dt}\Xi f$ belongs to $L^{q}\left([0,T] ; {\mathbb{H}}^{s}(\dom)\right)$ and
$\Xi f \in C\left([0,T] ; B^{s+2-\frac{2}{q}}_{2,q}\right)$.

Moreover $\Xi$ is  a linear and bounded operator between the above spaces, i.e. there exists $C=C_T>0$ such that   the following estimate is satisfied
\begin{equation}\label{eqn-3.8}
\begin{aligned}
    \Vert \Xi f\Vert_{L^{q}\left([0,T] ; {\mathbb{H}}^{s+2}(\dom)\right)}
&+ \Vert \frac{d}{dt}\Xi f \Vert_{L^{q}\left([0,T] ; {\mathbb{H}}^{s}(\dom)\right)}
+\Vert \Xi f\Vert_{C([0,T] ; B^{s+2-\frac{2}{q}}_{2,q})}
\\
&\leq C\Vert f\Vert_{L^{q}\left([0,T]; {\mathbb{H}}^{s}(\dom)\right)}, \;\; f\in L^{q}\left([0,T] ; {\mathbb{H}}^{s}(\dom)\right).
\end{aligned}
\end{equation}
\end{theorem}

The following Lemma is a corollary to
\autoref{thm-maximal regularity}.
\begin{cor}\label{lemma-lplq}
Let $\zb$ satisfy
\[
\frac{d}{dt}\zb(t)+A \zb(t)= \fb(t), \;\; \zb(0)=0
\]
i.e. $\zb(t)=\int_{0}^{t}e^{-(t-s)A} \fb(s)\, ds$, $t \geq 0$. Then for every $T>0$, the following holds:
\begin{trivlist}
\item[(i)]
For every $\fb \in L^{\frac{4}{3}}([0,T]; \mathbb{L}^{2}(\dom))$, we have
\begin{align}
\zb \in C([0,T]; \mathbb{H}^\frac{1}{2}(\dom)),
\end{align}
and there exists a constant $c=c_T>0$ such that
\begin{align}\label{eqn-maximal}
\Vert  \zb  \Vert _{C([0,T]; \mathbb{H}^{\frac{1}{2}}(\dom))\cap L^{2}([0,T]; \mathbb{H}^{\frac{3}{2}}(\dom))}\leq c \Vert   \fb \Vert _{L^{\frac{4}{3}}([0,T]; \mathbb{L}^{2}(\dom))}.
\end{align}
\item[(ii)]
For every $\fb \in L^{2}([0,T]; \mathbb{H}^{1}(\dom))$, we have
\begin{align}
\zb \in C([0,T]; \mathbb{H}^{2}(\dom)),
\end{align}
and  there exists $c=c_T>0$ such that
\begin{align}\label{eqn-infty H2}
\Vert  \zb  \Vert _{C([0,T]; \mathbb{H}^{2}(\dom))}\leq c\Vert \fb\Vert_{ L^{2}([0,T]; \mathbb{H}^{1}(\dom))}.
\end{align}
\end{trivlist}
\end{cor}
This follows from \autoref{thm-maximal regularity}, the embedding $B^{1/2}_{2,4/3}\hookrightarrow\mathbb H^{1/2}$ and interpolation; the details are given in Appendix~\ref{app-linear-heat-proof}.

\subsection{Local existence of solutions to (\ref{eqn-Burgers-intro.4})}\label{subsection-3.1}

Assume that the sequence of random variables $\{r_{\kb}^j\}_{(k,j)\in \mathcal{K}}$ satisfies \autoref{A1}.
For a given element $\ub_{0}\in \mathbb{L}^{2}(\dom)$, we define
a random variable
\[
\overrightarrow{\ub_{0}}: \Omega  \to    \mathbb{L}^{2}(\dom)
\]
called a "randomization" of $u_0$, using a formula similar to
$(\ref{eqn-randomization-u})$, i.e.
\begin{align}\label{eqn-randomization of u_0}
\overrightarrow{\ub_{0}}(\omega) :=\sum_{(\kb,j)\in \mathcal{K}}r_{\kb}^j( \omega)\ub_{0,\kb}^j e_k^j.
\end{align}

In the next result we will use the following notation
\begin{equation}\label{eqn-K_T}
K_T := \left( \sum_{(\kb,j)\in \mathcal{K}} \vert \ub_{0,\kb}^j \vert ^2(1-\exxp{-2|\kb|^{2}T})
\right)^\frac{1}{2}.
\end{equation}

\begin{prop}\label{prop-ouz1} Assume that a sequence  $\{r_{\kb}^j\}_{(k,j)\in \mathcal{K}}$
satisfies  \autoref{A1} and $\ub_{0}\in \mathbb{L}^{2}(\dom)$. Let
$\overrightarrow{\ub_{0}}\in \mathbb{L}^2(\Omega; \mathbb{L}^{2}(\dom))$
be the randomization of $\ub_{0}$ defined by \autoref{def-randomization-u}.
Then there exists a  constant  $c>0$  such that for every  $T \in (0,1]$ the following assertions hold
\begin{align}\label{eq-623-11}
  \Vert e^{-\cdot A}\overrightarrow{\ub_{0}}\Vert_{L^{2}(\Omega; L^{4}([0,T];\mathbb{L}^{6}(\dom)))}
&\leq c\Vert \ub_{0}\Vert_{L^2}T^{\frac{1}{4}},\\
\label{eq-623-1}
\Vert \nabla e^{-\cdot A}\overrightarrow{\ub_{0}}\Vert_{L^{2}(\Omega; L^{2}([0,T];\mathbb{L}^{3}(\dom)))}
&\leq  c K_T,
\end{align}
where $K_T$ is defined in \eqref{eqn-K_T}.
\end{prop}
\begin{proof}
Let $E_1=L^4([0,T];\mathbb L^6(\dom))$ and $E_2=L^2([0,T];\mathbb L^3(\dom))$.
Recall that $\overrightarrow{\ub_0}$ is defined by \eqref{eqn-randomization of u_0}. Since $Ae_\kb^j=|\kb|^2e_\kb^j$, by \eqref{eqn-vector-valued-second-moment} in \autoref{lem-moments of rvs}, we infer that there exists a constant $c=c(M_2)>0$, independent of $T$ and $\ub_0$, such that
\begin{align}\label{eqn-3.261}
\Vert e^{-\cdot A}\overrightarrow{\ub_0}\Vert_{L^2(\Omega;E_1)}
&=\Big\Vert\sum_{(\kb,j)\in\mathcal K}r_\kb^j\ub_{0,\kb}^je^{-|\kb|^2\cdot}e_\kb^j\Big\Vert_{L^2(\Omega;E_1)}\leq c\Big(\sum_{(\kb,j)\in\mathcal K}|\ub_{0,\kb}^j|^2\Vert e^{-|\kb|^2\cdot}e_\kb^j\Vert_{E_1}^2\Big)^{\frac12}.
\end{align}

In view of \eqref{eq-ONB-real} and \eqref{eq-ONB-real-2}, we have $\Vert e_\kb^j\Vert_{L^6}\leq c$. Since $e^{-|\kb|^2t}\leq1$, we infer that
\begin{align}\label{eqn-Minkowski-3}
\sum_{(\kb,j)\in\mathcal K}|\ub_{0,\kb}^j|^2\Vert e^{-|\kb|^2\cdot}e_\kb^j\Vert_{E_1}^2\leq cT^{\frac12}\sum_{(\kb,j)\in\mathcal K}|\ub_{0,\kb}^j|^2=cT^{\frac12}\Vert\ub_0\Vert_{L^2}^2.
\end{align}
Therefore, \eqref{eqn-3.261} and \eqref{eqn-Minkowski-3} imply \eqref{eq-623-11}.

Next, we prove \eqref{eq-623-1}. In view of \eqref{eq-ONB-real} and \eqref{eq-ONB-real-2}, we infer that there exists a constant $c>0$, independent of $(\kb,j)\in\mathcal K$, such that $\Vert\nabla e_\kb^j\Vert_{L^3}\leq c|\kb|$.
Therefore,
\begin{align}
\Vert e^{-|\kb|^2\cdot}\nabla e_\kb^j\Vert_{E_2}^2
&=\int_0^T\Vert e^{-|\kb|^2t}\nabla e_\kb^j\Vert_{L^3}^2\,dt=\Vert\nabla e_\kb^j\Vert_{L^3}^2\int_0^Te^{-2|\kb|^2t}\,dt\nonumber\\
&\leq c|\kb|^2\int_0^Te^{-2|\kb|^2t}\,dt
=\frac{c}{2}(1-e^{-2|\kb|^2T}).
\end{align}

Since $0\leq 1-e^{-2|\kb|^2T}\leq1$ for every $(\kb,j)\in\mathcal K$, by the definition of $K_T$ and the Parseval identity \eqref{eqn-Parseval-real}, we infer that
\begin{align}
K_T^2
&=\sum_{(\kb,j)\in\mathcal K}|\ub_{0,\kb}^j|^2(1-e^{-2|\kb|^2T})\leq\sum_{(\kb,j)\in\mathcal K}|\ub_{0,\kb}^j|^2=\Vert\ub_0\Vert_{L^2}^2<\infty.
\end{align}
Therefore, by the above estimate, we deduce that
\begin{align}
\sum_{(\kb,j)\in\mathcal K}|\ub_{0,\kb}^j|^2\Vert e^{-|\kb|^2\cdot}\nabla e_\kb^j\Vert_{E_2}^2
&\leq c\sum_{(\kb,j)\in\mathcal K}|\ub_{0,\kb}^j|^2(1-e^{-2|\kb|^2T})=cK_T^2<\infty.
\end{align}

Therefore, by \eqref{eqn-vector-valued-second-moment}, we infer that
\begin{align}\label{eqn-Minkowski-4}
\Vert\nabla e^{-\cdot A}\overrightarrow{\ub_0}\Vert_{L^2(\Omega;E_2)}^2
&=\Big\Vert\sum_{(\kb,j)\in\mathcal K}r_\kb^j\ub_{0,\kb}^je^{-|\kb|^2\cdot}\nabla e_\kb^j\Big\Vert_{L^2(\Omega;E_2)}^2\leq cM_2\sum_{(\kb,j)\in\mathcal K}|\ub_{0,\kb}^j|^2\Vert e^{-|\kb|^2\cdot}\nabla e_\kb^j\Vert_{E_2}^2\leq cM_2K_T^2.
\end{align}
which implies \eqref{eq-623-1} holds.
The proof is complete.
\end{proof}

To proceed to the next proposition, we begin with the following lemma.
\begin{lemma}\label{lemma-moment -alpha}
Assume that $b$ is the given constant in \eqref{eqn-Burgers-intro}.
If
\[\alpha(t):=\exp(bB(t)),\;\; t>0,\]
then for all $m\in\mathbb R$ and $T>0$ the following moment estimates hold
\begin{align}\label{eq-moment -alpha}
\mathbb{E}\left[\sup\limits_{t\in [0,T]}\alpha^{m}(t)\right]\leq 2\exp{\Big(m^{2}b^{2}T\Big)}.
\end{align}
\end{lemma}
\begin{proof}
Since $e^ {\frac{1}{2} m^2 b^2 T}  < \infty$, by the Novikov condition\cite[Corollary III.5.13]{Karatzas+Shreve_1991} we have
\begin{align*}
\exp{\Big(mbB(t)-\frac{m^{2}}{2}b^{2}t\Big)} \mbox{ is a martingale}.
\end{align*}
Then by Jensen's inequality and the Doob maximum inequality\cite[Theorem I.3.8]{Karatzas+Shreve_1991}, we infer that
\begin{equation}\label{eqn-3.25}
 \begin{aligned}
\mathbb{E}\left[\sup\limits_{t\in [0,T]}\alpha^{m}(t)\right]&=\mathbb{E}\sup\limits_{t\in [0,T]}\exp{\Big(mbB(t)\Big)}\\
&\leq \mathbb{E}\sup\limits_{t\in [0,T]}\exp{\Big(mbB(t)-\frac{m^{2}}{2}b^{2}t\Big)}\exp{\Big(\frac{m^{2}}{2}b^{2}T\Big)}\\
&\leq \Bigg{(}\mathbb{E}\sup\limits_{t\in [0,T]}\exp{\Big(2mbB(t)-m^{2}b^{2}t\Big)}\Bigg{)}^{\frac{1}{2}}\exp{\Big(\frac{m^{2}}{2}b^{2}T\Big)}\\
&\leq 2\Bigg{(}\mathbb{E}\exp{\Big(2mbB(T)-m^{2}b^{2}T\Big)}\Bigg{)}^{\frac{1}{2}}\exp{\Big(\frac{m^{2}}{2}b^{2}T\Big)}.
\end{aligned}
\end{equation}
Note that $B(T)\sim N(0, T)$, by moment generating function for normal distribution we infer
\begin{align}\label{eqn-3.26}
\mathbb{E}\left[\exp\left(2mbB(T) - m^2b^2T\right)\right]
= \exp\left( \frac{1}{2} \cdot 4m^2b^2T  \right) \exp(- m^2b^2T)
= \exp\left(m^2b^2T\right).
\end{align}
Therefore, \eqref{eqn-3.25} and \eqref{eqn-3.26} implies \eqref{eq-moment -alpha} holds.
We complete the proof.
\end{proof}
\begin{prop}\label{prop-ouz2}
Suppose that \autoref{A2} holds. Define a process $\wb$ by the following formula
\begin{align}\label{def-w}
\wb(t) := \int_{0}^{t} e^{(t-s)\Delta}\alpha^{-1}(s)\fb(s)\,ds, \;\;\; t \in [0,\infty),
\end{align}
Then, for every $T>0$, $\wb\in L^6(\Omega;C([0,T];\mathbb H^2(\dom)))$. In particular,
\[
\wb \in L^{6}(\Omega; L^{4}([0,T];\mathbb{L}^{6}(\dom))) \cap L^{6}(\Omega; L^{2}([0,T];\mathbb{H}^{1,3}(\dom))),
\]
where the space $\mathbb{H}^{1,3}(\dom)$ is defined by
\begin{align}\label{eq-space H13}
\mathbb{H}^{1,3}(\dom):=\{\ub\in \mathbb{H}^{1}(\dom): \Lambda \ub \in \mathbb{L}^{3}(\dom)\}
\end{align}
and $\Vert \ub\Vert_{H^{1,3}}:=\Vert \ub\Vert_{L^2}+ \Vert \Lambda\ub\Vert_{L^3}$.
On the torus, by the boundedness of the periodic Riesz transforms on $L^3(\dom)$, this norm is equivalent to $\Vert\ub\Vert_{L^2}+\Vert\nabla\ub\Vert_{L^3}$

Moreover, there exists $c>0$ such that  for every $T\in (0,1]$,
\begin{align}\label{eq-646-2}
&\Vert \wb\Vert_{L^{6}(\Omega; L^{4}([0,T];\mathbb{L}^{6}(\dom)))}\leq cT^{\frac{1}{4}},\\
\label{eq-623-2}
&\Vert \nabla\wb\Vert_{L^{6}(\Omega; L^{2}([0,T];\mathbb{L}^{3}(\dom)))}\leq cT^{\frac{1}{2}}.
\end{align}
\end{prop}
\begin{proof}[Proof of \autoref{prop-ouz2}]
Fix $T>0$. By \autoref{A2}, $\fb\in L^2(0,T;\mathbb H^1(\dom))$. Apply \eqref{eqn-infty H2} pathwise to $\alpha^{-1}\fb$ and take sixth moments. By \eqref{eq-moment -alpha} with $m=-6$,
\[
\mathbb E\|\wb\|_{C([0,T];\mathbb H^2)}^6
\leq c_T^6\|\fb\|_{L^2(0,T;\mathbb H^1)}^6
      \mathbb E\sup_{0\leq s\leq T}\alpha^{-6}(s)
\leq 2c_T^6 e^{36b^2T}\|\fb\|_{L^2(0,T;\mathbb H^1)}^6.
\]
The Sobolev embeddings $\mathbb H^2\hookrightarrow\mathbb L^6$ and $\mathbb H^2\hookrightarrow\mathbb H^{1,3}$ give the asserted mixed norms. For $0<T\leq1$, restricting the $C([0,1];\mathbb H^2)$ bound to $[0,T]$ gives the factors $T^{1/4}$ and $T^{1/2}$ in \eqref{eq-646-2} and \eqref{eq-623-2}, with a constant independent of $T$.

For the corresponding Fourier estimates on $0<T\leq1$, let
\begin{align}
\phi_{\kb}^{\jb}(t):=\int_{0}^{t}\exxp{|\kb|^{2}s}\alpha^{-1}(s)\hat{\fb}_{\kb}^j(s)\,ds,\;\; t\in [0,T],\;\;(\kb,j)\in \mathcal{K}.
\end{align}
  Notice that for each $(\kb,j)\in \mathcal{K}$, we have
\begin{align}\label{eq-psi}
\psi_{\kb}^{\jb}(t):=\int_{0}^{t}\exxp{-|\kb|^{2}(t-s)}\alpha^{-1}(s)\hat{\fb}_{\kb}^\jb(s)e_\kb^\jb\,ds= \exxp{-|\kb|^{2}t}e_\kb^\jb\phi_{\kb}^{\jb}(t),\;\;t\in [0,T].
\end{align}
According to \autoref{lemma-moment -alpha} and the H\"{o}lder inequality, for $p\geq 1$, we infer that for $t\in [0,T](\subset [0,1])$, there exists a constant $c=c(p)$ such that for $(\kb,j)\in K_0\times\{1,2,3,4,5,6\}$
\begin{align*}
\mathbb{E}|\phi_{\kb}^{\jb}(t)|^{2p}&\leq \Big{(}\int_{0}^{t}  \exxp{|\kb|^{2}s}|\hat{\fb}_{\kb}^j(s)|\,ds \Big{)}^{2p}\mathbb{E}
\left(
\sup_{s\in[0,1]}|\alpha^{-1}(s)|^{2p}
\right)\\
&\leq \Big{(}\int_{0}^{t}  \exxp{|\kb|^{2}s}|\hat{\fb}_{\kb}^j(s)|\,ds \Big{)}^{2p}2\exxp{4p^{2}b^2}\\
&\leq 2\exxp{4p^{2}b^2}\Big{(}\int_{0}^{t} \exxp{2|\kb|^{2}s} \,ds \Big{)}^{p}\Big{(}\int_{0}^{t}  |\hat{\fb}_{\kb}^j(s)|^{2} \,ds\Big{)}^{p}\\
&\leq \frac{\exxp{2p|\kb|^{2}t+4p^2b^2}}{2^{p-1}|\kb|^{2p}}\Big{(}\int_{0}^{t}  |\hat{\fb}_{\kb}^j(s)|^{2} \,ds\Big{)}^{p}.
\end{align*}
Then by \eqref{eq-psi}, it further yields that
\begin{align*}
\Vert \psi_{\kb}^{\jb}(t)\Vert_{L^{2p}(\Omega)}^{2p} =\mathbb{E}|\exxp{-|\kb|^{2}t}e_k^j\phi_{\kb}^{\jb}(t)|^{2p} \leq \exxp{-2p|\kb|^{2}t}\mathbb{E}|\phi_{\kb}^{\jb}(t)|^{2p}\leq 2^{1-p} {|\kb|^{-2p}}\exxp{4p^2b^2}\Big{(}\int_{0}^{t}  |\hat{\fb}_{\kb}^j(s)|^{2}\, ds\Big{)}^{p}.
\end{align*}
When $p=3$ we get there exists $c=c(b)>0$

such that
\begin{align}\label{eq-moment estimate z}
\Vert \psi_{\kb}^{\jb}(t)\Vert_{L^{6}(\Omega)}\leq c {|\kb|^{-1}}\Big{(}\int_{0}^{t}  |\hat{\fb}_{\kb}^j(s)|^{2} \,ds\Big{)}^{\frac{1}{2}},\;\; k\in K_0.
\end{align}
For $(\kb,j)\in\{0\}\times\{1,2,3\}$, $t\in(0,1]$, by \autoref{lemma-moment -alpha} and the H\"{o}lder inequality we infer that there exists $c=c(b)>0$ such that
\begin{equation}\label{eq-moment estimate z-0}
\begin{aligned}
\Vert \psi_{0}^\jb(t)\Vert_{L^{6}(\Omega)}=&\Vert \int_{0}^{t}\alpha^{-1}(s)\hat{\fb}_{0}^j(s)e_0^\jb\,ds\Vert_{L^{6}(\Omega)}\leq \left(\mathbb{E}|\int_{0}^{t}\alpha^{-1}(s)\hat{\fb}_{0}^j(s)\,ds|^6\right)^\frac{1}{6}\\
&\leq \int_{0}^{t}  |\hat{\fb}_{\kb}^j(s)|\,ds \left(\mathbb{E} ( |\alpha^{-6}(t)|_{L^{\infty}([0,1])})\right)^\frac{1}{6}\\
&\leq 2^\frac{1}{6}e^{6b^2}\int_{0}^{t}  |\hat{\fb}_{\kb}^j(s)|\,ds \leq 2^\frac{1}{6}e^{6b^2}t^\frac{1}{2}\Big{(}\int_{0}^{t}  |\hat{\fb}_{\kb}^j(s)|^2\,ds \Big{)}^\frac{1}{2}\\
&\leq  c\Big{(}\int_{0}^{t}  |\hat{\fb}_{\kb}^j(s)|^2\,ds \Big{)}^\frac{1}{2}.
\end{aligned}
\end{equation}
Under \autoref{A2},
\begin{align}
 \fb \in L^{2}([0,1]; \mathbb{H}^{\delta}(\dom)), \;\; \delta>\frac{3}{2}.
\end{align}
Applying Minkowski's inequality, inequalities \eqref{eq-moment estimate z}, \eqref{eq-moment estimate z-0} and Young's inequality yields
\begin{align*}
&\Vert \wb\Vert_{L^{6}(\Omega; L^{4}([0,T];\mathbb{L}^{6}(\dom)))}=\Big\Vert \sum_{(\kb,j)\in \mathcal{K}}\psi_{\kb}^{\jb}\Big\Vert_{L^{6}(\Omega; L^{4}([0,T];\mathbb{L}^{6}(\dom)))}
\leq \sum_{(\kb,j)\in \mathcal{K}}\Big\Vert \psi_{\kb}^{\jb}\Big\Vert_{L^{4}([0,T];\mathbb{L}^{6}(\dom;L^{6}(\Omega)))}\\
&\leq  c\sum_{(\kb,j)\in K_0\times\{1,2,3,4,5,6\}} \Big\Vert   {|\kb|^{-1}}\Big{(}\int_{0}^{1}  |\hat{\fb}_{\kb}^\jb(s)|^{2} \,ds\Big{)}^{\frac{1}{2}}  \Big\Vert_{L^{4}([0,T];\mathbb{L}^{6}(\dom))}+c\sum_{j=1}^3 \Big\Vert   \Big{(}\int_{0}^{1}  |\hat{\fb}_{0}^\jb(s)|^{2} \,ds\Big{)}^{\frac{1}{2}}  \Big\Vert_{L^{4}([0,T];\mathbb{L}^{6}(\dom))}\\       &=cT^{\frac{1}{4}}\sum_{(\kb,j)\in K_0\times\{1,2,3,4,5,6\}}|\kb|^{-1}\Big{(}\int_{0}^{1}  |\hat{\fb}_{\kb}^\jb(s)|^{2} \,ds\Big{)}^{\frac{1}{2}}+cT^{\frac{1}{4}}\sum_{j=1}^3  \Big{(}\int_{0}^{1}  |\hat{\fb}_{0}^\jb(s)|^{2} \,ds\Big{)}^{\frac{1}{2}}  \\
&\leq  cT^{\frac{1}{4}}\sum_{\kb\in K_0}\frac{|\kb|^{-2}}
{1+|\kb|^{2\delta}}+cT^{\frac{1}{4}}
\sum_{(\kb,j)\in K_0\times\{1,2,3,4,5,6\}}\int_{0}^{1} (1+|\kb|^{2\delta}) |\hat{\fb}_{\kb}^\jb(s)|^{2} \,ds+cT^{\frac{1}{4}}\sum_{j=1}^3  \int_{0}^{1}  |\hat{\fb}_{0}^\jb(s)|^{2} \,ds+cT^{\frac{1}{4}}\\
&=  cT^{\frac{1}{4}}\sum_{\kb\in K_0}\frac{|\kb|^{-2}}
{1+|\kb|^{2\delta}}+cT^{\frac{1}{4}}\Vert \fb \Vert_{L^2([0,1];\mathbb{H}^\delta(\dom))}^2+cT^{\frac{1}{4}}
\leq cT^{\frac{1}{4}}.
\end{align*}

According to Minkowski's inequality, inequalities \eqref{eq-moment estimate z} and \eqref{eq-moment estimate z-0}, since $\delta>\frac{3}{2}$ we have
\begin{align*}
\Vert \nabla\wb\Vert_{L^{6}(\Omega; L^{2}([0,T];\mathbb{L}^{3}(\dom)))}&\leq c\Big\Vert \sum_{(\kb,j)\in \mathcal{K}}|\kb|\psi_{\kb}^{\jb}\Big\Vert_{L^{6}(\Omega; L^{2}([0,T];\mathbb{L}^{3}(\dom)))}
\leq c\sum_{(\kb,j)\in \mathcal{K}}\Big\Vert |\kb|\psi_{\kb}^{\jb}\Big\Vert_{L^{2}([0,T];\mathbb{L}^{3}(\dom;L^{6}(\Omega)))}\\
&\leq c\sum_{(\kb,j)\in \mathcal{K}}\Big\Vert \Big{(}\int_{0}^{t}  |\hat{\fb}^j_{\kb}|^{2} \,ds\Big{)}^{\frac{1}{2}}\Big\Vert_{L^{2}([0,T];\mathbb{L}^{3}(\dom))}\\
&\leq cT^{\frac{1}{2}}\sum_{(\kb,j)\in \mathcal{K}}\left(\frac{1}{1+|\kb|^{2\delta}}\right)^\frac{1}{2}\Big{(}\int_{0}^{1}  (1+|\kb|^{2\delta})|\hat{\fb}^j_{\kb}(s)|^{2} \,ds\Big{)}^{\frac{1}{2}}\\
&\leq cT^{\frac{1}{2}}\sum_{\kb\in K}\frac{1}{1+|\kb|^{2\delta}}+cT^{\frac{1}{2}}\sum_{(\kb,j)\in \mathcal{K}}\int_{0}^{1}  (1+|\kb|^{2\delta})|\hat{\fb}^j_{\kb}(s)|^{2} \,ds\\
&= cT^{\frac{1}{2}}\sum_{\kb\in K}\frac{1}{1+|\kb|^{2\delta}}+cT^{\frac{1}{2}}\Vert \fb\Vert_{L^2([0,1];\mathbb{H}^\delta(\dom))}^2\leq cT^{\frac{1}{2}}.
\end{align*}
The proof is complete.
\end{proof}

\begin{cor}\label{cor-w-continuous-L2}
Assume \autoref{A2}. Then
\[
\wb \in C([0,\infty);\mathbb{H}^2(\dom))\cap L^{4}_{\loc}([0,\infty); \mathbb{L}^{6}(\dom))\cap L^{2}_{\loc}([0,\infty); \mathbb{H}^{1,3}(\dom)),
\;\; \mathbb{P}\text{-a.s.}.
\]
\end{cor}
\begin{proof}
By (ii) of \autoref{lemma-lplq} we infer that
\[
\wb\in C([0,T];\mathbb{H}^{2}(\dom )),\;\;\mathbb{P}\text{-a.s.}.
\]
In view of \autoref{prop-ouz2}, the proof of \autoref{cor-w-continuous-L2} is complete.
\end{proof}

\begin{prop}\label{prop-ouz}
Let $\ub_{0}\in \mathbb{L}^{2}(\dom)$ and $\overrightarrow{\ub_{0}}\in \mathbb{L}^2(\Omega; \mathbb{L}^{2}(\dom))$
be the randomization of $\ub_{0}$ defined by \autoref{def-randomization-u}. Assume that \autoref{A1} and \autoref{A2} hold. Then there exists a constant $c>0$, independent of $\lambda>0$ and $T\in(0,1]$, such that the process

\begin{align*}
\mathbf{z}(t, \overrightarrow{\ub_{0}})&=e^{-tA}\overrightarrow{\ub_{0}}+\wb(t,\xb),\;\;t\geq 0
\end{align*}
 satisfies that for every $T\in (0,1]$,
\begin{align}\label{eq-646}
&\Vert \mathbf{z}(t, \overrightarrow{\ub_{0}})\Vert_{L^{2}(\Omega; L^{4}([0,T];\mathbb{L}^{6}(\dom)))}\leq c\Vert \ub_{0}\Vert_{L^2}T^{\frac{1}{4}}+cT^{\frac{1}{4}},\\
\label{eq-623}
&\Vert \nabla\mathbf{z}(t, \overrightarrow{\ub_{0}})\Vert_{L^{2}(\Omega; L^{2}([0,T];\mathbb{L}^{3}(\dom)))}\leq
c K_T +cT^{\frac{1}{2}},
\end{align}
where $K_T$ is defined in \eqref{eqn-K_T}.
Moreover, for every $\lambda>0$,
\begin{align}
\label{eq-combination}
&\mathbb{P}(\mathbb{S}_{\lambda, T, \ub_{0}})\leq c\lambda^{-2}\Big{(}\Vert \ub_{0}\Vert_{L^2}T^{\frac{1}{4}}+ T^{\frac{1}{4}}+T^{\frac{1}{2}}+ K_T \Big{)}^{2},
\end{align}
 where
\begin{align}\label{eq-S}
\mathbb{S}_{\lambda, T, \ub_{0}}=\{\omega\in \Omega: \Vert \mathbf{z}(t, \overrightarrow{\ub_{0}})\Vert_{L^{4}([0,T];\mathbb{L}^{6}(\dom))}
+ \Vert \nabla\mathbf{z}(t, \overrightarrow{\ub_{0}})\Vert_{L^{2}([0,T];\mathbb{L}^{3}(\dom))}\geq \lambda\}.
\end{align}
\end{prop}
\begin{proof}
By the definition of $ \mathbf{z}(t, \overrightarrow{\ub_{0}})$, in view of inequalities \eqref{eq-623-11} and \eqref{eq-646-2} we infer
\begin{align*}
\Vert \mathbf{z}(t, \overrightarrow{\ub_{0}})\Vert_{L^{2}(\Omega; L^{4}([0,T];\mathbb{L}^{6}(\dom)))}&\leq \Vert e^{-\cdot A}\overrightarrow{\ub_{0}}\Vert_{L^{2}(\Omega; L^{4}([0,T];\mathbb{L}^{6}(\dom)))}+\Vert\wb\Vert_{L^{2}(\Omega; L^{4}([0,T];\mathbb{L}^{6}(\dom)))}\\
&\leq  c\Vert \ub_{0}\Vert_{L^2}T^{\frac{1}{4}}+cT^{\frac{1}{4}},
\end{align*}
which implies \eqref{eq-646} holds. Moreover, by inequalities \eqref{eq-623-1} and \eqref{eq-623-2} we infer
\begin{align*}
\Vert \nabla \mathbf{z}(t, \overrightarrow{\ub_{0}})\Vert_{L^{2}(\Omega; L^{2}([0,T];\mathbb{L}^{3}(\dom)))}&\leq \Vert\nabla e^{-\cdot A}\overrightarrow{\ub_{0}}\Vert_{L^{2}(\Omega; L^{2}([0,T];\mathbb{L}^{3}(\dom)))}+\Vert\nabla \wb\Vert_{L^{2}(\Omega; L^{2}([0,T];\mathbb{L}^{3}(\dom)))}\\
&\leq c K_T+cT^{\frac{1}{2}}.
\end{align*}
 Hence, \eqref{eq-623} holds. By inequalities $(\ref{eq-646})$ and $(\ref{eq-623})$, we infer that there exists a constant  $c>0$ such that
\begin{equation}
\begin{aligned}
&\mathbb{E}\!\left[\left(\Vert \mathbf{z}(t, \overrightarrow{\ub_{0}})\Vert_{L^{4}([0,T];\mathbb{L}^{6}(\dom))}
+ \Vert \nabla\mathbf{z}(t, \overrightarrow{\ub_{0}})\Vert_{L^{2}([0,T];\mathbb{L}^{3}(\dom))}\right)^2\right]\\
&= \left\Vert \Vert \mathbf{z}(t, \overrightarrow{\ub_{0}})\Vert_{L^{4}([0,T];\mathbb{L}^{6}(\dom))}
+ \Vert \nabla\mathbf{z}(t, \overrightarrow{\ub_{0}})\Vert_{L^{2}([0,T];\mathbb{L}^{3}(\dom))}\right\Vert_{L^2(\Omega)}^2\\
&\leq \left(\Vert \mathbf{z}(t, \overrightarrow{\ub_{0}})\Vert_{L^{2}(\Omega; L^{4}([0,T];\mathbb{L}^{6}(\dom)))}+\Vert \nabla\mathbf{z}(t, \overrightarrow{\ub_{0}})\Vert_{L^{2}(\Omega; L^{2}([0,T];\mathbb{L}^{3}(\dom)))}\right)^2\\
&\leq c\Big{(}\Vert \ub_{0}\Vert_{L^2}T^{\frac{1}{4}}+ T^{\frac{1}{4}}+T^{\frac{1}{2}}+ K_T \Big{)}^{2}
\end{aligned}
\end{equation}
Finally, by Chebyshev's inequality we infer
\begin{equation}
\begin{aligned}
\mathbb{P}(\mathbb{S}_{\lambda, T, \ub_{0}})&\leq \lambda^{-2}\mathbb{E}\!\left[\left(\Vert \mathbf{z}(t, \overrightarrow{\ub_{0}})\Vert_{L^{4}([0,T];\mathbb{L}^{6}(\dom))}
+ \Vert \nabla\mathbf{z}(t, \overrightarrow{\ub_{0}})\Vert_{L^{2}([0,T];\mathbb{L}^{3}(\dom))}\right)^2\right]\\
&\leq c\lambda^{-2}\Big{(}\Vert \ub_{0}\Vert_{L^2}T^{\frac{1}{4}}+ T^{\frac{1}{4}}+T^{\frac{1}{2}}+K_T\Big{)}^{2}.
\end{aligned}
\end{equation}
The proof is complete.
\end{proof}

By the above assumptions,
$\mathbf z(t,\overrightarrow{\ub_0})
=e^{-tA}\overrightarrow{\ub_0}+\wb(t)$,
where $\wb$ is defined in \eqref{def-w}, is $\mathbb F$-adapted.

\begin{definition}[\textbf{Local weak solutions to (\ref{eqn-Burgers-intro.4})}]\label{def-3.1}
Suppose $\ub_{0}$ belongs to $\mathbb{L}^{2}(\dom)$, $\overrightarrow{\ub_{0}}$ is its randomization.
\par
(i)A local weak pathwise solution to (\ref{eqn-Burgers-intro.4}) is a pair $(\vb, \tau) $, where $\tau:\Omega \to (0,\infty)$ is
a stopping time, $\vb: \Omega \times [0,\tau] \to \mathbb{H}^{\frac{1}{2}}(\dom)$ is a continuous $(\mathscr F_t)_{t\ge0}$-adapted process, such that almost surely
\[\vb\in C([0,\tau]; \mathbb{H}^{\frac{1}{2}}(\dom))\cap L^{2}([0,\tau]; \mathbb{H}^{\frac{3}{2}}(\dom)),\;\;\vb(0)=0\]
and
\[
\partial_{t} \vb\in L^{1}([0,\tau]; \mathbb{H}^{-\frac{1}{2}}(\dom)),
\]
moreover, there exists $\Omega_0\in\mathscr F$ with $\mathbb P(\Omega_0)=1$ such that, for every $\omega\in\Omega_0$, for almost every $t\in[0,\tau(\omega)]$ and every
$\eta\in\mathbb H^{\frac12}(\dom)$,
\begin{align*}
&\Big{\langle}\partial_{t} \vb(t), \eta\Big{\rangle}+\Big{\langle}\Lambda^{\frac{3}{2}} \vb(t), \Lambda^{\frac{1}{2}}\eta\Big{\rangle}+\alpha(t)\Big{\langle}((\vb(t)+\mathbf{z}(t,\overrightarrow{\ub_{0}} ))\cdot \nabla) (\vb(t)+\mathbf{z}(t,\overrightarrow{\ub_{0}} )), \eta\Big{\rangle}=0.
\end{align*}

(ii) Weak pathwise solutions of (\ref{eqn-Burgers-intro.4}) are said to be pathwise unique up to a stopping time $\tau>0$ if given any pair of solutions $(\vb^{1}, \tau), (\vb^{2}, \tau)$ which coincide at $t=0$ on the event $\tilde{\Omega}=\{ \vb^{1}(0)=\vb^{2}(0)\}\subset \Omega$, then
 \[
 \mathbb{P}(\Vert \vb^{1}(t)-\vb^{2}(t) \Vert_{L^2}\mathds1_{\tilde{\Omega}}=0; t \in [0, \tau]  )=1.
 \]
\end{definition}

\begin{definition}[\textbf{Local mild solutions to (\ref{eqn-Burgers-intro.4})}]\label{def-3.2}
Suppose $\ub_{0}$ belongs to $\mathbb{L}^{2}(\dom)$ and $\overrightarrow{\ub_{0}}$ is its randomization.
\par
(i) A local mild pathwise solution to (\ref{eqn-Burgers-intro.4}) is a pair $(\vb, \tau) $ where $\tau: \Omega \to (0,\infty)$ is a stopping time with respect to $(\mathscr{F}_t)_{t\geq0}$
and
$\vb: [0,\tau] \times \Omega \mapsto \mathbb{H}^{\frac{1}{2}}(\dom)$
is an $(\mathscr{F}_t)$-adapted process
such that  for almost all $\omega\in \Omega$,
\[\vb\in C([0,\tau]; \mathbb{H}^{\frac{1}{2}}(\dom))\cap L^{2}([0,\tau]; \mathbb{H}^{\frac{3}{2}}(\dom)) \]
and
\begin{align*}
&\vb(t)=-\int_{0}^{t}\alpha(s)e^{(t-s)\Delta} \Big{(}((\vb(s)+\mathbf{z}(s, \overrightarrow{\ub_{0}}))\cdot \nabla) (\vb(s)+\mathbf{z}(s, \overrightarrow{\ub_{0}}))\Big{)}\,ds, \;\; t\in [0, \tau].
\end{align*}
(ii) Mild solutions of (\ref{eqn-Burgers-intro.4}) are said to be pathwise unique up to a stopping time $\tau>0$ if given any pair of solutions $(\vb^{1}, \tau), (\vb^{2}, \tau)$ which coincide at $t=0$ on the event $\tilde{\Omega}=\{ \vb^{1}(0)=\vb^{2}(0)\}\subset \Omega$, then
 \[
 \mathbb{P}(\Vert \vb^{1}(t)-\vb^{2}(t) \Vert_{L^2}\mathds{1}_{\tilde{\Omega}}=0; t \in [0, \tau]  )=1.
 \]
\end{definition}

\begin{definition}
\label{def-M^2X_T}
Define the space $X_T$ by
\[X_T:= C([0,T]; \mathbb{H}^{\frac{1}{2}}(\dom)) \cap L^{2}([0,T]; \mathbb{H}^{\frac{3}{2}}(\dom)),\]
endowed with the intersection  norm, see  \eqref{eqb-X cap Y}.
Denote by $\mathbb{M}^2\left(X_T\right)$, the Banach space of all (equivalence classes of) $\mathbb{H}^{\frac{1}{2}}(\dom)$-valued $\mathbb{F}$-progressively measurable processes $u$ whose trajectories almost surely belong to $X_T$, and satisfying
\begin{align}\label{M^2 X_T}
\Vert u\Vert_{\mathbb{M}^2\left(X_T\right)}^2:=\mathbb{E}\left(\sup _{t \in[0, T]}\Vert u(t) \Vert_{H^\frac12}^2+\int_0^T\Vert u(t)\Vert_{H^\frac32}^2\, d t\right)<\infty.
\end{align}
The norm in \eqref{M^2 X_T} is equivalent to the norm of $L^2(\Omega;X_T)$, and the progressively measurable processes form a closed subspace of this space.
Indeed, any convergent sequence has a subsequence converging almost surely in $X_T$, whose pointwise limit admits a progressively measurable version.
\end{definition}

To prove the existence and uniqueness of the local mild solution to
\eqref{eqn-Burgers-intro.4}, we apply the Banach fixed-point theorem
in a closed ball of the space of adapted processes whose trajectories belong to $X_T$.

\begin{theorem}\label{thm-local l2}
Assume that  $\ub_{0}\in \mathbb{L}^{2}(\dom)$ and $\overrightarrow{\ub_{0}}\in\mathbb{L}^2(\Om;\mathbb{L}^{2}(\dom))$ be defined  in Definition \ref{def-randomization-u}, using a sequence $\{r_{\kb}^j\}_{(\kb,j)\in\mathcal K}$ satisfying \autoref{A1}, and assume that $f$ satisfies \autoref{A2}.
Then there exists a local mild solution $(\vb, \tau) $ to (\ref{eqn-Burgers-intro.4}) according to Definition \ref{def-3.2}.
Any two such local solutions are indistinguishable on the intersection of their existence intervals.

\end{theorem}
\begin{proof}[Proof of \autoref{thm-local l2}]
Recall that $\mathbf{z}(t, \overrightarrow{\ub_{0}}):=e^{-tA}\overrightarrow{\ub_{0}}+\wb(t)$, where $\wb$ is defined by \eqref{def-w}, satisfies the following Cauchy problem:

\begin{align}
&d\mathbf{z}(t, \overrightarrow{\ub_{0}})-\Delta\mathbf{z}(t, \overrightarrow{\ub_{0}})dt=\alpha^{-1}(t)\fb(t)dt, \ \ \mathrm{on}\ [0,\infty)\times \dom,\\
&\mathbf{z}(0, \overrightarrow{\ub_{0}})=\overrightarrow{\ub_{0}}\in \mathbb{L}^{2}(\dom).\nonumber
\end{align}

Let $\theta:[0,\infty)\longrightarrow[0,1]$ be a non-increasing smooth function with compact support, satisfying $\theta=1$ on $[0,1]$, $\theta=0$ on $[2,\infty)$, and $\Vert\theta'\Vert_{L^\infty([0,\infty))}\leq2$.
For $R>0$, set $\theta_R(r):=\theta(r/R)$ for $r\geq0$, so that
\begin{align*}
 \theta_R(r)=
 \begin{cases}
 1,\qquad 0\leq r\leq R, \\
 0,\qquad r\geq 2R.
 \end{cases}
\end{align*}
We have the following easy lemma about $\theta_R$, see e.g. \cite[Lemma 4.3]{Brz+Millet_2014}.
\begin{lemma}\label{lemma-theta}
For every $x,y\in[0,\infty)$,
\begin{align}\label{eqn-theta-differential}
\vert\theta_R(x)-\theta_R(y)\vert \leq \frac{2}{R}\vert x-y\vert.
\end{align}
\end{lemma}

Define
\begin{align}\label{eq-Xzt}
\Vert \vb\Vert_{X_t} &:= \Vert \vb\Vert_{L^\infty([0,t];\mathbb H^{\frac12}(\dom))} + \Vert \vb\Vert_{L^2([0,t];\mathbb H^{\frac32}(\dom))}, \;\;t\geq 0  \\
\Vert \mathbf z\Vert_{Y_t} &:= \Vert \mathbf z\Vert_{L^4([0,t];\mathbb L^6(\dom))} + \Vert \nabla\mathbf z\Vert_{L^2([0,t];\mathbb L^3(\dom))},\;\;t\geq 0.
\end{align}

Consider the following truncated equation:
\begin{align}\label{eq-truncated}
\vb_R^N(t) = -\int_0^t \theta_R^2\big( \varepsilon_s\big) \alpha(s){1}_{\{\Vert\alpha\Vert_{L^\infty([0,s])}\leq N\}}e^{(t-s)\Delta} \Big( (\vb_R^N(s)+\mathbf z(s))\cdot\nabla \Big) (\vb_R^N(s)+\mathbf z(s))\,ds,\;\;t\geq 0 ,
\end{align}
where
\begin{align}
    \varepsilon_t:=\Vert \vb_R^N\Vert_{X_t} +\Vert \mathbf z\Vert_{Y_t},\;\;t\geq 0.
\end{align}
Note that
\begin{align}\label{eqn-truncated-vz}
\theta_R(\varepsilon_s)=\theta_R\big( \Vert \vb_R^N\Vert_{X_s} + \Vert \mathbf z\Vert_{Y_s} \big) =\begin{cases}
1,\;\;0\leq \Vert \vb_R^N\Vert_{X_s} + \Vert \mathbf z\Vert_{Y_s}\leq R\\
0,\;\; \Vert \vb_R^N\Vert_{X_s} + \Vert \mathbf z\Vert_{Y_s}\geq 2R.
\end{cases}
\end{align}

Define the map $L_R^N:\mathbb{M}^2\left(X_T\right)\to \mathbb{M}^2\left(X_T\right)$ by

\begin{align*}
    L_R^N(\vb_R^N):=-\int_0^t \theta_R^2\big( \varepsilon_s \big) \alpha(s){1}_{\{\Vert\alpha\Vert_{L^\infty([0,s])}\leq N\}}e^{(t-s)\Delta} \Big( (\vb_R^N(s)+\mathbf z(s))\cdot\nabla \Big) (\vb_R^N(s)+\mathbf z(s))\,ds.
\end{align*}

By the Sobolev embedding and interpolation inequalities, for every $g\in X_t$,
\begin{align}\label{eqn-3.52}
\Vert\nabla g\Vert_{L^2([0,t];\mathbb L^3(\dom))}
\leq c\Vert g\Vert_{X_t},
\end{align}
and
\begin{align}\label{eqn-3.54}
\Vert g\Vert_{L^4([0,t];\mathbb L^6(\dom))}
\leq c\Vert g\Vert_{L^\infty([0,t];\mathbb H^{\frac12})}^{\frac12}
\Vert g\Vert_{L^2([0,t];\mathbb H^{\frac32})}^{\frac12}
\leq c\Vert g\Vert_{X_t}.
\end{align}
Hence, $\Vert g\Vert_{Y_t}\leq c\Vert g\Vert_{X_t}$, and the H\"older inequality gives
\begin{align}\label{eqn-g nablah}
\Vert(g\cdot\nabla)h\Vert_{L^{\frac43}([0,t];\mathbb L^2(\dom))}
&\leq\Vert g\Vert_{L^4([0,t];\mathbb L^6)}
\Vert\nabla h\Vert_{L^2([0,t];\mathbb L^3)}
\leq\Vert g\Vert_{Y_t}\Vert h\Vert_{Y_t},
\end{align}
and, for $g,h\in X_t$,
\begin{align}\label{eqn-g nablah-1}
\Vert(g\cdot\nabla)h\Vert_{L^{\frac43}([0,t];\mathbb L^2(\dom))}
\leq c\Vert g\Vert_{X_t}\Vert h\Vert_{X_t}.
\end{align}

Since $\mathbf z$ and $\alpha$ are progressively measurable,
$L_R^N$ maps $\mathbb M^2(X_T)$ into itself.
At $t=0$, the running norms are defined by continuity; thus
$\Vert v\Vert_{X_0}=\Vert v(0)\Vert_{\mathbb H^{\frac12}(\dom)}$
and $\Vert\mathbf z\Vert_{Y_0}=0$.
The estimates proved in Appendix~\ref{app-truncated-map-estimates} yield
constants $c_1,c_2>0$, independent of $N,R$ and $T\in(0,1]$, such that
\begin{align}\label{eq-map-1}
\Vert L_R^N(\vb_R^N)\Vert_{\mathbb M^2(X_T)}
\leq c_1N\big(R^2+R\Vert\vb_R^N\Vert_{\mathbb M^2(X_T)}\big),
\end{align}
and
\begin{align}\label{eq-difference-1}
\Vert L_R^N(\vb_R^{N,1})-L_R^N(\vb_R^{N,2})\Vert_{\mathbb M^2(X_T)}
\leq c_2NR\Vert\vb_R^{N,1}-\vb_R^{N,2}\Vert_{\mathbb M^2(X_T)}.
\end{align}

Next we will use the contraction principle to consider the local existence of mild solution $\vb_R^N$.
Let us choose and fix $R$ such that
\[
0<R<\min\{\frac{1}{2c_1N},\frac{1}{2c_2N}\},
\]
then $R$ satisfies
\begin{align}\label{eq-lambda}
2c_1NR^2\leq R\;\; \text{and}\;\; c_2NR\leq \frac12.
\end{align}

Therefore, estimates \eqref{eq-map-1} and \eqref{eq-difference-1} imply that, for every  $\Vert \vb_R^{N,i}\Vert_{\mathbb{M}^2\left(X_{T}\right)}\leq R$, $i=1,2$
\begin{align*}
&\Vert L_R^N(\vb_R^{N,i})\Vert_{\mathbb{M}^2\left(X_{T}\right)}\leq 2c_1NR^2\leq R,\\
&\Vert L_R^N(\vb_R^{N,1})- L_R^N(\vb_R^{N,2}) \Vert_{\mathbb{M}^2\left(X_{T}\right)}\leq c_2NR\Vert \vb_R^{N,1}-\vb_R^{N,2}\Vert_{\mathbb{M}^2\left(X_{T}\right)}\leq \frac12\Vert \vb_R^{N,1}-\vb_R^{N,2}\Vert_{\mathbb{M}^2\left(X_{T}\right)},
\end{align*}
which implies the map $L$ is contractive on the ball of radius $R$ centered at zero in the space $\mathbb{M}^2\left(X_{T}\right)$.
Therefore, by the Banach Fixed-Point Theorem, the map $L$ has a unique fixed point $\vb_R^N\in \mathbb{M}^2\left(X_{T}\right)$, such that for $0\leq t\leq T$, $\vb_R^N$ satisfies truncated equation \eqref{eq-truncated}.

Define stopping times $\tau_R^N$ and $\sigma_N$ by
\begin{align*}
&\tau_R^N
:=
\inf\left\{
t\in [0,T]:
\Vert \vb_R^N\Vert_{X_t} +
\Vert \mathbf z\Vert_{Y_t}
\geq R
\right\}\wedge T,
\\
&\sigma_N
:=
\inf\left\{
t\in [0,T]:
\Vert\alpha\Vert_{L^\infty([0,t])}
\geq N
\right\}\wedge T.
\end{align*}
Note that if $\left\{
t\in [0,T]:
\Vert \vb_R^N\Vert_{X_t} +
\Vert \mathbf z\Vert_{Y_t}
\geq R
\right\}=\emptyset$, then $\tau_R^N=\infty\wedge T=T$ and if $\left\{
t\in [0,T]:
\Vert\alpha\Vert_{L^\infty([0,t])}
\geq N
\right\}=\emptyset$, then $\sigma_N=T$.
Define $\hat\tau:=\tau_R^N\wedge\sigma_N$, then $\hat\tau$ is still a stopping time.
Note that $\theta_R\big( \Vert \vb_R^N\Vert_{X_t} + \Vert \mathbf z\Vert_{Y_t}\big)=1$ on $[0,\tau_R^N]$. Define $\vb=\vb_R^N$ on $[0,\hat\tau]$, then $\vb$ satisfies
\begin{align}\label{eq-mild v}
\vb(t)
=
-\int_0^t
\alpha(s)e^{(t-s)\Delta}
\Big((\vb(s)+\mathbf z(s))\cdot\nabla\Big)
(\vb(s)+\mathbf z(s))\,ds .
\end{align}

In the following we aim to show that $\hat\tau>0,\;\;\mathbb{P}-\text{a.s.}$.

Define
\begin{align}\label{eq-Sc}
\mathbb{S}^{c} = \bigcup_{n=1}^{\infty} \mathbb{S}_{\lambda, \frac{1}{n}, \ub_{0}}^{c}.
\end{align}

By \eqref{eqn-K_T} and the dominated convergence theorem, we infer that
\begin{align}
K_T^2=\sum_{(\kb,j)\in\mathcal K}|\ub_{0,\kb}^j|^2(1-e^{-2|\kb|^2T})\to0,\qquad T\downarrow0,
\end{align}
where we have used $\sum_{(\kb,j)\in\mathcal K}|\ub_{0,\kb}^j|^2=\Vert\ub_0\Vert_{L^2}^2<\infty$. Therefore, by \eqref{eq-combination}, we infer that
\begin{align}
\mathbb P(\mathbb S_{\lambda,\frac1n,\ub_0})\to0,\qquad n\to\infty.
\end{align}
Since $\mathbb S_{\lambda,\frac1n,\ub_0}^c\uparrow\mathbb S^c$ as $n\to\infty$, by the continuity of probability measures from below, we deduce that
\begin{align}
\mathbb P(\mathbb S^c)=\lim_{n\to\infty}\mathbb P(\mathbb S_{\lambda,\frac1n,\ub_0}^c)=1.
\end{align}

For arbitrary $\omega \in \mathbb{S}^{c}$, there exists $n\geq 1$ such that $\omega \in \mathbb{S}_{\lambda, \frac{1}{n}, \ub_{0}}^{c}$.
By \autoref{prop-ouz},
for almost every $\omega\in\Omega$, there exists $T(\omega)>0$ such that
\[
\mathbf z(\cdot,\omega)\in L^4([0,T(\omega)];\mathbb L^6(\dom)),\;\;
\nabla\mathbf z(\cdot,\omega)\in L^2([0,T(\omega)];\mathbb L^3(\dom)).
\]
Therefore,
\[
\lim_{t\downarrow0}
\left(
\Vert \vb_R^N\Vert_{X_t}
+
\Vert \mathbf z\Vert_{Y_t}
\right)=\Vert \vb_R^N\Vert_{X_0}
+
\Vert \mathbf z\Vert_{Y_0}=0,\;\;\mathbb{P}-\text{a.s.},
\]
which implies $\tau_R^N>0$ almost surely.

Moreover, let us choose and fix $N>|\alpha(0)|=1$, we infer that $\sigma_N>0$ almost surely. Hence,
\[
\hat\tau=\tau_R^N\wedge\sigma_N>0,\;\;\mathbb{P}-\text{a.s.}.
\]

Since $\vb_R^N\in \mathbb{M}^2\left(X_{T}\right)$ is $\mathscr{F}_t$-adapted, $\vb=\vb_R^N$ on $[0,\hat\tau]$ and $\hat\tau$ is a stopping time, $\vb(\cdot\wedge \hat\tau)$ is $\mathscr{F}_t$-adapted.
Thus, $(\vb,\hat\tau)$ is a local mild pathwise solution to
$(\ref{eqn-Burgers-intro.4})$ in the sense of \autoref{def-3.2}.

It remains to prove pathwise uniqueness.
Let $(\vb^{1},\tau_{1})$ and $(\vb^{2},\tau_{2})$ be two local mild
solutions to \eqref{eqn-Burgers-intro.4} with $\vb^{1}(0)=\vb^{2}(0)=0$. Set
\[
\tau^\ast:=\tau_{1}\wedge\tau_{2},
\]
and
\[
\vb^\ast:=\vb^{1}-\vb^{2}\mbox{ on } [0,\tau^\ast].
\]
In the following we aim to prove that
\[
\vb^\ast(t)=0,\qquad 0\leq t\leq\tau^\ast,
\quad \mathbb P\text{-a.s.}.
\]

For every $t\in [0,\tau^\ast]$, since $\vb^{1}$ and $\vb^{2}$ satisfy \eqref{eq-mild v}, we infer that
\begin{align}\label{eq-pathwise-difference}
\vb^\ast(t)
=&-\int_{0}^{t}\alpha(s)e^{(t-s)\Delta}
\Big[
(\vb^\ast(s)\cdot\nabla)(\vb^{1}(s)+\mathbf z(s))
+((\vb^{2}(s)+\mathbf z(s))\cdot\nabla)\vb^\ast(s)
\Big]\,ds.
\end{align}

Let us choose and fix $\omega$  and $T<\tau^\ast(\omega)$. Then
\[
\vb^1,\vb^2\in X_T,
\]
and
\[
\mathbf z\in L^4([0,T];\mathbb L^6(\dom)),
\qquad
\nabla\mathbf z\in L^2([0,T];\mathbb L^3(\dom)).
\]
Moreover, by the continuity of $\alpha$,
\[
\Vert\alpha\Vert_{L^\infty([0,T])}<\infty.
\]

By \eqref{eqn-3.52} and \eqref{eqn-3.54}, we infer that
\[
\vb^1,\vb^2\in L^4([0,T];\mathbb L^6(\dom)),
\qquad
\nabla\vb^1,\nabla\vb^2
\in L^2([0,T];\mathbb L^3(\dom)).
\]

Hence, there exists $t_0\in(0,T]$ sufficiently small such that
\begin{align}\label{eq-small-pathwise}
&c\Vert\alpha\Vert_{L^\infty([0,T])}
\Big(
\Vert\vb^1+\mathbf z\Vert_{L^4([0,t_0];\mathbb L^6)}
+\Vert\nabla(\vb^1+\mathbf z)\Vert_{L^2([0,t_0];\mathbb L^3)}
\nonumber\\
&\qquad\qquad+
\Vert\vb^2+\mathbf z\Vert_{L^4([0,t_0];\mathbb L^6)}
+\Vert\nabla(\vb^2+\mathbf z)\Vert_{L^2([0,t_0];\mathbb L^3)}
\Big)
<\frac12.
\end{align}

By \autoref{lemma-lplq}, the H\"older inequality,
\eqref{eqn-3.52} and \eqref{eqn-3.54}, it follows from \eqref{eq-pathwise-difference} that
\begin{align*}
\Vert\vb^\ast\Vert_{X_{t_0}}
\leq&
c\Vert\alpha\Vert_{L^\infty([0,T])}
\Big(
\Vert\vb^1+\mathbf z\Vert_{L^4([0,t_0];\mathbb L^6)}
+\Vert\nabla(\vb^1+\mathbf z)\Vert_{L^2([0,t_0];\mathbb L^3)}
\\
&\qquad+
\Vert\vb^2+\mathbf z\Vert_{L^4([0,t_0];\mathbb L^6)}
+\Vert\nabla(\vb^2+\mathbf z)\Vert_{L^2([0,t_0];\mathbb L^3)}
\Big)
\Vert\vb^\ast\Vert_{X_{t_0}}
\\
\leq&\frac12\Vert\vb^\ast\Vert_{X_{t_0}},
\end{align*}
which implies
\[
\vb^\ast=0\qquad\mbox{on }[0,t_0].
\]

Set
\[
\sigma:=\sup\left\{
s<\tau^\ast(\omega):
\vb^1(t)=\vb^2(t)\mbox{ for all }t\in[0,s]
\right\}.
\]
By the continuity of $\vb^1$ and $\vb^2$, we have $\vb^1(\sigma)=\vb^2(\sigma)$ for every $\sigma<\tau^\ast(\omega)$. Applying the above local uniqueness argument with initial time $\sigma$, we can extend the equality to $[\sigma,\sigma+h]$ for some $h>0$, which contradicts the definition of $\sigma$. Therefore,
\[
\sigma=\tau^\ast(\omega).
\]
Hence,
\[
\vb^1(t)=\vb^2(t),\;\;t\in[0,\tau^\ast(\omega)).
\]
If $\tau^\ast(\omega)<\infty$, the equality at $\vb^1(\tau^\ast)=\vb^2(\tau^\ast)$ follows from continuity.
This completes the proof of \autoref{thm-local l2}.

\end{proof}

\begin{lemma}\label{eq-equivalence}
Under the conditions of \autoref{thm-local l2}, the existence of local mild solutions and local weak solutions to (\ref{eqn-Burgers-intro.4}) are equivalent.
\end{lemma}
\begin{proof}
By \eqref{eqn-g nablah} and \eqref{eqn-g nablah-1}, the nonlinear term belongs to $L^{\frac43}([0,\tau];\mathbb L^2(\dom))$.
The semigroup differentiation formula in \autoref{semigroup-derivative} therefore gives the weak formulation for every mild solution.
Conversely, subtracting the corresponding Duhamel integral from a weak solution gives a solution of the heat equation with zero initial value, and hence the difference vanishes.
The details are given in Appendix~\ref{app-mild-weak-proof}.
\end{proof}

In view of \autoref{thm-local l2} and \autoref{eq-equivalence}, we establish the existence of local weak solutions to (\ref{eqn-Burgers-intro.4}).
\begin{theorem}\label{thm-local l2-weak}
 Assume that  $\ub_{0}\in \mathbb{L}^{2}(\dom)$ and $\overrightarrow{\ub_{0}}\in\mathbb{L}^2(\Omega;\mathbb{L}^{2}(\dom))$ be defined  in Definition \ref{def-randomization-u}, using a sequence $\{r_{\kb}^j\}_{(\kb,j)\in\mathcal K}$ satisfying \autoref{A1}, and assume that $f$ satisfies \autoref{A2}. Then there exists a local weak solution $(\vb, \tau) $ to (\ref{eqn-Burgers-intro.4}) according to Definition \ref{def-3.1}.
Any two such local solutions are indistinguishable on the intersection of their existence intervals.
\end{theorem}

\section{Global existence and uniqueness of the solutions to (\ref{ee-o}) and (\ref{eqn-Burgers-intro}) }\label{section-4}
In this section, we will pursue the global well-posedness of random Burgers equation (\ref{ee-o}).

\subsection{Local existence and regularities of solutions to (\ref{ee-o})}\label{subsection-({ee-o})}

Let us begin with the definition of local weak and local strong solutions.

\begin{definition}[\textbf{Local weak solutions to (\ref{ee-o})}]\label{def-3.1o}
Assume that  $\ub_{0}\in \mathbb{L}^{2}(\dom)$ and $\overrightarrow{\ub_{0}}$ be its randomization defined in Definition \ref{def-randomization-u}.
\begin{trivlist}
\item[(i)] A local weak pathwise solution to  equation  (\ref{ee-o}) with initial data $\tilde\vb(0)=\overrightarrow{\ub_{0}}$ is a pair $(\tilde{\vb}, \tau) $ where $\tau:\Omega \to (0,\infty)$ is
a stopping time, $\tilde{\vb}: \Omega \times [0,\tau] \to \mathbb{L}^{2}(\dom)$ is a continuous $(\mathscr F_t)_{t\ge0}$-adapted process such that almost surely
\begin{align}\label{eqn-tilde v}
\tilde{\vb}\in  L^{2}([0,\tau]; \mathbb{H}^{1,3}(\dom))\cap  L^4([0,\tau]; \mathbb{L}^{6}(\dom))\cap C([0,\tau]; \mathbb{L}^{2}(\dom)),
\end{align}
where $\mathbb{H}^{1,3}(\dom)$ is defined in \eqref{eq-space H13}, and
\begin{equation}\label{eqn-d_t tilde v}
\partial_t \tilde{\vb}
\in L^1([0,\tau];\mathbb{H}^{-1}(\dom)),
\end{equation}
moreover, for almost every $t\in [0,\tau]$ and for all $\eta\in \mathbb{H}^{1}(\dom)$,
\begin{align}\label{eqn-local weak solution-tilde v}
&\Big{\langle} \frac{d\tilde{\vb}}{dt}, \eta\Big{\rangle}+\Big{\langle}\Lambda \tilde{\vb}, \Lambda\eta\Big{\rangle}+\alpha(t)\Big{\langle}((\tilde{\vb}+\wb)\cdot \nabla) (\tilde{\vb}+\wb), \eta\Big{\rangle}=0,\ \mathbb{P}\text{-a.s.}.
\end{align}
\item[(ii)] Weak pathwise solutions of (\ref{ee-o}) are said to be pathwise unique up to a random positive time $\tau>0$ if given any pair of solutions $(\tilde{\vb}^{1}, \tau), (\tilde{\vb}^{2}, \tau)$ which coincide at $t=0$ on the event $\tilde{\Omega}=\{ \tilde{\vb}^{1}(0)=\tilde{\vb}^{2}(0)\}\subset \Omega$, then
 \[
 \mathbb{P}(\Vert \tilde{\vb}^{1}(t)-\tilde{\vb}^{2}(t) \Vert_{L^2}\mathds{1}_{\tilde{\Omega}}=0; t \in [0, \tau]  )=1.
 \]
 \end{trivlist}
\end{definition}
\begin{definition}[\textbf{Local strong solutions to (\ref{ee-o})}]\label{def-3.1ostr}
Assume the initial data $\tilde{\vb}_0\in \mathbb{H}^{i}(\dom)$, $i=1, 2$.
A local pathwise strong solution  to (\ref{ee-o})  is a pair $(\tilde{\vb}, \tau )$, where $\tau: \Omega \to (0,\infty)$ is a positive
stopping time and $\tilde{\vb}: \Omega \times [0,\tau] \to \mathbb{H}^{i}(\dom)$ is a continuous process such that almost surely
\[\tilde{\vb} \in C([0,\tau]; \mathbb{H}^{i}(\dom))\cap L^2([0,\tau]; \mathbb{H}^{i+1}(\dom)),\]

and is adapted, satisfies \eqref{eqn-tilde v}-\eqref{eqn-local weak solution-tilde v}.

The definition of uniqueness of $(\tilde{\vb}, \tau )$ is the same as that in \autoref{def-3.1o}.
\end{definition}
From \autoref{thm-local l2} and \autoref{eq-equivalence}, we have
\begin{theorem}\label{thm-local weak solution}
Assume that $\ub_{0}\in \mathbb{L}^{2}(\dom)$ and let
$\overrightarrow{\ub_{0}}\in \mathbb{L}^2(\Omega; \mathbb{L}^{2}(\dom))$
be the randomization of $\ub_{0}$ defined in \autoref{def-randomization-u}, using a sequence $\{r_{\kb}^j\}_{(\kb,j)\in\mathcal K}$ satisfying \autoref{A1}, and assume that $f$ satisfies \autoref{A2}. Then, there exists a local weak solution $(\tilde{\vb}, \tau) $ to (\ref{ee-o}) in the sense of \autoref{def-3.1o}.
Any two such local solutions with the same initial value are indistinguishable on the intersection of their existence intervals.
\end{theorem}
\begin{proof}

By  \autoref{thm-local l2},  there exists a local weak solution $(\vb, \tau)$   to (\ref{eqn-Burgers-intro.4}) in the sense of \autoref{def-3.1}, i.e. the equation with initial data $0$, such that almost surely
\begin{align}
    &\vb \in C([0,\tau]; \mathbb{H}^{\frac{1}{2}}(\dom))\cap L^{2}([0,\tau]; \mathbb{H}^{\frac{3}{2}}(\dom)), \\
    &\partial_{t} \vb\in L^{1}([0,\tau]; \mathbb{H}^{-\frac{1}{2}}(\dom)).
    \label{eqn-d_tv}
\end{align}

By \eqref{eqn-3.52}, \eqref{eqn-3.54} and the embeddings $\mathbb H^{1/2}(\dom)\hookrightarrow\mathbb L^3(\dom)\hookrightarrow\mathbb L^2(\dom)$, we obtain
\begin{align}\label{eq-v-CL2-L2H13-L4L6}
\vb \in  C([0,\tau]; \mathbb{L}^{2}(\dom))\cap L^{2}([0,\tau]; \mathbb{H}^{1,3}(\dom))\cap  L^{4}([0,\tau]; \mathbb{L}^{6}(\dom)),\;\;\mathbb{P}\text{-a.s.}
\end{align}
and
\begin{align}\label{eq-lamv-L2H1/2-L2L3}
\Lambda\vb \in L^{2}([0,\tau]; \mathbb{H}^{\frac{1}{2}}(\dom))\subset L^{2}([0,\tau]; \mathbb{L}^{3}(\dom)),\;\;\mathbb{P}\text{-a.s.}.
\end{align}

Since $\ub_{0}\in \mathbb{L}^2(\dom)$, by \autoref{rem- randomization} we have
$\overrightarrow{\ub_{0}}\in \mathbb{L}^2(\dom)$, $\mathbb{P}$-a.s..

In view of \autoref{prop-ouz1},  almost surely
\begin{align}\label{eq-etDelta}
e^{ \cdot \,\Delta } \overrightarrow{\ub_{0}}\in C([0,\tau]; \mathbb{L}^{2}(\dom))\cap  L^{4}([0,\tau]; \mathbb{L}^{6}(\dom))\cap L^{2}([0,\tau]; \mathbb{H}^{1,3}(\dom)),
\end{align}
where $\mathbb{H}^{1,3}(\dom)$ is defined in \eqref{eq-space H13}.

 Let  us introduce an auxiliary process
\begin{align}\label{eq-tidev-v}
\tilde{\vb}(t):=\vb(t)+e^{t\,\Delta } \overrightarrow{\ub_{0}},\;\; t\in [0,\tau].
\end{align}
In view of \eqref{eq-v-CL2-L2H13-L4L6} and \eqref{eq-etDelta}, we get \eqref{eqn-tilde v} holds.

We will show that property \eqref{eqn-d_t tilde v} holds. Let us choose and fix $\omega\in \mathbb{S}^{c}$, where $\mathbb{S}^{c}$ is defined in \autoref{thm-local l2}.
In view of definition \eqref{eq-tidev-v} of $\tilde{\vb}$
\begin{align}\label{eq-dt-tildev}
\partial_t \tilde{\vb}
= \partial_t \vb + \partial_t\bigl(e^{t\Delta}\overrightarrow{\ub_0}\bigr)
= \partial_t \vb + \Delta e^{t\Delta}\overrightarrow{\ub_0}.
\end{align}
Since $\mathbb{H}^{-\frac12}(\dom)\hookrightarrow \mathbb{H}^{-1}(\dom)$, from \eqref{eqn-d_tv} we infer that
\[
\partial_t \vb \in L^1([0,\tau];\mathbb{H}^{-1}(\dom)).
\]
For the second term in \eqref{eq-dt-tildev}, for every $\varphi\in \mathbb{H}^1(\dom)$ with
$\Vert\varphi\Vert_{H^1}\leq 1$,
\[
\langle \Delta e^{t\Delta}\overrightarrow{\ub_0},\varphi \rangle
= -\langle \nabla e^{t\Delta}\overrightarrow{\ub_0}, \nabla\varphi \rangle, \;\; t >0.
\]
Hence, by the H\"older inequality
\[
\Vert\Delta e^{t\Delta}\overrightarrow{\ub_0}\Vert_{H^{-1}}
\leq \Vert\nabla e^{t\Delta}\overrightarrow{\ub_0}\Vert_{L^2} , \;\; t >0.
\]
Therefore, inequality \eqref{eq-623-1} in \autoref{prop-ouz1} implies
\[
\Delta e^{ \cdot \Delta}\overrightarrow{\ub_0}
\in L^1([0,\tau];\mathbb{H}^{-1}(\dom)).
\]
Hence, by \eqref{eq-dt-tildev}, we  infer that property \eqref{eqn-d_t tilde v} is true and
in order to conclude the proof  that $(\tilde{\vb}, \tau)$ is a local weak solution to (\ref{ee-o}) with $\tilde{\vb}_0=\overrightarrow{\ub_{0}}$ in the sense of \autoref{def-3.1o},  it is sufficient to show that $(\tilde{\vb}, \tau)$ satisfies \eqref{eqn-local weak solution-tilde v}. For this purpose let us choose and fix  $\eta\in \mathbb{H}^{1}(\dom)$.
Let us temporarily assume that $ \eta \in D(\Lambda^{2})$. Then, since $\Delta$ and $e^{t\,\Delta } $ are self-adjoint on $\mathbb{L}^2(\dom)$, in view of \eqref{eq-tidev-v} we infer that
\begin{align}\label{eq-4.2}
\langle  \Lambda^{\frac{3}{2}} \vb, \Lambda^{\frac{1}{2}}\eta \rangle= \langle \vb,\Lambda^{2} \eta \rangle=\langle \tilde{\vb}-e^{t\,\Delta }\overrightarrow{\ub_{0}}, \Lambda^{2}\eta \rangle=\langle  \Lambda \tilde{\vb},  \Lambda\eta\rangle - \langle \Lambda^{2} e^{t\,\Delta }\overrightarrow{\ub_{0}}, \eta \rangle,  \text{ a.e. }t\in(0,\tau).
\end{align}

Moreover,
\begin{equation}\label{eq-4.3}
\begin{aligned}
\langle \partial_{t}\vb, \eta\rangle &=\partial_{t}\langle \vb, \eta\rangle= \partial_{t}\langle \tilde{\vb}, \eta  \rangle -\partial_{t}\langle e^{t\,\Delta }\overrightarrow{\ub_{0}} , \eta \rangle= \langle \partial_{t}\tilde{\vb}, \eta  \rangle- \langle \overrightarrow{\ub_{0}} , \partial_{t}e^{t\,\Delta }\eta \rangle\\
&=\langle \partial_{t}\tilde{\vb}, \eta  \rangle- \langle \overrightarrow{\ub_{0}} , \Delta e^{t\,\Delta }\eta \rangle=\langle \partial_{t}\tilde{\vb}, \eta  \rangle+\langle \Lambda^{2} e^{t\,\Delta }\overrightarrow{\ub_{0}}, \eta \rangle,\;\text{ a.e. }t\in(0,\tau).
\end{aligned}
\end{equation}

In view of equalities \eqref{eq-z} and \eqref{eq-tidev-v}, we infer that
\begin{align}\label{eq-v-tildev}
\vb(t)+\mathbf{z}(t,\overrightarrow{\ub_{0}} )=\tilde{\vb}(t)+\wb(t),\;\; t\in [0,\tau].
\end{align}
Since $(\vb, \tau)$ is a local weak solution to (\ref{eqn-Burgers-intro.4}), i.e. for almost every $t \in [0,\tau]$
\begin{align*}
&\Big{\langle}\partial_{t} \vb(t), \eta\Big{\rangle}+\Big{\langle}\Lambda^{\frac{3}{2}} \vb(t), \Lambda^{\frac{1}{2}}\eta\Big{\rangle}+\alpha(t)\Big{\langle}((\vb(t)+\mathbf{z}(t,\overrightarrow{\ub_{0}} ))\cdot \nabla) (\vb(t)+\mathbf{z}(t,\overrightarrow{\ub_{0}} )), \eta\Big{\rangle}=0,
\end{align*}
therefore, by equalities \eqref{eq-4.2}, \eqref{eq-4.3} and \eqref{eq-v-tildev}, for $\eta\in D(\Lambda^{2})$, we infer
\begin{align*}
&\Big{\langle} \partial_{t}\tilde{\vb}, \eta\Big{\rangle}+\Big{\langle}\Lambda \tilde{\vb}, \Lambda\eta\Big{\rangle}+\Big{\langle}\alpha(t)((\tilde{\vb}+\wb)\cdot \nabla) (\tilde{\vb}+\wb), \eta\Big{\rangle}=0,\;\;a.e.t\in[0,\tau].
\end{align*}
Since
$D(\Lambda^{2})$ is dense in $\mathbb{H}^{1}(\dom)$, the above equality holds for $\eta \in \mathbb{H}^{1}(\dom)$. In view of \autoref{def-3.1o}, $\tilde{\vb}$ is a weak solution to (\ref{ee-o}).

The following is devoted to the uniqueness of local weak solution $(\tilde{\vb}, \tau) $. For $i=1,2$, let $\ub_{0,i}\in \mathbb{L}^{2}(\dom)$ and $\overrightarrow{\ub_{0,i}}$ be the randomization of $ \ub_{0,i}$ such that \autoref{A1} holds. Assume that $(\tilde{\vb}^{i}, \tau)$, $i=1,2$, are local weak solutions to (\ref{ee-o}) with initial data $\overrightarrow{\ub_{0,i}}$, respectively. Set
\begin{align}\label{eq-vec v}
\vec{\vb}= \tilde{\vb}^{1}- \tilde{\vb}^{2},\mbox{ on }[0,\tau].
\end{align}
In view of property \eqref{eqn-tilde v} and \eqref{eqn-d_t tilde v} in \autoref{def-3.1o}, $\vec{\vb}\in L^2([0,\tau];\mathbb{H}^1(\dom))\cap C([0,\tau];\mathbb{L}^2(\dom))$ and
$\partial_t\vec{\vb}\in L^1([0,\tau];\mathbb{H}^{-1}(\dom))$. Moreover, $\vec{\vb}$ satisfies the following equation on $[0,\tau]\times\dom$:
\begin{align}\label{eq-diff-v}
\partial_{t} \vec{\vb} -\Delta \vec{\vb} +\alpha (\vec{\vb}\cdot \nabla) (\tilde{\vb}^{1}+\wb ) +\alpha((\tilde{\vb}^{2}+\wb )\cdot \nabla) \vec{\vb}=0.
\end{align}

In view of \eqref{eqn-tilde v} and
\autoref{prop-ouz2}, by the H\"older inequality,
we infer that
\begin{align*}
&\alpha\Big[
(\vec{\vb}\cdot\nabla)(\tilde{\vb}^{1}+\wb)
+((\tilde{\vb}^{2}+\wb)\cdot\nabla)\vec{\vb}
\Big]\in
L^{\frac43}([0,\tau];\mathbb{L}^{2}(\dom))
\subset L^1([0,\tau];\mathbb{L}^{2}(\dom)).
\end{align*}
Let $P_m$ be the orthogonal Fourier projection
onto a finite-dimensional subspace of
$\mathbb{L}^{2}(\dom)$.

Applying $P_m$ to \eqref{eq-diff-v}, since $P_m$
commutes with $\Delta$, we infer that
\begin{align*}
&\partial_t P_m\vec{\vb}
-\Delta P_m\vec{\vb}+\alpha P_m\Big[
(\vec{\vb}\cdot\nabla)(\tilde{\vb}^{1}+\wb)
+((\tilde{\vb}^{2}+\wb)\cdot\nabla)\vec{\vb}
\Big]=0.
\end{align*}
Taking the inner product with $P_m\vec{\vb}$ in
$\mathbb{L}^2(\dom)$ and integrating by parts,
we infer that
\begin{align*}
&\frac12\frac{d}{dt}
\Vert P_m\vec{\vb}(t)\Vert_{L^2}^2
+\Vert\nabla P_m\vec{\vb}(t)\Vert_{L^2}^2=-\alpha(t)\Big\langle
(\vec{\vb}(t)\cdot\nabla)
(\tilde{\vb}^{1}(t)+\wb(t))+
((\tilde{\vb}^{2}(t)+\wb(t))\cdot\nabla)
\vec{\vb}(t),P_m\vec{\vb}(t)
\Big\rangle,
\end{align*}
for almost every $t\in[0,\tau]$.
Consequently, integrating over $[0,t]$, we obtain
\begin{align*}
\frac12\Vert P_m\vec{\vb}(t)\Vert_{L^2}^2
+\int_0^t\Vert\nabla P_m\vec{\vb}(s)\Vert_{L^2}^2\,ds
=&\frac12\Vert P_m\vec{\vb}(0)\Vert_{L^2}^2
-\int_0^t\alpha(s)\Big\langle
(\vec{\vb}(s)\cdot\nabla)
(\tilde{\vb}^{1}(s)+\wb(s))\nonumber\\
&+
((\tilde{\vb}^{2}(s)+\wb(s))\cdot\nabla)
\vec{\vb}(s),P_m\vec{\vb}(s)
\Big\rangle\,ds.
\end{align*}

Since
\[
P_m\vec{\vb}\rightarrow\vec{\vb}
\quad\mbox{in }
C([0,\tau];\mathbb{L}^{2}(\dom))
\cap L^2([0,\tau];\mathbb{H}^{1}(\dom)),
\]
letting $m\rightarrow\infty$, we deduce that
\begin{align*}
&\frac12\Vert\vec{\vb}(t)\Vert_{L^2}^{2}
+\int_0^t\Vert\vec{\vb}(s)\Vert_1^2\,ds=\frac12\Vert\vec{\vb}(0)\Vert_{L^2}^{2}-\int_0^t\alpha(s)
\Big\langle
(\vec{\vb}\cdot\nabla)(\tilde{\vb}^{1}+\wb)
+((\tilde{\vb}^{2}+\wb)\cdot\nabla)\vec{\vb},
\vec{\vb}
\Big\rangle\,ds.
\end{align*}
Hence, the corresponding differential
energy identity holds for almost every
$t\in[0,\tau]$.

Therefore,
\begin{align*}
\frac{1}{2}\frac{d}{dt}\Vert \vec{\vb}\Vert_{L^2}^{2}+\Vert \vec{\vb}\Vert_{1}^{2}&\leq \alpha\Vert \nabla (\tilde{\vb}^{1}+\wb )\Vert_{L^3}\Vert \vec{\vb}\Vert_{L^2}\Vert \vec{\vb}\Vert_{L^6}+\alpha\Vert \tilde{\vb}^{2}+\wb\Vert_{L^6}\Vert \nabla \vec{\vb}\Vert_{L^2} \Vert \vec{\vb}\Vert_{L^3}\\
&\leq c\alpha\Vert \nabla (\tilde{\vb}^{1}+\wb )\Vert_{L^3}\Vert \vec{\vb}\Vert_{L^2}\Vert \vec{\vb}\Vert_{H^1}+c\alpha\Vert \tilde{\vb}^{2}+\wb\Vert_{L^6}\Vert \nabla \vec{\vb}\Vert_{L^2}\Vert \vec{\vb}\Vert_{L^2}^{\frac{1}{2}}\Vert \vec{\vb}\Vert_{H^{1}}^{\frac{1}{2}}.
\end{align*}
Recall that $\Vert\cdot\Vert_{H^1}^2=\Vert\cdot\Vert_{L^2}^2+\Vert\cdot\Vert_{1}^2$, see \eqref{eq-homogeneous norm}. By the Young inequality, we infer that there exists a constant $\eps \in (0,1]$ such that
\begin{align*}
\frac{1}{2}\frac{d}{dt}\Vert \vec{\vb}\Vert_{L^2}^{2}+\Vert \vec{\vb}\Vert_{1}^{2}
&\leq c\alpha\Vert \nabla (\tilde{\vb}^{1}+\wb )\Vert_{L^3}\Vert \vec{\vb}\Vert_{L^2}(\Vert \vec{\vb}\Vert_{L^2}+\Vert \vec{\vb}\Vert_{1})+c\alpha\Vert \tilde{\vb}^{2}+\wb\Vert_{L^6}\Vert \nabla \vec{\vb}\Vert_{L^2}\Vert \vec{\vb}\Vert_{L^2}^{\frac{1}{2}}(\Vert \vec{\vb}\Vert_{L^2}^{\frac{1}{2}}+\Vert \vec{\vb}\Vert_{1}^{\frac{1}{2}})\\
&\leq \eps\Vert \vec{\vb}\Vert_{1}^{2}+ c\alpha\Vert \nabla (\tilde{\vb}^{1}+\wb )\Vert_{L^3}\Vert \vec{\vb}\Vert_{L^2}^{2}+c\alpha^2\Vert \nabla (\tilde{\vb}^{1}+\wb )\Vert_{L^3}^2\Vert \vec{\vb}\Vert_{L^2}^{2}+c\alpha^{4}\Vert \tilde{\vb}^{2}+\wb\Vert_{L^6}^{4}\Vert \vec{\vb}\Vert_{L^2}^{2}\\
&\;\;+c\alpha^{2}\Vert \tilde{\vb}^{2}+\wb\Vert_{L^6}^{2}\Vert \vec{\vb}\Vert_{L^2}^{2}\\
&=\eps\Vert \vec{\vb}\Vert_{1}^{2}+g\Vert \vec{\vb}\Vert_{L^2}^{2},
\end{align*}
where the function $g:[0,\tau]\to \mathbb{R}$ is defined by
\[g(t):= c\alpha(t)\Vert \nabla (\tilde{\vb}^{1}+\wb )(t)\Vert_{L^3}+c\alpha^2(t)\Vert \nabla (\tilde{\vb}^{1}+\wb )(t)\Vert_{L^3}^2+c\alpha^{4}(t)\Vert \tilde{\vb}^{2}(t)+\wb(t)\Vert_{L^6}^{4}+c\alpha^{2}(t)\Vert \tilde{\vb}^{2}(t)+\wb(t)\Vert_{L^6}^{2}.\]
Therefore,
\[
\frac{d}{dt}\Vert \vec{\vb}(t)\Vert_{L^2}^{2}\leq 2g(t)\Vert \vec{\vb}(t)\Vert_{L^2}^{2}.
\]
In view of property $(\ref{eqn-tilde v})$ and \autoref{prop-ouz2}, we infer that $\int_{0}^{\tau} g(t)\,dt<\infty $. By the Gronwall inequality,
\begin{align*}
\Vert \vec{\vb}(t)\Vert_{L^2}^{2}\leq \Vert \vec{\vb}(0)\Vert_{L^2}^{2}\exp\Big{(}2\int_{0}^{\tau}g(t)\,dt   \Big{)} .
\end{align*}
Since $\vec{\vb}(0)=\overrightarrow{\ub_{0,1}}- \overrightarrow{\ub_{0,2}}$, we infer that if $\overrightarrow{\ub_{0,1}}= \overrightarrow{\ub_{0,2}}$, then for $t\in [0,\tau]$
\[
\Vert\tilde{\vb}^{1}(t)- \tilde{\vb}^{2}(t)\Vert_{L^2}^{2}=\Vert \vec{\vb}(t)\Vert_{L^2}^{2}=0,
\]
which implies the uniqueness of the local weak solutions to (\ref{ee-o}).

The proof of \autoref{thm-local weak solution} is complete.
\end{proof}

Assume \autoref{A2} holds. In view of \autoref{cor-w-continuous-L2},
\[
\wb \in C([0,\infty);\mathbb{H}^2(\dom))\cap L^{4}_{\loc}([0,\infty); \mathbb{L}^{6}(\dom))\cap L^{2}_{\loc}([0,\infty); \mathbb{H}^{1,3}(\dom)),
\;\; \mathbb{P}\text{-a.s.}.
\]
Motivated by the local $H^1$ argument in \cite[Proposition 3.1]{Dong+Guo+Wu+Zhou_2023}, we prove the following result for \eqref{ee-o}.
\begin{prop}\label{prop-tildev-local strong}
 Assume that $\tilde\vb_{0}\in \mathbb{H}^{i}(\dom)$, $i=1,2$ and \autoref{A2} holds. Then there exists a unique local strong solution  $(\tilde{\vb}, \tau_i)$  to (\ref{ee-o}) in the sense of \autoref{def-3.1ostr} such that almost surely
\[
\tilde{\vb} \in C([0,\tau_i]; \mathbb{H}^{i}(\dom))\cap L^2([0,\tau_i]; \mathbb{H}^{i+1}(\dom)).
\]
Moreover, the local strong pathwise solution $\tilde\vb$ to (\ref{ee-o}) is Lipschitz continuous with respect to the initial data $\tilde{\vb}_0$ in $\mathbb{H}^i(\dom)$.

For every $T,M>0$, there exists
\[
\delta_0:= \delta_0\left(i,M,\Vert\alpha\Vert_{L^\infty([0,T+1])},
\Vert\wb\Vert_{C([0,T+1];\mathbb{H}^2(\dom))}
\right)>0,\;\;\mathbb{P}\text{-a.s.}
\]
such that, for every $t_0\in[0,T]$ and
$\vb^0\in\mathbb H^i(\dom)$ satisfying $\Vert \vb^0\Vert_{H^i}\leq M$,
equation \eqref{ee-o} with $\tilde{\vb}(t_0)=\vb^0$ admits a unique strong solution on $[t_0,t_0+\delta_0]$.
Moreover, there exists a constant $L\geq1$, depending only on $i$, $M$, $\Vert\alpha\Vert_{L^\infty([0,T+1])}$ and $\Vert\wb\Vert_{C([0,T+1];\mathbb H^2(\dom))}$, such that the solutions $\tilde{\vb}^{\,1}$ and $\tilde{\vb}^{\,2}$ starting at $t_0$ from $a_1,a_2\in\mathbb H^i(\dom)$ with $\Vert a_1\Vert_{H^i},\Vert a_2\Vert_{H^i}\leq M$ satisfy
\begin{align}\label{eq-uniform-local-Lipschitz}
\Vert\tilde{\vb}^{\,1}-\tilde{\vb}^{\,2}\Vert_{C([t_0,t_0+\delta_0];\mathbb H^i(\dom))}
\leq L\Vert a_1-a_2\Vert_{H^i}.
\end{align}
The constants $\delta_0\in(0,1]$ and $L$ can be chosen independently of $t_0\in[0,T]$ and of $a_1,a_2$ in that ball.

\end{prop}

\begin{proof}
Following the proof of
\cite[Proposition 3.1 and Lemma 3.1]{Dong+Guo+Wu+Zhou_2023},
we infer the local existence, uniqueness and Lipschitz
continuity with respect to the initial data for
equation \eqref{ee-o}.

Indeed, for $i=1,2$, let $I=[0,h]$, $h\in(0,1]$, and set
\[
X_I^i:=C(I;\mathbb H^i(\dom))
\cap L^2(I;\mathbb H^{i+1}(\dom)).
\]
For every $u\in X_I^i$, let us define
\[
F(u):=\alpha((u+\wb)\cdot\nabla)(u+\wb).
\]
By the H\"older inequality and the Sobolev embedding theorem,
we infer that there exists a constant $c>0$ such that
\begin{align*}
&\Vert F(u)-F(v)\Vert_{L^2(I;\mathbb H^{i-1}(\dom))}\leq ch^{\gamma_i}\Vert\alpha\Vert_{L^\infty(I)}
\Big(\Vert u\Vert_{X_I^i}+\Vert v\Vert_{X_I^i}
+2\Vert\wb\Vert_{C(I;\mathbb H^2(\dom))}\Big)
\Vert u-v\Vert_{X_I^i},
\end{align*}
where $\gamma_1=\frac14$ and $\gamma_2=\frac12$.

Therefore, by \autoref{thm-maximal regularity},
choosing $h>0$ sufficiently small and applying the
Banach fixed point theorem, we obtain a unique local
strong solution to \eqref{ee-o}.
Moreover, the same estimate implies the local Lipschitz
continuity with respect to the initial data in
$\mathbb H^i(\dom)$.

Let us choose and fix $T,M>0$ and $t_0\in[0,T]$. By the change of
variable $t=t_0+s$, equation \eqref{ee-o} can be written as
\begin{align*}
&\partial_s\tilde{\vb}
-\Delta\tilde{\vb}+\alpha(t_0+s)
\left((\tilde{\vb}+\wb(t_0+s))\cdot\nabla\right)(\tilde{\vb}+\wb(t_0+s))=0,\\
&\tilde{\vb}(0)=\vb^0.
\end{align*}
For $s\in[0,1]$, we infer that
\begin{align*}
&\sup_{t_0\in[0,T]}
\Vert\alpha(t_0+\cdot)\Vert_{L^\infty([0,1])}
\leq
\Vert\alpha\Vert_{L^\infty([0,T+1])},\\
&\sup_{t_0\in[0,T]}
\Vert\wb(t_0+\cdot)\Vert_{C([0,1];\mathbb{H}^2(\dom))}
\leq
\Vert\wb\Vert_{C([0,T+1];\mathbb{H}^2(\dom))}.
\end{align*}
Therefore,
we infer that for every $t_0\in[0,T]$ and $\Vert\vb^0\Vert_{\mathbb{H}^i(\dom)}\leq M$,  there exists
\[
\delta_0
=
\delta_0\left(
i,M,
\Vert\alpha\Vert_{L^\infty([0,T+1])},
\Vert\wb\Vert_{C([0,T+1];\mathbb{H}^2(\dom))}
\right)>0
\]
such that
\[
\tilde{\vb}
\in
C([t_0,t_0+\delta_0];\mathbb H^i(\dom))
\cap
L^2([t_0,t_0+\delta_0];\mathbb H^{i+1}(\dom)).
\]

In the above fixed-point argument, the radius of the ball in $X_I^i$, where $I=[t_0,t_0+\delta_0]$, can be chosen using only $i$, $M$ and the two coefficient bounds on $[0,T+1]$. Decreasing $\delta_0\in(0,1]$ if necessary, the contraction constant can therefore be taken to be at most $\frac12$, uniformly in $t_0$ and in the initial data in the ball of radius $M$. Subtracting the two mild equations and applying the heat semigroup estimate and the above nonlinear estimate, we infer that
\[
\Vert\tilde{\vb}^{\,1}-\tilde{\vb}^{\,2}\Vert_{X_I^i}
\leq C\Vert a_1-a_2\Vert_{H^i}+\frac12\Vert\tilde{\vb}^{\,1}-\tilde{\vb}^{\,2}\Vert_{X_I^i}.
\]
Absorbing the last term yields \eqref{eq-uniform-local-Lipschitz} with $L=\max\{1,2C\}$, independently of $t_0$ and the two initial data in that ball.

Notice that $\delta_0$ is obtained for each fixed
$\omega$ and depends only on $i$, $M$,
$\Vert\alpha\Vert_{L^\infty([0,T+1])}$
and $\Vert\wb\Vert_{C([0,T+1];\mathbb H^2(\dom))}$.
Therefore, $\delta_0$ is a pathwise lower bound
for the existence time and is not necessarily
a stopping time with respect to the filtration
$\{\mathscr F_t\}_{t\geq0}$.
The proof of \autoref{prop-tildev-local strong} is complete.
\end{proof}

\begin{prop}\label{prop-tildev-adapted}
For $i=1,2$, assume that $\tilde\vb_{0}$ is $\mathscr F_0$-measurable and belongs to $\mathbb{H}^{i}(\dom)$ almost surely  and \autoref{A2} holds. Then the local strong solution $(\tilde{\vb}, \tau_i)$ is adapted with respect to $(\mathscr F_t)_{t\ge0}$.

Moreover, let $\sigma$ be an almost surely finite $(\mathscr F_t)_{t\geq0}$-stopping time and let $a$ be an $\mathscr F_\sigma$-measurable $\mathbb H^i(\dom)$-valued random variable.
Then there exist a stopping time $\tau>\sigma$ almost surely and an adapted strong solution on $[\sigma,\tau]$ with value $a$ at $\sigma$, obtained by replacing $\alpha(s),\wb(s)$ by $\alpha(\sigma+s),\wb(\sigma+s)$ and using the filtration $\mathscr G_s:=\mathscr F_{\sigma+s}$.
\end{prop}

\begin{proof}
Fix $i=1,2$, define
\[
X_t^i
:=
C([0,t];\mathbb H^i(\dom))
\cap
L^2([0,t];\mathbb H^{i+1}(\dom)),
\;\; t>0,
\]
endowed with the norm
\[
\|u\|_{X_t^i}
:=
\|u\|_{C([0,t];\mathbb H^i(\dom))}
+
\|u\|_{L^2([0,t];\mathbb H^{i+1}(\dom))}.
\]

For $N\in\mathbb N$, define a stopping time
\[
\rho_N
:=
\inf\left\{
t\geq0:
\alpha(t)+\alpha(t)^{-1}+\|\wb(t)\|_{\mathbb H^2}
\geq N
\right\}\wedge N.
\]
Note that
\[
\rho_N\uparrow\infty,
\qquad \mathbb P\text{-a.s.}
\]

Recall that for fixed $R>0$, the function $\theta_R(r)=\theta(r/R)$ defined above satisfies
\begin{align*}
\theta_R(r)=\begin{cases}
    1,\qquad 0\leq r\leq R,\\
    0,\qquad r\geq 2R.
\end{cases}
\end{align*}

For fixed $R,N\in\mathbb N$, consider the approximated equation
\begin{align}
&\tilde{\vb}^{R,N,0}(t)
=
e^{t\Delta}\tilde{\vb}_0,\label{eq-approximated v0}\\
&\tilde{\vb}^{R,N,n+1}(t)
=
e^{t\Delta}\tilde{\vb}_0-
\int_0^t
e^{(t-s)\Delta}
\theta_R\left(
\|\tilde{\vb}^{R,N,n}\|_{X_s^i}
\right)
\alpha(s\wedge\rho_N)\nonumber\\
&\;\;\;\;\;\times
\Big(
(\tilde{\vb}^{R,N,n}(s)+\wb(s\wedge\rho_N))
\cdot\nabla
\Big)
(\tilde{\vb}^{R,N,n}(s)+\wb(s\wedge\rho_N))
\,ds .\label{eq-approximated vn}
\end{align}
We first prove that every $\tilde{\vb}^{R,N,n}$ is adapted. The
initial approximation is adapted since $\tilde{\vb}_0$ is
$\mathscr F_0$-measurable. Suppose that
$\tilde{\vb}^{R,N,n}$ is adapted. Then the process
\[
s\longmapsto
\|\tilde{\vb}^{R,N,n}\|_{X_s^i}
\]
is adapted and continuous. Since $\alpha(\cdot\wedge\rho_N)$ and
$\wb(\cdot\wedge\rho_N)$ are progressively measurable, the integrand
in the definition of $\tilde{\vb}^{R,N,n+1}$ is progressively
measurable. Consequently, its Bochner integral over $[0,t]$ is
$\mathscr F_t$-measurable, and hence
$\tilde{\vb}^{R,N,n+1}(t)$ is $\mathscr F_t$-measurable. Therefore,
by induction, $\tilde{\vb}^{R,N,n}$ is adapted for every
$n\in\mathbb N$.

{
Let us choose and fix $R,N\in\mathbb N$, $T>0$ and $i=1,2$.
In view of the definition of $\rho_N$, we infer that
\[
\Vert\alpha(\cdot\wedge\rho_N)\Vert_{L^\infty([0,T])}
+\Vert\wb(\cdot\wedge\rho_N)\Vert_{C([0,T];\mathbb H^2(\dom))}
\leq 2N.
\]

For every $I=[t_0,t_0+h]\subset[0,T]$, $h\in(0,1]$, recall that $X_I^i:=C(I;\mathbb H^i(\dom))\cap
L^2(I;\mathbb H^{i+1}(\dom))$.
Let $G_{R,N}(u)$ denote the truncated nonlinear term
in \eqref{eq-approximated vn}, i.e.
\begin{align}\label{eq-GRN}
G_{R,N}(u)(t)
:=\theta_R\left(\Vert u\Vert_{X_t^i}\right)
\alpha(t\wedge\rho_N)
\left((u(t)+\wb(t\wedge\rho_N))\cdot\nabla\right)
(u(t)+\wb(t\wedge\rho_N)),
\quad t\in[0,T].
\end{align}

Following the proof of \autoref{prop-tildev-local strong},
let $u,v\in X_T^i$ coincide on $[0,t_0]$.
By the Lipschitz continuity of $\theta_R$, we infer that
\begin{align*}
&\left|\theta_R\left(\Vert u\Vert_{X_t^i}\right)
-\theta_R\left(\Vert v\Vert_{X_t^i}\right)\right|\leq \Vert\theta_R'\Vert_{L^\infty([0,\infty))}
\Vert u-v\Vert_{X_t^i}
\leq c_R\Vert u-v\Vert_{X_I^i},
\qquad t\in I.
\end{align*}

Moreover, since $\theta_R(r)=0$ for $r\geq2R$,
we have $\Vert u\Vert_{X_t^i}<2R$ whenever
$\theta_R(\Vert u\Vert_{X_t^i})\neq0$.
The same property holds for $v$.
}
{
Let us choose and fix $R,N\in\mathbb N$ and $T>0$.
For $I=[t_0,t_0+h]\subset[0,T]$, $h\in(0,1]$,
let $u,v\in X_T^i$ coincide on $[0,t_0]$.
By the Lipschitz continuity of $\theta_R$, we infer that
\begin{align*}
&\left|\theta_R(\Vert u\Vert_{X_t^i})
-\theta_R(\Vert v\Vert_{X_t^i})\right|
\leq \frac{c}{R}\Vert u-v\Vert_{X_I^i},
\qquad t\in I.
\end{align*}
Moreover, if $\theta_R(\Vert u\Vert_{X_t^i})\neq0$ then $\Vert u\Vert_{X_t^i}<2R$.

Following the proof of \autoref{prop-tildev-local strong},
by considering the cut-off exit times and applying
the H\"older inequality and the Sobolev embedding theorem,
we infer that
\begin{align}\label{eq-adapted-cutoff-Lipschitz}
&\Vert G_{R,N}(u)-G_{R,N}(v)\Vert_{L^2(I;\mathbb H^{i-1}(\dom))}
\leq C_{R,N}h^{\gamma_i}\Vert u-v\Vert_{X_I^i},
\end{align}
where $\gamma_1=\frac14$ and $\gamma_2=\frac12$.

By \autoref{thm-maximal regularity} and
\eqref{eq-adapted-cutoff-Lipschitz}, we infer that
\begin{align*}
\Vert\Phi(u)-\Phi(v)\Vert_{X_I^i}
&\leq C_0C_{R,N}h^{\gamma_i}
\Vert u-v\Vert_{X_I^i}\leq\frac12\Vert u-v\Vert_{X_I^i},
\end{align*}
provided that $C_0C_{R,N}h^{\gamma_i}\leq\frac12$.

Therefore, by the Banach fixed point theorem and
successive continuation, we obtain a unique solution
\[
\tilde{\vb}^{R,N}\in
C([0,T];\mathbb H^i(\dom))
\cap L^2([0,T];\mathbb H^{i+1}(\dom)).
\]
Since the Picard approximations are adapted and
converge to the corresponding fixed point in $X_I^i$,
we infer that $\tilde{\vb}^{R,N}$ is adapted on $[0,T]$.
}

Define
\[
\tau_{R,N}
:=
\rho_N\wedge
\inf\left\{
t\geq0:
\|\tilde{\vb}^{R,N}\|_{X_t^i}\geq R
\right\}.
\]
The process
\[
t\longmapsto
\|\tilde{\vb}^{R,N}\|_{X_t^i}
\]
is adapted, continuous and nondecreasing. Therefore,
$\tau_{R,N}$ is an $(\mathscr F_t)_{t\geq0}$-stopping time. Note that $\theta_R\left(
\|\tilde{\vb}^{R,N}\|_{X_s^i}
\right)=1$ on
$[0,\tau_{R,N})$. Hence,
$\tilde{\vb}^{R,N}$ satisfies equation \eqref{ee-o} on
$[0,\tau_{R,N})$.

Set
\[
\sigma_k:=\tau_{k,k},
\qquad
\tilde{\vb}^{k}:=\tilde{\vb}^{k,k},
\qquad k\in\mathbb N.
\]
By pathwise uniqueness in
\autoref{prop-tildev-local strong},
\[
\sigma_k\leq\sigma_{k+1}
\]
and
\[
\tilde{\vb}^{k+1}(t)
=
\tilde{\vb}^{k}(t),
\qquad
0\leq t<\sigma_k,
\quad\mathbb P\text{-a.s.}
\]
Define
\[
\widehat{\tau}_i
:=
\lim_{k\to\infty}\sigma_k
\]
and
\[
\tilde{\vb}(t)
:=
\tilde{\vb}^{k}(t),
\qquad 0\leq t<\sigma_k.
\]

For every fixed $t\geq0$,
\[
\tilde{\vb}(t){1}_{\{t<\widehat{\tau}_i\}}
=
\lim_{k\to\infty}
\tilde{\vb}^{k}(t){1}_{\{t<\sigma_k\}},
\qquad \mathbb P\text{-a.s.}
\]
For every $k\in\mathbb N$, the random variable $\tilde{\vb}^{k}(t){1}_{\{t<\sigma_k\}}$ is $\mathscr F_t$-measurable.

Hence
$\tilde{\vb}(t){1}_{\{t<\widehat{\tau}_i\}}$ is
$\mathscr F_t$-measurable. Thus, $(\tilde{\vb},\widehat{\tau}_i)$ is adapted for $i=1,2$.
By pathwise uniqueness in \autoref{prop-tildev-local strong}, this
adapted solution agrees with the local strong solution on their common
interval.

It remains to justify the construction after a stopping time $\sigma$. Let $\sigma$ be an almost surely finite $(\mathscr F_t)$-stopping time and let $a$ be an $\mathscr F_\sigma$-measurable $\mathbb H^i(\dom)$-valued random variable.
For $N\geq3$, set
\[
E_N:=\left\{\Vert a\Vert_{H^i}+\alpha(\sigma)+\alpha(\sigma)^{-1}+\Vert\wb(\sigma)\Vert_{H^2}\leq N\right\}\in\mathscr F_\sigma,
\qquad E_2:=\varnothing,
\]
and
\[
\rho_N^\sigma:=\inf\left\{s\geq0:\alpha(\sigma+s)+\alpha(\sigma+s)^{-1}+\Vert\wb(\sigma+s)\Vert_{H^2}\geq2N\right\}\wedge1.
\]
By continuity and adaptedness of $\alpha$ and $\wb$, $\rho_N^\sigma$ is an $(\mathscr F_{\sigma+s})$-stopping time, and
\[
E_N\uparrow\Omega,\qquad \rho_N^\sigma>0\ \text{on }E_N,\qquad \mathbb P\text{-a.s.}
\]
Define
\begin{align*}
\alpha_N^\sigma(s)&:=\mathds1_{E_N}\alpha\bigl(\sigma+(s\wedge\rho_N^\sigma)\bigr)+\mathds1_{E_N^c},\\
\wb_N^\sigma(s)&:=\mathds1_{E_N}\wb\bigl(\sigma+(s\wedge\rho_N^\sigma)\bigr).
\end{align*}
Since $E_N\in\mathscr F_\sigma$, $\alpha_N^\sigma$ and $\wb_N^\sigma$ are continuous and $(\mathscr F_{\sigma+s})$-adapted, with
\[
\Vert\mathds1_{E_N}a\Vert_{H^i}\leq N,
\qquad
\sup_{s\geq0}\left(\alpha_N^\sigma(s)+(\alpha_N^\sigma(s))^{-1}+\Vert\wb_N^\sigma(s)\Vert_{H^2}\right)\leq2N.
\]
Therefore, by \autoref{prop-tildev-local strong} and the preceding Picard iteration, there exist a deterministic $\delta_N=\delta_N(i,N)\in(0,1]$, independent of $\sigma$, and an $(\mathscr F_{\sigma+s})$-adapted solution $U_N\in X_{\delta_N}^i$ satisfying
\begin{align*}
U_N(s)
={}&e^{-sA}\mathds1_{E_N}a\\
&-\int_0^s\alpha_N^\sigma(r)e^{-(s-r)A}\bigl((U_N(r)+\wb_N^\sigma(r))\cdot\nabla\bigr)(U_N(r)+\wb_N^\sigma(r))\,dr,
\qquad s\in[0,\delta_N].
\end{align*}
On $E_N$, we have
\[
U_N(0)=a,\qquad
\alpha_N^\sigma(s)=\alpha(\sigma+s),\qquad
\wb_N^\sigma(s)=\wb(\sigma+s),
\qquad 0\leq s\leq\rho_N^\sigma.
\]

Set
\[
\tau:=\sigma+\sum_{N=3}^{\infty}\mathds1_{E_N\setminus E_{N-1}}(\delta_N\wedge\rho_N^\sigma).
\]
For every $s\geq0$,
\[
\{\tau-\sigma>s\}
=
\bigcup_{N=3}^{\infty}\left((E_N\setminus E_{N-1})\cap\{\delta_N\wedge\rho_N^\sigma>s\}\right)
\in\mathscr F_{\sigma+s}.
\]
Hence, $\tau$ is an $(\mathscr F_t)$-stopping time and $\tau>\sigma$, $\mathbb P$-a.s..
Define
\[
\tilde{\vb}(\sigma+s,\omega):=U_N(s,\omega),
\qquad
0\leq s\leq\tau(\omega)-\sigma(\omega),\quad \omega\in E_N\setminus E_{N-1}.
\]
Since $E_N\setminus E_{N-1}\in\mathscr F_\sigma$, we infer that $\tilde{\vb}$ is adapted and satisfies \eqref{ee-o} on $[\sigma,\tau]$, with
\[
\tilde{\vb}(\sigma)=a,\qquad
\tilde{\vb}\in C([\sigma,\tau];\mathbb H^i(\dom))\cap L^2([\sigma,\tau];\mathbb H^{i+1}(\dom)),
\qquad\mathbb P\text{-a.s.}
\]
Uniqueness follows from \autoref{prop-tildev-local strong}. The proof of \autoref{prop-tildev-adapted} is complete.

\end{proof}

\subsection{Global existence and regularities of solutions to (\ref{ee-o}) and (\ref{eqn-Burgers-intro})}\label{subsection-({ee-o2})}
In this subsection, we will extend the local existence of weak solutions of (\ref{ee-o}) to global existence. We begin with the definition of maximal and global solutions.

\begin{definition}[\textbf{Maximal and global weak solutions to (\ref{ee-o})}] \label{def-2.2}
Let $\ub_{0}\in \mathbb{L}^{2}(\dom)$ and $\overrightarrow{\ub_{0}}$ be its randomization defined in Definition \ref{def-randomization-u}.
 \begin{enumerate}

\item Let $\xi$ be an $(\mathscr F_t)_{t\geq0}$-stopping time taking values in $(0,\infty]$.
We say that $(\tilde{\vb},\xi)$ is a maximal weak pathwise solution to \eqref{ee-o} if there exists an increasing sequence of stopping times $\{\tau_n\}_{n\in\mathbb N}$ such that
\[
\tau_n\uparrow\xi,\qquad \mathbb P\text{-a.s.},
\]
and $(\tilde{\vb},\tau_n)$ is a local weak pathwise solution to \eqref{ee-o} for every $n\in\mathbb N$.
Moreover, for any local weak pathwise solution $(\widehat{\vb},\hat\xi)$ to \eqref{ee-o} satisfying
\[
\widehat{\vb}(t)=\tilde{\vb}(t),\qquad
0\leq t<\xi\wedge\hat\xi,\quad \mathbb P\text{-a.s.},
\]
we have
\[
\mathbb P(\hat\xi>\xi)=0.
\]

\item If $(\tilde{\vb}, \xi)$ is a maximal weak pathwise solution and $\xi=\infty$ $\mathbb{P}$-a.s., then we say the solution $\tilde{\vb}$ is global.
 \end{enumerate}
\end{definition}

\begin{definition}[\textbf{Maximal and global strong solutions to (\ref{ee-o})}] \label{def-2.2 strong}
 Assume the initial data $\tilde{\vb}_0\in \mathbb{H}^{i}(\dom)$, $i=1, 2$.
 \begin{enumerate}

\item Let $\xi_i$ be an $(\mathscr F_t)_{t\geq0}$-stopping time taking values in $(0,\infty]$.
We say that $(\tilde{\vb},\xi_i)$ is a maximal strong pathwise solution to \eqref{ee-o} if there exists an increasing sequence of stopping times $\{\tau_{i,n}\}_{n\in\mathbb N}$ such that
\[
\tau_{i,n}\uparrow\xi_i,\qquad \mathbb P\text{-a.s.},
\]
and $(\tilde{\vb},\tau_{i,n})$ is a local strong pathwise solution to \eqref{ee-o} for every $n\in\mathbb N$.
Moreover, for any local strong pathwise solution $(\widehat{\vb},\hat\xi_i)$ to \eqref{ee-o} satisfying
\[
\widehat{\vb}(t)=\tilde{\vb}(t),\qquad
0\leq t<\xi_i\wedge\hat\xi_i,\quad
\mathbb P\text{-a.s.},
\]
we have
\[
\mathbb P(\hat\xi_i>\xi_i)=0.
\]

\item If $(\tilde{\vb}, \xi_i)$ is a maximal strong pathwise solution and $\xi_i=\infty$ $\mathbb{P}$-a.s., then we say the solution $\tilde{\vb}$ is global.
 \end{enumerate}
\end{definition}

\begin{prop}\label{prop-maximal-tildev}
Assume that $\ub_{0}\in\mathbb{L}^{2}(\dom)$ and let
$\overrightarrow{\ub_{0}}\in\mathbb{L}^{2}(\Omega;\mathbb{L}^{2}(\dom))$
be the randomization of $\ub_{0}$ defined by
\autoref{def-randomization-u}, using a sequence $\{r_{\kb}^j\}_{(\kb,j)\in\mathcal K}$ satisfying \autoref{A1},
and assume that $f$ satisfies \autoref{A2}.
Then there exists a unique maximal weak pathwise solution
$(\tilde{\vb},\xi)$ to \eqref{ee-o}
in the sense of \autoref{def-2.2}.
\end{prop}

\begin{proof}
In view of \autoref{thm-local weak solution},
there exists a unique local weak pathwise solution
to \eqref{ee-o} with initial data
$\overrightarrow{\ub_{0}}$.

Let $\mathcal{T}$ be the collection of all stopping times
$\tau$ such that there exists a local weak pathwise
solution to \eqref{ee-o} on $[0,\tau]$
with initial data $\overrightarrow{\ub_{0}}$.

Let us choose and fix $\tau_1,\tau_2\in\mathcal{T}$,
and let $(\tilde{\vb}^{1},\tau_1)$ and
$(\tilde{\vb}^{2},\tau_2)$ be the corresponding
local weak solutions.

By the pathwise uniqueness in
\autoref{thm-local weak solution}, we infer that
\[
\tilde{\vb}^{1}(t)=\tilde{\vb}^{2}(t),
\qquad 0\leq t\leq\tau_1\wedge\tau_2,
\quad\mathbb{P}\text{-a.s.}
\]
Therefore, the two solutions can be consistently
pasted on $[0,\tau_1\vee\tau_2]$, which implies
\[
\tau_1\vee\tau_2\in\mathcal{T}.
\]
Hence, $\mathcal{T}$ is upward directed.

Let us define
\[
\xi:=\operatorname*{ess\,sup}_{\tau\in\mathcal{T}}\tau.
\]
Since $\mathcal{T}$ is upward directed,
there exists an increasing sequence
$\{\tau_n\}_{n\in\mathbb N}\subset\mathcal{T}$
such that
\[
\tau_n\uparrow\xi,\qquad\mathbb{P}\text{-a.s.}
\]
In particular, $\xi$ is a stopping time.

For every $n\in\mathbb N$, let
$(\tilde{\vb}^{n},\tau_n)$ be the corresponding
local weak solution.
By pathwise uniqueness, we infer that
\[
\tilde{\vb}^{n+1}(t)=\tilde{\vb}^{n}(t),
\qquad 0\leq t\leq\tau_n,
\quad\mathbb{P}\text{-a.s.}
\]
Consequently, we can define
\[
\tilde{\vb}(t):=\tilde{\vb}^{n}(t),
\qquad 0\leq t\leq\tau_n.
\]
Thus, $(\tilde{\vb},\xi)$ is a weak pathwise
solution on $[0,\xi)$, localized by
$\{\tau_n\}_{n\in\mathbb N}$.

We claim that $(\tilde{\vb},\xi)$ is maximal.
Indeed, suppose that there exists a local weak
pathwise solution $(\widehat{\vb},\hat\xi)$
to \eqref{ee-o} satisfying
\[
\widehat{\vb}(t)=\tilde{\vb}(t),
\qquad 0\leq t<\xi\wedge\hat\xi,
\quad\mathbb{P}\text{-a.s.}
\]
Since $\hat\xi\in\mathcal{T}$,
by the definition of $\xi$, we infer that
\[
\hat\xi\leq\xi,\qquad\mathbb{P}\text{-a.s.}
\]
which implies
\[
\mathbb P(\hat\xi>\xi)=0.
\]
Therefore, $(\tilde{\vb},\xi)$ is a maximal
weak pathwise solution in the sense of
\autoref{def-2.2}.

Uniqueness follows from the pathwise uniqueness
of local weak solutions.
The proof of \autoref{prop-maximal-tildev}
is complete.
\end{proof}

Moreover, the local strong solution to \eqref{ee-o}
can be extended to a unique maximal strong solution.

\begin{prop}\label{prop-maximal global-tildev}
For $i=1,2$, assume that
$\tilde{\vb}_0\in\mathbb H^i(\dom)$
and \autoref{A2} holds.
Then there exists a unique maximal strong
pathwise solution $(\tilde{\vb},\xi_i)$
to \eqref{ee-o} in the sense of
\autoref{def-2.2 strong}.
Moreover,
\begin{align}\label{eq-veps-H1-blowup}
\limsup_{t\uparrow\xi_i}
\Vert\tilde{\vb}(t)\Vert_{\mathbb H^i(\dom)}
=\infty
\quad\text{on }\{\xi_i<\infty\},
\quad\mathbb P\text{-a.s.}
\end{align}
The same existence, uniqueness and blow-up alternative hold with any deterministic initial time $r\geq0$ and any
$\mathscr F_r$-measurable $\mathbb H^i(\dom)$-valued initial value $a$, with $t\uparrow\xi_i$ taken over $t\in[r,\xi_i)$.
\end{prop}

\begin{proof}
Let us choose and fix $i=1,2$.
In view of \autoref{prop-tildev-local strong}
and \autoref{prop-tildev-adapted},
there exists a unique adapted local strong
solution to \eqref{ee-o}.

Following the proof of \autoref{prop-maximal-tildev},
by the pathwise uniqueness of local strong solutions,
we infer that there exists a unique maximal strong
pathwise solution $(\tilde{\vb},\xi_i)$
to \eqref{ee-o}, together with an increasing
sequence of stopping times
$\{\tau_{i,n}\}_{n\in\mathbb N}$ satisfying
\[
\tau_{i,n}\uparrow\xi_i,\qquad
\mathbb P\text{-a.s.}
\]

In the following, we aim to prove
\eqref{eq-veps-H1-blowup}.
Suppose by contradiction that
\[
\mathbb P\left(
\xi_i<\infty,\;
\sup_{0\leq t<\xi_i}
\Vert\tilde{\vb}(t)\Vert_{\mathbb H^i(\dom)}
<\infty
\right)>0.
\]

For every $K\in\mathbb N$, let us set
\[
\rho_K:=\inf\left\{t\geq0:
\alpha(t)+\alpha^{-1}(t)
+\Vert\wb(t)\Vert_{\mathbb H^2(\dom)}
\geq K\right\}\wedge K.
\]
Since $\alpha$ and $\wb$ are continuous,
we infer that
\[
\rho_K\uparrow\infty,\qquad
\mathbb P\text{-a.s.}
\]

Consequently, there exist $T,M,K\in\mathbb N$
such that
\[
\mathbb P(E_{T,M,K})>0,
\]
where
\[
E_{T,M,K}:=
\left\{\xi_i\leq T,\;
\sup_{0\leq t<\xi_i}
\Vert\tilde{\vb}(t)\Vert_{\mathbb H^i(\dom)}
\leq M,\;
\rho_K>T+1\right\}.
\]

By the uniform local existence estimate in
\autoref{prop-tildev-local strong},
there exists
\[
\delta_0=\delta_0(i,M,K,T)\in(0,1]
\]
such that the equation with coefficients stopped
at $\rho_K$ admits a unique strong solution
on $[t_0,t_0+\delta_0]$ whenever
$t_0\leq T\wedge\rho_K$ and
$\Vert\tilde{\vb}(t_0)\Vert_{\mathbb H^i(\dom)}
\leq M$.
The bounded stopping-time restart assertion in
\autoref{prop-tildev-adapted} gives the same deterministic lower bound
for adapted restarts at stopping times.

For every $n\in\mathbb N$, let
$(\tilde{\vb}^{n},\tau_{i,n})$
be the corresponding local strong solution
and define
\[
A_n:=
\left\{\tau_{i,n}\leq T,\;
\tau_{i,n}<\rho_K,\;
\Vert\tilde{\vb}^{n}(\tau_{i,n})
\Vert_{\mathbb H^i(\dom)}\leq M\right\}
\in\mathscr F_{\tau_{i,n}}.
\]
Applying the adapted local strong construction
at $\tau_{i,n}$ on $A_n$, we obtain a local
strong solution extending
$\tilde{\vb}^{n}$ up to
\[
\hat\tau_{i,n}:=
\begin{cases}
(\tau_{i,n}+\delta_0)\wedge\rho_K,
&\text{on }A_n,\\
\tau_{i,n},
&\text{on }A_n^c.
\end{cases}
\]

For every $\omega\in E_{T,M,K}$,
we have $\omega\in A_n$ for every $n\in\mathbb N$.
Since
\[
\tau_{i,n}\uparrow\xi_i
\quad\text{and}\quad
\rho_K>T+1,
\]
we infer that, for sufficiently large $n$,
\[
\hat\tau_{i,n}
=\tau_{i,n}+\delta_0>\xi_i
\quad\text{on }E_{T,M,K}.
\]
Hence,
\[
\mathbb P\left(
\bigcup_{n\in\mathbb N}
\{\hat\tau_{i,n}>\xi_i\}
\right)
\geq\mathbb P(E_{T,M,K})>0.
\]
This contradicts the maximality of
$(\tilde{\vb},\xi_i)$.

Therefore,
\[
\mathbb P\left(
\xi_i<\infty,\;
\sup_{0\leq t<\xi_i}
\Vert\tilde{\vb}(t)\Vert_{\mathbb H^i(\dom)}
<\infty
\right)=0.
\]
In view of the continuity of $\tilde{\vb}$
on every compact interval contained in
$[0,\xi_i)$, we infer that
\eqref{eq-veps-H1-blowup} holds.

For initial time $r$ and random initial value $a$, use the restart assertion of \autoref{prop-tildev-adapted} in the filtration $\mathscr G_s=\mathscr F_{r+s}$ to obtain the local adapted solutions, and take the essential supremum of their lifetimes as above. For the blow-up alternative, repeat the preceding argument with $\sup_{r\leq t<\xi_i}$ in place of $\sup_{0\leq t<\xi_i}$ and choose $T\geq r$. Retain the same physical-time coefficient exit time $\rho_K$. The bounded restart assertion of \autoref{prop-tildev-adapted} then justifies the same deterministic $\delta_0$ and pasting at $\tau_{i,n}$ on $A_n$. No moment assumption on $a$ is needed,
since its norm is localized on $\mathscr F_r$-measurable events.

The proof of
\autoref{prop-maximal global-tildev}
is complete.
\end{proof}

\begin{cor}
\label{cor-ae-H2}
Assume that $u_0 \in \mathbb{L}^{2}(\dom )$.
Let $\overrightarrow{\ub_{0}} \in  \mathbb{L}^{2}(\Omega; \mathbb{L}^{2}(\dom ))$ be the randomization of $\ub_0$ given by Definition \ref{def-randomization-u}, using a sequence $\{r_{\kb}^j\}_{(\kb,j)\in\mathcal K}$ satisfying \autoref{A1}, and assume that $f$ satisfies \autoref{A2}.
Let $(\tilde{\vb},\xi)$ be the unique maximal weak solution to equation (\ref{ee-o}) in the sense of \autoref{def-2.2}.

Then, for almost every $\omega\in\Omega$ and every $\varepsilon\in(0,\xi(\omega))$, there exists $t^\ast\in(0,\varepsilon)$ such that
\[
\tilde{\vb}(t^\ast,\omega)\in\mathbb H^2(\dom).
\]
\end{cor}

\begin{proof}
For every $\tau<\xi$, we infer that
\[
\tilde{\vb}\in L^2([0,\tau];\mathbb H^{1,3}(\dom)),\;\;\mathbb P\text{-a.s.}.
\]
Since $\mathbb H^{1,3}(\dom)\hookrightarrow\mathbb H^1(\dom)$, for almost every $\omega\in\Omega$,
\[
\tilde{\vb}(t,\omega)\in\mathbb H^1(\dom),\;\;\text{for a.e. }t\in(0,\xi(\omega)).
\]

Let us choose and fix  $\varepsilon\in(0,\xi(\omega))$. There exists  $t_0\in\left(0,\frac{\varepsilon}{2}\right)$ such that
\[
\tilde{\vb}(t_0,\omega)\in\mathbb H^1(\dom).
\]
By \autoref{prop-tildev-local strong} and uniqueness, for some
$\delta>0$,
\[
\tilde{\vb}(\cdot,\omega)
\in
L^2\left(
[t_0,t_0+\delta];
\mathbb H^2(\dom)
\right).
\]
Hence
\[
\tilde{\vb}(t,\omega)\in\mathbb H^2(\dom),\;\;\text{for a.e. }t\in\left(t_0,t_0+\min\left\{\delta,\frac{\varepsilon}{2}\right\}\right).
\]
Therefore, there exists
\[
t^\ast
\in
\left(
t_0,
t_0+\min\left\{\delta,\frac{\varepsilon}{2}\right\}
\right)
\subset(0,\varepsilon)
\]
such that
\[
\tilde{\vb}(t^\ast,\omega)\in\mathbb H^2(\dom).
\]
The proof of \autoref{cor-ae-H2} is complete.

\end{proof}

To prove the global existence of $\tilde{\vb}$ to (\ref{ee-o}), we smooth the coefficients of the equation. Recall that $\left(\hat{ \fb}_\kb^\jb\right)_{ (k,j) \in \mathcal{K}}$ is the Fourier coefficient of $f$ such that
\[
\fb(x) =  \sum_{(k,j) \in \mathcal{K}} \hat{ \fb}_\kb^\jb \eb_k^j(x),\;\;x\in\dom.
\]
For $m\in\mathbb{N}$, set
\begin{align*}
\wb_m(t,\xb):=\sum_{(\kb,j)\in \mathcal{K}, |\kb|\leq m}\int_{0}^{t}\exp(-|\kb|^{2}(t-s))\alpha^{-1}(s)\hat{\fb}_{\kb}(s)e_k^j\,ds, \;\; (t,\xb) \in [0,\infty)\times \dom.
\end{align*}

We first extend $\hat{\fb}_{\kb}^{\jb}$ and $B$ to $\mathbb{R}$ by setting
\[
\hat{\fb}_{\kb}^{\jb}(t)=0,\;\; B(t)=0,\;\;\mbox{ on } t<0,
\]
and use the same notation for these extensions.
We denote by $\hat{\fb}_{\kb}^{\jb,\eps}$ and $B^\eps$ the mollification
of $\hat{\fb}_{\kb}^{\jb}$ and $B$ on $\mathbb{R}$ respectively, i.e.,
\[
\hat{\fb}_{\kb}^{\jb,\eps}
=
\eta_\eps*\hat{\fb}_{\kb}^{\jb},\;\;B^\eps=\eta_\eps*B.
\]

Set
\[
\alpha^\eps:=e^{bB^\eps}.
\]
Thus,
$\alpha^\eps\in C^\infty([0,T];\mathbb{R})$,
where $C^\infty([0,T];\mathbb{R})$ is defined as follows:
\[
C^\infty([0,T]; \mathbb{R}) := \left\{ f : [0,T] \to \mathbb{R} \,\middle|\, \partial_t^k f \text{ exists and is continuous for all } k \in \mathbb{N} \right\}.
\]
For $(t,\xb) \in [0,\infty)\times \dom $, define
\begin{align}
\wb_m^{\eps}:=\wb_m^{\eps}(t,\xb)&=\sum_{(\kb,j)\in \mathcal{K}, |\kb|\leq m}\int_{0}^{t}\exp(-|\kb|^{2}(t-s))(\alpha^{\eps})^{-1}(s)\hat{ \fb}_{\kb}^{j,\eps}(s)\eb_k^j\,ds, \label{eq-wmesp}\\
\wb^{\eps}:=\wb^{\eps}(t,\xb)&= \sum_{(\kb,j)\in \mathcal{K}}\int_{0}^{t}\exp(-|\kb|^{2}(t-s))(\alpha^{\eps})^{-1}(s)\hat{ \fb}_{\kb}^{j,\eps}(s)\eb_k^j\,ds.\label{eq-wesp}
\end{align}
Obviously, $\wb_m^{\eps} \in C^{\infty}([0,T]\times \dom; \mathbb{R}^{3})$, where $C^{\infty}([0,T]\times \dom; \mathbb{R}^{3})$ is defined as follows:
\[
C^\infty([0,T] \times \dom; \mathbb{R}^3) := \left\{ f : [0,T] \times \dom \to \mathbb{R}^3 \,\middle|\,
\partial_t^k \partial_x^\alpha f \text{ exists and is continuous for all } k \in \mathbb{N}, \alpha \in \mathbb{N}^3 \right\}.
\]

Define
\[
\fb^{\eps}(x) =  \sum_{(k,j) \in \mathcal{K}} \hat{ \fb}_{\kb}^{\jb,\eps} \eb_k^j(x),\;\;x\in\dom.
\]
In view of \autoref{lem-mollification}, we have the following lemma.
\begin{lemma}\label{lem-uniform convergence}
Assume that \autoref{A2} holds. For every $T>0$, let $\eps\rightarrow 0$, then we have the following result
\begin{trivlist}
\item[(i)]

\begin{align}\label{eq-uniform convergence1}
&\Vert\alpha^\eps-\alpha\Vert_{L^\infty([0,T])}
+
\Vert(\alpha^\eps)^{-1}-\alpha^{-1}\Vert_{L^\infty([0,T])}
\rightarrow 0,
\qquad \mathbb P\text{-a.s.}
\end{align}
\item[(ii)]
\begin{align}\label{eq-uniform convergence3}
\Vert\fb^\eps-\fb\Vert_{L^2([0,T];\mathbb H^\delta(\dom))}
\to 0,
\qquad \mathbb P\text{-a.s.}
\end{align}
Moreover,
\begin{align}\label{eq-f-eps-uniform}
\sup_{\eps\in(0,1)}
\Vert\fb^\eps\Vert_{L^2([0,T];\mathbb H^\delta(\dom))}
<\infty.
\end{align}
\item[(iii)] Let $\eps\to0$ and $m\to\infty$,
\begin{align}\label{eq-uniform convergence2}
\sup_{t\in[0,T]}
\Vert\wb_m^\eps(t)-\wb(t)\Vert_{H^2}
\rightarrow0,
\qquad\mathbb P\text{-a.s.}.
\end{align}
\end{trivlist}
\end{lemma}
The proof is given in \autoref{app-proof-uniform-convergence}.

To prove the global existence of $\tilde{\vb}$, for $\eps>0$ and $m\in\mathbb{N}$, we need to introduce an auxiliary function $\vb^{\eps,m}$, which satisfies the following equation:
\begin{align}\label{ee-epsilon}
&\partial_{t} \vb^{\eps,m}-\Delta \vb^{\eps,m}+\alpha^{\eps}(t)((\vb^{\eps,m}+\wb_m^{\eps})\cdot \nabla) (\vb^{\eps,m}+\wb_m^{\eps})=0,\ \mathrm{on}\ [0,\infty)\times \dom,\\
&\vb^{\eps,m}(0)=\tilde\vb_{0}(\xb),\quad \xb=(x_{1},x_{2},x_{3})\in \dom.\nonumber
\end{align}
Fix a coefficient path outside the null set on which the preceding regularity assertions may fail. The time regularization uses coefficient values at later times, so we treat \eqref{ee-epsilon} as a deterministic parabolic equation on this path. Its local and maximal existence times are pathwise lifetimes. Applying the deterministic fixed-point estimates in the proof of \autoref{prop-tildev-local strong} to the smooth coefficients $\alpha^\eps$ and $\wb_m^\eps$ gives the following local result.

\begin{lemma}\label{lem-vepsm-local strong}
 Assume that $\tilde\vb_{0}\in \mathbb{H}^{1}(\dom)$ and \autoref{A2} holds, then there exists a unique local strong solution  $(\vb^{\eps,m}, \tau^{\eps,m})$  to (\ref{ee-epsilon})  such that
\[
\vb^{\eps,m}\in C([0, \tau^{\eps,m}]; \mathbb{H}^{1}(\dom)) \cap L^{2}([0, \tau^{\eps,m}]; \mathbb{H}^{2}(\dom)),\;\;\mathbb{P}\text{-a.s.}.
\]
Moreover, the local strong pathwise solution $\vb^{\eps,m}$ to (\ref{ee-epsilon}) is Lipschitz continuous with respect to the initial data $\tilde{\vb}_0$.

\end{lemma}

Moreover, the local strong solution to (\ref{ee-epsilon}) can be extended to a unique maximal strong solution, denoted by $(\vb^{\eps,m}, \xi^{\eps,m})$.
\begin{prop}\label{prop-maximal global-veps}
Assume that $\tilde\vb_{0}\in \mathbb{H}^{i}(\dom)$, $i=1,2$ and \autoref{A2} holds. Problem (\ref{ee-epsilon}) has a unique maximal strong solution $(\vb^{\eps,m}, \xi^{\eps,m}_i)$.
\end{prop}

The Burgers equation does not preserve the spatial mean, even for initial data of zero mean. The following estimate controls the full Sobolev norm by the homogeneous seminorm and the initial data. The mean estimates used in its proof are given in \autoref{app-mean-control}.

\begin{cor}\label{cor-dominate L2}
Assume that $\tilde{\vb}_{0}\in \mathbb{H}^{1}(\dom)$ and \autoref{A2} holds. Let  $(\vb^{\eps,m}, \xi^{\eps,m})$ be the maximal strong solutions of (\ref{ee-epsilon}). If for $ s>0$ and $t\in [0, \xi^{\eps,m})$,
\[
\Lambda^s \vb^{\eps,m}(t)\in\mathbb{L}^2(\dom).
\]
Then, we have
\begin{align}\label{eq-2.8}
 &\Vert \vb^{\eps,m}(t)\Vert_{s}\leq \Vert \vb^{\eps,m}(t)\Vert_{H^{s}}\\
 &\leq 2\Big (\Vert \vb^{\eps,m}(t)\Vert_{s}+\frac{1}{(2\pi)^\frac{3}{2}}\int_{0}^{t}\alpha^{\eps}(s)\left (\Vert \vb^{\eps,m}(s)\Vert_{\frac{1}{2}}^{2}+\Vert \wb_m^{\eps}(s)\Vert_{\frac{1}{2}}^{2}\right )\,ds\Big )+\frac{1}{(2\pi)^\frac{3}{2}}\Vert\tilde{\vb}_{0}\Vert_{L^1}.\nonumber
\end{align}

Let $s\geq \frac12$ and $0<T<\xi^{\eps,m}$. If $\sup_{t\in[0,T]}\Vert\vb^{\eps,m}(t)\Vert_s<\infty$, then there is a finite constant $C_T$ such that

\begin{align}\label{eqn-H^s bound}
\sup_{t\in[0,T]}
\Vert \vb^{\eps,m}(t)\Vert_{H^s}
\leq C_T \left(1+\sup_{t\in[0,T]}\Vert \vb^{\eps,m}(t)\Vert_s+\sup_{t\in[0,T]}\Vert \vb^{\eps,m}(t)\Vert_s^2 \right).
\end{align}
The constant $C_T$ can be chosen to depend only on $T$, $\Vert\tilde{\vb}_0\Vert_{L^1}$, $\Vert\alpha^\eps\Vert_{L^\infty([0,T])}$ and $\Vert\wb_m^\eps\Vert_{C([0,T];\mathbb{H}^{1/2}(\dom))}$.
\end{cor}
The proof is given in \autoref{app-mean-control}.

The Galerkin approximation and the regularity estimates are established in \autoref{app-galerkin-regularity}. In particular, they yield the following classical regularity.

\begin{lemma}\label{eq-classical solution}
Assume that $\tilde\vb_{0}\in \mathbb{H}^{1}(\dom)$ and \autoref{A2} holds, the maximal strong solution $(\vb^{\eps,m}, \xi^{\eps,m})$ to (\ref{ee-epsilon}) is also classical on $(0, \xi^{\eps,m})$.
\end{lemma}
The proof is given in \autoref{app-proof-classical-solution}.
\begin{prop}
	\label{prop-6.1.2} Assume that \autoref{A2} holds. Let $(\vb^{\eps,m},\xi^{\eps,m} )$ be the maximal strong solution to (\ref{ee-epsilon}) with initial data $\tilde\vb_{0}\in \mathbb{H}^{1}(\dom)$ in \autoref{prop-maximal global-veps}. If the initial data $ \tilde\vb_{0}\in \mathbb{H}^{2}(\dom)$, then for arbitrary $T>0$,  the solution $\vb^{\eps,m}$ satisfies
	\begin{align}\label{eq-epsilon estimates}
\sup_{m\in\mathbb N}\sup\limits_{\eps\in (0,1)}	\sup\limits_{t\in \lbrack 0,\xi^{\eps,m} \wedge T )}\Vert\vb^{\eps,m}(t)\Vert_{L^\infty } &\leq C(T,\Vert\tilde\vb_{0}\Vert_{L^\infty}, \Vert \fb\Vert_{L^{2}([0,T+1]; \mathbb{H}^{\delta}(\dom))},  \Vert\alpha^{-1}\Vert_{L^{\infty}([0,T+1])})<\infty,\;\mathbb{P}\text{-a.s.}.
	\end{align}

\end{prop}
\begin{proof}
For $m\in \mathbb{N}$, similar to \eqref{eqn-P_n},  we set
 \begin{align*}
&\fb_{m}:=\fb_{m}(t,\xb)\!=\! \sum_{(\kb,j)\in \mathcal{K}, |\kb|\leq m} \hat{\fb}^j_{\kb}(t) \eb_k^j,\;\;\; \\
&\fb_{m}^{\eps}:=\fb_{m}^{\eps}(t,\xb)\!=\! \sum_{(\kb,j)\in \mathcal{K}, |\kb|\leq m} \hat{\fb}^{j,\eps}_{\kb}(t) \eb_k^j.
\end{align*}

Recall that $\wb_m^\eps$ satisfies
\[
\partial_t\wb_m^\eps-\Delta\wb_m^\eps
=
(\alpha^\eps)^{-1}\fb_m^\eps,
\qquad
\wb_m^\eps(0)=0.
\]
Define $ \tilde\ub^{\eps,m}=\vb^{\eps,m}+\wb_m^\eps$, then (\ref{ee-epsilon}) is equivalent to
\begin{align}\label{ee-5.2}
&\partial_{t} \tilde\ub^{\eps,m}-\Delta \tilde\ub^{\eps,m}+\alpha^{\eps}(\tilde\ub^{\eps,m}\cdot \nabla) \tilde\ub^{\eps,m}=(\alpha^{\eps})^{-1}\fb_{m}^{\eps},\ \mathrm{on}\ (0,\xi^{\eps,m})\times \dom,\\
&\tilde\ub^{\eps,m}(0)=\tilde\vb_{0} \in \mathbb{H}^{2}(\dom).\nonumber
\end{align}
According to  \autoref{eq-classical solution} the local strong solution $\tilde\ub^{\eps,m}$ is classical on $(0, \xi^{\eps,m})$.

Let $\tilde\ub^{\eps,m}=(\tilde u_1^{\eps,m}, \tilde u_2^{\eps,m}, \tilde u_3^{\eps,m})$.
By \eqref{ee-5.2}, for $i=1,2,3$,
\[
\partial_t\tilde u_i^{\eps,m}
-\Delta\tilde u_i^{\eps,m}
+\alpha^\eps
(\tilde\ub^{\eps,m}\cdot\nabla)
\tilde u_i^{\eps,m}
=
(\alpha^\eps)^{-1}f_{m,i}^\eps.
\]
Set
\[
\eta(t)
:=
(\alpha^\eps)^{-1}(t)
\Vert\fb_m^\eps(t)\Vert_{L^\infty}.
\]
For $i=1,2,3$, define
\[
\tilde u_{i,+}^{\eps,m}(t,\xb)
:=
\tilde u_i^{\eps,m}(t,\xb)
-\int_0^t\eta(s)\,ds.
\]
Then
\begin{align*}
&\partial_t\tilde u_{i,+}^{\eps,m}
-\Delta\tilde u_{i,+}^{\eps,m}
+\alpha^\eps
(\tilde\ub^{\eps,m}\cdot\nabla)
\tilde u_{i,+}^{\eps,m}=(\alpha^\eps)^{-1}f_{m,i}^\eps-\eta(t)
\leq0.
\end{align*}
Applying the maximum principle
\cite[Theorem 3.1.1]{Topping_2006}, we deduce
\[
\tilde u_i^{\eps,m}(t,\xb)
\leq
\Vert\tilde\vb_0\Vert_{L^\infty}
+\int_0^t\eta(s)\,ds.
\]

Similarly, define
\[
\tilde u_{i,-}^{\eps,m}(t,\xb)
:=
\tilde u_i^{\eps,m}(t,\xb)
+\int_0^t\eta(s)\,ds.
\]
Then
\begin{align*}
&\partial_t\tilde u_{i,-}^{\eps,m}
-\Delta\tilde u_{i,-}^{\eps,m}
+\alpha^\eps
(\tilde\ub^{\eps,m}\cdot\nabla)
\tilde u_{i,-}^{\eps,m}=(\alpha^\eps)^{-1}f_{m,i}^\eps+\eta(t)
\geq0.
\end{align*}
Applying the maximum principle again, we infer that
\[
\tilde u_i^{\eps,m}(t,\xb)
\geq
-\Vert\tilde\vb_0\Vert_{L^\infty}
-\int_0^t\eta(s)\,ds.
\]
Therefore,
\[
|\tilde u_i^{\eps,m}(t,\xb)|
\leq
\Vert\tilde\vb_0\Vert_{L^\infty}
+
\int_0^t
(\alpha^\eps)^{-1}(s)
\Vert\fb_m^\eps(s)\Vert_{L^\infty}\,ds,
\]
which implies there exists a constant $c>0$ such that
\begin{align}\label{eq-epsilon estimates2}
\sup_{t\in[0,\xi^{\eps,m}\wedge T)}
\Vert\tilde\ub^{\eps,m}(t)\Vert_{L^\infty}
\leq
c\Vert\tilde\vb_0\Vert_{L^\infty}
+
c\int_0^T
(\alpha^\eps)^{-1}(s)
\Vert\fb_m^\eps(s)\Vert_{L^\infty}\,ds.
\end{align}

In view of \autoref{lem-uniform convergence},
\begin{align*}
&\Vert(\alpha^{\eps})^{-1}-\alpha^{-1}\Vert_{L^\infty([0,T])}\to 0,\;\;\mbox{ as }\eps\to 0,\\
&\Vert\fb^\eps-\fb
\Vert_{L^2([0,T];\mathbb H^\delta(\dom))}
\to0,
\;\;\mbox{ as }\eps\to0.
\end{align*}

By the Sobolev embedding theorem,
\[
\mathbb H^\delta(\dom)
\hookrightarrow
\mathbb L^\infty(\dom),\;\;\delta>\frac32.
\]
Therefore, by the H\"older inequality,
\begin{align}\label{eq-eta-uniform}
\int_0^T\eta(s)\,ds
&=
\int_0^T
(\alpha^\eps)^{-1}(s)
\Vert\fb_m^\eps(s)\Vert_{L^\infty}\,ds
\leq
cT^{\frac12}
\Vert(\alpha^\eps)^{-1}\Vert_{L^\infty([0,T])}
\Vert\fb_m^\eps\Vert_{L^2([0,T];\mathbb H^\delta(\dom))}.
\end{align}
Since $\fb_m^\eps$ is the Fourier projection of $\fb^\eps$, we infer that
\[
\Vert\fb_m^\eps(t)\Vert_{H^\delta}
\leq
\Vert\fb^\eps(t)\Vert_{H^\delta},
\]
which implies
\[
\Vert\fb_m^\eps\Vert_{L^2([0,T];\mathbb H^\delta(\dom))}
\leq
\Vert\fb^\eps\Vert_{L^2([0,T];\mathbb H^\delta(\dom))}
\leq
\Vert\fb\Vert_{L^2([0,T+1];\mathbb H^\delta(\dom))}.
\]
Consequently,
\[
\sup_{\eps\in(0,1)}
\sup_{m\in\mathbb N}
\int_0^T\eta(s)\,ds
\leq
C\left(
T,
\Vert\fb\Vert_{L^2([0,T+1];\mathbb H^\delta(\dom))},
\Vert\alpha^{-1}\Vert_{L^\infty([0,T+1])}
\right)
<\infty,
\quad\mathbb P\text{-a.s.}
\]
Thus, by \eqref{eq-epsilon estimates2}
\[
\sup_{\eps\in(0,1)}
\sup_{m\in\mathbb N}
\sup_{t\in[0,\xi^{\eps,m}\wedge T)}
\Vert\tilde\ub^{\eps,m}(t)\Vert_{L^\infty}
<\infty,
\quad\mathbb P\text{-a.s.}
\]

In view of \autoref{lemma-lplq} and the Sobolev embedding inequality, we infer that there exists a constant $c>0$ such that
\begin{equation}\label{eq-epsilon estimates1}
\begin{aligned}
\sup_{\eps\in(0,1)}
\sup_{m\in\mathbb N}
\sup_{t\in[0,\xi^{\eps,m}\wedge T)}
\Vert\wb_m^\eps(t)\Vert_{L^\infty}
&\leq
c
\sup_{\eps\in(0,1)}
\sup_{m\in\mathbb N}
\sup_{t\in[0,\xi^{\eps,m}\wedge T)}
\Vert\wb_m^\eps(t)\Vert_{H^2}
\\
&\leq
c
\sup_{\eps\in(0,1)}
\sup_{m\in\mathbb N}
\Vert(\alpha^\eps)^{-1}\fb_m^\eps
\Vert_{L^2([0,T];\mathbb H^1(\dom))}
\\
&\leq
c
\sup_{\eps\in(0,1)}
\Vert(\alpha^\eps)^{-1}\Vert_{L^\infty([0,T])}
\sup_{\eps\in(0,1)}
\sup_{m\in\mathbb N}
\Vert\fb_m^\eps
\Vert_{L^2([0,T];\mathbb H^1(\dom))}
\\
&\leq
c
\sup_{\eps\in(0,1)}
\Vert(\alpha^\eps)^{-1}\Vert_{L^\infty([0,T])}
\sup_{\eps\in(0,1)}
\Vert\fb^\eps
\Vert_{L^2([0,T];\mathbb H^\delta(\dom))}
\\
&\leq
C\left(T,
\Vert\fb\Vert_{L^2([0,T+1];\mathbb H^\delta(\dom))},
\Vert\alpha^{-1}\Vert_{L^\infty([0,T+1])}
\right)
<\infty,
\quad\mathbb P\text{-a.s.}
\end{aligned}
\end{equation}

Since $\mathbb{H}^2(\dom)\hookrightarrow \mathbb{L}^{\infty}(\dom)$ and $\tilde\vb_{0}\in \mathbb{H}^{2}(\dom)$, we infer that $\Vert\tilde\vb_{0}\Vert_{L^\infty}<\infty$.
Note that $\vb^{\eps,m}=\tilde\ub^{\eps,m}-\wb_m^\eps$, in view of \eqref{eq-epsilon estimates2} and \eqref{eq-epsilon estimates1}, we deduce
\begin{equation}
\begin{aligned}
&\sup_{\eps\in(0,1)}
\sup_{m\in\mathbb N}
\sup_{t\in[0,\xi^{\eps,m}\wedge T)}
\Vert\vb^{\eps,m}(t)\Vert_{L^\infty}\\
\leq &
\sup_{\eps\in(0,1)}
\sup_{m\in\mathbb N}
\sup_{t\in[0,\xi^{\eps,m}\wedge T)}
\Vert\tilde\ub^{\eps,m}(t)\Vert_{L^\infty}
+
\sup_{\eps\in(0,1)}
\sup_{m\in\mathbb N}
\sup_{t\in[0,\xi^{\eps,m}\wedge T)}
\Vert\wb_m^\eps(t)\Vert_{L^\infty}
\\
\leq &
C\left(
T,
\Vert\tilde\vb_0\Vert_{L^\infty},
\Vert\fb\Vert_{L^2([0,T+1];\mathbb H^\delta(\dom))},
\Vert\alpha^{-1}\Vert_{L^\infty([0,T+1])}
\right)
<\infty,
\quad\mathbb P\text{-a.s.}
\end{aligned}
\end{equation}
We complete the proof of \autoref{prop-6.1.2}.
\end{proof}

\begin{prop}\label{prop-global classical}
Assume that \autoref{A2} holds. The maximal strong solution $(\vb^{\eps,m},\xi^{\eps,m} )$ to (\ref{ee-epsilon}) with initial data $\tilde\vb_{0}\in \mathbb{H}^{1}(\dom)$ is global.
\end{prop}
\begin{proof}
It suffices to prove $\xi^{\eps,m}=\infty$, $\mathbb{P}-$a.s.. Since $(\vb^{\eps,m},\xi^{\eps,m} )$ is the maximal strong solution to (\ref{ee-epsilon}) with initial data $\tilde\vb_{0}\in \mathbb{H}^{1}(\dom)$ and satisfies for every $\tau^{\eps,m}\in (0,\xi^{\eps,m})$
\[
\vb^{\eps,m}\in C([0, \tau^{\eps,m}]; \mathbb{H}^{1}(\dom)) \cap L^{2}([0, \tau^{\eps,m}]; \mathbb{H}^{2}(\dom)),\;\;\mathbb{P}\text{-a.s.}.
\]
Therefore, there exists $\tau^\ast\in (0,\xi^{\eps,m})$, such that $ \vb^{\eps,m}(\tau^\ast) \in \mathbb{H}^{2}(\dom)$, $\mathbb{P}-$ a.s.. In the following, we regard $ \tau^\ast$ as the initial time for the system (\ref{ee-epsilon}).
Multiplying (\ref{ee-epsilon}) with $\Lambda^{2}\vb^{\eps,m}$ and integrating over $\dom$ yields, for almost every $t\in(\tau^\ast,\xi^{\eps,m})$,
\begin{align*}
	\frac{1}{2}\partial_{t}\Vert \vb^{\eps,m}\Vert_{1}^{2}+\Vert \vb^{\eps,m}\Vert_{2}^{2} &\leq
    \vert\alpha^\varepsilon(t)\vert\Vert\vb^{\eps,m}+\wb_m^{\eps}\Vert_{L^\infty}\Vert  \vb^{\eps,m} + \wb_m^{\eps} \Vert_{1}\Vert \vb^{\eps,m}\Vert_{2}.
\end{align*}

Since $\vb^{\eps,m}(\tau^\ast)\in\mathbb H^2(\dom)$ and $\wb_m^\eps\in C([0,T];\mathbb H^2(\dom))$,
following the proof of \autoref{prop-6.1.2}, by the maximum principle, we infer that for every $T>\tau^\ast$,
\[\sup\limits_{t\in [\tau^\ast,\xi^{\eps,m}\wedge T)} \Vert \vb^{\eps,m}+\wb_m^{\eps} \Vert_{L^\infty}\leq \sup\limits_{t\in [\tau^\ast,\xi^{\eps,m}\wedge T)} \Big{(}\Vert \vb^{\eps,m}\Vert_{L^\infty}+\Vert\wb_m^{\eps} \Vert_{L^\infty}\Big{)}<\infty, \mathbb{P}\text{-a.s.}.\]
Therefore, by the Young inequality and the Gronwall inequality, we infer that
\[
\sup_{t\in[\tau^\ast,\xi^{\eps,m}\wedge T)}
\Vert\vb^{\eps,m}(t)\Vert_1<\infty,
\quad\mathbb P\text{-a.s.}
\]
Moreover, since
$\wb_m^\eps\in C([0,T];\mathbb H^2(\dom))$,
by the Sobolev embedding theorem, we infer that
\[
\sup_{t\in[\tau^\ast,\xi^{\eps,m}\wedge T)}
\Vert\vb^{\eps,m}(t)\Vert_{L^\infty}<\infty,
\quad\mathbb P\text{-a.s.}
\]
Recall that $\Vert\vb^{\eps,m}(t)\Vert_{\mathbb H^1(\dom)}^2
=\Vert\vb^{\eps,m}(t)\Vert_{L^2}^2
+\Vert\vb^{\eps,m}(t)\Vert_1^2$. Therefore, we infer that
\[
\sup_{t\in[\tau^\ast,\xi^{\eps,m}\wedge T)}
\Vert\vb^{\eps,m}(t)\Vert_{\mathbb H^1(\dom)}
<\infty,\quad\mathbb P\text{-a.s.}
\]
Since $\vb^{\eps,m}\in
C([0,\tau^\ast];\mathbb H^1(\dom))$,
if $\xi^{\eps,m}<\infty$, choosing
$T>\xi^{\eps,m}$ yields
\[
\sup_{t\in[0,\xi^{\eps,m})}
\Vert\vb^{\eps,m}(t)\Vert_{\mathbb H^1(\dom)}
<\infty.
\]
In view of \autoref{prop-tildev-local strong} and \eqref{ee-epsilon}, the solution can be extended
beyond $\xi^{\eps,m}$, which contradicts its maximality.
Therefore,
\begin{align}\label{eq-xiepsiloninfty}
	\xi^{\eps,m}=\infty,\ \ \  \mathbb{P}\text{-a.s.}.
\end{align}
The proof of \autoref{prop-global classical} is complete.
\end{proof}

Then we have the following important Lemma.
\begin{lemma}\label{lem-z-lim}
Assume that $\tilde\vb_0\in\mathbb{H}^2(\dom)$ and \autoref{A2} holds. For $m\in\mathbb{N}$, let $\{\vb^{\eps_{k},m}\}_{k\in \mathbb{N}}$ be an arbitrary subsequence of the global strong solutions $\{ \vb^{\eps,m}\}_{\eps> 0}$ to (\ref{ee-epsilon}) such that
\[
\eps_{k}\to 0,\;\;\mbox{ as } k\to \infty.
\]
Let $(\tilde\vb,\xi_2)$ be the maximal strong solution to the equation \eqref{ee-o} in the sense of \autoref{def-2.2 strong}.
Then, almost surely, for every $S\in (0,\xi_2)$,
\begin{align}\label{eq-z convergence H2}
\sup\limits_{t\in [0, S]}\Vert \vb^{\eps_{k},m}(t)-\tilde{\vb}(t)\Vert_{H^2}^{2}\rightarrow 0, \ \ \  \ as\ k\rightarrow \infty\ \mathrm{and}\ m\rightarrow \infty.
\end{align}

\end{lemma}
\begin{proof}

Let us choose and fix $T>0$. In the following, we consider $\omega\in\{\xi_2>T\}$. Firstly, we aim to show that for almost every $\omega\in\{\xi_2>T\}$, there exists $K_0(\omega),M_0(\omega)\in\N$ such that
\begin{align}\label{eq-veps-H23-bound}
\sup_{k\geq K_0(\omega)}\sup_{m\geq M_0(\omega)}\left[\sup\limits_{t\in [0, T]}\Vert \vb^{\eps_{k},m}(t)\Vert_{H^{2}}+\int_0^{T} \Vert \vb^{\eps_{k},m}(t)\Vert_{H^{3}}^{2}\,dt\right]<\infty.
\end{align}
Multiplying (\ref{ee-epsilon}) with $-\Delta \vb^{\eps_{k},m}$ and integrating over $\dom$, the H\"older inequality, the Young inequality and the triangle inequality yields,
\begin{align*}
	\frac{1}{2}\partial_{t} \Vert \vb^{\eps_{k},m}\Vert_{1}^{2} + \Vert  \vb^{\eps_{k},m}\Vert_{2}^{2} &= -\alpha^{\eps_{k}}\langle ((\vb^{\eps_{k},m}+\wb_m^{\eps_{k}} )\cdot \nabla) (\vb^{\eps_{k},m}+\wb_m^{\eps_{k}} ),  -\Delta \vb^{\eps_{k},m} \rangle\\
    &\leq \alpha^{\eps_{k}} \Vert\vb^{\eps_{k},m}+\wb_m^{\eps_{k}} \Vert_{L^\infty}\Vert \vb^{\eps_{k},m}+\wb_m^{\eps_{k}} \Vert_{1}\Vert \vb^{\eps_{k},m}\Vert_{2}\\
    &\leq (\alpha^{\eps_{k}})^2 \Vert\vb^{\eps_{k},m}+\wb_m^{\eps_{k}} \Vert_{L^\infty}^2(\Vert \vb^{\eps_{k},m}\Vert_{1}^2+\Vert\wb_m^{\eps_{k}} \Vert_{1}^2)+\frac{1}{2}\Vert \vb^{\eps_{k},m}\Vert_{2}^2,
\end{align*}
which implies
\begin{align}\label{eqn-v_meps}
	\partial_{t} \Vert \vb^{\eps_{k},m}\Vert_{1}^{2} + \Vert  \vb^{\eps_{k},m}\Vert_{2}^{2} \leq 2(\alpha^{\eps_{k}})^2 \Vert\vb^{\eps_{k},m}+\wb_m^{\eps_{k}} \Vert_{L^\infty}^2(\Vert \vb^{\eps_{k},m}\Vert_{1}^2+\Vert\wb_m^{\eps_{k}} \Vert_{1}^2).
\end{align}
By (\ref{eq-uniform convergence2}) in \autoref{lem-uniform convergence}, we infer that
\[
\sup\limits_{t\in [0, T]}\Vert \wb_m^{\eps_{k}} (t)-  \wb(t)\Vert_{H^2}\rightarrow 0,\  \mathbb{P}-\mathrm{a.s.},\ \mathrm{as}\ k\rightarrow \infty\ \mathrm{and} \ m\rightarrow \infty,
\]
In view of \autoref{cor-w-continuous-L2}, $\wb\in C([0,T]; \mathbb{H}^{2}(\dom))$.

Therefore, for almost every $\omega\in\Omega$, there exists $K_0(\omega),M_0(\omega)\in\N$ such that
\[
\sup_{k\geq K_0(\omega)}\sup_{m\geq M_0(\omega)}\sup_{t\in[0,T]}
\Vert\wb_m^{\eps_k}(t,\omega)\Vert_{H^2}
\leq \sup_{t\in[0,T]}\Vert\wb(t,\omega)\Vert_{H^2}+1<\infty.
\]

Therefore, in view of (\ref{eqn-v_meps}), \autoref{prop-6.1.2} and  using the Gronwall inequality we infer that
\begin{align*}
\sup_{k\geq K_0(\omega)}\sup_{m\geq M_0(\omega)}\left[\sup\limits_{t\in [0 ,T]}\Vert \vb^{\eps_{k},m}(t)\Vert_{1}+\int_{0}^{T} \Vert \vb^{\eps_{k},m}(t)\Vert_{2}^{2}\,dt\right]<\infty, \mathbb{P}\text{-a.s.}.
\end{align*}
By \autoref{prop-6.1.2} and \autoref{prop-global classical}, increasing $K_0$ so that $\eps_k<1$ for $k\geq K_0$ gives
\begin{equation}\label{eqn-approx-L2-uniform}
\sup_{k\geq K_0,\,m\geq M_0}\sup_{0\leq t\leq T}\Vert\vb^{\eps_k,m}(t)\Vert_{L^2}^2
\leq |\dom|\sup_{k\geq K_0,\,m\geq M_0}\sup_{0\leq t\leq T}\Vert\vb^{\eps_k,m}(t)\Vert_{L^\infty}^2<\infty .
\end{equation}
For $j=1,2$, the norm identity \eqref{eq-homogeneous norm} gives
\begin{align}\label{eqn-approx-full-Sobolev}
&\sup_{0\leq t\leq T}\Vert h(t)\Vert_{H^j}^2+\int_0^T\Vert h(t)\Vert_{H^{j+1}}^2\,dt\nonumber\\
&\qquad\leq (1+T)\sup_{0\leq t\leq T}\Vert h(t)\Vert_{L^2}^2
+\sup_{0\leq t\leq T}\Vert h(t)\Vert_j^2+\int_0^T\Vert h(t)\Vert_{j+1}^2\,dt .
\end{align}
Applying this inequality with $j=1$ to $h=\vb^{\eps_k,m}$ yields
\begin{align*}
\sup_{k\geq K_0(\omega)}\sup_{m\geq M_0(\omega)}\left[\sup\limits_{t\in [0 ,T]}\Vert \vb^{\eps_{k},m}(t)\Vert_{H^{1}}+\int_{0}^{T} \Vert \vb^{\eps_{k},m}(t)\Vert_{H^{2}}^{2}\,dt\right]<\infty, \mathbb{P}\text{-a.s.}.
\end{align*}

We record the energy identity used below.
If $h\in C(I;\mathbb H^2)\cap L^2(I;\mathbb H^3)$ solves $\partial_t h-\Delta h=F$ with $F\in L^2(I;\mathbb H^1)$ on a compact interval $I$, then $\Lambda^2h\in L^2(I;\mathbb H^1)$ and $\partial_t\Lambda^2h\in L^2(I;\mathbb H^{-1})$.
The Hilbert space chain rule and the Fourier multiplier identity \eqref{eq-Lambda-1} give, for almost every $t\in I$,
\begin{equation}\label{eqn-H2-energy-pairing}
\frac12\frac{d}{dt}\Vert h\Vert_2^2+\Vert h\Vert_3^2
=\langle\Lambda F,\Lambda^3h\rangle_{L^2}.
\end{equation}
The product rule and $\mathbb H^2\hookrightarrow\mathbb L^\infty\cap\mathbb W^{1,6}$ give
\begin{equation}\label{eqn-product-H2-H1}
\Vert(a\cdot\nabla)b\Vert_{H^1}\leq C\Vert a\Vert_{H^2}\Vert b\Vert_{H^2}.
\end{equation}
For each fixed $k,m$, the required regularity of $\vb^{\eps_k,m}$ holds near $0$ by \autoref{prop-maximal global-veps} with $i=2$, and on compact positive-time intervals by \autoref{prop-global classical} and the $\mathbb H^4$ regularity in Appendix~\ref{app-proof-classical-solution}.
Writing $V=\vb^{\eps_k,m}+\wb_m^{\eps_k}$ and using \eqref{eqn-H2-energy-pairing} in \eqref{ee-epsilon}, we obtain
\begin{align*}
\frac12\partial_t\Vert\vb^{\eps_k,m}\Vert_2^2+\Vert\vb^{\eps_k,m}\Vert_3^2
&=-\alpha^{\eps_k}\langle\Lambda[(V\cdot\nabla)V],\Lambda^3\vb^{\eps_k,m}\rangle_{L^2}\\
&\leq C\alpha^{\eps_k}\big(\Vert\nabla V\Vert_{L^4}^2+\Vert V\Vert_{L^\infty}\Vert V\Vert_2\big)\Vert\vb^{\eps_k,m}\Vert_3 .
\end{align*}

Recall the periodic Gagliardo--Nirenberg inequality
\begin{align*}
\Vert\nabla h\Vert_{L^4}^2
\leq
c\Vert h\Vert_{L^\infty}\Vert h\Vert_{H^2}
+
c\Vert h\Vert_{L^2}^2.
\end{align*}
Therefore, we infer that
\begin{align*}
\frac12\partial_t\Vert\vb^{\eps_k,m}\Vert_2^2
+\Vert\vb^{\eps_k,m}\Vert_3^2
&\leq
c\alpha^{\eps_k}
\Vert\vb^{\eps_k,m}+\wb_m^{\eps_k}\Vert_{L^\infty}
\Vert\vb^{\eps_k,m}+\wb_m^{\eps_k}\Vert_{H^2}
\Vert\vb^{\eps_k,m}\Vert_3
\\
&\quad+
c\alpha^{\eps_k}
\Vert\vb^{\eps_k,m}+\wb_m^{\eps_k}\Vert_{L^2}^2
\Vert\vb^{\eps_k,m}\Vert_3.
\end{align*}
Hence, by the Young inequality,
\begin{align}
\partial_t\Vert\vb^{\eps_k,m}\Vert_2^2
+\Vert\vb^{\eps_k,m}\Vert_3^2
&\leq
c(\alpha^{\eps_k})^2
\Vert\vb^{\eps_k,m}+\wb_m^{\eps_k}\Vert_{L^\infty}^2
\Vert\vb^{\eps_k,m}+\wb_m^{\eps_k}\Vert_{H^2}^2
\nonumber\\
&\quad+
c(\alpha^{\eps_k})^2
\Vert\vb^{\eps_k,m}+\wb_m^{\eps_k}\Vert_{L^2}^4.
\label{eq-h2-pre}
\end{align}

By \autoref{prop-6.1.2}, we have
\[
\sup_{k\geq K_0(\omega)}
\sup_{m\geq M_0(\omega)}
\sup_{t\in[0,T]}
\Vert\vb^{\eps_k,m}(t)\Vert_{L^\infty}
<\infty.
\]
Moreover, by \eqref{eq-uniform convergence2} and
$\wb\in C([0,T];\mathbb H^2(\dom))$,
\[
\sup_{k\geq K_0(\omega)}
\sup_{m\geq M_0(\omega)}
\sup_{t\in[0,T]}
\Vert\wb_m^{\eps_k}(t)\Vert_{H^2}
<\infty.
\]
Since $\mathbb H^2(\dom)\hookrightarrow\mathbb L^\infty(\dom)$,
we infer that
\[
\sup_{k\geq K_0(\omega)}
\sup_{m\geq M_0(\omega)}
\sup_{t\in[0,T]}
\Vert\vb^{\eps_k,m}(t)+\wb_m^{\eps_k}(t)\Vert_{L^\infty}
<\infty.
\]
In view of \eqref{eq-uniform convergence1},
\[
\sup_{k\geq K_0(\omega)}
\sup_{t\in[0,T]}
|\alpha^{\eps_k}(t)|
<\infty.
\]
Therefore, by \eqref{eq-h2-pre} we infer that
\[
\partial_t\Vert\vb^{\eps_k,m}\Vert_2^2
+\Vert\vb^{\eps_k,m}\Vert_3^2
\leq
C_T(\omega)
\left(
1+\Vert\vb^{\eps_k,m}\Vert_2^2
\right).
\]
Thus, by the Gronwall inequality,
\[
\sup_{k\geq K_0(\omega)}
\sup_{m\geq M_0(\omega)}
\left[
\sup_{t\in[0,T]}
\Vert\vb^{\eps_k,m}(t)\Vert_2^2
+
\int_0^T
\Vert\vb^{\eps_k,m}(t)\Vert_3^2\,dt
\right]
<\infty.
\]

Combining this estimate with \eqref{eqn-approx-L2-uniform} and \eqref{eqn-approx-full-Sobolev} for $j=2$ proves \eqref{eq-veps-H23-bound}.

Define a process
\[
y^{\eps_{k},m}:= \vb^{\eps_{k},m}-\tilde{\vb},\;\;k,m\in\mathbb{N},\;\;t\in [0,\xi_2).
\]
Then, we have on $[0, T]$
\begin{align*}
\partial_{t}y^{\eps_{k},m}-\Delta y^{\eps_{k},m}\!\!& +\!\![\alpha^{\eps_{k}}-\alpha]((\vb^{\eps_{k},m}+\wb_m^{\eps_{k}} )\cdot \nabla) (\vb^{\eps_{k},m}+\wb_m^{\eps_{k}} )+\alpha[ y^{\eps_{k},m}+(\wb_m^{\eps_{k}} -\wb) ]\cdot \nabla (\vb^{\eps_{k},m}+\wb_m^{\eps_{k}} )\\
&+\alpha (( \tilde{\vb}+\wb)\cdot \nabla) [ y^{\eps_{k},m}+(\wb_m^{\eps_{k}} -\wb) ]=0.
\end{align*}

Multiplying $y^{\eps_{k},m}$ and integrating over $\dom$ yields,
\begin{align}\label{eq-difference L2}
\frac{1}{2} \partial_{t} \Vert y^{\eps_{k},m}\Vert_{L^2}^{2}+\Vert y^{\eps_{k},m}\Vert_{1}^{2}
&=(\alpha-\alpha^{\eps_{k}})\langle ((\vb^{\eps_{k},m}+\wb_m^{\eps_{k}} )\cdot \nabla) (\vb^{\eps_{k},m}+\wb_m^{\eps_{k}} ),  y^{\eps_{k},m} \rangle\\
&-\alpha\langle ([ y^{\eps_{k},m}+(\wb_m^{\eps_{k}} -\wb) ]\cdot \nabla) (\vb^{\eps_{k},m}+\wb_m^{\eps_{k}} ),    y^{\eps_{k},m} \rangle\nonumber\\
&-\alpha\langle (( \tilde{\vb}+\wb)\cdot \nabla) [ y^{\eps_{k},m}+(\wb_m^{\eps_{k}} -\wb) ], y^{\eps_{k},m} \rangle\nonumber\\
&:=J_{1}+J_{2}+J_{3}.\nonumber
\end{align}
Let us choose and fix $\delta\in(0,\frac12)$.
By the H\"{o}lder inequality, the Sobolev inequality and the Young inequality we get
\begin{align*}
\vert J_{1}\vert
&\leq |\alpha-\alpha^{\eps_{k}}|\Vert\vb^{\eps_{k},m}+\wb_m^{\eps_{k}} \Vert_{L^\infty}
\Vert\nabla( \vb^{\eps_{k},m}+\wb_m^{\eps_{k}} )\Vert_{L^2}\Vert y^{\eps_{k},m}\Vert_{L^2} \\
&\leq c|\alpha-\alpha^{\eps_{k}}|\Vert\vb^{\eps_{k},m}+\wb_m^{\eps_{k}} \Vert_{H^2}^2 \Vert y^{\eps_{k},m}\Vert_{L^2}\\
&\leq c|\alpha-\alpha^{\eps_{k}}|^{2}+c\Vert\vb^{\eps_{k},m}+\wb_m^{\eps_{k}} \Vert_{H^2}^4 \Vert y^{\eps_{k},m}\Vert_{L^2}^{2}.
\end{align*}
Similarly, since $H^2(\dom)\hookrightarrow W^{1,6}(\dom)$ we have
\begin{align*}
J_{2}
&\leq \alpha \Vert y^{\eps_{k},m}+(\wb_m^{\eps_{k}} -\wb)\Vert_{L^2}\Vert\nabla( \vb^{\eps_{k},m}+\wb_m^{\eps_{k}} )\Vert_{L^3}
\Vert y^{\eps_{k},m}\Vert_{L^6}\\
&\leq \alpha \left(\Vert y^{\eps_{k},m}\Vert_{L^2}+\Vert\wb_m^{\eps_{k}} -\wb\Vert_{L^2}\right)\Vert \vb^{\eps_{k},m}+\wb_m^{\eps_{k}} \Vert_{H^2}
\Vert y^{\eps_{k},m}\Vert_{H^1}\\
&\leq \alpha \Vert \vb^{\eps_{k},m}+\wb_m^{\eps_{k}} \Vert_{H^2} \left(\Vert y^{\eps_{k},m}\Vert_{L^2}+\Vert\wb_m^{\eps_{k}} -\wb\Vert_{L^2}\right)\left(\Vert y^{\eps_{k},m}\Vert_{L^2}+\Vert y^{\eps_{k},m}\Vert_{1}\right)
\\
&\leq C_{T,\omega,\delta}\big(\Vert y^{\eps_k,m}\Vert_{L^2}^2+\Vert\wb_m^{\eps_k}-\wb\Vert_{L^2}^2\big)+\delta\Vert y^{\eps_k,m}\Vert_1^2 .
\end{align*}
Here $C_{T,\omega,\delta}$ is uniform for $k\geq K_0$, $m\geq M_0$, by \eqref{eq-veps-H23-bound}, \eqref{eq-uniform convergence2} and the boundedness of $\alpha$ on $[0,T]$.
For the last term we have
\begin{align*}
J_{3}
&\leq\alpha\Vert  \tilde{\vb}+\wb\Vert_{L^\infty}\Vert  y^{\eps_{k},m}+(\wb_m^{\eps_{k}} -\wb)\Vert_{1}\Vert y^{\eps_{k},m}\Vert_{L^2}\\
&\leq\alpha\Vert  \tilde{\vb}+\wb\Vert_{H^2}\left(\Vert  y^{\eps_{k},m}\Vert_{1}+\Vert\wb_m^{\eps_{k}} -\wb\Vert_{1}\right)\Vert y^{\eps_{k},m}\Vert_{L^2}\\
&\leq c\alpha^2\Vert  \tilde{\vb}+\wb\Vert_{H^2}^2\Vert y^{\eps_{k},m}\Vert_{L^2}^2+\delta\Vert y^{\eps_{k},m}\Vert_{1}^{2}+c\Vert \wb_m^{\eps_{k}} -\wb\Vert_{1}^{2},
\end{align*}
Combine the above estimates, we infer that
\begin{align*}
\frac{1}{2} \partial_{t} \Vert y^{\eps_{k},m}\Vert_{L^2}^{2}+(1-2\delta)\Vert y^{\eps_{k},m}\Vert_{1}^{2}\leq c\Vert y^{\eps_{k},m}\Vert_{L^2}^2+c|\alpha(t)-\alpha^{\eps_{k}}(t)|^{2}+c\Vert \wb_m^{\eps_{k}} -\wb\Vert_{H^1}^{2}
\end{align*}
where we have used $\sup\limits_{t\in [0, T]}\Vert \tilde{\vb}+\wb\Vert_{H^{2}(\dom)}<\infty,\, \mathbb{P}\text{-a.s.}. $
Hence, by the Gronwall inequality, we deduce
\begin{align}\label{eq-w convergence l2}
\sup\limits_{t\in [0, T]}\Vert y^{\eps_{k},m}(t)\Vert_{L^2}^{2}&\leq c\Big{(}\sup\limits_{t\in [0, T]}\Vert \wb_m^{\eps_{k}}(t) -\wb(t)\Vert_{H^1}^{2}+ \sup\limits_{t\in [0, T]}|\alpha^{\eps_{k}}-\alpha|^{2}\Big{)}\nonumber\\
&\rightarrow 0, \ \ \  \ as\ k\rightarrow \infty,\ \mathrm{and}\ m\rightarrow \infty, \ \mathbb{P}\text{-a.s.}.
\end{align}

On $[0,T]\subset[0,\xi_2)$ both solutions lie in $C([0,T];\mathbb H^2)\cap L^2([0,T];\mathbb H^3)$, and \eqref{eqn-product-H2-H1} places each nonlinear term in $L^2([0,T];\mathbb H^1)$.
Set $D=y^{\eps_k,m}+\wb_m^{\eps_k}-\wb$ and recall $V=\vb^{\eps_k,m}+\wb_m^{\eps_k}$.
Applying \eqref{eqn-H2-energy-pairing} to the difference gives
\begin{align}\label{eq-difference h2}
&\frac12\partial_t\Vert y^{\eps_k,m}\Vert_2^2+\Vert y^{\eps_k,m}\Vert_3^2\nonumber\\
&=(\alpha-\alpha^{\eps_k})\langle\Lambda[(V\cdot\nabla)V],\Lambda^3y^{\eps_k,m}\rangle_{L^2}\nonumber\\
&\quad-\alpha\langle\Lambda[(D\cdot\nabla)V],\Lambda^3y^{\eps_k,m}\rangle_{L^2}\nonumber\\
&\quad-\alpha\langle\Lambda[((\tilde\vb+\wb)\cdot\nabla)D],\Lambda^3y^{\eps_k,m}\rangle_{L^2}\\
&=:I_1+I_2+I_3 .\nonumber
\end{align}
Let us choose and fix $\delta_1\in(0,\frac13)$. By the H\"{o}lder inequality, the Sobolev inequality and the Young inequality we get
\begin{align*}
I_{1}
&\leq c\sup\limits_{t\in [0, T]}|\alpha(t)-\alpha^{\eps_k}(t)|\Vert y^{\eps_{k},m}\Vert_{3}(\Vert \vb^{\eps_{k},m}+\wb_m^{\eps_{k}} \Vert_{H^2}^{2}+\Vert\vb^{\eps_{k},m}+\wb_m^{\eps_{k}} \Vert_{L^\infty}
\Vert \vb^{\eps_{k},m}+\wb_m^{\eps_{k}} \Vert_{2} )\\
&\leq c\sup\limits_{t\in [0, T]}|\alpha(t)-\alpha^{\eps_k}(t)|^{2}(\Vert \vb^{\eps_{k},m}+\wb_m^{\eps_{k}} \Vert_{H^2}^{4}+\Vert\vb^{\eps_{k},m}+\wb_m^{\eps_{k}} \Vert_{L^\infty}^{2}
\Vert \vb^{\eps_{k},m}+\wb_m^{\eps_{k}} \Vert_{2}^{2} )+\delta_1 \Vert y^{\eps_{k},m}\Vert_{3}^{2}\\
&\leq c\sup\limits_{t\in [0, T]}|\alpha(t)-\alpha^{\eps_k}(t)|^{2}+\delta_1 \Vert y^{\eps_{k},m}\Vert_{3}^{2}.
\end{align*}
Similarly, we have
\begin{align*}
I_{2}
&\leq c\alpha\Vert y^{\eps_{k},m}\Vert_{3}\Big{[} \Vert \nabla[y^{\eps_{k},m}+(\wb_m^{\eps_{k}} -\wb)]\Vert_{L^3} \Vert \nabla(\vb^{\eps_{k},m}+\wb_m^{\eps_{k}}) \Vert_{L^6}
+\Vert y^{\eps_{k},m}+(\wb_m^{\eps_{k}} -\wb)\Vert_{L^\infty}\Vert \vb^{\eps_{k},m}+\wb_m^{\eps_{k}} \Vert_{2}\Big{]}\\
&\leq c\alpha\Vert y^{\eps_{k},m}\Vert_{3} \Vert y^{\eps_{k},m}+(\wb_m^{\eps_{k}} -\wb)\Vert_{H^2} \Vert \vb^{\eps_{k},m}+\wb_m^{\eps_{k}} \Vert_{H^2}\\
&\leq c\alpha\Vert y^{\eps_{k},m}\Vert_{3} \left[\Vert y^{\eps_{k},m}\Vert_{2}+\Vert y^{\eps_{k},m}\Vert_{L^2}+\Vert(\wb_m^{\eps_{k}} -\wb)\Vert_{H^2} \right]\Vert \vb^{\eps_{k},m}+\wb_m^{\eps_{k}} \Vert_{H^2}\\
&\leq c(\Vert y^{\eps_{k},m}\Vert_{2}^{2}+\Vert y^{\eps_{k},m}\Vert_{L^2}^{2}+\Vert \wb_m^{\eps_{k}} -\wb\Vert_{H^2}^{2})+\delta_1 \Vert y^{\eps_{k},m}\Vert_{3}^{2},
\end{align*}
and
\begin{align*}
I_{3}
&\leq\alpha\Vert y^{\eps_{k},m}\Vert_{3}\Vert  \tilde{\vb}+\wb\Vert_{H^2}\Vert  y^{\eps_{k},m}+(\wb_m^{\eps_{k}} -\wb)\Vert_{H^2}\\
&\leq\alpha\Vert y^{\eps_{k},m}\Vert_{3}\Vert  \tilde{\vb}+\wb\Vert_{H^2}\left[\Vert y^{\eps_{k},m}\Vert_{2}+\Vert y^{\eps_{k},m}\Vert_{L^2}+\Vert(\wb_m^{\eps_{k}} -\wb)\Vert_{H^2} \right]\\
&\leq c(\Vert y^{\eps_{k},m}\Vert_{2}^{2}+\Vert y^{\eps_{k},m}\Vert_{L^2}^2+\Vert \wb_m^{\eps_{k}} -\wb\Vert_{H^2}^{2})+\delta_1 \Vert y^{\eps_{k},m}\Vert_{3}^{2},
\end{align*}
By (\ref{eq-difference h2}) and estimates of $I_{1}, I_{2}, I_{3}$ we infer that
\[
\frac{1}{2} \partial_{t} \Vert y^{\eps_{k},m}\Vert_{2}^{2}+(1-3\delta_1)\Vert y^{\eps_{k},m}\Vert_{3}^{2}\leq c(\Vert y^{\eps_{k},m}\Vert_{2}^{2}+\Vert y^{\eps_{k},m}\Vert_{L^2}^2+\Vert \wb_m^{\eps_{k}} -\wb\Vert_{H^2}^{2}+\sup\limits_{t\in [0, T]}|\alpha(t)-\alpha^{\eps_k}(t)|^{2})
\]
where we have used $\sup\limits_{t\in [0, T]}\Vert \tilde{\vb}+\wb\Vert_{H^{2}(\dom)}<\infty,\, \mathbb{P}\text{-a.s.}. $
Therefore, by the Gronwall inequality we deduce
\begin{align}
\sup\limits_{t\in [0, T]}\Vert y^{\eps_{k},m}(t)\Vert_{2}^{2}&\leq c \Big{(}\sup\limits_{t\in [0, T]}\Vert \wb_m^{\eps_{k}}(t) -\wb(t)\Vert_{H^2}^{2}+ \sup\limits_{t\in [0, T]}|\alpha(t)-\alpha^{\eps_{k}}(t)|^{2}+\sup\limits_{t\in [0, T]}\Vert y^{\eps_{k},m}(t)\Vert_{L^2}^{2}\Big{)}\nonumber\\
&\rightarrow 0, \ \ \  \ as\ k\rightarrow \infty,\ \mathrm{and}\ m\rightarrow \infty, \ \mathbb{P}\text{-a.s.},
\end{align}

Combining the above estimates, we infer that for every fixed $T>0$,
\[
\sup_{t\in[0,T]}
\Vert\vb^{\eps_k,m}(t)-\tilde{\vb}(t)\Vert_{H^2}^2
\to 0
\quad\text{on }\{\xi_2>T\},
\quad\mathbb P\text{-a.s.},\;\;\mbox{ as } k\rightarrow\infty, m\rightarrow\infty.
\]

Since $\mathbb Q_+$ is countable, there exists $\Omega_0\subset\Omega$ with $\mathbb P(\Omega_0)=1$
such that the above convergence holds for every $T\in\mathbb Q_+$ satisfying $T<\xi_2(\omega)$, for every $\omega\in\Omega_0$.

Let us choose and fix $\omega\in\Omega_0$ and $S\in(0,\xi_2(\omega))$. There exists $T\in\mathbb Q_+$ such that $S<T<\xi_2(\omega)$. Therefore,
\begin{align*}
&\sup_{t\in[0,S]}
\Vert\vb^{\eps_k,m}(t)-\tilde{\vb}(t)\Vert_{H^2}^2\leq
\sup_{t\in[0,T]}
\Vert\vb^{\eps_k,m}(t)-\tilde{\vb}(t)\Vert_{H^2}^2
\to 0, \;\;\mbox{ as } k\rightarrow\infty, m\rightarrow\infty.
\end{align*}
Hence, \eqref{eq-z convergence H2} holds.
The proof of \autoref{lem-z-lim} is complete.

\end{proof}

From the uniform convergence of $\vb^{\eps_{k},m}$ to $\tilde{\vb}$ we established in \autoref{lem-z-lim}, we will obtain the uniform a priori estimates for $\tilde{\vb}$.
\begin{prop}\label{prop-global weak}
Assume that $\tilde\vb_0\in\mathbb{H}^2(\dom)$ and \autoref{A2} holds.
Let $(\tilde\vb,\xi_2)$ be the maximal strong solution to the equation \eqref{ee-o} in the sense of \autoref{def-2.2 strong}. Then for every $T>0$
\begin{align*}
	\sup\limits_{t\in \lbrack 0,\xi_2\wedge T )}\Vert\tilde{\vb}\Vert_{L^\infty}&\leq C(T,\Vert\tilde\vb_0\Vert_{L^\infty}, \Vert \fb\Vert_{L^{2}([0,T+1]; \mathbb{H}^\delta(\dom)) }, \Vert\alpha^{-1}\Vert_{L^{\infty}([0,T+1])})<\infty,\;\mathbb{P}\text{-a.s.}.
\end{align*}

\end{prop}

\begin{proof}
Let $\{\vb^{\eps_{k},m}\}_{k\in \mathbb{N}}$ be an arbitrary subsequence of the global strong solutions $\{ \vb^{\eps,m}\}_{\eps> 0}$ to (\ref{ee-epsilon}) such that
\[
\eps_{k}\to 0,\;\;\mbox{ as } k\to\infty.
\]

Let us choose and fix $T>0$. Since $\eps_k\to0$ as $k\to\infty$, let us choose and fix $k$ sufficiently large, such that $\eps_k\in(0,1)$. By \autoref{prop-6.1.2} and \eqref{eq-xiepsiloninfty}, we infer that
\begin{align*}
\sup\limits_{m\in\mathbb N}\sup\limits_{t\in[0,T]}\Vert\vb^{\eps_k,m}(t)\Vert_{L^\infty}\leq C\left(T,\Vert\tilde\vb_0\Vert_{L^\infty},\Vert\fb\Vert_{L^2([0,T+1];\mathbb H^\delta(\dom))},\Vert\alpha^{-1}\Vert_{L^\infty([0,T+1])}\right)=:C_\omega(T)<\infty,
\;\;\mathbb P\text{-a.s.}
\end{align*}

Let us choose and fix $\tau\in(0,\xi_2\wedge T)$. In view of \autoref{lem-z-lim}, we infer that
\begin{align}\label{eq-tilde v uniform}
\sup\limits_{t\in[0,\tau]}
\Vert\tilde\vb(t)\Vert_{L^\infty}
&\leq
\sup\limits_{t\in[0,\tau]}
\Vert\tilde\vb(t)-\vb^{\eps_k,m}(t)\Vert_{L^\infty}
+
\sup\limits_{t\in[0,\tau]}
\Vert\vb^{\eps_k,m}(t)\Vert_{L^\infty}
\nonumber\\
&\leq
c\sup\limits_{t\in[0,\tau]}
\Vert\tilde\vb(t)-\vb^{\eps_k,m}(t)\Vert_{H^2}
+
C_\omega(T).
\end{align}
Letting $m\to\infty$ and $k\to\infty$, we infer that
\[
\sup\limits_{t\in[0,\tau]}
\Vert\tilde\vb(t)\Vert_{L^\infty}
\leq C_\omega(T),
\qquad \mathbb P\text{-a.s.}
\]
Since $C_\omega(T)$ is independent of $\tau$, and
$\tau\in(0,\xi_2\wedge T)$ is arbitrary, we deduce that
\[
\sup\limits_{t\in[0,\xi_2\wedge T)}
\Vert\tilde\vb(t)\Vert_{L^\infty}
\leq C_\omega(T)<\infty,
\qquad \mathbb P\text{-a.s.}
\]
Since $T>0$ is arbitrary, the proof of \autoref{prop-global weak} is complete.
\end{proof}

\begin{prop}\label{prop-pathwise-H2-continuation}
Assume that \autoref{A2} holds. Let us choose and fix $\omega\in\Omega$ outside a $\mathbb P$-null set, $t_0\geq0$ and $\tilde{\vb}_0\in\mathbb H^2(\dom)$. Then there exists a unique maximal pathwise strong solution $\tilde{\vb}$ to \eqref{ee-o} on $[t_0,\zeta)$ with $\tilde{\vb}(t_0)=\tilde{\vb}_0$ such that
\[
\tilde{\vb}\in C([t_0,\zeta);\mathbb H^2(\dom))\cap L^2_{\mathrm{loc}}([t_0,\zeta);\mathbb H^3(\dom)).
\]
Moreover, if $\zeta<\infty$, then
\begin{align}\label{eq-pathwise-H2-blowup}
\limsup_{t\uparrow\zeta}\Vert\tilde{\vb}(t)\Vert_{\mathbb H^2}=\infty.
\end{align}
\end{prop}

\begin{proof}
Since $\alpha(\cdot,\omega)$ is continuous and $\wb(\cdot,\omega)\in C([0,\infty);\mathbb H^2(\dom))$, by the local existence and pathwise uniqueness established in \autoref{prop-tildev-local strong}, we infer that there exists a unique maximal pathwise strong solution $\tilde{\vb}$ to \eqref{ee-o} on $[t_0,\zeta)$ with $\tilde{\vb}(t_0)=\tilde{\vb}_0$.

It remains to prove \eqref{eq-pathwise-H2-blowup}. Suppose by contradiction that $\zeta<\infty$ and
\[
M:=\sup_{t\in[t_0,\zeta)}\Vert\tilde{\vb}(t)\Vert_{\mathbb H^2}<\infty.
\]
Since $\Vert\tilde{\vb}(t_1)\Vert_{\mathbb H^2}\leq M$ for every $t_1\in[t_0,\zeta)$, by \autoref{prop-tildev-local strong} we infer that there exists $\delta_0:=\delta_0\left(M,\Vert\alpha\Vert_{L^\infty([0,\zeta+1])}, \Vert\wb\Vert_{C([0,\zeta+1];\mathbb{H}^2(\dom))} \right)>0$, such that there exists a unique  strong solution  to \eqref{ee-o} on $[t_1,t_1+\delta_0]$ with initial data $\tilde{\vb}(t_1)$ at time $t_1$. Let us choose and fix $t_1\in[t_0,\zeta)$ sufficiently close to $\zeta$ such that $t_1+\delta_0>\zeta$. By pathwise uniqueness, this solution extends $\tilde{\vb}$ beyond $\zeta$, which contradicts the maximality of $\zeta$. Therefore, \eqref{eq-pathwise-H2-blowup} holds. The proof is complete.
\end{proof}

\begin{theorem}\label{thm-global weak transformed equation}

Assume that $u_0 \in \mathbb{L}^{2}(\dom)$.
Let $\overrightarrow{\ub_{0}} \in  \mathbb{L}^{2}(\Omega; \mathbb{L}^{2}(\dom))$  be a randomization of $\ub_0$ defined by   \eqref{eqn-randomization of u_0}, using a sequence $\{r_{\kb}^j\}_{(\kb,j)\in\mathcal K}$ satisfying \autoref{A1}, and assume that $f$ satisfies \autoref{A2}.
Then there is a unique global weak solution $\tilde{\vb}$ to (\ref{ee-o}).
Moreover, there exists an event of probability one on which, for every $0<t_0<T<\infty$,
\[
\tilde{\vb}\in C([t_0,T];\mathbb H^2(\dom))\cap L^2([t_0,T];\mathbb H^3(\dom)).
\]

\end{theorem}
\begin{proof}[Proof of Theorem \ref{thm-global weak transformed equation}]

Let us choose and fix $u_0 \in \mathbb{L}^{2}(\dom)$ and its randomization of $\ub_0$  defined by   \eqref{eqn-randomization of u_0}.
Let $(\tilde{\vb}, \xi)$ be the maximal weak solution to (\ref{ee-o}) which exists by \autoref{prop-maximal-tildev}.

In order to prove that the existence of $\tilde{\vb}$ is global, it suffices to prove $\xi=\infty$, $\mathbb{P}-$a.s..
Let us choose and fix $\omega\in\Omega$ outside a null set such that
all the pathwise properties used below hold. By \autoref{cor-ae-H2}, there exists $\tau^\ast\in(0,\xi(\omega))$ such that
\[
\tilde{\vb}(\tau^\ast,\omega)\in\mathbb H^2(\dom).
\]

In the following, we regard $\tau^\ast$ as the initial time for the system \eqref{ee-o}. By \autoref{prop-pathwise-H2-continuation}, there exists a unique maximal pathwise strong solution on $[\tau^\ast,\xi_2)$, which we still denote by $\tilde{\vb}$. By pathwise uniqueness, this strong solution agrees with the original weak solution on their common interval.

Multiplying (\ref{ee-o}) with $\Lambda^{2}\tilde{\vb}$ and integrating over $\dom$, by the Young inequality, there exists $C>0$ such that
\begin{align*}
	\frac{1}{2}\partial_{t}\Vert \tilde{\vb}\Vert_{1}^{2}+\Vert \tilde{\vb}\Vert_{2}^{2} &\leq \alpha\Vert \tilde{\vb}+\wb\Vert_{L^\infty}\Vert  \tilde{\vb} + \wb\Vert_{1}\Vert \tilde{\vb}\Vert_{2}\leq \frac12\Vert \tilde{\vb}\Vert_{2}^2+C\alpha^2\Vert \tilde{\vb}+\wb\Vert_{L^\infty}^2\Vert  \tilde{\vb} + \wb\Vert_{1}^2.
\end{align*}
Fix a sample path, write $s=\tau^\ast$, and choose $T>s$ and $s<S<\xi_2\wedge T$.
Set $U=\tilde\vb+\wb$ and $g=\alpha^{-1}\fb$. Adding the equation for $\wb$ to \eqref{ee-o} gives
\[
\partial_tU-\Delta U+\alpha(U\cdot\nabla)U=g,
\qquad U(s)=\tilde\vb(s)+\wb(s)\in\mathbb H^2 .
\]
By \autoref{A2} and $\delta>3/2$, $g\in L^2([s,T];\mathbb L^\infty)$.
Since $U\in C([s,S];\mathbb H^2)$, the equation gives $\partial_tU\in L^2([s,S];\mathbb L^2)$.
For each component $j$, let
\[
M_j(t)=\Vert U_j(s)\Vert_{L^\infty}+\int_s^t\Vert g_j(r)\Vert_{L^\infty}\,dr,
\qquad q=(U_j-M_j)_+ .
\]
Then $q(s)=0$ and $\nabla q=\mathbf 1_{\{U_j>M_j\}}\nabla U_j$.
The positive-part chain rule permits testing the component equation with $q$.
Since $(g_j-M_j')q\leq0$, the Cauchy--Schwarz and Young inequalities give
\[
\frac12\frac{d}{dt}\Vert q\Vert_{L^2}^2+\Vert\nabla q\Vert_{L^2}^2
\leq \Vert\alpha U\Vert_{L^\infty}\Vert\nabla q\Vert_{L^2}\Vert q\Vert_{L^2}
\leq \frac12\Vert\nabla q\Vert_{L^2}^2+\frac12\Vert\alpha U\Vert_{L^\infty}^2\Vert q\Vert_{L^2}^2 .
\]
The coefficient is bounded on $[s,S]$, so Gronwall's inequality gives $q=0$ there.
Applying the same argument to $-U_j$ with force $-g_j$ yields $|U_j(t)|\leq M_j(t)$ almost everywhere in space.
Consequently,
\[
\sup_{s\leq t\leq S}\Vert U(t)\Vert_{L^\infty}
\leq \sqrt{3}\left(\Vert U(s)\Vert_{L^\infty}+\int_s^T\Vert g(r)\Vert_{L^\infty}\,dr\right).
\]
The bound is independent of $S$.
Taking the supremum over $S<\xi_2\wedge T$ and using $\wb\in C([0,T];\mathbb H^2)$ from \autoref{cor-w-continuous-L2} gives
\begin{align}\label{eqn-vw-Linfty}
\sup\limits_{t\in [\tau^\ast,\xi_2 \wedge T)}  (\Vert \tilde{\vb}\Vert_{L^\infty}+\Vert\wb\Vert_{L^\infty})<\infty, \mathbb{P}\text{-a.s.}.
\end{align}
The Gronwall inequality implies that
\begin{align}\label{eqn-H1H2}
\sup\limits_{t\in [\tau^\ast,\xi_2 \wedge T)}  \Vert  \tilde{\vb}(t)\Vert_{1}+\int_{\tau^\ast}^{\xi_2 \wedge T}\Vert  \tilde{\vb}(t)\Vert_{2}^{2}\,dt  <\infty,\ \ \ \ \mathbb{P}\text{-a.s.}.
\end{align}
On each $[s,S]$, \eqref{eqn-product-H2-H1} gives $(U\cdot\nabla)U\in C([s,S];\mathbb H^1)$.
Applying \eqref{eqn-H2-energy-pairing} to \eqref{ee-o} and using $\Vert\nabla U\Vert_{L^4}\leq C\Vert U\Vert_2$ on the torus, we obtain
\begin{align*}
\frac12\partial_t\Vert\tilde\vb\Vert_2^2+\Vert\tilde\vb\Vert_3^2
&=-\alpha\langle\Lambda[(U\cdot\nabla)U],\Lambda^3\tilde\vb\rangle_{L^2}\\
&\leq C\alpha\big(\Vert U\Vert_{L^\infty}\Vert U\Vert_2+\Vert U\Vert_2^2\big)\Vert\tilde\vb\Vert_3\\
&\leq \frac12\Vert\tilde\vb\Vert_3^2+C\alpha^2\big(\Vert U\Vert_{L^\infty}^2+\Vert U\Vert_2^2\big)\Vert U\Vert_2^2 .
\end{align*}
Note that $\wb\in C([0,T];\mathbb H^2)$ and  $\alpha^2(\Vert\tilde\vb+\wb\Vert_{L^\infty}^2+
\Vert\tilde\vb+\wb\Vert_2^2)$ is integrable on $[\tau^\ast,\xi_2\wedge T)$. In view of \eqref{eqn-vw-Linfty} and \eqref{eqn-H1H2}, by the Gronwall inequality we infer that
\[
\sup\limits_{t\in [\tau^\ast,\xi_2 \wedge T)}  \Vert  \tilde{\vb}(t)\Vert_{2} <\infty,\ \ \ \ \mathbb{P}\text{-a.s.}.
\]

Since $\Vert\tilde\vb\Vert_{H^2}^2
\leq |\dom|\Vert\tilde\vb\Vert_{L^\infty}^2+
\Vert\tilde\vb\Vert_2^2$, if $\xi_2<\infty$, taking $T>\xi_2$ gives
\[
\sup_{t\in[\tau^\ast,\xi_2)}
\Vert\tilde{\vb}(t,\omega)\Vert_{H^2}<\infty,
\]
which contradicts the blow-up alternative \eqref{eq-pathwise-H2-blowup} in
\autoref{prop-pathwise-H2-continuation}.

Hence
\[
\xi_2=\infty.
\]

{
It remains to prove that the pathwise continuation yields a global adapted weak solution. For every $r\in\mathbb Q_+$, define
\[
E_r:=\{r<\xi,\ \tilde{\vb}(r)\in\mathbb H^2(\dom)\}.
\]
Let $V_r$ be $\tilde{\vb}(r)$ on $\{r<\xi\}$ and zero elsewhere. Since $\xi$ is a stopping time and the solution on its lifetime is adapted, $V_r$ is $\mathscr F_r$-measurable in $\mathbb L^2(\dom)$. If $P_N$ denotes the Fourier cutoff at frequency $N$, then
\[
E_r=\{r<\xi\}\cap\bigl\{\sup_{N\geq1}\|P_N V_r\|_{H^2}<\infty\bigr\}\in\mathscr F_r.
\]
By \autoref{cor-ae-H2}, for almost every $\omega\in\Omega$, there exists $\tau^\ast\in(0,\xi(\omega))$ such that $\tilde{\vb}(\tau^\ast)\in\mathbb H^2(\dom)$. In view of \autoref{prop-pathwise-H2-continuation} and the pathwise uniqueness of weak solutions, there exists $\delta>0$ such that $\tilde{\vb}(t)\in\mathbb H^2(\dom)$ for every $t\in[\tau^\ast,\tau^\ast+\delta]\cap[0,\xi(\omega))$. Therefore, by the density of $\mathbb Q_+$ in $(0,\infty)$, we deduce that
\[
\mathbb P\left(\bigcup_{r\in\mathbb Q_+}E_r\right)=1.
\]

Now fix $r\in\mathbb Q_+$ and put $\tilde{\vb}_r=\mathds1_{E_r}V_r$. The finite-dimensional variables $\mathds1_{E_r}P_N V_r$ are $\mathscr F_r$-measurable and converge to $\tilde{\vb}_r$ in $\mathbb H^2(\dom)$ on every path. Thus $\tilde{\vb}_r$ is measurable as an $\mathbb H^2(\dom)$-valued initial condition.
Taking $r$ as the initial time for equation \eqref{ee-o}, by the shifted version of \autoref{prop-maximal global-tildev}, we obtain a unique maximal adapted strong solution $(\tilde{\vb}^{\,r},\xi_r)$ with $\tilde{\vb}^{\,r}(r)=\tilde{\vb}_r$. The above pathwise estimates apply on $[r,\xi_r\wedge T)$, with constants allowed to depend on $\tilde{\vb}_r(\omega)$ and the coefficient path. On $\{\xi_r<\infty\}$, choosing an integer $T>\xi_r$ therefore gives a bounded $H^2$ norm on $[r,\xi_r)$, contradicting the shifted adapted blow-up alternative \eqref{eq-veps-H1-blowup}. Hence $\xi_r=\infty$ almost surely. Moreover, by the pathwise uniqueness of weak solutions, we deduce that
\[
\tilde{\vb}^{\,r}(t)=\tilde{\vb}(t),\qquad r\leq t<\xi(\omega),\quad\text{on }E_r.
\]

In the following, we aim to prove that $\xi=\infty$, $\mathbb P$-a.s.. Assume by contradiction that $\mathbb P(\xi<\infty)>0$. Since $\mathbb Q_+$ is countable and $\mathbb P(\bigcup_{r\in\mathbb Q_+}E_r)=1$, there exist $r\in\mathbb Q_+$ and $T>r$ such that
\[
\mathbb P(E_r\cap\{\xi<T\})>0.
\]
Let $\{\tau_n\}_{n\in\mathbb N}$ be the increasing sequence of stopping times in \autoref{prop-maximal-tildev} satisfying $\tau_n\uparrow\xi$, $\mathbb P$-a.s.. For every $n\in\mathbb N$, set
\[
F_n:=E_r\cap\{\tau_n\geq r\}.
\]
Since $E_r\in\mathscr F_r$ and $\tau_n$ is a stopping time, we infer that $F_n\in\mathscr F_r$. Moreover, since $\tau_n\uparrow\xi$ and $r<\xi$ on $E_r$, we have
\[
E_r=\bigcup_{n=1}^\infty F_n.
\]
Therefore, there exists $n\in\mathbb N$ such that
\[
\mathbb P(F_n\cap\{\xi<T\})>0.
\]

Let us choose and fix such an $n$. Define
\[
\widehat\tau_n:=
\begin{cases}
\tau_n\vee T,&\text{on }F_n,\\
\tau_n,&\text{on }F_n^c.
\end{cases}
\]
Since $F_n\in\mathscr F_r$, $\tau_n\geq r$ on $F_n$ and $T>r$, we infer that $\widehat\tau_n$ is a stopping time. Moreover, $\widehat\tau_n\geq\tau_n$ almost surely.

On $F_n$, we define $\widehat{\vb}=\tilde{\vb}$ on $[0,r]$ and $\widehat{\vb}=\tilde{\vb}^{\,r}$ on $[r,\widehat\tau_n]$. On $F_n^c$, we define $\widehat{\vb}=\tilde{\vb}$ on $[0,\tau_n]$. Since $F_n\in\mathscr F_r$ and $\tilde{\vb}^{\,r}(r)=\tilde{\vb}(r)$ on $F_n$, by the adaptedness and pathwise uniqueness of the corresponding solutions, we infer that $(\widehat{\vb},\widehat\tau_n)$ is a local weak pathwise solution to \eqref{ee-o} satisfying
\[
\widehat{\vb}(t)=\tilde{\vb}(t),\qquad 0\leq t<\xi\wedge\widehat\tau_n,\quad\mathbb P\text{-a.s.}
\]
However, on $F_n\cap\{\xi<T\}$, we have
\[
\widehat\tau_n=\tau_n\vee T=T>\xi.
\]
Therefore,
\[
\mathbb P(\widehat\tau_n>\xi)\geq\mathbb P(F_n\cap\{\xi<T\})>0,
\]
which contradicts the maximality of $(\tilde{\vb},\xi)$ in \autoref{def-2.2}. Hence, we deduce that
\[
\xi=\infty,\qquad\mathbb P\text{-a.s.}
\]
For positive rational numbers $r<R$, choose the $\mathbb H^2$ time in \autoref{cor-ae-H2} before $r$ and apply the above global strong continuation and pathwise uniqueness to obtain $\tilde{\vb}\in C([r,R];\mathbb H^2(\dom))\cap L^2([r,R];\mathbb H^3(\dom))$, almost surely.
Taking a countable intersection of these full-probability events and then restricting the time intervals gives the stated regularity simultaneously for every $0<t_0<T<\infty$.
Finally, the uniqueness follows from the pathwise uniqueness established in \autoref{thm-local weak solution}. The proof is complete.
}

\end{proof}

\subsection{Backward uniqueness}\label{subsec-backward-uniqueness}

\begin{theorem}\label{thm-global tilde v weak}
Assume that the sequences of random variables $\{r_{\kb}^{j,1}\}_{(\kb,j)\in\mathcal K}$  and $\{r_{\kb}^{j,2}\}_{(\kb,j)\in\mathcal K}$ satisfy \autoref{A1}, and let
$\ub_0^1,\ub_0^2\in\mathbb L^2(\dom)$.
Let $\overrightarrow{\ub_0^1}$ and $\overrightarrow{\ub_0^2}$ be their randomizations defined by \eqref{eqn-randomization of u_0}. Assume that both random initial values are $\mathscr F_0$-measurable, but they are not necessarily independent of each other.

Under \autoref{A2}, let $\tilde{\vb}^1$ and $\tilde{\vb}^2$ be two weak pathwise solutions to \eqref{ee-o} in the sense of \autoref{def-2.2}, with initial data $\overrightarrow{\ub_0^1}$ and $\overrightarrow{\ub_0^2}$, respectively. Assume that the two solutions are driven by the same Brownian motion $B$, with the same coefficients $b$ and $f$.

Then for fixed $T>0$, there is an event $\Omega_T$ such that $\mathbb{P} (\Omega_T)=1$ and for every $\omega\in\Omega_T$
\begin{equation}\label{eq-backwarduniqueness}
\tilde{\vb}^{1}(T)=\tilde{\vb}^{2}(T)
\quad\Longrightarrow\quad
\tilde{\vb}^{1}(t)=\tilde{\vb}^{2}(t)\quad\text{for every }t\in[0,T],
\end{equation}
which is so called backward uniqueness of the weak solution $\tilde{\vb}$ to (\ref{ee-o}).
\end{theorem}
\begin{proof}
Let $\tilde{\vb}^{1}$ and $\tilde{\vb}^{2}$ be two weak solutions to (\ref{ee-o}). Then we have for $i=1,2$
\begin{equation}
	\partial _{t}\tilde{\vb}^{i}-\Delta \tilde{\vb}^{i}+\alpha (t)((\tilde{\vb}^{i}+\wb
	)\cdot \nabla) (\tilde{\vb}^{i}+\wb)=0,\ \mathrm{on}\ [0,\infty)\times
	\dom.
\end{equation}
Denote $\tilde{\tilde{\vb}}:=\tilde{\vb}^{1}-\tilde{\vb}^{2}$, then we get
\begin{equation}\label{ee-6.1.3}
	\partial _{t}\tilde{\tilde{\vb}}-\Delta \tilde{\tilde{\vb}}+\alpha (t)(\tilde{\tilde{\vb}}\cdot \nabla) (\tilde{\vb}^{1}+
	\wb)+\alpha (t)((\tilde{\vb}^{2}+\wb)\cdot \nabla) \tilde{\tilde{\vb}}=0\ \mathrm{on}\ [0,\infty)\times
	\dom.
\end{equation}

Let us choose and fix $\omega\in\Omega$ outside a null set and suppress $\omega$
from the notation. Let $T>0$ and $\eta\in(0,T)$.
By \autoref{cor-ae-H2}, for $i=1,2$, there exist $t_i^\ast\in(0,\eta)$ such that
\[
\tilde{\vb}^{i}(t_i^\ast)\in\mathbb H^2(\dom).
\]
In view of \autoref{thm-global weak transformed equation},
$\tilde{\vb}^{i}$ is a strong solution on $[t_i^\ast,T]$.
Set
\[
t^\ast:=\max\{t_1^\ast,t_2^\ast\}.
\]
Then
\[
t^\ast<\eta,
\qquad
\tilde{\vb}^{i}(t^\ast)\in\mathbb H^2(\dom),
\quad i=1,2,
\]
and hence
\[
\tilde{\tilde{\vb}}(t^\ast)\in\mathbb H^2(\dom).
\]
We now consider \eqref{ee-6.1.3} on $[t^\ast,T]$.

For $t\in [t^\ast, T]$, we define the linear operator $M(t)$ by
\begin{align*}
M(t)x&:=\alpha (t)(x\cdot \nabla) (\tilde{\vb}^{1}(t)+\wb(t))+\alpha
(t)((\tilde{\vb}^{2}(t)+\wb(t))\cdot \nabla) x,\;\;x\in \mathbb{H}^1(\dom).
\end{align*}
and set
\[
\phi_{\eps} (t):=\frac{\langle -\Delta\tilde{\tilde{\vb}}(t),\tilde{\tilde{\vb}}(t)\rangle}{\Vert \tilde{\tilde{\vb}}(t)\Vert_{L^2}^{2}+\eps},\;\;\eps>0.
\]

On $[t^\ast,T]$, $\tilde{\tilde{\vb}}\in C([t^\ast,T];\mathbb H^2(\dom))\cap L^2([t^\ast,T];\mathbb H^3(\dom))$ and $\wb\in C([0,T];\mathbb H^2(\dom))$. Hence $M\tilde{\tilde{\vb}}$ and $\partial_t\tilde{\tilde{\vb}}=\Delta\tilde{\tilde{\vb}}-M\tilde{\tilde{\vb}}$ belong to $L^2([t^\ast,T];\mathbb L^2(\dom))$.
In particular, the numerator of $\phi_\eps$ is absolutely continuous with derivative $2\langle-\Delta\tilde{\tilde{\vb}},\partial_t\tilde{\tilde{\vb}}\rangle$ almost everywhere. Indeed, the identity holds for the finite Fourier projections, and their Laplacians and time derivatives converge in $L^2([t^\ast,T];\mathbb L^2)$, so the products converge in $L^1$. The denominator is absolutely continuous and bounded below by $\eps$. Thus the following differential identities hold almost everywhere. From (\ref{ee-6.1.3}) we infer that for every $\eps>0$
\begin{align}\label{eq-dtlog}
\frac{d}{dt}\log (\Vert \tilde{\tilde{\vb}}(t)\Vert_{L^2}^{2}+\eps) &=2\frac{\langle\tilde{\tilde{\vb}}(t),\tilde{\tilde{\vb}}^{\prime}(t)\rangle}{\Vert \tilde{\tilde{\vb}}(t)\Vert_{L^2}^{2}+\eps} =2\frac{\langle\tilde{\tilde{\vb}}(t),\Delta \tilde{\tilde{\vb}}(t)-M(t)\tilde{\tilde{\vb}}(t)\rangle}{\Vert \tilde{\tilde{\vb}}(t)\Vert_{L^2}^{2}+\eps} =-2\phi_{\eps}(t)-2\frac{\langle\tilde{\tilde{\vb}}(t),M(t)\tilde{\tilde{\vb}}(t)\rangle}{\Vert \tilde{\tilde{\vb}}(t)\Vert_{L^2}^{2}+\eps}.
\end{align}

Firstly, we aim to estimate the first term  in \eqref{eq-dtlog}.
Define
\begin{align*}
F(t):=\partial_t\tilde{\tilde{\vb}}(t)-\Delta\tilde{\tilde{\vb}}(t),\;\;t\in[t^\ast,T].
\end{align*}

Therefore, we infer that
\begin{align*}
\phi ^{\prime}_{\eps}(t) &=\frac{2\langle-\Delta \tilde{\tilde{\vb}}(t),\tilde{\tilde{\vb}}^{\prime}(t)\rangle}{\Vert \tilde{\tilde{\vb}}(t)\Vert_{L^2}^{2}+\eps}-2\langle\tilde{\tilde{\vb}}(t),\tilde{\tilde{\vb}}^{\prime}(t)\rangle\frac{\langle-\Delta \tilde{\tilde{\vb}}(t),\tilde{\tilde{\vb}}(t)\rangle}{(\Vert \tilde{\tilde{\vb}}(t)\Vert_{L^2}^{2}+\eps)^{2}} \\
&=\frac{2\langle-\Delta \tilde{\tilde{\vb}}(t),\Delta \tilde{\tilde{\vb}}(t)+F(t)\rangle}{\Vert \tilde{\tilde{\vb}}(t)\Vert_{L^2}^{2}+\eps}
-2\frac{\langle-\Delta\tilde{\tilde{\vb}}(t),\tilde{\tilde{\vb}}(t)\rangle\langle\tilde{\tilde{\vb}}(t),\Delta \tilde{\tilde{\vb}}(t)
+F(t)\rangle}{(\Vert \tilde{\tilde{\vb}}(t)\Vert_{L^2}^{2}+\eps)^{2}} \\
&=-\frac{2\Vert \Delta \tilde{\tilde{\vb}}(t)\Vert_{L^2}^{2}}{\Vert \tilde{\tilde{\vb}}(t)\Vert_{L^2}^{2}+\eps}-\frac{2\langle\Delta \tilde{\tilde{\vb}}(t),F(t)\rangle}{\Vert \tilde{\tilde{\vb}}(t)\Vert_{L^2}^{2}+\eps}+\frac{2(\langle\Delta\tilde{\tilde{\vb}}(t),\tilde{\tilde{\vb}}(t)\rangle\langle\tilde{\tilde{\vb}}(t),F(t)\rangle+\langle\Delta \tilde{\tilde{\vb}}(t),\tilde{\tilde{\vb}}(t)\rangle^{2})}{(\Vert \tilde{\tilde{\vb}}(t)\Vert_{L^2}^{2}+\eps)^{2}} \\
&=-\frac{2\Vert \Delta \tilde{\tilde{\vb}}(t)\Vert_{L^2}^{2}}{\Vert \tilde{\tilde{\vb}}(t)\Vert_{L^2}^{2}+\eps}-\frac{2\langle\Delta \tilde{\tilde{\vb}}(t),F(t)\rangle}{\Vert \tilde{\tilde{\vb}}(t)\Vert_{L^2}^{2}+\eps}+\frac{2(\langle\Delta\tilde{\tilde{\vb}}(t),\tilde{\tilde{\vb}}(t)\rangle\langle\tilde{\tilde{\vb}}(t),F(t)\rangle+\langle\Delta \tilde{\tilde{\vb}}(t),\tilde{\tilde{\vb}}(t)\rangle^{2}+\frac{1}{4}\langle\tilde{\tilde{\vb}}(t),F(t)\rangle^{2})}{(\Vert \tilde{\tilde{\vb}}(t)\Vert_{L^2}^{2}+\eps)^{2}}\\
&\;\;\;\;-\frac{1}{2}\frac{\langle\tilde{\tilde{\vb}}(t),F(t)\rangle^{2}}{(\Vert \tilde{\tilde{\vb}}(t)\Vert_{L^2}^{2}+\eps)^{2}} \\
&=-\frac{2\Vert \Delta \tilde{\tilde{\vb}}(t)\Vert_{L^2}^{2}}{\Vert \tilde{\tilde{\vb}}(t)\Vert_{L^2}^{2}+\eps}-\frac{2\langle\Delta \tilde{\tilde{\vb}}(t),F(t)\rangle}{\Vert \tilde{\tilde{\vb}}(t)\Vert_{L^2}^{2}+\eps}+\frac{2\langle\Delta \tilde{\tilde{\vb}}(t)+\frac{1}{2}F(t),\tilde{\tilde{\vb}}(t)\rangle^{2}}{(\Vert \tilde{\tilde{\vb}}(t)\Vert_{L^2}^{2}+\eps)^{2}}-\frac{1}{2}\frac{\langle\tilde{\tilde{\vb}}(t),F(t)\rangle^{2}}{(\Vert \tilde{\tilde{\vb}}(t)\Vert_{L^2}^{2}+\eps)^{2}}.
\end{align*}
By the H\"{o}lder inequality, we deduce
\begin{align}\label{eqn-dt-phieps}
\phi^{\prime}_{\eps}(t) &\leq -\frac{2\Vert \Delta \tilde{\tilde{\vb}}(t)\Vert_{L^2}^{2}}{\Vert \tilde{\tilde{\vb}}(t)\Vert_{L^2}^{2}+\eps}-\frac{2\langle\Delta \tilde{\tilde{\vb}}(t),F(t)\rangle}{\Vert \tilde{\tilde{\vb}}(t)\Vert_{L^2}^{2}+\eps}+\frac{2\left\Vert \Delta \tilde{\tilde{\vb}}(t)+\frac{1}{2}F(t)\right\Vert _{L^2}^{2}}{\Vert \tilde{\tilde{\vb}}(t)\Vert_{L^2}^{2}+\eps}-\frac{1}{2}\frac{\langle\tilde{\tilde{\vb}}(t),F(t)\rangle^{2}}{(\Vert \tilde{\tilde{\vb}}(t)\Vert_{L^2}^{2}+\eps)^{2}}\nonumber \\
&=-\frac{2\Vert \Delta \tilde{\tilde{\vb}}(t)\Vert_{L^2}^{2}}{\Vert \tilde{\tilde{\vb}}(t)\Vert_{L^2}^{2}+\eps}-\frac{2\langle\Delta \tilde{\tilde{\vb}}(t),F(t)\rangle}{\Vert \tilde{\tilde{\vb}}(t)\Vert_{L^2}^{2}+\eps}+2\frac{\left\Vert \Delta\tilde{\tilde{\vb}}(t)\right\Vert_{L^2}^{2}+\langle\Delta \tilde{\tilde{\vb}}(t),F(t)\rangle+\frac{1}{4}\Vert F(t)\Vert_{L^2}^{2}}{\Vert \tilde{\tilde{\vb}}(t)\Vert_{L^2}^{2}+\eps}-\frac{1}{2}\frac{\langle\tilde{\tilde{\vb}}(t),F(t)\rangle^{2}}{(\Vert \tilde{\tilde{\vb}}(t)\Vert_{L^2}^{2}+\eps)^{2}} \nonumber\\
&=\frac{1}{2}\frac{\left\Vert F(t)\right\Vert_{L^2}^{2}}{\Vert \tilde{\tilde{\vb}}(t)\Vert_{L^2}^{2}+\eps}-\frac{1}{2}\frac{\langle\tilde{\tilde{\vb}}(t),F(t)\rangle^{2}}{(\Vert \tilde{\tilde{\vb}}(t)\Vert_{L^2}^{2}+\eps)^{2}} \leq \frac{\left\Vert F(t)\right\Vert_{L^2}^{2}}{2\Vert \tilde{\tilde{\vb}}(t)\Vert_{L^2}^{2}+2\eps}.
\end{align}

By virtue of the H\"{o}lder
inequality and the Sobolev embedding theorem, there exists $C=C(\sup_{[t^\ast,T]}\vert \alpha\vert)\geq 0$ such that
\begin{equation}\label{eqn-Mv-L2}
\begin{aligned}
	\Vert M(t)\tilde{\tilde{\vb}}(t)\Vert_{L^2} &\leq \Vert \alpha (t)(\tilde{\tilde{\vb}}(t)\cdot \nabla)
	(\tilde{\vb}^{1}(t)+\wb(t))\Vert _{L^2}+\Vert \alpha (t)((\tilde{\vb}^{2}(t)+\wb(t)
	)\cdot \nabla) \tilde{\tilde{\vb}}(t)\Vert _{L^2} \\
	&\leq C\left( \Vert (\tilde{\tilde{\vb}}(t)\cdot \nabla) (\tilde{\vb}^{1}(t)+\wb(t))\Vert
	_{L^2}+\Vert ((\tilde{\vb}^{2}(t)+\wb(t))\cdot \nabla) \tilde{\tilde{\vb}}(t)\Vert_{L^2}\right)  \\
	&\leq C\left( \Vert  \tilde{\tilde{\vb}}(t)\Vert_{L^3}\Vert\nabla (\tilde{\vb}^{1}(t)+\wb(t))\Vert_{L^6}+\Vert \tilde{\vb}^{2}(t)+\wb(t)\Vert_{L^{\infty}}\Vert\tilde{\tilde{\vb}}(t)\Vert_{1}\right)  \\
	&\leq C\left( \Vert\tilde{\tilde{\vb}}(t)\Vert_{H^{1}}\Vert\tilde{\vb}^{1}(t)+\wb(t)\Vert_{H^2}+\Vert
	\tilde{\vb}^{2}(t)+\wb(t)\Vert_{H^2}\Vert\tilde{\tilde{\vb}}(t)\Vert_{1}\right)  \\
	&\leq C\left( \Vert\tilde{\vb}^{1}(t)+\wb(t)\Vert_{H^2}+\Vert\tilde{\vb}^{2}(t)+\wb(t)\Vert_{H^{2}}\right) \Vert\tilde{\tilde{\vb}}(t)\Vert_{H^{1}} \\
	&\leq C\left( \Vert\tilde{\vb}^{1}(t)\Vert_{H^2}+\Vert\tilde{\vb}^{2}(t)\Vert_{H^{2}}+\Vert \wb(t)\Vert_{H^{2}}\right) \Vert\tilde{\tilde{\vb}}(t)\Vert_{H^{1}}.
\end{aligned}
\end{equation}

For clarity we introduce the following auxiliary notation
\[
H(t) \coloneqq \Vert\tilde{\vb}^{1}(t)\Vert_{H^2}^{2}+\Vert\tilde{\vb}^{2}(t)\Vert_{H^{2}}^{2}+\Vert \wb(t)\Vert_{H^{2}}^{2}.
\]

We observe that $F(t)=-M(t)\tilde{\tilde{\vb}}(t)$, by \eqref{eqn-dt-phieps} and \eqref{eqn-Mv-L2} we infer
\begin{align*}
\phi^\prime_{\eps}(t) &\leq
\frac{\left\Vert M(t)\tilde{\tilde{\vb}}(t)\right\Vert_{L^2}^{2}}{2\Vert \tilde{\tilde{\vb}}(t)\Vert_{L^2}^{2}+2\eps} \leq \frac{CH(t) \Vert\tilde{\tilde{\vb}}(t)\Vert_{H^{1}}^{2}}{\Vert \tilde{\tilde{\vb}}(t)\Vert_{L^2}^{2}+\eps} =\frac{CH(t)(\Vert\tilde{\tilde{\vb}}(t)\Vert_{L^{2}}^{2}+\Vert\tilde{\tilde{\vb}}(t)\Vert_{1}^{2})}{\Vert \tilde{\tilde{\vb}}(t)\Vert_{L^2}^{2}+\eps}\leq CH(t) (\phi_{\eps}(t)+1).
\end{align*}

By the Gronwall inequality, we infer that there exists a constant $C>0$ such that
\begin{equation}\label{eqn-phieps}
\begin{aligned}
\phi_{\eps} (t) &\leq \Big(\phi_{\eps} (t^\ast)+C\int_{t^\ast}^{t}H(s)\,ds\Big)e^{C\int_{t^\ast}^{t}H(s)\,ds} .
\end{aligned}
\end{equation}
Assume that $\Vert \tilde{\tilde{\vb}}(t^\ast)\Vert_{L^2}\neq 0$. Define
\[
\phi(t^\ast) :=\frac{\langle -\Delta\tilde{\tilde{\vb}}(t^\ast),\tilde{\tilde{\vb}}(t^\ast)\rangle}{\Vert \tilde{\tilde{\vb}}(t^\ast)\Vert_{L^2}^{2}}.
\]
Since $\phi_{\eps} (t^\ast)\leq \phi (t^\ast)$, there exists a constant $C_0>0$ such that
\begin{equation}\label{eqn-phieps-2}
\begin{aligned}
\phi_{\eps} (t)  &\leq \Big{(}\phi (t^\ast)+C\int_{t^\ast}^{T}H(s)\,ds\Big{)}e^{C\int_{t^\ast}^{T}H(s)\,ds} =: C_0.
\end{aligned}
\end{equation}

For the second term in \eqref{eq-dtlog}, by the H\"older inequality and inequality (\ref{eqn-Mv-L2}) we infer
\begin{equation}\label{eqn-474}
\begin{aligned}
\frac{\vert\langle\tilde{\tilde{\vb}}(t),M(t)\tilde{\tilde{\vb}}(t)\rangle\vert}{\Vert \tilde{\tilde{\vb}}(t)\Vert_{L^2}^{2}+\eps} &\leq \frac{\Vert \tilde{\tilde{\vb}}(t)\Vert_{L^2}\left\Vert M(t)\tilde{\tilde{\vb}}(t)\right\Vert _{L^2}}{\Vert \tilde{\tilde{\vb}}(t)\Vert_{L^2}^{2}+\eps} \\
&\leq C\frac{\Vert \tilde{\tilde{\vb}}(t)\Vert_{H^{1}}^{2}\left(\Vert\tilde{\vb}^{1}(t)\Vert_{H^2}+\Vert\tilde{\vb}^{2}(t)\Vert_{H^{2}}+\Vert \wb(t)\Vert_{H^{2}}\right) }{\Vert \tilde{\tilde{\vb}}(t)\Vert_{L^2}^{2}+\eps} \\
&=C\frac{(\Vert \tilde{\tilde{\vb}}(t)\Vert_{L^2}^{2}+\Vert \tilde{\tilde{\vb}}(t)\Vert_{1}^{2})\left(\Vert\tilde{\vb}^{1}(t)\Vert_{H^2}+\Vert\tilde{\vb}^{2}(t)\Vert_{H^{2}}+\Vert\wb(t)\Vert_{H^{2}}\right) }{\Vert \tilde{\tilde{\vb}}(t)\Vert_{L^2}^{2}+\eps} \\
&\leq C\left(\Vert\tilde{\vb}^{1}(t)\Vert_{H^2}+\Vert\tilde{\vb}^{2}(t)\Vert_{H^{2}}+\Vert \wb(t)\Vert_{H^{2}}\right) +C\phi_{\eps} (t)\left(\Vert\tilde{\vb}^{1}(t)\Vert_{H^2}+\Vert\tilde{\vb}^{2}(t)\Vert_{H^{2}}+\Vert \wb(t)\Vert_{H^{2}}\right)  \\
&\leq C(C_0+1)\left( \Vert\tilde{\vb}^{1}(t)\Vert_{H^2}+\Vert\tilde{\vb}^{2}(t)\Vert_{H^{2}}+\Vert \wb(t)\Vert_{H^{2}}\right) .
\end{aligned}
\end{equation}

Hence, in view of \eqref{eqn-phieps} and \eqref{eqn-474}, equality \eqref{eq-dtlog} implies
\begin{align*}
\frac{d}{dt}\log (\Vert \tilde{\tilde{\vb}}(t)\Vert_{L^2}^{2}+\eps)
&\geq -C[\Vert\tilde{\vb}^{1}(t)\Vert_{H^2}+\Vert\tilde{\vb}^{2}(t)\Vert_{H^{2}}+\Vert \wb(t)\Vert_{H^{2}}+1].
\end{align*}
Integrating in time yields
\begin{align*}
\log (\left\Vert \tilde{\tilde{\vb}}\left( T\right) \right\Vert _{L^2}^{2}+\eps)
\geq -C\left( T-t^\ast
\right) -C \int_{t^\ast}^{T}(\Vert\tilde{\vb}^{1}(t)\Vert_{H^2}+\Vert\tilde{\vb}^{2}(t)\Vert_{H^{2}}+\Vert \wb(t)\Vert_{H^{2}})\,dt+\log (\left\Vert \tilde{\tilde{\vb}}(t^\ast) \right\Vert _{L^2}^{2}+\eps).
\end{align*}

Set $D_\eps(t)=\Vert\tilde{\tilde{\vb}}(t)\Vert_{L^2}^2+\eps$ and $K=C(T -t^\ast)+C\int_{t^\ast}^T(\Vert\tilde\vb^1\Vert_{H^2}+\Vert\tilde\vb^2\Vert_{H^2}+\Vert\wb\Vert_{H^2})\,dt<\infty$. Then
\[
D_\eps(T)\geq e^{-K}D_\eps(t^\ast).
\]
If $\tilde{\tilde{\vb}}(T)=0$, letting $\eps\downarrow0$ contradicts the assumption that $\tilde{\tilde{\vb}}(t^\ast)\neq0$, Hence
\[
\left\Vert \tilde{\tilde{\vb}}(t^\ast)\right\Vert_{L^2}=0.
\]
By the pathwise uniqueness established in \autoref{thm-global weak transformed equation},
\[
\left\Vert \tilde{\tilde{\vb}}(t)\right\Vert_{L^2}=0,\;\;t\in[t^\ast,T].
\]

Let us choose and fix $t\in(0,T]$ and $\eta\in(0,t)$. Since $t^\ast(\eta)<\eta<t\leq T$, we infer that
\[
\left\Vert \tilde{\tilde{\vb}}(t)\right\Vert_{L^2}=0.
\]

Hence, $\tilde{\tilde{\vb}}(t)=0$ on $(0,T]$. Since
$\tilde{\vb}^{1},\tilde{\vb}^{2}\in
C([0,T];\mathbb{L}^{2}(\dom))$,
we have
$\tilde{\tilde{\vb}}\in
C([0,T];\mathbb{L}^{2}(\dom))$.
Therefore,
\[
\Vert\tilde{\tilde{\vb}}(0)\Vert_{L^2}
=
\lim_{t\downarrow0}
\Vert\tilde{\tilde{\vb}}(t)\Vert_{L^2}
=0.
\]
Thus the backward uniqueness holds on $[0,T]$.
The proof of \autoref{thm-global tilde v weak} is complete.
\end{proof}

\subsection{Proof of the main result}\label{subsec-main-proof}

In the proof below we will use the following shortcut notation.
\begin{equation}\label{eqn-Z_T}
 \mathcal Z_T \coloneqq   C([0,T];\mathbb{L}^{2}(\dom))
\cap L^{2}([0,T];\mathbb{H}^{1,3}(\dom))
\cap L^{4}([0,T];\mathbb{L}^{6}(\dom)).
\end{equation}

\begin{proof}[Proof of Theorem \ref{thm-main result}]
By the scaling described in the Introduction, it is sufficient
to prove the result for $\nu=1$.
Let us consider the process $\wb$ defined by formula \eqref{def-w}  in Proposition \ref{prop-ouz2}, and $\tilde \vb$ is the global weak solution to \eqref{ee-o} in the sense of \autoref{def-2.2}.
Fix $T>0$.
In view of \autoref{cor-w-continuous-L2} and \autoref{thm-global weak transformed equation}, for almost all $\omega\in\Omega$,
\begin{align*}
\tilde{\vb}+\wb\in \mathcal Z_T.
\end{align*}
Since the function $[0,\infty) \ni t \mapsto \alpha(t)=e^{bB(t)} \in \mathbb{R}$ is continuous,  almost surely, for almost all $\omega\in\Omega$,
\[
\ub(t,\omega)=\alpha(t,\omega)\big(\tilde{\vb}(t,\omega)+\wb(t,\omega)\big)
\in \mathcal Z_T .
\]
Taking a countable intersection of the full-probability events corresponding to $T\in\mathbb N$, we infer that $\ub\in Z_\infty$ on one event of probability one.

Let us choose and fix $\eta\in \mathbb{H}^1(\dom)$. Since $\wb$ satisfies the following equation
\begin{equation}\label{eq-w-linear}
\partial_t \wb(t) = \Delta \wb(t) + \alpha^{-1}(t)\fb(t),
\quad \; \text{a.e.}t\geq 0.
\end{equation}
we infer that for a.e.  $t\geq 0$,
\begin{equation}\label{eq-w-weak}
\big\langle \partial_t\wb(t),\eta\big\rangle
+\big\langle \Lambda \wb(t),\Lambda\eta\big\rangle
=\big\langle \alpha^{-1}(t)\fb(t),\eta\big\rangle, \;\; \mathbb{P}\mbox{-a.s.}
\end{equation}
Since $\tilde{\vb}$ is the global weak pathwise solution to \eqref{ee-o}, $\tilde{\vb}$ satisfies for a.e. $t\in [0, T]$,
\begin{equation}\label{eq-vtilde-weak}
\big\langle \partial_t\tilde{\vb}(t),\eta\big\rangle
+\big\langle \Lambda \tilde{\vb}(t),\Lambda\eta\big\rangle
+\alpha(t)\big\langle (\tilde{\vb}+\wb)(t)\cdot\nabla(\tilde{\vb}+\wb)(t),\eta\big\rangle=0.
\end{equation}
Adding \eqref{eq-w-weak} and \eqref{eq-vtilde-weak} yields that for a.e. $t\in [0, T]$,
\begin{equation}\label{eq-v-weak}
\big\langle \partial_t(\tilde{\vb}+\wb)(t),\eta\big\rangle
+\big\langle \Lambda (\tilde{\vb}+\wb)(t),\Lambda\eta\big\rangle
+\alpha(t)\big\langle (\tilde{\vb}+\wb)(t)\cdot\nabla(\tilde{\vb}+\wb)(t),\eta\big\rangle
=\big\langle \alpha^{-1}(t)\fb(t),\eta\big\rangle .
\end{equation}
Since $\ub=\alpha(\tilde{\vb}+\wb)$, using the Stratonovich product rule, i.e. $d\alpha=b\alpha\circ dB$, we infer that
\begin{equation}\label{eq-strat-product}
\begin{aligned}
d\langle \ub,\eta\rangle&=d\big(\alpha\langle (\tilde{\vb}+\wb),\eta\rangle\big)
=\alpha\,d\langle (\tilde{\vb}+\wb),\eta\rangle+\langle (\tilde{\vb}+\wb),\eta\rangle \circ d\alpha\\
&=\alpha\,d\langle (\tilde{\vb}+\wb),\eta\rangle + b\alpha\langle (\tilde{\vb}+\wb),\eta\rangle\circ dB
=\alpha\langle \partial_t(\tilde{\vb}+\wb),\eta\rangle\,dt + b\langle \ub,\eta\rangle\circ dB.
\end{aligned}
\end{equation}
Multiplying \eqref{eq-v-weak} by $\alpha(t)$, we infer that for a.e. $t\in [0, T]$,
\begin{align}\label{eqn-489}
\alpha(t)\langle \partial_t(\tilde{\vb}+\wb)(t),\eta\rangle
= -\alpha(t)\langle \Lambda(\tilde{\vb}+\wb)(t),\Lambda\eta\rangle
-\alpha(t)^2\langle (\tilde{\vb}+\wb)(t)\cdot\nabla(\tilde{\vb}+\wb)(t),\eta\rangle
+\langle \fb(t),\eta\rangle.
\end{align}
Since $\ub=\alpha(\tilde{\vb}+\wb)$ and $\Lambda(\alpha(\tilde{\vb}+\wb))=\alpha\Lambda(\tilde{\vb}+\wb)$, we deduce
\begin{align}\label{eqn-490}
\alpha\langle \Lambda(\tilde{\vb}+\wb),\Lambda\eta\rangle=\langle \Lambda\ub,\Lambda\eta\rangle,
\qquad
\alpha^2\langle (\tilde{\vb}+\wb)\cdot\nabla(\tilde{\vb}+\wb),\eta\rangle
=\langle (\ub\cdot\nabla)\ub,\eta\rangle.
\end{align}
Substituting \eqref{eqn-489} and \eqref{eqn-490} into \eqref{eq-strat-product} yields, for a.e. $t\in [0, T]$
\[
d\langle \ub(t),\eta\rangle
=
\Big(
-\langle \Lambda\ub(t),\Lambda\eta\rangle
-\langle (\ub\cdot\nabla)\ub(t),\eta\rangle
+\langle \fb(t),\eta\rangle
\Big)\,dt
+ b\langle \ub(t),\eta\rangle\circ dB(t).
\]
Note that $\ub(0)=\overrightarrow{\ub_0}$, the above equality means rigorously that for every $t \in [0,T]$ and $\eta\in \mathbb{H}^1(\dom)$,

\begin{equation}\label{eq-Stratonovich}
\begin{aligned}
        \langle\ub(t),\eta\rangle
&=\langle\overrightarrow{\ub_0},\eta\rangle
-\int_0^t\langle \Lambda\ub(s),\Lambda\eta\rangle\,ds
-\int_0^t\langle (\ub(s)\cdot\nabla)\ub(s),\eta\rangle\,ds
\\
&+\int_0^t\langle \fb(s),\eta\rangle\,ds
+ b\int_0^t\langle \ub(s),\eta\rangle\circ dB(s), \;\; \mathbb{P}\mbox{-a.s.}.
\end{aligned}
\end{equation}

In view of \eqref{eq-v-weak}, we infer that for every $\eta\in\mathbb H^1(\dom)$, the function $t\mapsto\langle\tilde{\vb}(t)+\wb(t),\eta\rangle$ is absolutely continuous on $[0,T]$. Since $\alpha$ is a continuous semimartingale, by the integration by parts formula, we infer that $\langle\ub(t),\eta\rangle=\alpha(t)\langle\tilde{\vb}(t)+\wb(t),\eta\rangle$ is a continuous semimartingale. Therefore, by \eqref{eq-Stratonovich}, we deduce that $\ub$ is a global weak solution to \eqref{eqn-Burgers-intro} in the sense of \autoref{def-solution}.

It remains to prove uniqueness. Let $\ub^1$ and $\ub^2$ be two weak solutions to \eqref{eqn-Burgers-intro} with the same initial data. For $i=1,2$, set
\[
\tilde{\vb}^{\,i}:=\alpha^{-1}\ub^i-\wb.
\]

Let us choose and fix $\eta\in\mathbb H^1(\dom)$. By \autoref{def-solution}, the processes $\langle\ub^i,\eta\rangle$, $i=1,2$, are continuous semimartingales. Since $d\alpha^{-1}=-b\alpha^{-1}\circ dB$, by the Stratonovich product rule, we infer that
\begin{align*}
d\langle\alpha^{-1}\ub^i,\eta\rangle
&=\alpha^{-1}\circ d\langle\ub^i,\eta\rangle+\langle\ub^i,\eta\rangle\circ d\alpha^{-1}=\alpha^{-1}\circ d\langle\ub^i,\eta\rangle-b\alpha^{-1}\langle\ub^i,\eta\rangle\circ dB.
\end{align*}
Substituting the weak formulation \eqref{eq-global-weak-solution} into the above equality, we infer that
\begin{align*}
d\langle\alpha^{-1}\ub^i,\eta\rangle
={}&-\alpha^{-1}\langle\Lambda\ub^i,\Lambda\eta\rangle\,dt-\alpha^{-1}\langle(\ub^i\cdot\nabla)\ub^i,\eta\rangle\,dt+\alpha^{-1}\langle\fb,\eta\rangle\,dt.
\end{align*}
Integrate this identity and subtract the time integral of \eqref{eq-w-weak}. For $i=1,2$, take a common event of probability one on which the resulting scalar identities hold at all nonnegative rational times for a countable dense subset of $\mathbb H^1(\dom)$, and on which $\ub^i,\wb\in\mathcal Z_T$ for every positive integer $T$. On this event $\tilde{\vb}^{\,i}\in\mathcal Z_T$ because $\alpha^{-1}$ is bounded on $[0,T]$. Write $\mathbf y^i=\tilde{\vb}^{\,i}+\wb$. By spatial and temporal H\"older inequalities,
\[
\|(\mathbf y^i\cdot\nabla)\mathbf y^i\|_{L^{4/3}(0,T;\mathbb L^2)}
\leq \|\mathbf y^i\|_{L^4(0,T;\mathbb L^6)}
       \|\nabla\mathbf y^i\|_{L^2(0,T;\mathbb L^3)}.
\]
Also $\Delta\tilde{\vb}^{\,i}\in L^2(0,T;\mathbb H^{-1})$. The drift is therefore Bochner integrable in $\mathbb H^{-1}(\dom)$ on each finite interval. Continuity in time of its indefinite integral and of $\tilde{\vb}^{\,i}$, followed by density in the test function, gives for all $t\geq0$
\[
\tilde{\vb}^{\,i}(t)=\overrightarrow{\ub_0}+
\int_0^t\!\big[\Delta\tilde{\vb}^{\,i}(s)
-\alpha(s)(\mathbf y^i(s)\cdot\nabla)\mathbf y^i(s)\big]\,ds
\qquad\text{in }\mathbb H^{-1}(\dom).
\]
Consequently $\partial_t\tilde{\vb}^{\,i}\in L^1_{\loc}([0,\infty);\mathbb H^{-1})$ and \eqref{ee-o} holds almost everywhere with the prescribed initial value. The process $\wb$ is adapted by \eqref{def-w}, whose integral up to $t$ uses only the coefficient path on $[0,t]$. Thus $\alpha^{-1}\ub^i-\wb$ is continuous and adapted in $\mathbb L^2$, hence progressively measurable. Both processes satisfy \autoref{def-2.2}.
By the uniqueness in \autoref{thm-global weak transformed equation}, we deduce
\[
\tilde{\vb}^{1}=\tilde{\vb}^{2},\;\; \mathbb{P}\text{-a.s.},
\]
which implies
\[
\ub^1=\alpha(\tilde{\vb}^{1}+\wb)=\alpha(\tilde{\vb}^{2}+\wb)=\ub^2,\;\; \mathbb{P}\text{-a.s.}
\]
Hence, the weak solution to \eqref{eqn-Burgers-intro} is unique.

{
It remains to prove continuous dependence. Let $\ub_{0,n}\to\ub_0$ in $\mathbb L^2(\dom)$, and let $\ub_n$ and $\ub$ be the corresponding global weak solutions constructed using the same random multipliers, $B$ and $\fb$. Denote by $\tilde{\vb}_n$, $\tilde{\vb}$ and $\vb_n$, $\vb$ the corresponding transformed solutions and mild solutions, respectively. Recall that $\mathbf z_n=e^{-\cdot A}\overrightarrow{\ub_{0,n}}+\wb$ and $\mathbf z=e^{-\cdot A}\overrightarrow{\ub_0}+\wb$.

Let us choose and fix $T>0$. By \eqref{eq-randomization-L2}, \eqref{eq-623-11} and \eqref{eq-623-1}, we infer that there exists $c>0$ such that
\begin{align*}
&\mathbb E\Vert\overrightarrow{\ub_{0,n}}-\overrightarrow{\ub_0}\Vert_{L^2}^2
+\mathbb E\Vert\mathbf z_n-\mathbf z\Vert_{Y_{T\wedge1}}^2
\leq c\Vert\ub_{0,n}-\ub_0\Vert_{L^2}^2\to0,\qquad n\to\infty.
\end{align*}
Hence, from any subsequence we may extract a further subsequence, still indexed by $n$, such that
\[
\sum_{n=1}^{\infty}\left(\Vert\overrightarrow{\ub_{0,n}}-\overrightarrow{\ub_0}\Vert_{L^2}^2+\Vert\mathbf z_n-\mathbf z\Vert_{Y_{T\wedge1}}^2\right)<\infty,\qquad\mathbb P\text{-a.s.}
\]
Let us choose and fix $\omega\in\Omega$ outside the corresponding null sets. In the following, all constants may depend on $\omega$, but are independent of $n$.

We first obtain an estimate on a common small-time interval.
For $s\in(0,T\wedge1)$, $g\in X_s$ and $h\in Y_s$, set
\[
\Phi_h(g)(t):=-\int_0^t\alpha(r)e^{(t-r)\Delta}\big((g(r)+h(r))\cdot\nabla\big)(g(r)+h(r))\,dr,
\qquad 0\leq t\leq s.
\]
By \autoref{lemma-lplq}, \eqref{eqn-3.52}, \eqref{eqn-3.54} and \eqref{eqn-g nablah}, we infer that
\begin{align*}
\Vert\Phi_h(g)\Vert_{X_s}
&\leq c\big(\Vert g\Vert_{X_s}+\Vert h\Vert_{Y_s}\big)^2,\\
\Vert\Phi_{h_1}(g_1)-\Phi_{h_2}(g_2)\Vert_{X_s}
&\leq c\big(\Vert g_1\Vert_{X_s}+\Vert g_2\Vert_{X_s}+\Vert h_1\Vert_{Y_s}+\Vert h_2\Vert_{Y_s}\big)\\
&\quad\times\big(\Vert g_1-g_2\Vert_{X_s}+\Vert h_1-h_2\Vert_{Y_s}\big).
\end{align*}
Here $c$ depends on $\Vert\alpha\Vert_{L^\infty([0,T\wedge1])}$ but not on $s$; for the maximal regularity estimate, extend the forcing by zero to $[0,1]$.
Choose $R>0$ such that $8cR\leq1$.
Whenever $\Vert h\Vert_{Y_s}\leq R$, the map $\Phi_h$ maps
\[
\mathcal B_s:=\{g\in X_s:g(0)=0,\ \Vert g\Vert_{X_s}\leq R\}
\]
into itself and has contraction constant at most $\frac12$.
Since $\Vert\mathbf z\Vert_{Y_t}\to0$ as $t\to0+$, choose $s\in(0,T\wedge1)$ such that $\Vert\mathbf z\Vert_{Y_s}\leq R/2$.
By $\mathbf z_n\to\mathbf z$ in $Y_{T\wedge1}$, we also have $\Vert\mathbf z_n\Vert_{Y_s}\leq R$ for all sufficiently large $n$.
Thus, $\Phi_{\mathbf z}$ and $\Phi_{\mathbf z_n}$ have fixed points in the same ball $\mathcal B_s$.
Adding the corresponding randomized heat evolutions and applying the pathwise uniqueness argument in the proof of \autoref{thm-local weak solution}, we identify these fixed points with the remainders $\vb$ and $\vb_n$ of the global solutions on $[0,s]$.
Consequently,
\begin{align*}
\Vert\vb_n-\vb\Vert_{X_s}
&\leq\Vert\Phi_{\mathbf z_n}(\vb_n)-\Phi_{\mathbf z_n}(\vb)\Vert_{X_s}
+\Vert\Phi_{\mathbf z_n}(\vb)-\Phi_{\mathbf z}(\vb)\Vert_{X_s}\\
&\leq\frac12\Vert\vb_n-\vb\Vert_{X_s}+4cR\Vert\mathbf z_n-\mathbf z\Vert_{Y_s}.
\end{align*}
Absorbing the first term on the right-hand side, we infer that
\[
\Vert\vb_n-\vb\Vert_{X_s}\leq8cR\Vert\mathbf z_n-\mathbf z\Vert_{Y_s}\to0,\qquad n\to\infty.
\]
The time $s$ is chosen on the fixed sample path and is used only for this comparison; it need not be a stopping time.
Since $\tilde{\vb}_n=\vb_n+e^{-\cdot A}\overrightarrow{\ub_{0,n}}$ and $\tilde{\vb}=\vb+e^{-\cdot A}\overrightarrow{\ub_0}$, by the heat semigroup energy estimate, we deduce that
\begin{align*}
&\Vert\tilde{\vb}_n-\tilde{\vb}\Vert_{C([0,s];\mathbb L^2(\dom))}
+\Vert\tilde{\vb}_n-\tilde{\vb}\Vert_{L^2([0,s];\mathbb H^1(\dom))}\leq c\left(\Vert\overrightarrow{\ub_{0,n}}-\overrightarrow{\ub_0}\Vert_{L^2}+\Vert\mathbf z_n-\mathbf z\Vert_{Y_s}\right)\to0.
\end{align*}
Since the right-hand side is square-summable with respect to $n$, by Tonelli's theorem, there exists $t_0\in(0,s)$ such that
\[
\tilde{\vb}_n(t_0)\to\tilde{\vb}(t_0)\quad\mbox{in }\mathbb H^1(\dom).
\]

By \autoref{cor-ae-H2} and the proof of \autoref{thm-global weak transformed equation}, we infer that $\tilde{\vb}\in C([t_0,T];\mathbb H^1(\dom))$. Set
\[
M:=\sup_{t\in[t_0,T]}\Vert\tilde{\vb}(t)\Vert_{H^1}.
\]
By \autoref{prop-tildev-local strong} and \eqref{eq-uniform-local-Lipschitz}, let $h\in(0,1]$ and $L\geq1$ be a common existence time and Lipschitz constant for initial data in the $\mathbb H^1(\dom)$-ball of radius $M+1$, uniformly over starting times in $[t_0,T]$. Define
\[
N:=\left\lceil\frac{T-t_0}{h}\right\rceil,
\qquad t_j:=\min\{t_0+jh,T\},\quad j=0,\ldots,N.
\]
Let $n$ be sufficiently large that $\Vert\tilde{\vb}_n(t_0)-\tilde{\vb}(t_0)\Vert_{H^1}\leq L^{-N}$. Applying the uniform local estimate successively on $[t_j,t_{j+1}]$ and using pathwise uniqueness, we infer by induction that
\[
\Vert\tilde{\vb}_n-\tilde{\vb}\Vert_{C([t_j,t_{j+1}];\mathbb H^1(\dom))}
\leq L^{j+1}\Vert\tilde{\vb}_n(t_0)-\tilde{\vb}(t_0)\Vert_{H^1}\leq1,
\quad j=0,\ldots,N-1.
\]
Indeed, the bound on each preceding interval gives $\Vert\tilde{\vb}_n(t_j)\Vert_{H^1}\leq M+1$, so the same constants $h$ and $L$ apply on the next interval. Consequently,
\[
\Vert\tilde{\vb}_n-\tilde{\vb}\Vert_{C([t_0,T];\mathbb H^1(\dom))}
\leq L^N\Vert\tilde{\vb}_n(t_0)-\tilde{\vb}(t_0)\Vert_{H^1}\to0,
\qquad n\to\infty,
\]
where $L^N$ is independent of $n$.

Combining this with the convergence on $[0,s]$, and recalling that $\ub_n=\alpha(\tilde{\vb}_n+\wb)$ and $\ub=\alpha(\tilde{\vb}+\wb)$, we infer that
\[
\Vert\ub_n-\ub\Vert_{C([0,T];\mathbb L^2(\dom))}
\leq\Vert\alpha\Vert_{L^\infty([0,T])}\Vert\tilde{\vb}_n-\tilde{\vb}\Vert_{C([0,T];\mathbb L^2(\dom))}\to0.
\]
Since every subsequence admits a further subsequence converging almost surely, we conclude that $\ub_n\to\ub$ in probability in $C([0,T];\mathbb L^2(\dom))$.
}

It remains to prove the backward uniqueness of (\ref{eqn-Burgers-intro}). Fix $T>0$ and a path outside the null sets in \autoref{thm-global tilde v weak}.
Let $\ub^{1}$ and $\ub^{2}$ be the two weak solutions specified in the statement of the theorem, corresponding to the random initial values $\overrightarrow{\ub_{0}^1}$ and $\overrightarrow{\ub_{0}^2}$, respectively.
Recall that $\ub=\alpha\tilde{\vb}+\alpha\wb$. If $ \ub^{1}(T)= \ub^{2}(T)$, we infer that $\tilde{\vb}^{1}(T)= \tilde{\vb}^{2}(T)$. In view of \autoref{thm-global tilde v weak}, we deduce
\[
{\tilde{\vb}}^{1}(t)={\tilde{\vb}}^{2}(t), \mbox{ for }t\in [0, T], \mathbb{P}\text{-a.s.},
\]
therefore
\[
\ub^{1}(t)=\ub^{2}(t), \mbox{ for }t\in [0, T], \mathbb{P}\text{-a.s.},
\]
completes the proof of our main result \autoref{thm-main result}.
\end{proof}

\noindent{\bf  Acknowledgements}
We would also like to thank Professor Tusheng Zhang for numerous inspiring discussions. This work was partially supported by the National Key R\&D Program of China (Grant No. 2020YFA0712700); the Fundamental Research Funds for the Central Universities (Grant Nos. 2025CDJZKPT-09 and 2025CDJ-LAISYB-001); the NNSF of China (Grant Nos. 12090014, 11931004, 12471138 and 12371146); the Key Laboratory of Random Complex Structures and Data Science, Academy of Mathematics and Systems Science, Chinese Academy of Sciences (Grant No. 2008DP173182); the Natural Science Foundation Project of CQ (Grant No. CSTB2025NSCQ-GPX0722); and the Key Laboratory of Nonlinear Analysis and its Applications (Chongqing University), Ministry of Education.

\appendix
\renewcommand{\theequation}{\Alph{section}.\arabic{equation}}

\section{Fourier series on $\dom$}\label{sec-Fourier}
All Fourier expansions in this appendix are understood in $L^2$; scalar bases are applied componentwise to vector fields.

We will frequently work with the complex-valued analogs of these spaces, denoted by $\mathbb{L}^2(\dom; \mathbb{C}^3)$ which
is a complex Hilbert space with the inner product
\begin{equation}\label{eqn-L^2-complex inner product}
\langle \ub, \vb \rangle := \int_{\dom} \ub(x) \cdot \overline{\vb(x)} \, dx.
\end{equation}

Let $e_k(x)=(2\pi)^{-3/2}e^{ik\cdot x}$. Then
\[
\Vert e_k\Vert_{L^2}=1,
\qquad
\langle e_n,e_k\rangle=(2\pi)^{-3}\int_{\dom}e^{i(n-k)\cdot x}\,dx=0,
\quad n\ne k.
\]
The sequence $(e_k)_{k\in\mathbb Z^3}$ is an orthonormal basis of $L^2(\dom;\mathbb C)$. Consequently, the functions
\begin{equation}\label{eq-ONB-complex}
\begin{aligned}
     \tilde{\eb}_k^1(x)&:=(2\pi)^{-\frac32}\left(e^{ik\cdot x},0,0\right), \;\; x \in \dom,
     \\
     \tilde{\eb}_k^2(x)&:=(2\pi)^{-\frac32}\left(0,e^{ik\cdot x},0\right), \;\; x \in \dom,
     \\
     \tilde{\eb}_k^3(x)&:=(2\pi)^{-\frac32}\left(0,0,e^{ik\cdot x}\right), \;\; x \in \dom,
\end{aligned}
\end{equation}
form an orthonormal basis of $\mathbb L^2(\dom;\mathbb C^3)$.

\begin{cor}[Parseval identity, {\cite[Theorem 4.18]{Rudin_1987-RCA}}]\label{cor-Parseval}
For every $\ub\in\mathbb L^2(\dom;\mathbb C^3)$,
\begin{equation}\label{eqn-u-Fourier expansion}
\ub=\sum_{(k,j)\in\mathbb Z^3\times\{1,2,3\}}\tilde{\ub}_k^j\tilde{\eb}_k^j
=\sum_{k\in\mathbb Z^3}\check{\ub}_k e_k,
\end{equation}
where
\begin{equation}\label{eqn-u_k}
\tilde{\ub}_k^j=\langle\ub,\tilde{\eb}_k^j\rangle
=(2\pi)^{-3/2}\int_{\dom}u^j(x)e^{-ik\cdot x}\,dx,
\qquad j=1,2,3,
\end{equation}
and $\check{\ub}_k=(\tilde{\ub}_k^1,\tilde{\ub}_k^2,\tilde{\ub}_k^3)$. Moreover,
\begin{equation}\label{eqn-Parseval formula}
\Vert\ub\Vert_{L^2}^2
=\sum_{(k,j)\in\mathbb Z^3\times\{1,2,3\}}|\tilde{\ub}_k^j|^2
=\sum_{k\in\mathbb Z^3}|\check{\ub}_k|^2.
\end{equation}
\end{cor}

\begin{remark}\label{rem-Parseval}
For the unnormalized expansion $\ub=\sum_k\hat{\ub}_k e^{ik\cdot x}$, the coefficients satisfy
\[
\hat{\ub}_k=(2\pi)^{-3}\int_{\dom}\ub(x)e^{-ik\cdot x}\,dx,
\qquad
\check{\ub}_k=(2\pi)^{3/2}\hat{\ub}_k.
\]
Hence $\Vert\ub\Vert_{L^2}^2=(2\pi)^3\sum_k|\hat{\ub}_k|^2$.
\end{remark}

\begin{remark}
Assume $\ub\in \mathbb{L}^2(\dom;\mathbb{C}^3)$. The following are equivalent:
\begin{trivlist}
\item[(i)] $\ub(x)\in \mathbb{R}^3$ for almost every $x\in\dom$.
\item[(ii)] For $j=1,2,3$, in $L^2(\dom;\mathbb C)$,
\[
\sum_{k\in \mathbb{Z}^3} \tilde{\ub}_k^j e^{ik\cdot x}=\sum_{k\in \mathbb{Z}^3} \overline{\tilde{\ub}_{-k}^j} e^{ik\cdot x}
\]
\item[(iii)] For all $k\in \mathbb{Z}^3$ and $j=1,2,3$
\[
\tilde{\ub}_k^j =\overline{\tilde{\ub}_{-k}^j} .
\]
\end{trivlist}
The equivalence to (iii) is a consequence of the uniqueness of the Fourier series expansion.
\end{remark}
\begin{theorem}\label{thm-Fourier-real}
\begin{align}
\mathbb{L}^2(\dom;\mathbb{R}^{3})&=
\Bigl\{ \ub \in \mathbb{L}^2(\dom; \mathbb{C}^{3}) :
\overline{\tilde{\ub}^j_{\kb}}=\tilde{\ub}^j_{-\kb}, \;\; \mbox{ for all } \kb \in \mathbb{Z}^3, j=1,2,3. \Bigr\}.
\end{align}
Moreover, there exist $\ub_\kb^j\in\mathbb R$, $(\kb,j)\in\mathcal K$, such that
\begin{align}\label{eq-u-Fourier-real}
\ub = \sum_{j=1}^3 \ub_0^j \eb_0^j+ \sum_{j=1}^3 \sum_{\kb\in K_0} \ub_\kb^j\eb_\kb^{j}+ \sum_{j=1}^3 \sum_{\kb\in K_0} \ub_\kb^{j+3}\eb_\kb^{j+3}.
\end{align}
\end{theorem}
\begin{proof}[Proof of \eqref{eq-u-Fourier-real}]
Define
$\mathbb{Z}^3_0=\mathbb{Z}^3\setminus \{(0,0,0)\}$. Note that
\begin{equation}\label{eqn-u-Fourier expansion1}
\begin{aligned}
\ub = \sum_{(\kb,j) \in \mathbb{Z}^3 \times \{1,2,3\}} \tilde{\ub}_\kb^j \tilde{\eb}_\kb^j=\sum_{j=1}^3 \tilde{\ub}_0^j \tilde{\eb}_0^j+ \sum_{j=1}^3 (\sum_{\kb\in \mathbb{Z}_0^3} \tilde{\ub}_\kb^j \tilde{\eb}_\kb^j).
\end{aligned}
\end{equation}
Since
\[
\kb\in K_0 \mbox{ if and only if } -\kb \in \mathbb{Z}^3_0 \setminus K_0.
\]
Since $\overline{\tilde{\ub}^j_{\kb}}=\tilde{\ub}^j_{-\kb}$,
\[
\tilde{\ub}_{-\kb}^j \tilde{\eb}_{-\kb}^j=\overline{\tilde{\ub}_\kb^j} \tilde{\eb}_{-\kb}^j=\overline{\tilde{\ub}_\kb^j \tilde{\eb}_\kb^j},
\]
we get
\begin{equation}\label{eqn-u-Fourier expansion2}
\begin{aligned}
\ub &= \sum_{j=1}^3 \tilde{\ub}_0^j \tilde{\eb}_0^j+ \sum_{j=1}^3 (\sum_{\kb\in K_0} \tilde{\ub}_\kb^j \tilde{\eb}_\kb^j+\sum_{\kb\in \mathbb{Z}^3_0 \setminus K_0} \tilde{\ub}_\kb^j \tilde{\eb}_\kb^j)\\
&= \sum_{j=1}^3 \tilde{\ub}_0^j \tilde{\eb}_0^j+ \sum_{j=1}^3 \sum_{\kb\in K_0} (\tilde{\ub}_\kb^j \tilde{\eb}_\kb^j+ \tilde{\ub}_{-\kb}^j \tilde{\eb}_{-\kb}^j)\\
&= \sum_{j=1}^3 \tilde{\ub}_0^j \tilde{\eb}_0^j+ \sum_{j=1}^3 \sum_{\kb\in K_0} (\tilde{\ub}_\kb^j \tilde{\eb}_\kb^j+ \overline{\tilde{\ub}_\kb^j \tilde{\eb}_\kb^j})\\
&= \sum_{j=1}^3 \tilde{\ub}_0^j \tilde{\eb}_0^j+ 2(2\pi)^{-\frac32}\sum_{\kb\in K_0} (\Re(\tilde{\ub}_\kb^1 \eb^{i\kb\cdot \xb}),\Re(\tilde{\ub}_\kb^2 \eb^{i\kb\cdot \xb}),\Re(\tilde{\ub}_\kb^3 \eb^{i\kb\cdot \xb})).
\end{aligned}
\end{equation}
Note that
\[
\Re(\tilde{\ub}_\kb^j \eb^{i\kb\cdot \xb})=\Re(\tilde{\ub}_\kb^j)\cdot \cos \kb \xb-\Im(\tilde{\ub}_\kb^j)\cdot \sin \kb \xb\;\; j=1,2,3.
\]
Therefore,
\begin{equation}\label{eqn-u-Fourier expansion3}
\begin{aligned}
\ub
&= \sum_{j=1}^3 \tilde{\ub}_0^j \eb_0^j- \sum_{j=1}^3 \sum_{\kb\in K_0} \sqrt{2}\Im(\tilde{\ub}_\kb^j)\eb_\kb^{j}+ \sum_{j=1}^3 \sum_{\kb\in K_0} \sqrt{2}\Re(\tilde{\ub}_\kb^j)\eb_\kb^{j+3}.
\end{aligned}
\end{equation}
For $j=1,2,3$, set $\ub_0^j:=\tilde{\ub}_0^j$. For $k\in K_0 $, let
\begin{align}
\ub_\kb^j:=\begin{cases}
-\sqrt{2}\Im(\tilde{\ub}_\kb^j)\;\; j=1,2,3\\
\sqrt{2}\Re(\tilde{\ub}_\kb^{j-3}),\;\;j=4,5,6.
\end{cases}
\end{align}
where $\tilde{\ub}_\kb^j$ is defined in \eqref{eqn-u_k}. Therefore, we get
\begin{align*}
\ub = \sum_{j=1}^3 \ub_0^j \eb_0^j+ \sum_{j=1}^3 \sum_{\kb\in K_0} \ub_\kb^j\eb_\kb^{j}+ \sum_{j=1}^3 \sum_{\kb\in K_0} \ub_\kb^{j+3}\eb_\kb^{j+3}= \sum_{(k,j) \in \mathcal{K}} u_k^j \eb_k^j.
\end{align*}
\end{proof}

\section{Homogeneous Sobolev Spaces on $\tor$}
The definition of homogeneous Sobolev space on $\tor$ for $s\geq0$ presented below is adapted from \cite[Section 1.7.1]{RRS}.

For every $u \in \mathbb{H}^s(\tor)$, $s\geq 0$, the expression
\[
\Vert u\Vert_{\dot{H}^s}^2:=(2 \pi)^3 \sum_{k \in \dot{\mathbb{Z}}^3}|k|^{2 s}\left|\hat{u}_k\right|^2
\]
where we exclude $k=(0,0,0)$ from the summation by setting
\[
\dot{\mathbb{Z}}^3:=\left\{k \in \mathbb{Z}^3:|k| \neq 0\right\},
\]
is clearly well defined.
\begin{definition}
The space $\dot{\mathbb{L}}^2(\tor)$ consists of all functions $u \in \mathbb{L}^2(\tor)$ that have zero mean,
\[
\int_{\tor} u(x) \, d x=0.
\]
It is easy to see that $\dot{\mathbb{L}}^2(\tor)$ is given by the collection of all (real) functions with Fourier expansion in which $\hat{u}_0=0$, i.e. all real functions of the form
\[
\sum_{k \in \dot{\mathbb{Z}}^3} \hat{u}_k \mathrm{e}^{\mathrm{i} k \cdot x} \quad \text { with } \quad \sum_{k \in \dot{\mathbb{Z}}^3}\left|\hat{u}_k\right|^2<\infty .
\]
\end{definition}
\begin{definition}
The homogeneous Sobolev space $\mathbb{\dot{H}}^s(\tor)$ is defined for $s \geq 0$ as
\[
\mathbb{\dot{H}}^s(\tor):=\mathbb{H}^s(\tor) \cap \dot{\mathbb{L}}^2(\tor) .
\]
The nesting of homogeneous spaces, $\mathbb{\dot{H}}^{s_2}(\tor) \subset \mathbb{\dot{H}}^{s_1}(\tor)$ if $s_2>s_1$, holds in the periodic case.

For $s>0$, denote $(\dot{\mathbb{H}}^s(\tor))'$ by $\dot{\mathbb{H}}^{-s}(\tor)$, represented by zero-mean periodic distributions, with norm
\[
\Vert\vb\Vert_{\dot H^{-s}}^2:=(2\pi)^3
\sum_{k\in\dot{\mathbb Z}^3}|k|^{-2s}|\hat{\vb}_k|^2.
\]
\end{definition}
Indeed, suppose that $s>0$, $u \in \dot{\mathbb{H}}^s$, and $\vb$ is given by a formal expression
\[
\vb(x) = \sum_{k \in {\dot{\mathbb{Z}}^3}} \hat{\vb}_k e^{i k x} \quad \text{with} \quad \sum_{k \in {\dot{\mathbb{Z}}^3}} |k|^{-2s} |\hat{\vb}_k|^2 < \infty.
\]
Then we can attach meaning to the expression
\[
\langle \vb, u \rangle_{{\dot{\mathbb{H}}^{-s}} \times \dot{\mathbb{H}}^s} := {(2\pi)^3}\sum_{k \in {\dot{\mathbb{Z}}^3}} \hat{\vb}_k {\cdot}\hat{u}_{-k}.
\]

\section{Sobolev and heat semigroup tools}
\begin{definition}
If $u \in \mathbb{L}^2(\dom)$ and $\alpha\in \mathbb{N}^3$, then we say that the weak derivative $D^\alpha u$ exists if and only if there exists $\vb\in \mathbb{L}^2(\dom)$ such that for every $\phi\in C_0^\infty(\dom)$
\begin{align}\label{eq-weak derivative}
\int_{\dom}u(x) D^\alpha \phi(x)\,dx=(-1)^{\vert\alpha\vert}\int_{\dom}\vb(x) \phi(x)\,dx
\end{align}
where both sides of the above equality $\in \mathbb{R}^3$.
\end{definition}

\begin{definition}
If $n\in\mathbb{N}$, then $\mathbb{H}^n(\dom)=\mathbb{H}^n(\dom;\mathbb{R}^3)$ is the space of all $u\in\mathbb{L}^2(\dom)$ such that for every $|\alpha|\leq n$ the weak derivative $D^\alpha u$ exists and belongs to $\mathbb{L}^2(\dom)$.
\end{definition}
\begin{prop}
$\mathbb{H}^n(\dom)$ is a Hilbert space with inner product and norm defined as below:
\begin{align}\label{eq-inner product}
\langle u,\vb\rangle_{\mathbb{H}^n}=\sum_{|\alpha|\leq n}\int_{\dom}D^\alpha u(x)\cdot D^\alpha \vb(x)\,dx=\sum_{|\alpha|\leq n}\int_{\dom}\langle D^\alpha u(x), D^\alpha \vb(x)\rangle_{\mathbb{R}^3}\,dx.
\end{align}
\end{prop}

One can also define $\mathbb{H}^s(\dom)$ for $s>0$. Several equivalent definitions are available.

In view of \autoref{def-Laplace}, we define operator $\rA$  in the Hilbert space $\mathbb{L}^2(\dom)$ by
\begin{equation}\label{eqn-Laplace-3}
 \begin{aligned}
&D(\rA)=\mathbb{H}^2(\dom)\\
&\rA u := (-\Delta u_1, -\Delta u_2,-\Delta u_3) \end{aligned}
\end{equation}
with equivalent norms. $\rA$ is non-negative and self-adjoint.

\begin{definition}\label{def-graph norm}
If $X$ is a normed space, $\rA: D(\rA)\to X$ is a linear operator, where $D(\rA)\subset X$. The graph norm is defined by
\begin{align}\label{eq-graph norm1}
\vert x \vert_{D(\rA)}:=\vert x \vert+\vert \rA x \vert,\;\;x\in D(\rA).
\end{align}
If $X$ is a Hilbert space, then we can also use "graph" norm defined by
\begin{align}\label{eq-graph norm2}
\Vert x \Vert_{D(\rA)}^2:=\vert x \vert^2+\vert \rA x \vert^2,\;\;x\in D(\rA).
\end{align}
\end{definition}
\begin{remark}
If $X$ is Banach. $\rA$ is a closed operator if and only if $D(\rA)$ is endowed with a complete graph norm.
\end{remark}
\begin{lemma}
We can define many norms on $\mathbb{R}^2$. If $x=(x_1,x_2)\in \mathbb{R}^2$.
\begin{equation}
\begin{aligned}
&\vert x\vert_2:=(\sum_{i=1}^2 x_i^2)^{\frac{1}{2}},\\
&\vert x\vert_1:=\sum_{i=1}^2 \vert x_i\vert,\\
&\vert x\vert_\infty:=\max_{i=1,2} \vert x_i\vert.
\end{aligned}
\end{equation}
These three norms are equivalent.
\end{lemma}
\begin{cor}
If $X$ is a Hilbert space, then $\Vert \cdot \Vert_{D(\rA)}$ and $\vert \cdot \vert_{D(\rA)}$ are equivalent.
\end{cor}

\begin{prop}
The $\mathbb{H}^2$-norm is equivalent to the "graph" norm, i.e.  {we have}
\begin{align}
{\frac12}\Vert u \Vert_{D(\rA)}^2\leq \Vert u\Vert_{H^2}^2\leq {\frac32}\Vert u \Vert_{D(\rA)}^2,\;\; \forall u\in D(\rA)=\mathbb{H}^2(\dom).
\end{align}
\end{prop}
\begin{proof}
Note that
\begin{equation}
\begin{aligned}
&\Vert u \Vert_{D(\rA)}^2=\Vert u \Vert_{L^2}^2+\Vert \rA{u} \Vert_{L^2}^2=\int_{\dom}\vert u\vert^2+\vert\frac{\partial^2u}{\partial x_1^2}\vert^2+\vert\frac{\partial^2u}{\partial x_2^2}\vert^2+\vert\frac{\partial^2u}{\partial x_3^2}\vert^2\,dx
{+2\sum_{1\leq i<j\leq3}\Vert\partial_{ij}u\Vert_{L^2}^2},\\
&\Vert u\Vert_{H^2{,\mathrm{der}}}^2=\sum_{|\alpha|\leq 2}\int_{\dom}\vert D^\alpha u\vert^2\,dx.
\end{aligned}
\end{equation}
The identity for $\Vert\rA u\Vert_{L^2}^2$ follows by periodic integration by parts.
Since $\Vert\nabla u\Vert_{L^2}^2=\langle\rA u,u\rangle
\leq\frac12\Vert u\Vert_{D(\rA)}^2$, the claimed bounds follow.
\end{proof}
From now on, we use only the "graph" norm on $D(\rA)$, i.e.
\begin{align}\label{eq-graph norm3}
\Vert u \Vert_{D(\rA)}^2:=\Vert u \Vert_{L^2}^2+\Vert \rA u \Vert_{L^2}^2,\;\;u\in D(\rA).
\end{align}
Since the sequence $(\eb_\kb^j:(\kb,j) \in \mathcal{K})$ is an orthonormal basis of the real Hilbert space  $\mathbb{L}^2(\dom)$
\begin{align}\label{eqn-u-real1}
\ub(x) =  \sum_{(k,j) \in \mathcal{K}} u_k^j \eb_k^j(x),   \;\;\ub_\kb^j := \langle \ub , \eb_k^j \rangle.
\end{align}
and
\begin{equation}\label{eqn-Laplace1}
\rA e_k^j =\vert k \vert^2 e_k^j.
\end{equation}
Then we have
\begin{align}\label{eqn-Parseval-real1}
\Vert \ub\Vert_{L^2}^2 =  \sum_{(k,j) \in \mathcal{K}} |\ub_k^j|^2.
\end{align}
and
\begin{equation}\label{eq-Au norm}
\Vert \rA u \Vert_{L^2}^2=\sum_{(k,j) \in \mathcal{K}} |k|^4|\ub_k^j|^2.
\end{equation}
Hence,
\begin{equation}\label{eq-D(A) norm}
\begin{aligned}
\Vert u \Vert_{D(\rA)}^2=\sum_{(k,j) \in \mathcal{K}} |\ub_k^j|^2+\sum_{(k,j) \in \mathcal{K}} |k|^4|\ub_k^j|^2=\sum_{(k,j) \in \mathcal{K}} (1+|k|^4)|\ub_k^j|^2.
\end{aligned}
\end{equation}

With the Fourier norm \eqref{eqn-norm H^s}, $\Vert u\Vert_{H^2}^2=\Vert u\Vert_{D(\rA)}^2$.

The fractional Sobolev spaces $\mathbb{H}^s(\dom)$ can be defined {for $s>0$} as
\begin{align}\label{eq-fractional Sobolev spaces}
\mathbb{H}^s(\dom):=D(I+\rA^{\frac{s}{2}})
\end{align}
with the norm
\begin{align}\label{eq-fractional Sobolev norm}
\Vert u \Vert_{H^s}^2:= \sum_{(k,j) \in \mathcal{K}} (1+|k|^{2s})|\ub_k^j|^2
\end{align}
Equivalent definition:
\begin{align}\label{eq-fractional Sobolev spaces1}
\tilde{\mathbb{H}}^s(\dom):=D((I+\rA)^{\frac{s}{2}})
\end{align}
with the equivalent norm
\begin{align}\label{eq-fractional Sobolev norm1}
\Vert u \Vert_{\tilde{H}^s}^2:= \sum_{(k,j) \in \mathcal{K}} (1+|k|^4)^{\frac{s}{2}}|\ub_k^j|^2
\end{align}
The natural spectral norm is
\[
\Vert(I+\rA)^{s/2}u\Vert_{L^2}^2
=\sum_{(k,j)\in\mathcal K}(1+|k|^2)^s|\ub_k^j|^2.
\]
The equivalence follows from $1+r^4\leq(1+r^2)^2\leq2(1+r^4)$ for $r\geq0$.

\begin{prop}
$\Vert u \Vert_{H^s}$ and $\Vert u \Vert_{\tilde{H}^s}$ are equivalent.
\end{prop}
\begin{proof}
Let us choose and fix $s>0$. It is sufficient to show that there exist constants $C_1,C_2>0$ such that
\begin{align}\label{eqn-Hs norm equivalent}
C_1(1+|k|^{2s})\leq
(1+|k|^{4})^{\frac{s}{2}}
\leq  C_2(1+|k|^{2s}),
\;\;\forall k\in\mathbb Z^{3}.
\end{align}

Step 1: For $0\leq |k|\leq 1$.
We have $1\leq 1+|k|^{4}\leq 2$, then
\begin{align}\label{eqn-3.20}
  1 \leq (1+|k|^{4})^{\frac{s}{2}} \leq 2^{\frac{s}{2}}.
\end{align}
Moreover, since
$1\leq 1+|k|^{2s}\leq 2$ we have
\begin{align}\label{eqn-3.21}
\frac{1}{2}(1+|k|^{2s})\leq 1 \leq 1+|k|^{2s}
\end{align}
Combine \eqref{eqn-3.20} and \eqref{eqn-3.21} we infer that
\[
  \frac{1}{2}(1+|k|^{2s})
\leq  (1+|k|^{4})^{\frac{s}{2}}
\leq  2^{\frac{s}{2}}\,(1+|k|^{2s}).
\]

Step 2: For $|k|> 1$.
Since $|k|^{4}\leq 1+|k|^{4}\leq 2|k|^{4}$, we infer that
\begin{align}\label{eqn-3.22}
 |k|^{2s}\leq   (1+|k|^{4})^{\frac{s}{2}}\leq 2^{\frac{s}{2}}|k|^{2s}.
\end{align}
Moreover, we also have $1+|k|^{2s}\leq 2|k|^{2s}$, therefore
\begin{align}\label{eqn-3.23}
\frac{1}{2}(1+|k|^{2s})\leq   \frac{1}{2} 2|k|^{2s}=|k|^{2s}.
\end{align}
Combine \eqref{eqn-3.22} and \eqref{eqn-3.23} we infer that
\[
  \frac{1}{2}(1+|k|^{2s})\leq   (1+|k|^{4})^{\frac{s}{2}}\leq
  2^{\frac{s}{2}}(1+|k|^{2s}).
\]

Therefore \eqref{eqn-Hs norm equivalent} holds for $C_1=\frac{1}{2}$ and $C_2=2^{\frac{s}{2}}$.

\end{proof}

\subsection{Heat-semigroup estimates}\label{app-heat-semigroup-proofs}
We first state the following two known results. The first theorem can be found in \cite[Theorem from section 7.3.2]{Arendt_2002}. The second is a consequence of \cite[Theorem 6.13]{Pazy_1983}.

For $1<p<\infty$, let $\rA_p$ be the realization of $\rA$ on
$\mathbb{L}^p(\dom)$, that is,
\[
D(\rA_p)=\{u\in\mathbb{L}^p(\dom):-\Delta u\in\mathbb{L}^p(\dom)\},
\;\;
\rA_p u=(-\Delta u_1,-\Delta u_2,-\Delta u_3).
\]
Then $-\rA_p$ generates the analytic semigroup
$(e^{-t\rA_p})_{t\geq0}$ on $\mathbb{L}^p(\dom)$.
For $\omega>0$ and $\alpha>0$, we define
\[
D\left((\omega+\rA_p)^\alpha\right)
:=
\left\{
u\in \mathbb{L}^p(\dom):
(\omega+\rA_p)^\alpha u\in \mathbb{L}^p(\dom)
\right\},
\]
equipped with the graph norm
\[
\Vert u\Vert_{D((\omega+\rA_p)^\alpha)}
=
\Vert u\Vert_{L^p}
+
\Vert(\omega+\rA_p)^\alpha u\Vert_{L^p}.
\]

Since $\rA_p$ is the periodic Laplace operator on $\dom$, the
fractional domain is characterized by the Bessel potential space:
\[
D\left((\omega+\rA_p)^\alpha\right)
=
\mathbb{H}^{2\alpha,p}(\dom)
\]
with equivalent norms. Hence, by the Sobolev embedding on $\dom$,
\[
\mathbb{H}^{2\alpha,p}(\dom)
\hookrightarrow
\mathbb{L}^q(\dom),
\qquad
\frac1q=\frac1p-\frac{2\alpha}{3},
\]
provided $0<\alpha<\frac{3}{2p}$. Therefore,
\[
D\left((\omega+\rA_p)^\alpha\right)
\hookrightarrow
\mathbb{L}^q(\dom),
\qquad
\alpha=\frac32\left(\frac1p-\frac1q\right).
\]

\begin{theorem}\label{thm-Arendt-2004}
Consider the following conditions:
\begin{trivlist}
\item[(i)] There exist $c>0$, $1 \leq p<q \leq \infty$ such that
\[
\Vert e^{-tA}f\Vert_{L^q} \leq c t^{-\frac{3}{2}\left(\frac{1}{p}-\frac{1}{q}\right)} \Vert f\Vert_{L^p},\;\;\text { for all } 0<t \leq 1 .
\]
\item[(ii)] There exists a constant $c>0$ such that for all $1 \leq p<q \leq \infty$,
\[
\Vert e^{-tA}f\Vert_{L^q} \leq c t^{-\frac{3}{2}\left(\frac{1}{p}-\frac{1}{q}\right)} \Vert f\Vert_{L^p},\;\; 0<t \leq 1 .
\]
\item[(iii)]There exist $1<p<\infty$ and $0<\alpha<\frac{3}{2 p}$ such that $D\left(\left(\omega+A_p\right)^\alpha\right) \subset \mathbb{L}^q(\dom)$, where $q$ is defined by $\alpha=\frac{3}{2}\left(\frac{1}{p}-\frac{1}{q}\right)$.
\item[(iv)]For all $0<\alpha<\frac{3}{2 p}$ and all $1<p<\infty$ one has $D\left(\left(\omega+A_p\right)^\alpha\right) \subset \mathbb{L}^q(\dom)$, where $q$ is defined by $\alpha=\frac{3}{2}\left(\frac{1}{p}-\frac{1}{q}\right)$.
\end{trivlist}
Then all four assertions are equivalent.
\end{theorem}
\begin{theorem}\label{thm-Pazy}
Assume that $p\in[1,\infty)$. Then for every $\gamma\geq0$ and $T>0$ there exists $C=C(A,T,p,\gamma)>0$ such that
\begin{align}\label{eqn-2.39}
\left\Vert \Lambda^\gamma e^{-\frac{t}{2}A} h\right\Vert_{L^p}
\leq C t^{-\frac{\gamma}{2}}\Vert h\Vert_{L^p}, \;\; h \in\mathbb{L}^p(\dom), \;\;t \in (0,T].
\end{align}
\end{theorem}
\begin{proof}
Write $A=-\Delta$ and $h_0=t/2$. When $\gamma=0$ the claim follows from heat-semigroup contraction. For $\beta=\gamma/2>0$, choose an integer $m>\beta$. The kernel of $A^m e^{-vA}$ is the periodization of $(-\Delta)^m$ applied to the Euclidean heat kernel at time $v$. Gaussian scaling gives an $L^1(\mathbb R^3)$ norm of $C_m v^{-m}$, which also bounds the $L^1(\dom)$ norm of the periodization. Young's inequality therefore gives $\|A^m e^{-vA}\|_{L^p\to L^p}\leq C_m v^{-m}$ for $v>0$.
On trigonometric polynomials, integration on each Fourier mode gives
\[
A^\beta e^{-h_0 A}h=
\frac{1}{\Gamma(m-\beta)}
\int_0^\infty r^{m-\beta-1}A^m e^{-(h_0+r)A}h\,dr .
\]
Both sides vanish on the constant mode. The integral is bounded in $L^p$ by
\[
C\|h\|_{L^p}\int_0^\infty r^{m-\beta-1}(h_0+r)^{-m}\,dr
=C' h_0^{-\beta}\|h\|_{L^p}.
\]
Density in $L^p(\dom)$ proves \eqref{eqn-2.39}.
\end{proof}
\begin{proof}[Proof of \autoref{lem3.3-1}]
Let us choose and fix $p \in (1, \infty)$,  $q \in [p, \infty)$ and $r\in (p,\infty]$ satisfying
condition \eqref{eqn-pqr}, $\gamma\geq 0$ and $h \in \mathbb{L}^p(\dom)$. Then, for  $0<t\leq 1$,  we have
 \begin{align}\label{eqn-semigroup-101}
\Lambda^\gamma e^{-tA} h= \Lambda^\gamma e^{-\frac{t}{2}A}  e^{-\frac{t}{2}A} h =\Lambda^\gamma e^{-\frac{t}{2}A}g
\end{align}

By \autoref{thm-Pazy}, we infer that there exists $C>0$ independent of $g$ such that
\begin{align}\label{eqn-2.38-1}
\left\Vert \Lambda^\gamma e^{-\frac{t}{2}A} g\right\Vert_{L^q}
\leq C (\frac{t}{2})^{-\frac{\gamma}{2}}\Vert g\Vert_{L^q}=
C' t^{-\frac{\gamma}{2}}\Vert g\Vert_{L^q}
, \;\; t \in (0,1].
\end{align}

For $p<q$, \autoref{thm-Arendt-2004} gives the following estimate. For $p=q$, we have $r=\infty$ and the same bound follows from the $L^p$ contraction of the periodic heat semigroup. In both cases there is a constant $C$ independent of $h$ such that
\begin{equation}\label{eqn-2.38}
\begin{aligned}
\Vert g\Vert_{L^q}= \left\Vert  e^{-\frac{t}{2}A} h\right\Vert_{L^q}\leq C
(\frac{t}{2})^{-\frac{3}{2r}}\left\Vert  h\right\Vert_{L^p}
=C{''}t^{-\frac{3}{2r}}\left\Vert  h\right\Vert_{L^p}.
\end{aligned}
\end{equation}

Therefore, by \eqref{eqn-semigroup-101},  \eqref{eqn-2.38-1} and \eqref{eqn-2.38} we
deduce that
\begin{align}
\left\Vert \Lambda^\gamma e^{-tA} h\right\Vert_{L^q}
\leq C t^{-\frac{\gamma}{2}-\frac{3}{2 r}}\Vert h\Vert_{L^p}.
\end{align}
The proof is complete.
    \end{proof}

\begin{proof}[Proof of \autoref{lem-heat semigroup on homogenous spaces}]
Let us choose and fix $g \in {\mathbb{H}}^{s_1}(\dom)$. By assumption we infer that $\gamma \geq 0$. Hence
by \eqref{eq-homogeneous norm}, \eqref{eq-Lambda-1} and \autoref{lem3.3-1}, and the fact that the operators $e^{-rA}$ commute with fractional powers of $\Lambda$, there exists $C>0$ such that for every $0<r\leq 1$
\begin{equation*}
\begin{aligned}
\left\Vert  e^{-rA} g \right\Vert_{H^{s_2}}&\leq \big\Vert  \Lambda^{s_2} e^{-rA}g\big\Vert_{L^2}+\big\Vert  e^{-rA} g\big\Vert_{L^2}= \big\Vert  \Lambda^{\gamma} e^{-rA}\Lambda^{s_1} g\big\Vert_{L^2}+\big\Vert  e^{-rA} g\big\Vert_{L^2}\\
&\leq C r^{-\frac{\gamma}{2}} \big\Vert  \Lambda^{s_1} g \big\Vert_{L^2}+\big\Vert  g\big\Vert_{L^2}\leq (Cr^{-\frac{\gamma}{2}}+1) \Vert  g \Vert_{H^{s_1}},
\end{aligned}
\end{equation*}
which implies (i) holds.

By part (i), for every $\varphi\in \mathbb{H}^{s_1}(\dom)$, there exists $C>0$ such that
\[
\Vert e^{-rA}\varphi\Vert_{H^{s_2}} \leq (Cr^{-\frac{\gamma}{2}}+1) \Vert  \varphi \Vert_{H^{s_1}}.
\]
Therefore, in view of \eqref{eq-negative norm}, by the H\"older inequality and the fact that the operators $e^{-rA}$ is self-adjoint, for every $g \in \mathbb{H}^{-s_2}(\dom)$, we have
\[
\Vert e^{-rA}g\Vert_{H^{-s_1}}
= \sup_{\Vert \varphi\Vert_{H^{s_1}}=1}|\langle g, e^{-rA}\varphi \rangle| \leq \sup_{\Vert \varphi\Vert_{H^{s_1}}=1} \Vert g\Vert_{H^{-s_2}}\Vert e^{-rA}\varphi\Vert_{H^{s_2}} \leq (Cr^{-\frac{\gamma}{2}}+1)\Vert g\Vert_{H^{-s_2}}.
\]
The proof is complete.
\end{proof}

\subsection{Consequences of maximal regularity}\label{app-linear-heat-proof}
\begin{proof}[Proof of \autoref{lemma-lplq}]
In view of \autoref{thm-maximal regularity}, let $q=\frac{4}{3}$ and $s=0$, we infer that
\[
\zb \in C\left([0,T] ; B^{\frac12}_{2,\frac43}\right) .
\]
and there exists $C=C_T>0$ such that
\begin{align}\label{eqn-3.8-1}
\Vert \zb\Vert_{C([0,T] ; B^{\frac12}_{2,\frac43})}
\leq C\Vert   \fb \Vert _{L^{\frac{4}{3}}([0,T]; \mathbb{L}^{2}(\dom))}, \;\; f\in L^{\frac{4}{3}}([0,T]; \mathbb{L}^{2}(\dom)),
\end{align}
By \cite[Theorem in Section 3.5.4 and Remark 4 in Section 3.5.1]{Schmeisser+Triebel_1987},
\[
B^{\frac12}_{2,\frac43}\hookrightarrow B^{\frac12}_{2,2}  =\mathbb{H}^\frac12(\dom).
\]
Therefore, \eqref{eqn-3.8-1} implies $\zb \in C([0,T]; \mathbb{H}^\frac{1}{2}(\dom))$ and
\begin{align}\label{eqn-3.8-2}
\Vert \zb\Vert_{C([0,T] ; \mathbb{H}^\frac12)}\leq C\Vert \zb\Vert_{C([0,T] ; B^{\frac12}_{2,\frac43})}
\leq C\Vert f\Vert_{L^{\frac{4}{3}}([0,T]; \mathbb{L}^{2}(\dom))}, \;\; f\in L^{\frac{4}{3}}([0,T]; \mathbb{L}^{2}(\dom)).
\end{align}

Moreover, by \autoref{thm-maximal regularity} we infer that $\zb \in L^{\frac{4}{3}}([0,T]; \mathbb{H}^{2}(\dom)) $.
By the Sobolev interpolation inequality (see \cite[Corollary 5.1]{Brezis+Mironescu_2018}) we infer that there exists a constant $C>0$ such that
\begin{align}\label{eqn-interpolation H3/2}
\Vert u\Vert_{H^\frac{3}{2}}\leq C\Vert u\Vert_{H^\frac{1}{2}}^\frac{1}{3}\Vert u\Vert_{H^2}^\frac{2}{3}.
\end{align}
Then by inequalities \eqref{eqn-interpolation H3/2}, \eqref{eqn-3.8-2} and \eqref{eqn-3.8} in \autoref{thm-maximal regularity} we deduce there exists $C=C_T>0$ such that
\begin{equation}\label{eqn-3.8-3}
\begin{aligned}
\Vert  \zb  \Vert _{L^{2}([0,T]; \mathbb{H}^{\frac{3}{2}}(\dom))}&=\left(\int_0^T \Vert \zb(t)\Vert_{H^\frac{3}{2}}^2\,dt\right)^\frac{1}{2}\leq C\left(\int_0^T \Vert \zb\Vert_{H^\frac{1}{2}}^\frac{2}{3}\Vert \zb\Vert_{H^2}^\frac{4}{3}\,dt\right)^\frac{1}{2}\\
&\leq C \Vert \zb \Vert_{L^{\infty}([0,T]; \mathbb{H}^{\frac{1}{2}}(\dom))}^\frac{1}{3}\left(\int_0^T \Vert \zb \Vert_{H^2}^\frac{4}{3}\,dt\right)^\frac{1}{2}\\
& =C \Vert \zb \Vert_{L^{\infty}([0,T]; \mathbb{H}^{\frac{1}{2}}(\dom))}^\frac{1}{3}\Vert\zb \Vert_{L^{\frac{4}{3}}([0,T]; \mathbb{H}^2(\dom))}^\frac{2}{3}\\
&\leq C\Vert f\Vert_{L^{\frac{4}{3}}([0,T]; \mathbb{L}^{2}(\dom))}, \;\; f\in L^{\frac{4}{3}}([0,T]; \mathbb{L}^{2}(\dom)).
\end{aligned}
\end{equation}
Therefore, by \eqref{eqn-3.8-2} and \eqref{eqn-3.8-3} we infer that  (i) holds.

To prove (ii), by \cite[Theorem in Section 3.5.4 and Theorem 1 in Section 3.6.1]{Schmeisser+Triebel_1987}, we infer
\[
\left(\mathbb{H}^1(\dom),\mathbb{H }^3(\dom)\right)_{\frac12,2}=B^{2}_{2,2}=\mathbb{H}^2 (\dom).
\]
In view of \autoref{thm-maximal regularity}, let $q=2$ and $s=1$, we infer that
$\zb \in C\left([0,T] ; \mathbb{H}^2 (\dom)\right)$ and there exists $c=c_T>0$ such that
\begin{align}\label{eqn-3.9-2}
\Vert \zb\Vert_{C([0,T] ; \mathbb{H}^2(\dom))}
\leq c\Vert f\Vert_{L^2([0,T]; \mathbb{H}^{1}(\dom))}, \;\; f\in L^2([0,T]; \mathbb{H}^{1}(\dom)).
\end{align}
The proof is complete.
\end{proof}

\begin{lemma}\label{semigroup-derivative}
Let $T>0$ and assume that $\ub\in L^{\frac43}([0,T];\mathbb{L}^{2}(\dom))$.
Then, for almost every $t\in[0,T]$,
\begin{align}\label{eq-differential}
\frac{d}{dt}\int_{0}^{t}e^{(t-s)\Delta}\ub(s)\,ds
=
\ub(t)
+
\Delta\int_{0}^{t}e^{(t-s)\Delta}\ub(s)\,ds
\quad\text{in }\mathbb{L}^{2}(\dom).
\end{align}
Moreover,
\begin{align}\label{eq-differential-2}
\frac{d}{dt}\int_{0}^{t}e^{(t-s)\Delta}\ub(s)\,ds
\in
L^{1}([0,T];\mathbb{H}^{-\frac12}(\dom)).
\end{align}
\end{lemma}

\begin{proof}
Let us set
\[
U(t):=\int_0^t e^{(t-s)\Delta}\ub(s)\,ds,
\qquad t\in[0,T].
\]
By \autoref{thm-maximal regularity} with
$q=\frac43$ and $s=0$, we infer that
\begin{align}\label{eq-U-maximal}
U\in L^{\frac43}([0,T];\mathbb H^2(\dom)),
\qquad
\frac{dU}{dt}\in L^{\frac43}([0,T];\mathbb L^2(\dom)).
\end{align}

Let us choose and fix $\eta\in D(\Lambda^2)$ and
$\varphi\in C_c^\infty((0,T))$.
By Fubini's theorem and integration by parts,
we infer that
\begin{align*}
-\int_0^T\langle U(t),\eta\rangle\varphi'(t)\,dt
&=\int_0^T\langle\ub(t),\eta\rangle\varphi(t)\,dt
+\int_0^T\langle U(t),\Delta\eta\rangle\varphi(t)\,dt\\
&=\int_0^T\langle\ub(t)+\Delta U(t),\eta\rangle
\varphi(t)\,dt,
\end{align*}
where the last equality follows from the
self-adjointness of $\Delta$.

On the other hand, by \eqref{eq-U-maximal},
we have
\[
-\int_0^T\langle U(t),\eta\rangle\varphi'(t)\,dt
=\int_0^T\left\langle\frac{dU}{dt}(t),\eta\right\rangle
\varphi(t)\,dt.
\]
Since $D(\Lambda^2)$ is dense in $\mathbb L^2(\dom)$
and $\frac{dU}{dt},\ub,\Delta U$ belong to
$L^{\frac43}([0,T];\mathbb L^2(\dom))$,
by the uniqueness of distributional derivatives,
we infer that
\[
\frac{dU}{dt}(t)=\ub(t)+\Delta U(t)
\quad\text{in }\mathbb L^2(\dom),
\quad\text{for a.e. }t\in[0,T].
\]
Hence, \eqref{eq-differential} holds.

Finally, since
$\mathbb L^2(\dom)\hookrightarrow
\mathbb H^{-\frac12}(\dom)$,
by the H\"older inequality and \eqref{eq-U-maximal},
we infer that
\begin{align*}
\left\Vert\frac{dU}{dt}\right\Vert_
{L^1([0,T];\mathbb H^{-\frac12}(\dom))}
&\leq cT^{\frac14}
\left\Vert\frac{dU}{dt}\right\Vert_
{L^{\frac43}([0,T];\mathbb L^2(\dom))}
<\infty.
\end{align*}
Therefore, \eqref{eq-differential-2} holds.
The proof is complete.
\end{proof}

\section{Mollification and coefficient approximation}\label{sec-mollification}
In this subsection, we introduce a mollifier that will allow us to build smooth approximations to given functions.
\begin{definition}
If $U\subset \mathbb{R}$ is open and $\eps>0$, we write
\[
U_\eps := \{\, x\in U \mid \mathrm{dist}(x,\partial U)>\eps \,\}.
\]
\end{definition}
\begin{definition}
\item Define the standard mollifier $\eta\in C^\infty(\mathbb{R})$ by
\[
\eta(x):=
\begin{cases}
C\exp\!\left(\dfrac{1}{|x|^2-1}\right), & |x|<1,\\[1ex]
0, & |x|\ge 1,
\end{cases}
\]
where the constant $C>0$ is chosen such that
\[
\int_{\mathbb{R}} \eta(x)\,dx = 1.
\]

\item For each $\eps>0$, set
\[
\eta_\eps(x):=\frac{1}{\eps}\eta\!\left(\frac{x}{\eps}\right).
\]
The functions $\eta_\eps$ are $C^\infty$ and satisfy
\[
\int_{\mathbb{R}} \eta_\eps(x)\,dx = 1,
\;\;
\mathrm{spt}(\eta_\eps)=[-\eps,\eps].
\]
\end{definition}
\begin{definition}\label{def-mollification}
If $f:U\to\mathbb{R}$ is locally integrable, define its mollification by
\[
f^\eps := \eta_\eps * f
\qquad \text{in } U_\eps .
\]
That is,
\[
f^\eps(x)
= \int_U \eta_\eps(x-y)f(y)\,dy
= \int_{-\eps}^\eps \eta_\eps(y)f(x-y)\,dy,
\qquad x\in U_\eps.
\]
\end{definition}

\begin{lemma}\label{lem-mollification}
We have the following result.
\begin{itemize}
\item[(i)] $f^\eps\in C^\infty(U_\eps)$.
\item[(ii)] $f^\eps \to f$ a.e.\ as $\eps\to 0$.
\item[(iii)] If $f\in C(U)$, then $f^\eps \to f$ uniformly on compact subsets of $U$.
\item[(iv)] If $1\leq p<\infty$ and $f\in L^p_{\loc}(U)$, then
      $f^\eps \to f$ in $L^p_{\loc}(U)$.
\end{itemize}
\end{lemma}

\subsection{Convergence of the regularized coefficients}\label{app-proof-uniform-convergence}

\begin{proof}[Proof of \autoref{lem-uniform convergence}]
Since for almost every $\omega$,
\[
B(t,\omega)\in C([0,T]).
\]
Hence, by Part (iii) in \autoref{lem-mollification}, we infer that
\[
\Vert B^\eps-B\Vert_{L^\infty([0,T])}\to 0.
\]
Note that the exponential function is locally Lipschitz, thus there exists $c>0$ such that
\[
\Vert\alpha^{\eps}-\alpha\Vert_{L^{\infty}([0,T])}=\Vert e^{bB^{\eps}} -e^{bB}\Vert_{L^{\infty}([0,T])}\leq c\Vert B^\eps-B\Vert_{L^\infty([0,T])}\to 0.
\]
Similarly, we infer that
\[
\Vert(\alpha^\eps)^{-1}-\alpha^{-1}\Vert_{L^\infty([0,T])}
\leq
c\Vert B^\eps-B\Vert_{L^\infty([0,T])}
\to 0.
\]
Thus, \eqref{eq-uniform convergence1} holds.

Recall that
\[
\int_{\mathbb{R}} \eta_\eps(s)\,ds=\int_{B(0,\eps)} \eta_\eps(s)\,ds = 1.
\]
Since $\fb^\eps$ is defined as the convolution of $\fb$ with the mollifier $\eta_\eps$, see \autoref{def-mollification}. For every $t \in [0,T]$,
\[
\fb^\eps(t)-\fb(t)=\int_{B(0,\eps)}\eta_\eps(s)\fb(t-s)\,ds-\fb(t)=\int_{B(0,\eps)}\eta_\eps(s)\bigl(\fb(t-s)-\fb(t)\bigr)\,ds.
\]
By Minkowski's inequality we infer that
\[
\Vert\fb^\eps - \fb\Vert_{L^2([0,T]; \mathbb{H}^\delta(\dom))} \leq \int_{B(0,\eps)}\eta_\eps(s)\,
\Vert\fb(\cdot-s)-\fb(\cdot)\Vert_{L^2([0,T];\mathbb H^\delta(\dom))}\,ds.
\]
In view of \autoref{A2}, $\fb \in \mathbb{L}^2_{\loc}([0, \infty); \mathbb{H}^\delta(\dom))$ with $\delta > \frac{3}{2}$. Therefore, $\fb(t) \in L^2([0,T];\mathbb H^\delta(\dom))$, which implies
\[
\Vert\fb(\cdot-s)-\fb(\cdot)\Vert_{L^2([0,T];\mathbb H^\delta(\dom))} \to 0,\;\mbox{as}\;s \to 0.
\]
Moreover, by the triangle inequality,
\[
\|\fb(\cdot-s)-\fb(\cdot)\|_{L^2([0,T];\mathbb H^\delta(\dom))}
\le 2\|\fb\|_{L^2([0,T+1];\mathbb H^\delta(\dom))}\qquad \text{for }|s|<1.
\]
By the translation continuity established above and the unit mass of $\eta_\eps$, we have
\[
\Vert\fb^\eps-\fb\Vert_{L^2([0,T];\mathbb H^\delta(\dom))}
\leq \sup_{|s|<\eps}\Vert\fb(\cdot-s)-\fb(\cdot)\Vert_{L^2([0,T];\mathbb H^\delta(\dom))}
\longrightarrow0.
\]
Thus, \eqref{eq-uniform convergence3} holds.

Moreover, by Minkowski's inequality,
\begin{align*}
\Vert\fb^\eps\Vert_{L^2([0,T];\mathbb H^\delta(\dom))}
&\leq
\int_{B(0,\eps)}
\eta_\eps(s)
\Vert\fb(\cdot-s)\Vert_{L^2([0,T];\mathbb H^\delta(\dom))}
\,ds
\leq
c
\Vert\fb\Vert_{L^2([0,T+1];\mathbb H^\delta(\dom))},
\end{align*}
which implies \eqref{eq-f-eps-uniform} holds. Moreover, for $0<\eps<1$,
\[
\Vert(\alpha^\eps)^{-1}\Vert_{L^\infty([0,T])}
\leq \Vert\alpha^{-1}\Vert_{L^\infty([0,T+1])}.
\]

Recall that
\begin{align*}
&\wb_m^\eps(t)
=
\int_0^t
e^{(t-s)\Delta}
(\alpha^\eps)^{-1}(s)\fb_m^\eps(s)\,ds,\\
&\wb(t)
=
\int_0^t
e^{(t-s)\Delta}
\alpha^{-1}(s)\fb(s)\,ds.
\end{align*}
Therefore,
\begin{align*}
\wb_m^\eps(t)-\wb(t)
&=
\int_0^t
e^{(t-s)\Delta}
\left(
(\alpha^\eps)^{-1}(s)\fb_m^\eps(s)
-
\alpha^{-1}(s)\fb(s)
\right)\,ds.
\end{align*}
By \autoref{lemma-lplq}, we infer that
\begin{align}\label{eq-wm-w-estimate}
\sup_{t\in[0,T]}
\Vert\wb_m^\eps(t)-\wb(t)\Vert_{H^2}
&\leq
c
\Vert
(\alpha^\eps)^{-1}\fb_m^\eps-\alpha^{-1}\fb
\Vert_{L^2([0,T];\mathbb H^1(\dom))}.
\end{align}

Since $\fb_m^\eps=P_m\fb^\eps$, we infer that
\begin{align}\label{eq-force-joint-convergence}
\Vert
(\alpha^\eps)^{-1}\fb_m^\eps-\alpha^{-1}\fb
\Vert_{L^2([0,T];\mathbb H^1(\dom))}
&\leq
\Vert(\alpha^\eps)^{-1}\Vert_{L^\infty([0,T])}
\Vert
P_m(\fb^\eps-\fb)
\Vert_{L^2([0,T];\mathbb H^1(\dom))}
\nonumber\\
&\quad+
\Vert(\alpha^\eps)^{-1}\Vert_{L^\infty([0,T])}
\Vert
P_m\fb-\fb
\Vert_{L^2([0,T];\mathbb H^1(\dom))}
\nonumber\\
&\quad+
\Vert
(\alpha^\eps)^{-1}-\alpha^{-1}
\Vert_{L^\infty([0,T])}
\Vert\fb\Vert_{L^2([0,T];\mathbb H^1(\dom))}.
\end{align}

Note that
\[
\Vert
P_m(\fb^\eps-\fb)
\Vert_{L^2([0,T];\mathbb H^1(\dom))}
\leq
\Vert
\fb^\eps-\fb
\Vert_{L^2([0,T];\mathbb H^1(\dom))}.
\]
Recall that $\delta>\frac32$, by \eqref{eq-uniform convergence3} we infer that
\[
\Vert
\fb^\eps-\fb
\Vert_{L^2([0,T];\mathbb H^1(\dom))}
\to0,
\qquad\text{as }\eps\to0.
\]

Moreover, since $\fb\in L^2([0,T];\mathbb H^\delta(\dom))\subset L^2([0,T];\mathbb H^1(\dom))$, we infer that
\[
\Vert
P_m\fb-\fb
\Vert_{L^2([0,T];\mathbb H^1(\dom))}
\to0,
\qquad\text{as }m\to\infty.
\]

In view of \eqref{eq-uniform convergence1},
\[
\Vert
(\alpha^\eps)^{-1}-\alpha^{-1}
\Vert_{L^\infty([0,T])}
\to0,
\qquad\text{as }\eps\to0,
\]
which implies there exists $\eps_0>0$ such that
\[
\sup_{\eps\in(0,\eps_0)}
\Vert(\alpha^\eps)^{-1}\Vert_{L^\infty([0,T])}
<\infty.
\]
Combining these estimates with
\eqref{eq-force-joint-convergence}, we infer that
\[
\Vert
(\alpha^\eps)^{-1}\fb_m^\eps-\alpha^{-1}\fb
\Vert_{L^2([0,T];\mathbb H^1(\dom))}
\to0, \;\;\mbox{ as }\eps\to0 \mbox{ and }m\to\infty.
\]

Therefore, by \eqref{eq-wm-w-estimate},
\[
\sup_{t\in[0,T]}
\Vert\wb_m^\eps(t)-\wb(t)\Vert_{H^2}
\to0, \;\;\mbox{ as }\eps\to0 \mbox{ and }m\to\infty,
\]
 which implies \eqref{eq-uniform convergence2} holds.
This completes the proof of \autoref{lem-uniform convergence}.
\end{proof}

\section{Auxiliary proofs for local existence}\label{sec-local-technical}

\subsection{Second-moment estimates}\label{app-second-moment-estimate}
\begin{proof}[Proof of \autoref{lem-moments of rvs}]
Let us first prove \eqref{eqn-boundedness in L^p}. Let $\{\mathcal K_N\}_{N\in\mathbb N}$ be an increasing sequence of finite subsets of $\mathcal K$ such that
\[
\bigcup_{N\in\mathbb N}\mathcal K_N=\mathcal K.
\]
For every $N\in\mathbb N$, define
\[
S_N:=\sum_{(\kb,j)\in\mathcal K_N}c_{\kb}^jr_{\kb}^j.
\]
By the independence and zero-mean property of $\{r_{\kb}^j\}_{(\kb,j)\in\mathcal K}$, for $(\kb,j)\neq(\ell,i)$, we have
\[
\mathbb E[r_{\kb}^jr_\ell^i]=\mathbb E[r_{\kb}^j]\mathbb E[r_\ell^i]=0.
\]
Therefore, we infer that
\begin{align}\label{eq-scalar-finite-second-moment}
\mathbb E|S_N|^2
&=\sum_{(\kb,j),(\ell,i)\in\mathcal K_N}c_{\kb}^j\overline{c_\ell^i}\mathbb E[r_{\kb}^jr_\ell^i]=\sum_{(\kb,j)\in\mathcal K_N}|c_{\kb}^j|^2\mathbb E|r_{\kb}^j|^2
\leq M_2\sum_{(\kb,j)\in\mathcal K_N}|c_{\kb}^j|^2.
\end{align}
For $M>N$, applying \eqref{eq-scalar-finite-second-moment} to the finite set $\mathcal K_M\setminus\mathcal K_N$, we infer that
\begin{align}\label{eq-scalar-Cauchy-second-moment}
\Vert S_M-S_N\Vert_{L^2(\Omega)}^2
&=\mathbb E\left|\sum_{(\kb,j)\in\mathcal K_M\setminus\mathcal K_N}c_{\kb}^jr_{\kb}^j\right|^2\leq M_2\sum_{(\kb,j)\in\mathcal K_M\setminus\mathcal K_N}|c_{\kb}^j|^2.
\end{align}
Since $\{c_{\kb}^j\}_{(\kb,j)\in\mathcal K}\in\ell^2(\mathcal K;\mathbb C)$ and $\mathcal K_N\uparrow\mathcal K$, we infer that
\[
\sum_{(\kb,j)\in\mathcal K\setminus\mathcal K_N}|c_{\kb}^j|^2\to0,\qquad N\to\infty.
\]
Therefore,
\[
\lim_{N\to\infty}\sup_{M>N}\Vert S_M-S_N\Vert_{L^2(\Omega)}=0.
\]
which implies $\{S_N\}_{N\in\mathbb N}$ is a Cauchy sequence in $L^2(\Omega)$.

Since $L^2(\Omega)$ is complete, there exists $S\in L^2(\Omega)$ such that
\[
\lim_{N\to\infty}\mathbb E|S_N-S|^2=0.
\]
Hence, the series $\sum_{(\kb,j)\in\mathcal K}c_{\kb}^jr_{\kb}^j$ converges to $S$ in $L^2(\Omega)$.
Letting $N\to\infty$ in \eqref{eq-scalar-finite-second-moment}, we obtain
\[
\Big\Vert\sum_{(\kb,j)\in\mathcal K}c_{\kb}^jr_{\kb}^j\Big\Vert_{L^2(\Omega)}^2
\leq M_2\sum_{(\kb,j)\in\mathcal K}|c_{\kb}^j|^2,
\]
which implies \eqref{eqn-boundedness in L^p} holds.

In the following, we prove \eqref{eqn-vector-valued-second-moment}. Let us choose and fix an arbitrary finite subset
\[
\mathcal K_N=\{(\kb_\ell,j_\ell):\ell=1,\ldots,N\}
\subset\mathcal K.
\]
For simplicity, we denote
\[
r_\ell:=r_{\kb_\ell}^{j_\ell},
\qquad
F_\ell:=F_{\kb_\ell}^{j_\ell},
\qquad \ell=1,\ldots,N.
\]

Let $\{r_\ell'\}_{\ell=1}^N$ be an independent copy of $\{r_\ell\}_{\ell=1}^N$, and let $\{\varepsilon_\ell\}_{\ell=1}^N$ be a sequence of independent symmetric Bernoulli random variables, satisfying
\[
\mathbb P(\varepsilon_\ell=1)
=
\mathbb P(\varepsilon_\ell=-1)
=
\frac12,\qquad \ell=1,\ldots,N.
\]
Moreover, $\{\varepsilon_\ell\}_{\ell=1}^N$ are independent of both $\{r_\ell\}_{\ell=1}^N$ and $\{r_\ell'\}_{\ell=1}^N$.
We denote by $\mathbb E'$ and $\mathbb E_\varepsilon$ the expectations with respect to $\{r_\ell'\}_{\ell=1}^N$ and $\{\varepsilon_\ell\}_{\ell=1}^N$, respectively.

Since $\mathbb E'r_\ell'=0$, by Jensen's inequality, we infer that
\begin{align}\label{eq-symmetrization-second-moment}
\mathbb E\Big\Vert\sum_{\ell=1}^Nr_\ell F_\ell\Big\Vert_E^2
&=\mathbb E\Big\Vert\mathbb E'\Big[\sum_{\ell=1}^N(r_\ell-r_\ell')F_\ell\Big]\Big\Vert_E^2\leq\mathbb E\mathbb E'\Big\Vert\sum_{\ell=1}^N(r_\ell-r_\ell')F_\ell\Big\Vert_E^2.
\end{align}
Since $\{r_\ell-r_\ell'\}_{\ell=1}^N$ is a sequence of independent symmetric random variables,

by the triangle inequality and the fact that $\{r_\ell\}_{\ell=1}^N$ and $\{r_\ell'\}_{\ell=1}^N$ have the same joint distribution, we infer that
\begin{align}\label{eq-symmetrization-bernoulli}
\mathbb E\mathbb E'\Big\Vert\sum_{\ell=1}^N(r_\ell-r_\ell')F_\ell\Big\Vert_E^2
&=\mathbb E\mathbb E'\mathbb E_\varepsilon\Big\Vert\sum_{\ell=1}^N\varepsilon_\ell(r_\ell-r_\ell')F_\ell\Big\Vert_E^2\nonumber\\
&\leq 2\mathbb E\mathbb E_\varepsilon\Big\Vert\sum_{\ell=1}^N\varepsilon_\ell r_\ell F_\ell\Big\Vert_E^2
+2\mathbb E'\mathbb E_\varepsilon\Big\Vert\sum_{\ell=1}^N\varepsilon_\ell r_\ell'F_\ell\Big\Vert_E^2\nonumber\\
&=4\mathbb E\mathbb E_\varepsilon\Big\Vert\sum_{\ell=1}^N\varepsilon_\ell r_\ell F_\ell\Big\Vert_E^2.
\end{align}
Therefore, inequalities \eqref{eq-symmetrization-bernoulli} and  \eqref{eq-symmetrization-second-moment} imply
\begin{align}\label{eq-symmetrization-combined}
\mathbb E\Big\Vert\sum_{\ell=1}^Nr_\ell F_\ell\Big\Vert_E^2
\leq 4\mathbb E\mathbb E_\varepsilon\Big\Vert\sum_{\ell=1}^N\varepsilon_\ell r_\ell F_\ell\Big\Vert_E^2.
\end{align}

Let $E=L^q([0,T];\mathbb L^p(\dom))$, where $(q,p)=(4,6)$ or $(q,p)=(2,3)$. Let us choose and fix an arbitrary sequence of real numbers $\{a_\ell\}_{\ell=1}^N$ and functions $\{F_\ell\}_{\ell=1}^N\subset E$. In the following, we aim to prove that there exists a constant $c>0$, independent of $N$, such that
\begin{align}\label{eq-rademacher-second-moment}
\Big(\mathbb E_\varepsilon\Big\Vert\sum_{\ell=1}^N\varepsilon_\ell a_\ell F_\ell\Big\Vert_E^2\Big)^{\frac12}
\leq c\Big(\sum_{\ell=1}^N|a_\ell|^2\Vert F_\ell\Vert_E^2\Big)^{\frac12}.
\end{align}

Firstly, by the H\"older inequality on the probability space, we infer that
\begin{align}\label{eq-rademacher-holder-2}
\mathbb E_\varepsilon\Big\Vert\sum_{\ell=1}^N\varepsilon_\ell a_\ell F_\ell\Big\Vert_E^2
&\leq\Big(\mathbb E_\varepsilon\Big\Vert\sum_{\ell=1}^N\varepsilon_\ell a_\ell F_\ell\Big\Vert_E^6\Big)^{\frac13}
\Big(\mathbb E_\varepsilon 1^{\frac32}\Big)^{\frac23}=\Big(\mathbb E_\varepsilon\Big\Vert\sum_{\ell=1}^N\varepsilon_\ell a_\ell F_\ell\Big\Vert_E^6\Big)^{\frac13}.
\end{align}

Since $q,p\leq6$, by the Minkowski inequality \cite[Lemma 3.3.1]{Kwapien+Woycznski_1992}, we infer that
\begin{align}\label{eq-rademacher-minkowski-1}
\Big(\mathbb E_\varepsilon\Big\Vert\sum_{\ell=1}^N\varepsilon_\ell a_\ell F_\ell\Big\Vert_E^6\Big)^{\frac16}&=\Big\Vert\sum_{\ell=1}^N\varepsilon_\ell a_\ell F_\ell\Big\Vert_{L^6(\Omega_\varepsilon;L^q([0,T];\mathbb L^p(\dom)))}\leq\Big\Vert\sum_{\ell=1}^N\varepsilon_\ell a_\ell F_\ell\Big\Vert_{L^q([0,T];L^6(\Omega_\varepsilon;\mathbb L^p(\dom)))}\nonumber\\
&\leq\Big\Vert\sum_{\ell=1}^N\varepsilon_\ell a_\ell F_\ell\Big\Vert_{L^q([0,T];\mathbb L^p(\dom;L^6(\Omega_\varepsilon)))}=\Big\Vert\Big(\mathbb E_\varepsilon\Big|\sum_{\ell=1}^N\varepsilon_\ell a_\ell F_\ell\Big|^6\Big)^{\frac16}\Big\Vert_E.
\end{align}

For almost every $(t,x)\in[0,T]\times\dom$, by the Khintchine inequality for independent symmetric Bernoulli random variables, see \cite[Theorem 1.1]{Astashkin_2024}, we infer that there exists a constant $c>0$, independent of $N$, such that
\begin{align}\label{eq-rademacher-khintchine}
&\Big(\mathbb E_\varepsilon\Big|\sum_{\ell=1}^N\varepsilon_\ell a_\ell F_\ell(t,x)\Big|^6\Big)^{\frac16}\leq c\Big(\sum_{\ell=1}^N|a_\ell F_\ell(t,x)|^2\Big)^{\frac12}.
\end{align}
The same inequality holds for finite-dimensional vector-valued and matrix-valued functions by applying the scalar Khintchine inequality to their components. Therefore, taking the $E$-norm on both sides of \eqref{eq-rademacher-khintchine}, we infer that
\begin{align}\label{eq-rademacher-khintchine-2}
\Big\Vert\Big(\mathbb E_\varepsilon\Big|\sum_{\ell=1}^N\varepsilon_\ell a_\ell F_\ell\Big|^6\Big)^{\frac16}\Big\Vert_E
\leq c\Big\Vert\Big(\sum_{\ell=1}^N|a_\ell F_\ell|^2\Big)^{\frac12}\Big\Vert_E.
\end{align}
Moreover, since $q,p\geq2$, we have $\frac q2,\frac p2\geq1$. Hence, by the triangle inequality in $L^{\frac q2}([0,T];L^{\frac p2}(\dom))$, we infer that
\begin{align}\label{eq-rademacher-minkowski-3}
&\Big\Vert\Big(\sum_{\ell=1}^N|a_\ell F_\ell|^2\Big)^{\frac12}\Big\Vert_E^2=\Big\Vert\sum_{\ell=1}^N|a_\ell F_\ell|^2\Big\Vert_{L^{\frac q2}([0,T];L^{\frac p2}(\dom))}\leq\sum_{\ell=1}^N|a_\ell|^2\Big\Vert|F_\ell|^2\Big\Vert_{L^{\frac q2}([0,T];L^{\frac p2}(\dom))}=\sum_{\ell=1}^N|a_\ell|^2\Vert F_\ell\Vert_E^2.
\end{align}

Combining \eqref{eq-rademacher-holder-2}, \eqref{eq-rademacher-minkowski-1}, \eqref{eq-rademacher-khintchine-2} and \eqref{eq-rademacher-minkowski-3}, we infer that
\begin{align*}
\Big(\mathbb E_\varepsilon\Big\Vert\sum_{\ell=1}^N\varepsilon_\ell a_\ell F_\ell\Big\Vert_E^2\Big)^{\frac12}
&\leq\Big(\mathbb E_\varepsilon\Big\Vert\sum_{\ell=1}^N\varepsilon_\ell a_\ell F_\ell\Big\Vert_E^6\Big)^{\frac16}\leq\Big\Vert\Big(\mathbb E_\varepsilon\Big|\sum_{\ell=1}^N\varepsilon_\ell a_\ell F_\ell\Big|^6\Big)^{\frac16}\Big\Vert_E\nonumber\\
&\leq c\Big\Vert\Big(\sum_{\ell=1}^N|a_\ell F_\ell|^2\Big)^{\frac12}\Big\Vert_E\leq c\Big(\sum_{\ell=1}^N|a_\ell|^2\Vert F_\ell\Vert_E^2\Big)^{\frac12},
\end{align*}
which implies \eqref{eq-rademacher-second-moment} holds.

For $E=\mathbb L^6(\dom)$, the same argument without the time variable gives
\begin{align*}
\Big(\mathbb E_\varepsilon\Big\Vert\sum_{\ell=1}^N\varepsilon_\ell a_\ell F_\ell\Big\Vert_{L^6}^2\Big)^{\frac12}
&\leq\Big(\mathbb E_\varepsilon\Big\Vert\sum_{\ell=1}^N\varepsilon_\ell a_\ell F_\ell\Big\Vert_{L^6}^6\Big)^{\frac16}
=\Big\Vert\Big(\mathbb E_\varepsilon\Big|\sum_{\ell=1}^N\varepsilon_\ell a_\ell F_\ell\Big|^6\Big)^{\frac16}\Big\Vert_{L^6}\\
&\leq c\Big\Vert\Big(\sum_{\ell=1}^N|a_\ell F_\ell|^2\Big)^{\frac12}\Big\Vert_{L^6}
\leq c\Big(\sum_{\ell=1}^N|a_\ell|^2\Vert F_\ell\Vert_{L^6}^2\Big)^{\frac12},
\end{align*}
where the last inequality follows from the triangle inequality in $L^3(\dom)$. Hence, \eqref{eq-rademacher-second-moment} holds for all three choices of $E$ in \autoref{lem-moments of rvs}.

Applying \eqref{eq-rademacher-second-moment} with $a_\ell=r_\ell$, by \eqref{eq-symmetrization-combined}, we infer that
\begin{align}\label{eq-vector-finite-second-moment-bound}
\mathbb E\Big\Vert\sum_{\ell=1}^Nr_\ell F_\ell\Big\Vert_E^2
&\leq c\mathbb E\sum_{\ell=1}^N|r_\ell|^2\Vert F_\ell\Vert_E^2\leq cM_2\sum_{\ell=1}^N\Vert F_\ell\Vert_E^2.
\end{align}

Finally, let $\{\mathcal K_N\}_{N\in\mathbb N}$ be an increasing sequence of finite subsets of $\mathcal K$ such that $\bigcup_{N\in\mathbb N}\mathcal K_N=\mathcal K$, and define
\[
X_N:=\sum_{(\kb,j)\in\mathcal K_N}r_{\kb}^jF_{\kb}^j.
\]
For $M>N$, by \eqref{eq-vector-finite-second-moment-bound}, we infer that
\begin{align}
\Vert X_M-X_N\Vert_{L^2(\Omega;E)}^2
&\leq cM_2\sum_{(\kb,j)\in\mathcal K_M\setminus\mathcal K_N}\Vert F_{\kb}^j\Vert_E^2\leq cM_2\sum_{(\kb,j)\in\mathcal K\setminus\mathcal K_N}\Vert F_{\kb}^j\Vert_E^2\to0,\qquad N\to\infty.
\end{align}
Since $L^2(\Omega;E)$ is complete, the sequence $\{X_N\}_{N\in\mathbb N}$ converges in $L^2(\Omega;E)$.
Letting $N\to\infty$ in \eqref{eq-vector-finite-second-moment-bound}, we infer that
\begin{align}
\Big\Vert\sum_{(\kb,j)\in\mathcal K}r_{\kb}^jF_{\kb}^j\Big\Vert_{L^2(\Omega;E)}
\leq cM_2^{\frac12}\Big(\sum_{(\kb,j)\in\mathcal K}\Vert F_{\kb}^j\Vert_E^2\Big)^{\frac12},
\end{align}
which implies \eqref{eqn-vector-valued-second-moment} holds. The proof is complete.

\end{proof}

\subsection{Estimates for the truncated map}\label{app-truncated-map-estimates}
Let $0<T\leq1$, $N,R>0$, and let $\alpha$ and $\mathbf z$ be the coefficient and linear part used in the proof of \autoref{thm-local l2}.
We use the norms in \eqref{eq-Xzt} and a cut-off $\theta_R$ satisfying \eqref{eqn-theta-differential} and \eqref{eqn-truncated-vz}.
For an arbitrary $\vb_R^N\in\mathbb M^2(X_T)$, set
\[
\varepsilon_t:=\Vert\vb_R^N\Vert_{X_t}+\Vert\mathbf z\Vert_{Y_t},
\qquad 0\leq t\leq T,
\]
and define
\[
L_R^N(\vb_R^N)(t):=-\int_0^t\theta_R^2(\varepsilon_s)\alpha(s)
\mathds1_{\{\Vert\alpha\Vert_{L^\infty([0,s])}\leq N\}}
e^{(t-s)\Delta}\big((\vb_R^N+\mathbf z)\cdot\nabla\big)(\vb_R^N+\mathbf z)(s)\,ds.
\]
We verify \eqref{eq-map-1} and \eqref{eq-difference-1}.
All maximal-regularity constants below are taken on $[0,1]$ and are therefore uniform for $T\leq1$.

In the following, we write $\mathbf{z}$ instead of $\mathbf{z}(t, \overrightarrow{\ub_{0}})$ and $\vb_R^N$ instead of $ \vb_R^N(t,\xb)$ for short.
Then inequalities \eqref{eqn-g nablah} and \eqref{eqn-truncated-vz} imply there exists a constant $c>0$ such that

\begin{equation}\label{eqn-v nablaz}
\begin{aligned}
\Vert\theta_R^2(\varepsilon_s)\left( \vb_R^N\cdot \nabla\right) \vb_R^N  \Vert_{L^{\frac{4}{3}}([0,T ]; \mathbb{L}^{2}(\dom))}&\leq\Vert \Big(\theta_R(\varepsilon_s)\vb_R^N\cdot \nabla \Big)\theta_R(\varepsilon_s)\vb_R^N  \Vert_{L^{\frac{4}{3}}([0,T]; \mathbb{L}^{2}(\dom))}\\
&\leq c\Vert \theta_R(\varepsilon_s)\vb_R^N\Vert_{L^{4}([0,T]; \mathbb{L}^{6}(\dom))}\Vert \nabla \vb_R^N\Vert_{L^{2}([0,T]; \mathbb{L}^{3}(\dom))}\leq 2cR\Vert \vb_R^N\Vert_{X_T},\\
\Vert \theta_R^2(\varepsilon_s)\Big(\vb_R^N\cdot \nabla\Big) \mathbf{z}  \Vert_{L^{\frac{4}{3}}([0,T]; \mathbb{L}^{2}(\dom))}&\leq \Vert \Big(\theta_R(\varepsilon_s)\vb_R^N\cdot \nabla\Big) \theta_R(\varepsilon_s)\mathbf{z}  \Vert_{L^{\frac{4}{3}}([0,T]; \mathbb{L}^{2}(\dom))}\leq 2cR \Vert \vb_R^N\Vert_{X_T},\\
\Vert \theta_R^2(\varepsilon_s)\Big(\mathbf{z}\cdot \nabla\Big) \vb_R^N  \Vert_{L^{\frac{4}{3}}([0,T]; \mathbb{L}^{2}(\dom))}&\leq \Vert \Big(\theta_R(\varepsilon_s)\mathbf{z}\cdot \nabla\Big) \theta_R(\varepsilon_s)\vb_R^N  \Vert_{L^{\frac{4}{3}}([0,T]; \mathbb{L}^{2}(\dom))}\leq 2cR\Vert \vb_R^N\Vert_{X_T},\\
\Vert \theta_R^2(\varepsilon_s)\Big(\mathbf{z}\cdot \nabla\Big) \mathbf{z}  \Vert_{L^{\frac{4}{3}}([0,T ]; \mathbb{L}^{2}(\dom))}&\leq \Vert \Big(\theta_R(\varepsilon_s)\mathbf{z}\cdot \nabla\Big) \theta_R(\varepsilon_s)\mathbf{z}  \Vert_{L^{\frac{4}{3}}([0,T]; \mathbb{L}^{2}(\dom))}\leq 4cR^{2}.
\end{aligned}
\end{equation}

By \autoref{lemma-lplq}, we conclude from \eqref{eqn-v nablaz} and \eqref{eqn-truncated-vz} that for every $\vb_R^N\in X_{T}$, there exists a constant $c_1>0$ independent of $N$ and $R$, such that
\begin{align}\label{eq-map}
\Vert L_R^N(\vb_R^N)\Vert_{X_{T }}
&\leq cN\Vert \theta_R^2(\varepsilon_s)\Big((\vb_R^N+\mathbf{z})\cdot \nabla\Big) (\vb_R^N+\mathbf{z})\Vert_{L^{\frac{4}{3}}([0,T ]; \mathbb{L}^{2}(\dom))}
\nonumber\\
&\leq cN\Vert\theta_R^2(\varepsilon_s)\Big( \vb_R^N\cdot \nabla\Big) \vb_R^N  \Vert_{L^{\frac{4}{3}}([0,T ]; \mathbb{L}^{2}(\dom))}+cN\Vert \theta_R^2(\varepsilon_s)\Big(\mathbf{z}\cdot \nabla\Big) \vb_R^N  \Vert_{L^{\frac{4}{3}}([0,T ]; \mathbb{L}^{2}(\dom))}\nonumber \\
&+cN\Vert \theta_R^2(\varepsilon_s)\Big(\vb_R^N\cdot \nabla\Big) \mathbf{z}  \Vert_{L^{\frac{4}{3}}([0,T ]; \mathbb{L}^{2}(\dom))}+cN\Vert \theta_R^2(\varepsilon_s)\Big(\mathbf{z}\cdot \nabla\Big) \mathbf{z}  \Vert_{L^{\frac{4}{3}}([0,T ]; \mathbb{L}^{2}(\dom))}\nonumber\\
&\leq c_1N(R^2+ R\Vert \vb_R^N\Vert_{X_{T }}).
\end{align}

For $\vb_R^{N,1} , \vb_R^{N,2}\in \mathbb{M}^2\left(X_{T}\right)$, denote $\varepsilon_t^1:=\Vert \vb_R^{N,1}\Vert_{X_t} +\Vert \mathbf z\Vert_{Y_t}$ and $\varepsilon_t^2:=\Vert \vb_R^{N,2}\Vert_{X_t} +\Vert \mathbf z\Vert_{Y_t}$. In view of \autoref{lemma-lplq}, by inequalities \eqref{eqn-g nablah} and \eqref{eqn-truncated-vz}, together with the Cauchy-Schwarz inequality, we infer that there exists a constant $c>0$ such that
\begin{align}\label{eq-difference}
&\Vert L_R^N(\vb_R^{N,1})- L_R^N(\vb_R^{N,2}) \Vert_{X_{T }}\\
&\leq cN\Vert \Big(\theta_R(\varepsilon_t^1)(\vb_R^{N,1}+\mathbf{z})\cdot \nabla\Big) \theta_R(\varepsilon_t^1)(\vb_R^{N,1}+\mathbf{z})-\Big(\theta_R(\varepsilon_t^2)(\vb_R^{N,2}+\mathbf{z})\cdot\nabla \Big)\theta_R(\varepsilon_t^2) (\vb_R^{N,2}+\mathbf{z})\Vert_{L^{\frac{4}{3}}([0,T ]; \mathbb{L}^{2}(\dom))}\nonumber\\
&\leq cN\Vert \Big((\theta_R(\varepsilon_t^1)(\vb_R^{N,1}+\mathbf{z})-\theta_R(\varepsilon_t^2)(\vb_R^{N,2}+\mathbf{z}))\cdot \nabla\Big) \theta_R(\varepsilon_t^1)(\vb_R^{N,1}+\mathbf{z})\nonumber\\
&\;\;
+\Big(\theta_R(\varepsilon_t^2)(\vb_R^{N,2}+\mathbf{z})\cdot \nabla\Big)(\theta_R(\varepsilon_t^1) (\vb_R^{N,1}+\mathbf{z})-\theta_R(\varepsilon_t^2) (\vb_R^{N,2}+\mathbf{z}))\Vert_{L^{\frac{4}{3}}([0,T ]; \mathbb{L}^{2}(\dom))}\nonumber\\
&\leq cN\Vert\Big( (\theta_R(\varepsilon_t^1)\vb_R^{N,1}-\theta_R(\varepsilon_t^2)\vb_R^{N,2})\cdot \nabla\Big) \theta_R(\varepsilon_t^1)\vb_R^{N,1}\Vert_{L^{\frac{4}{3}}([0,T ];  \mathbb{L}^{2}(\dom))}\nonumber\\
&+cN\Vert \Big((\theta_R(\varepsilon_t^1)\vb_R^{N,1}-\theta_R(\varepsilon_t^2)\vb_R^{N,2})\cdot \nabla\Big) \theta_R(\varepsilon_t^1)\mathbf{z}\Vert_{L^{\frac{4}{3}}([0,T ];  \mathbb{L}^{2}(\dom))}\nonumber\\
&+cN\Vert \Big((\theta_R(\varepsilon_t^1)-\theta_R(\varepsilon_t^2))\mathbf{z}\cdot \nabla\Big)\theta_R(\varepsilon_t^1)\vb_R^{N,1}\Vert_{L^{\frac{4}{3}}([0,T];  \mathbb{L}^{2}(\dom))}+cN\Vert \Big((\theta_R(\varepsilon_t^1)-\theta_R(\varepsilon_t^2))\mathbf{z}\cdot \nabla\Big)\theta_R(\varepsilon_t^1)\mathbf{z}\Vert_{L^{\frac{4}{3}}([0,T];  \mathbb{L}^{2}(\dom))}\nonumber\\
&+cN\Vert \Big(\theta_R(\varepsilon_t^2)\vb_R^{N,2}\cdot \nabla\Big)\Big(\theta_R(\varepsilon_t^1)\vb_R^{N,1}-\theta_R(\varepsilon_t^2)\vb_R^{N,2}\Big)\Vert_{L^{\frac{4}{3}}([0,T];  \mathbb{L}^{2}(\dom))}\nonumber\\
&+cN\Vert \Big(\theta_R(\varepsilon_t^2)\mathbf{z}\cdot \nabla\Big)\Big(\theta_R(\varepsilon_t^1)\vb_R^{N,1}-\theta_R(\varepsilon_t^2)\vb_R^{N,2}\Big)\Vert_{L^{\frac{4}{3}}([0,T];  \mathbb{L}^{2}(\dom))}\nonumber\\
&+cN\Vert \Big(\theta_R(\varepsilon_t^2)\vb_R^{N,2}\cdot \nabla\Big)\Big((\theta_R(\varepsilon_t^1)-\theta_R(\varepsilon_t^2))\mathbf{z}\Big)\Vert_{L^{\frac{4}{3}}([0,T ];  \mathbb{L}^{2}(\dom))}+cN\Vert \Big(\theta_R(\varepsilon_t^2)\mathbf{z}\cdot \nabla\Big)\Big((\theta_R(\varepsilon_t^1)-\theta_R(\varepsilon_t^2))\mathbf{z}\Big)\Vert_{L^{\frac{4}{3}}([0,T];  \mathbb{L}^{2}(\dom))}\nonumber\\
&\leq cN\Vert\theta_R(\varepsilon_t^1)\vb_R^{N,1}-\theta_R(\varepsilon_t^2)\vb_R^{N,2}\Vert_{X_T}\Vert \theta_R(\varepsilon_t^1)\vb_R^{N,1}\Vert_{X_T}+cN\Vert\theta_R(\varepsilon_t^1)\vb_R^{N,1}-\theta_R(\varepsilon_t^2)\vb_R^{N,2}\Vert_{X_T}\Vert \theta_R(\varepsilon_t^1)\mathbf{z}\Vert_{Y_T}\\
&+cN\Vert (\theta_R(\varepsilon_t^1)-\theta_R(\varepsilon_t^2))\mathbf{z}\Vert_{Y_T}\Vert\theta_R(\varepsilon_t^1)\vb_R^{N,1}\Vert_{X_T}+cN\Vert (\theta_R(\varepsilon_t^1)-\theta_R(\varepsilon_t^2))\mathbf{z}\Vert_{Y_T}\Vert\theta_R(\varepsilon_t^1)\mathbf{z}\Vert_{Y_T}\nonumber\\
&+cN\Vert \theta_R(\varepsilon_t^2)\vb_R^{N,2}\Vert_{X_T}\Vert\theta_R(\varepsilon_t^1)\vb_R^{N,1}-\theta_R(\varepsilon_t^2)\vb_R^{N,2}\Vert_{X_T}+cN\Vert \theta_R(\varepsilon_t^2)\mathbf{z}\Vert_{Y_T}\Vert\theta_R(\varepsilon_t^1)\vb_R^{N,1}-\theta_R(\varepsilon_t^2)\vb_R^{N,2}\Vert_{X_T}\nonumber\\
&+cN\Vert\theta_R(\varepsilon_t^2)\vb_R^{N,2}\Vert_{X_T}\Vert ((\theta_R(\varepsilon_t^1)-\theta_R(\varepsilon_t^2))\mathbf{z})\Vert_{Y_T}+cN\Vert\theta_R(\varepsilon_t^2)\mathbf{z}\Vert_{Y_T}\Vert ((\theta_R(\varepsilon_t^1)-\theta_R(\varepsilon_t^2))\mathbf{z})\Vert_{Y_T}.
\end{align}

Let us denote
\[
D:=\Vert\vb_R^{N,1}-\vb_R^{N,2}\Vert_{X_T}.
\]
By \eqref{eqn-theta-differential}, for every $t\in[0,T]$,
we infer that
\begin{align*}
&|\theta_R(\varepsilon_t^1)
-\theta_R(\varepsilon_t^2)|\leq \frac{2}{R}
|\varepsilon_t^1-\varepsilon_t^2|\leq \frac{2}{R}
\Vert\vb_R^{N,1}-\vb_R^{N,2}\Vert_{X_t}
\leq \frac{2D}{R}.
\end{align*}
Since $0\leq\theta_R\leq1$, we deduce that
\begin{align}\label{eq-theta-difference-bound}
\sup_{t\in[0,T]}
|\theta_R(\varepsilon_t^1)
-\theta_R(\varepsilon_t^2)|
\leq \min\left\{1,\frac{2D}{R}\right\}.
\end{align}

Moreover, since $\varepsilon_t^i$ is nondecreasing
and $\theta_R(r)=0$ for $r\geq2R$, let us set
\[
\rho_i:=
\inf\{t\in[0,T]:\varepsilon_t^i\geq2R\}
\wedge T,\qquad i=1,2,
\]
where $\inf\varnothing=\infty$, and
\[
\rho:=\rho_1\vee\rho_2.
\]
If $\rho=0$, then both $\theta_R(\varepsilon_t^1)$ and $\theta_R(\varepsilon_t^2)$ vanish for every $t\in[0,T]$, so the two cut-off difference bounds proved below hold trivially. We prove those bounds for $\rho>0$.
Then, for every $t\in[0,\rho)$, there exists
$i\in\{1,2\}$ such that
$\varepsilon_t^i<2R$. Hence,
\[
\Vert\mathbf z\Vert_{Y_t}<2R.
\]
On the other hand, if $\varepsilon_t^2<2R$,
then
\[
\Vert\vb_R^{N,2}\Vert_{X_t}<2R.
\]
Otherwise, we have $\varepsilon_t^1<2R$,
and therefore
\begin{align*}
\Vert\vb_R^{N,2}\Vert_{X_t}
&\leq
\Vert\vb_R^{N,1}\Vert_{X_t}
+\Vert\vb_R^{N,1}-\vb_R^{N,2}\Vert_{X_t}\leq 2R+D.
\end{align*}
Consequently, by the continuity of the running
norms, we infer that
\begin{align}\label{eq-running-norm-bounds}
\Vert\mathbf z\Vert_{Y_\rho}
&\leq 2R,\nonumber\\
\Vert\vb_R^{N,2}\Vert_{X_\rho}
&\leq 2R+D.
\end{align}

Since $\theta_R(\varepsilon_t^1)
=\theta_R(\varepsilon_t^2)=0$ for $t>\rho$,
by \eqref{eq-theta-difference-bound} and
\eqref{eq-running-norm-bounds}, we deduce that
\begin{align}
&\Vert
(\theta_R(\varepsilon_\cdot^1)
-\theta_R(\varepsilon_\cdot^2))
\mathbf z\Vert_{Y_T}\leq
\sup_{t\in[0,T]}
|\theta_R(\varepsilon_t^1)
-\theta_R(\varepsilon_t^2)|
\Vert\mathbf z\Vert_{Y_\rho}\leq
2R\min\left\{1,\frac{2D}{R}\right\}
\leq 4D.
\end{align}

Similarly, note that
\begin{align*}
&\theta_R(\varepsilon_t^1)\vb_R^{N,1}(t)
-\theta_R(\varepsilon_t^2)\vb_R^{N,2}(t)=\theta_R(\varepsilon_t^1)
(\vb_R^{N,1}(t)-\vb_R^{N,2}(t))+
(\theta_R(\varepsilon_t^1)
-\theta_R(\varepsilon_t^2))
\vb_R^{N,2}(t),
\end{align*}
we infer that
\begin{align}
\Vert
\theta_R(\varepsilon_\cdot^1)\vb_R^{N,1}
-\theta_R(\varepsilon_\cdot^2)\vb_R^{N,2}
\Vert_{X_T}&\leq D+
\sup_{t\in[0,T]}
|\theta_R(\varepsilon_t^1)
-\theta_R(\varepsilon_t^2)|
\Vert\vb_R^{N,2}\Vert_{X_\rho}\nonumber\\
&\leq
D+(2R+D)\min\left\{1,\frac{2D}{R}\right\}
\leq 6D.
\end{align}
Indeed, if $D\leq R/2$, then
\[
D+(2R+D)\frac{2D}{R}
=5D+\frac{2D^2}{R}\leq6D.
\]
If $D>R/2$, then
\[
D+(2R+D)=2R+2D\leq6D.
\]

Finally, if $\rho_i=0$, then $\theta_R(\varepsilon_t^i)=0$ for every $t\in[0,T]$, so the sum of the two cut-off norms below is zero. If $\rho_i>0$, since $\theta_R(\varepsilon_t^i)=0$ for $t\geq\rho_i$ when $\rho_i<T$, and $\varepsilon_t^i<2R$ for $t<\rho_i$, continuity of the running norms gives
\begin{align}
&\Vert
\theta_R(\varepsilon_\cdot^i)\vb_R^{N,i}
\Vert_{X_T}
+\Vert
\theta_R(\varepsilon_\cdot^i)\mathbf z
\Vert_{Y_T}\leq
\Vert\vb_R^{N,i}\Vert_{X_{\rho_i}}
+\Vert\mathbf z\Vert_{Y_{\rho_i}}
\leq2R,\qquad i=1,2.
\end{align}
Consequently, the sum of the two cut-off norms is at most $2R$ in both cases.

Therefore, in view of \eqref{eq-difference}, there exists a constant $c_2>0$ independent of $N$ and $R$, such that
\begin{align}\label{eq-difference1}
&\Vert L_R^N(\vb_R^{N,1})- L_R^N(\vb_R^{N,2}) \Vert_{X_{T }}
\leq c_2NR\Vert \vb_R^{N,1}-\vb_R^{N,2}\Vert_{X_{T }},
\end{align}
Taking second moments in \eqref{eq-map} and \eqref{eq-difference1}, and using the definition \eqref{M^2 X_T}, gives \eqref{eq-map-1} and \eqref{eq-difference-1}, after increasing the constants if necessary.

\subsection{Equivalence of mild and weak solutions}\label{app-mild-weak-proof}
We give the detailed proof of \autoref{eq-equivalence} under the assumptions of \autoref{thm-local l2}.
Throughout the proof, $\mathbf z(t)=e^{-tA}\overrightarrow{\ub_0}+\wb(t)$, as in \eqref{eq-z}.

\begin{proof}[Proof of \autoref{eq-equivalence}]
Firstly, we prove that the existence of local mild solutions to (\ref{eqn-Burgers-intro.4}) implies the existence of local weak solutions to (\ref{eqn-Burgers-intro.4}).
Assume $(\vb, \tau )$ to be the local mild solution to (\ref{eqn-Burgers-intro.4}). Then $ \vb\in C([0,\tau]; \mathbb{H}^{\frac{1}{2}}(\dom))\cap L^{2}([0,\tau]; \mathbb{H}^{\frac{3}{2}}(\dom))$. Let us choose and fix $\omega\in\Omega$. Note that $\mathbf{z}:= \mathbf{z}(\cdot, \overrightarrow{\ub_{0}})\in L^{4}([0,\tau]; \mathbb{L}^{6}(\dom))$ and  $\nabla \mathbf{z}\in L^{2}([0,\tau]; \mathbb{L}^{3}(\dom))$, see \autoref{prop-ouz}, so by the H\"{o}lder inequality, the Sobolev embedding theorem, see \cite[Theorem 5.4]{Adams_1975} and the Sobolev interpolation inequality, see \cite[Corollary 5.1]{Brezis+Mironescu_2018}, we infer that
\begin{align}\label{eq-energy estimate}
\Vert( (\vb+\mathbf{z})\cdot \nabla) (\vb+\mathbf{z})\Vert_{L^2}&\leq \Vert \vb+\mathbf{z}\Vert_{L^6}\Vert \nabla (\vb+\mathbf{z})\Vert_{L^3}
\end{align}

Taking the $L^{\frac43}([0,\tau])$-norm on both sides of
\eqref{eq-energy estimate} and using the H\"older inequality, we infer that
\begin{align}
&\Vert
((\vb+\mathbf z)\cdot\nabla)(\vb+\mathbf z)
\Vert_{L^{\frac43}([0,\tau];\mathbb L^2(\dom))}
\nonumber\\
&\leq
\Vert \vb+\mathbf z\Vert_{L^4([0,\tau];\mathbb L^6(\dom))}
\Vert \nabla(\vb+\mathbf z)\Vert_{L^2([0,\tau];\mathbb L^3(\dom))}
\nonumber\\
&\leq
\Big(
\Vert\vb\Vert_{L^4([0,\tau];\mathbb L^6(\dom))}
+
\Vert\mathbf z\Vert_{L^4([0,\tau];\mathbb L^6(\dom))}
\Big)
\Big(
\Vert\nabla\vb\Vert_{L^2([0,\tau];\mathbb L^3(\dom))}
+
\Vert\nabla\mathbf z\Vert_{L^2([0,\tau];\mathbb L^3(\dom))}
\Big).
\label{eq-nonlinear-L43-1}
\end{align}

Therefore,
\[
((\vb+\mathbf z)\cdot\nabla)(\vb+\mathbf z)
\in
L^{\frac43}([0,\tau];\mathbb L^2(\dom)).
\]
Since $\alpha(\cdot,\omega)$ is continuous on $[0,\tau]$, we have $\Vert\alpha\Vert_{L^\infty([0,\tau])}<\infty$.
Hence,
\begin{align}\label{eq-nonlinear-L43}
\alpha((\vb+\mathbf z)\cdot\nabla)(\vb+\mathbf z)
\in
L^{\frac43}([0,\tau];\mathbb L^2(\dom)).
\end{align}

Note that $(\vb, \tau )$ is a local mild solution to (\ref{eqn-Burgers-intro.4}). By \autoref{semigroup-derivative} and \eqref{eq-nonlinear-L43}, we infer that for almost every $t\in[0,\tau]$
\[
\partial_t\vb
=
\Delta\vb
-
\alpha(t)((\vb+\mathbf z)\cdot\nabla)(\vb+\mathbf z)
\quad\text{in }\mathbb L^2(\dom).
\]
Moreover,
\[
\partial_t\vb
\in
L^1([0,\tau];\mathbb H^{-\frac12}(\dom)).
\]

Therefore, for almost every $t\in [0, \tau]$ and for all $\eta\in D(\Lambda^{2})$,
\begin{align*}
&\Big{\langle} \partial_{t} \vb , \eta\Big{\rangle}+\Big{\langle}\Lambda^{\frac{3}{2}} \vb, \Lambda^{\frac{1}{2}}\eta\Big{\rangle}+\Big{\langle}\alpha(t)((\vb+\mathbf{z})\cdot \nabla) (\vb+\mathbf{z}), \eta\Big{\rangle}=0.
\end{align*}
Note that
$D(\Lambda^{2})$ is dense in $\mathbb{H}^{\frac{1}{2}}(\dom)$, therefore, for every $\eta \in \mathbb{H}^{\frac{1}{2}}(\dom)$, we have
\begin{align*}
&\Big{\langle} \partial_{t} \vb(t), \eta\Big{\rangle}+\Big{\langle}\Lambda^{\frac{3}{2}} \vb(t), \Lambda^{\frac{1}{2}}\eta\Big{\rangle}+\Big{\langle}\alpha(t)((\vb(t)+\mathbf{z}(t))\cdot \nabla) (\vb(t)+\mathbf{z}(t)), \eta\Big{\rangle}=0, a.e. t\in [0, \tau].
\end{align*}

Conversely, we show that
any local weak solution to (\ref{eqn-Burgers-intro.4}) is also a local mild solution.

Let us choose and fix $\omega\in\Omega_0$, where $\mathbb P(\Omega_0)=1$ as in \autoref{def-3.1}.
Assume $(\vb, \tau )$ to be the local weak solution to (\ref{eqn-Burgers-intro.4})
in the sense of  \autoref{def-3.1}, i.e. for almost every $t \in [0,\tau]$ and all $\eta\in \mathbb{H}^{\frac{1}{2}}(\dom)$,
\begin{align*}
&\Big{\langle}\partial_{t} \vb(t), \eta\Big{\rangle}+\Big{\langle}\Lambda^{\frac{3}{2}} \vb(t), \Lambda^{\frac{1}{2}}\eta\Big{\rangle}+\alpha(t)\Big{\langle}((\vb(t)+\mathbf{z}(t,\overrightarrow{\ub_{0}} ))\cdot \nabla) (\vb(t)+\mathbf{z}(t,\overrightarrow{\ub_{0}} )), \eta\Big{\rangle}=0.
\end{align*}
Let
\begin{align*}
L(t, \xb): = -\int_{0}^{t}e^{(t-s)\Delta}(\alpha(s)((\vb+\mathbf{z})\cdot \nabla) (\vb+\mathbf{z}))(s,\xb) \, ds, t\in [0, \tau], \xb\in \dom.
\end{align*}
In view of \autoref{semigroup-derivative} and \eqref{eq-nonlinear-L43}, taking derivative with respect to time $t$ yields that
\begin{align*}
\partial_{t} L=\Delta L-\alpha(t)((\vb+\mathbf{z})\cdot \nabla) (\vb+\mathbf{z})\text{ in } \mathbb{L}^{2}(\dom), \text{ for a.e. } t\ \in [0,\tau]
\end{align*}
Therefore
\begin{align*}
\partial_{t} (L-\vb)=\Delta L- (\partial_{t}\vb+\alpha(t)((\vb+\mathbf{z})\cdot \nabla) (\vb+\mathbf{z}))=\Delta (L-\vb ),\;\;\mbox{ in }\mathbb{H}^{-\frac{1}{2}}(\dom).
\end{align*}

Since $\vb(0)=0$ and
$\vb\in C([0,\tau];\mathbb H^{\frac12}(\dom))$, we infer that
\[
\lim_{t\to0+}\vb(t)=0
\quad\text{in }\mathbb L^2(\dom).
\]
Moreover, by the contraction property of the heat semigroup,
\begin{align*}
\Vert L(t)\Vert_{L^2}
&\leq
\int_0^t
\Vert
\alpha(s)((\vb+\mathbf z)\cdot\nabla)(\vb+\mathbf z)(s)
\Vert_{L^2}\,ds
\rightarrow 0,\;\; \mbox{ as }t\to0+.
\end{align*}
Hence,
\[
\lim_{t\to0+}(L(t)-\vb(t))=0
\quad\text{in }\mathbb L^2(\dom).
\]

By a standard uniqueness result for the heat flow, we infer that
\[
L=\vb,\;\;\mbox{ on }[0,\tau].
\]
Therefore,
\[
\vb(t)
=
-\int_0^t
e^{(t-s)\Delta}
\alpha(s)((\vb+\mathbf z)\cdot\nabla)(\vb+\mathbf z)(s)\,ds\in C([0,\tau]; \mathbb{H}^{\frac{1}{2}}(\dom))\cap
L^{2}([0,\tau]; \mathbb{H}^{\frac{3}{2}}(\dom)),
\]
which implies $(\vb,\tau)$ is a local mild solution to
(\ref{eqn-Burgers-intro.4}).

This completes the proof.
\end{proof}

\section{Technical estimates for the global existence argument}\label{sec-global-approximation}

\subsection{Control of the spatial mean}\label{app-mean-control}

\begin{lemma}\label{lemma-2.1-aux}
If $\ab,\bb \in \mathbb{H}^1(\dom) $ then
\begin{align}\label{eqn-aux1}
\left|\int_{\dom}(\ab(\xb)\cdot \nabla) \bb (\xb)\,d\xb \right| &\leq \frac{1}{2} \left( \Vert  \ab\Vert_{\frac{1}{2}}^2+\Vert  \bb\Vert_{\frac{1}{2}}^2\right).
\end{align}
In particular,
\begin{align}\label{eqn-aux2}
\left|\int_{\dom}(\ab(\xb)\cdot \nabla) \ab (\xb)\,d\xb \right| &\leq  \Vert  \ab\Vert_{\frac{1}{2}}^2.
\end{align}

\end{lemma}

\begin{proof}[Proof of Lemma \ref{lemma-2.1-aux}]
We first assume that $\ab,\bb\in\mathbb{H}^{\frac12}(\dom)$ are trigonometric polynomials.
For $r,q=1,2,3$, in view of \eqref{eqn-u-real} we infer that
\begin{align*}
a_r(\xb)&=(2\pi)^{-\frac32}\hat{\ab}_0^r+(4\pi^3)^{-\frac12}\sum_{\kb\in K_0}
\left(\hat{\ab}_{\kb}^{r}\sin(\kb\cdot\xb)+\hat{\ab}_{\kb}^{r+3}\cos(\kb\cdot\xb)\right),
\\
b_q(\xb)&=(2\pi)^{-\frac32}\hat{\bb}_0^q+(4\pi^3)^{-\frac12}\sum_{\kb\in K_0}\left(\hat{\bb}_{\kb}^{q}\sin(\kb\cdot\xb)+\hat{\bb}_{\kb}^{q+3}\cos(\kb\cdot\xb)\right).
\end{align*}
For $\kb=(k_1,k_2,k_3)\in K_0$, we denote by $k_r$the $r$-th component of $\kb$. Then
\begin{align*}
\partial_r b_q(\xb)
&=
(4\pi^3)^{-\frac12}
\sum_{\kb\in K_0}
k_r
\left(
\hat{\bb}_{\kb}^{q}\cos(\kb\cdot\xb)
-
\hat{\bb}_{\kb}^{q+3}\sin(\kb\cdot\xb)
\right).
\end{align*}

Since $[(\ab\cdot\nabla)\bb]_q=\sum_{r=1}^3 a_r\partial_r b_q$, we infer that
\begin{align*}
\left[
\int_{\dom}
(\ab(\xb)\cdot\nabla)\bb(\xb)\,d\xb
\right]_q
&=
\sum_{r=1}^3
\int_{\dom}
a_r(\xb)\partial_r b_q(\xb)\,d\xb.
\end{align*}
Note that for $\kb,\ell\in K_0$,
\begin{align*}
&\int_{\dom}\sin(\kb\cdot\xb)\,d\xb=\int_{\dom}\cos(\kb\cdot\xb)\,d\xb=0,\\
&\int_{\dom}\sin(\kb\cdot\xb)\cos(\ell\cdot\xb)\,d\xb=\int_{\dom}\cos(\kb\cdot\xb)\sin(\ell\cdot\xb)\,d\xb=0,\\
&\int_{\dom}\sin(\kb\cdot\xb)\sin(\ell\cdot\xb)\,d\xb=\int_{\dom}\cos(\kb\cdot\xb)\cos(\ell\cdot\xb)\,d\xb=
\begin{cases}
4\pi^3, & \kb=\ell,\\
0, & \kb\neq\ell.
\end{cases}
\end{align*}

Therefore,
\begin{align*}
&\left[\int_{\dom}(\ab(\xb)\cdot\nabla)\bb(\xb)\,d\xb\right]_q\\
&=\frac{1}{4\pi^3}\sum_{r=1}^3\sum_{\kb,\ell\in K_0}\ell_r\int_{\dom}\Big[\hat{\ab}_{\kb}^{r}\sin(\kb\cdot\xb)
+
\hat{\ab}_{\kb}^{r+3}\cos(\kb\cdot\xb)
\Big]\times
\Big[
\hat{\bb}_{\ell}^{q}\cos(\ell\cdot\xb)
-
\hat{\bb}_{\ell}^{q+3}\sin(\ell\cdot\xb)
\Big]\,d\xb
\\
&=
\frac{1}{4\pi^3}
\sum_{r=1}^3
\sum_{\kb,\ell\in K_0}
\ell_r
\Bigg[
\hat{\ab}_{\kb}^{r}\hat{\bb}_{\ell}^{q}
\int_{\dom}
\sin(\kb\cdot\xb)\cos(\ell\cdot\xb)\,d\xb-
\hat{\ab}_{\kb}^{r}\hat{\bb}_{\ell}^{q+3}
\int_{\dom}
\sin(\kb\cdot\xb)\sin(\ell\cdot\xb)\,d\xb
\\
&\qquad\qquad
+
\hat{\ab}_{\kb}^{r+3}\hat{\bb}_{\ell}^{q}
\int_{\dom}
\cos(\kb\cdot\xb)\cos(\ell\cdot\xb)\,d\xb-
\hat{\ab}_{\kb}^{r+3}\hat{\bb}_{\ell}^{q+3}
\int_{\dom}
\cos(\kb\cdot\xb)\sin(\ell\cdot\xb)\,d\xb
\Bigg]\\
&=\sum_{\kb\in K_0}
\sum_{r=1}^3
k_r
\left(
\hat{\ab}_{\kb}^{r+3}\hat{\bb}_{\kb}^{q}
-
\hat{\ab}_{\kb}^{r}\hat{\bb}_{\kb}^{q+3}
\right).
\end{align*}

Hence, by the triangle inequality and the Cauchy--Schwarz inequality,
\begin{align*}
&\left|
\int_{\dom}
(\ab(\xb)\cdot\nabla)\bb(\xb)\,d\xb
\right|
=
\left[
\sum_{q=1}^3
\left|
\sum_{\kb\in K_0}
\sum_{r=1}^3
k_r
\left(
\hat{\ab}_{\kb}^{r+3}\hat{\bb}_{\kb}^{q}
-
\hat{\ab}_{\kb}^{r}\hat{\bb}_{\kb}^{q+3}
\right)
\right|^2
\right]^{\frac12}
\\
&\leq
\sum_{\kb\in K_0}
\left[
\sum_{q=1}^3
\left|
\sum_{r=1}^3
k_r
\left(
\hat{\ab}_{\kb}^{r+3}\hat{\bb}_{\kb}^{q}
-
\hat{\ab}_{\kb}^{r}\hat{\bb}_{\kb}^{q+3}
\right)
\right|^2
\right]^{\frac12}
=
\sum_{\kb\in K_0}
\left[
\sum_{q=1}^3
\left|
\left(
\sum_{r=1}^3
k_r\hat{\ab}_{\kb}^{r+3}
\right)
\hat{\bb}_{\kb}^{q}
-
\left(
\sum_{r=1}^3
k_r\hat{\ab}_{\kb}^{r}
\right)
\hat{\bb}_{\kb}^{q+3}
\right|^2
\right]^{\frac12}
\\
&\leq
\sum_{\kb\in K_0}
\left|
\sum_{r=1}^3
k_r\hat{\ab}_{\kb}^{r+3}
\right|
\left(
\sum_{q=1}^3
|\hat{\bb}_{\kb}^{q}|^2
\right)^{\frac12}+
\sum_{\kb\in K_0}
\left|
\sum_{r=1}^3
k_r\hat{\ab}_{\kb}^{r}
\right|
\left(
\sum_{q=1}^3
|\hat{\bb}_{\kb}^{q+3}|^2
\right)^{\frac12}
\\
&\leq
\sum_{\kb\in K_0}
|\kb|
\Bigg[
\left(
\sum_{r=1}^3
|\hat{\ab}_{\kb}^{r+3}|^2
\right)^{\frac12}
\left(
\sum_{q=1}^3
|\hat{\bb}_{\kb}^{q}|^2
\right)^{\frac12}
+
\left(
\sum_{r=1}^3
|\hat{\ab}_{\kb}^{r}|^2
\right)^{\frac12}
\left(
\sum_{q=1}^3
|\hat{\bb}_{\kb}^{q+3}|^2
\right)^{\frac12}
\Bigg]\\
&\leq
\sum_{\kb\in K_0}
|\kb|
\left(
\sum_{j=1}^6
|\hat{\ab}_{\kb}^{j}|^2
\right)^{\frac12}
\left(
\sum_{j=1}^6
|\hat{\bb}_{\kb}^{j}|^2
\right)^{\frac12}.
\end{align*}

Using the Cauchy--Schwarz inequality again, we obtain
\begin{align*}
\left|
\int_{\dom}
(\ab(\xb)\cdot\nabla)\bb(\xb)\,d\xb
\right|
&\leq
\left(
\sum_{\kb\in K_0}
\sum_{j=1}^6
|\kb||\hat{\ab}_{\kb}^{j}|^2
\right)^{\frac12}
\left(
\sum_{\kb\in K_0}
\sum_{j=1}^6
|\kb||\hat{\bb}_{\kb}^{j}|^2
\right)^{\frac12}
\\
&=
\Vert\ab\Vert_{\frac12}
\Vert\bb\Vert_{\frac12}
\leq
\frac12
\left(
\Vert\ab\Vert_{\frac12}^2
+
\Vert\bb\Vert_{\frac12}^2
\right).
\end{align*}
Thus, \eqref{eqn-aux1} holds for trigonometric polynomials.
By the density of trigonometric polynomials in $\mathbb H^1(\dom)$ and the H\"older inequality, \eqref{eqn-aux1} holds for every $\ab,\bb\in\mathbb H^1(\dom)$.
Inequality \eqref{eqn-aux2} follows directly from
\eqref{eqn-aux1}.
\end{proof}

The next lemma follows by integrating \eqref{ee-epsilon} over $\dom$.
\begin{lemma}\label{lemma-2.1}
Assume that $\tilde{\vb}_{0}\in \mathbb{H}^{1}(\dom)$ and \autoref{A2} holds. Let  $(\vb^{\eps,m}, \xi^{\eps,m})$ be the maximal strong solutions of (\ref{ee-epsilon}). Then, for  every
$t\in [0,\xi^{\eps,m}) $,

\begin{align*}
\left | \int_{\dom}\vb^{\eps,m}(t,\xb)\,d\xb \right | \leq 2\int_{0}^{t}\alpha^\eps(s)\left ( \Vert \vb^{\eps,m}(s)\Vert_{\frac{1}{2}}^{2}+\Vert \wb_m^{\eps}(s)\Vert_{\frac{1}{2}}^{2} \right )\,ds+\Vert\tilde{\vb}_{0}\Vert_{L^1}.
\end{align*}
\end{lemma}

\begin{proof}
Note that
\begin{align}\label{eq-ekj}
\int_\dom e_\kb^\jb\,dx=
0,\;\;k\neq 0.
\end{align}
Therefore,  for almost every $t\in [0,\xi^{\eps,m}) $
\[
\int_{\dom}\Delta\vb^{\eps,m}(t,\xb)\,d\xb=-\sum_{(\kb,\jb)\in\mathcal K}
|\kb|^2(\vb^{\eps,m})_\kb^\jb(t)
\int_{\dom}e_\kb^\jb(\xb)\,d\xb=0.
\]
Therefore, by integration by parts, we infer that
\[
\int_{\dom}\Delta\vb^{\eps,m}(t,\xb)\,d\xb=0,
\quad\text{for a.e. }t\in[0,\xi^{\eps,m}).
\]

Since $\partial_t\vb^{\eps,m}\in L^2([0,T];\mathbb L^2(\dom))$ for every $T<\xi^{\eps,m}$, we infer that $t\mapsto\int_\dom\vb^{\eps,m}(t,\xb)\,d\xb$ is absolutely continuous.
Therefore, integrating \eqref{ee-epsilon} over $\dom$, by \autoref{lemma-2.1-aux} and the Cauchy-Schwarz inequality, we infer that for almost every $t\in[0,\xi^{\eps,m})$,

\begin{align}
\left|\frac{d}{dt}\int_{\dom} \vb^{\eps,m}(t,\xb)\,d\xb\right|&= \left|\int_{\dom} \Delta\vb^{\eps,m}(t,\xb)\,d\xb-\alpha^\eps(t)\int_{\dom}((\vb^{\eps,m}+\wb_m^{\eps})\cdot \nabla) (\vb^{\eps,m}+\wb_m^{\eps}) (t,\xb)\,d\xb\right|\nonumber\\
&\leq \alpha^\eps(t) \Vert \vb^{\eps,m}(t)+\wb_m^{\eps}(t)\Vert_{\frac{1}{2}}^{2} \nonumber\\
&\leq 2\alpha^\eps(t)\Big{(}\Vert \vb^{\eps,m}(t)\Vert_{\frac{1}{2}}^{2}+\Vert \wb_m^{\eps}(t)\Vert_{\frac{1}{2}}^{2} \Big{)}.
\end{align}
By  integration of the above estimate with respect to $t$, we get
\begin{align}
\left | \int_{\dom}\vb^{\eps,m}(t,\xb)\,d\xb \right |&\leq  \left | \int_0^t \frac{d}{ds}\int_{\dom} \vb^{\eps,m}(s,\xb)\,d\xb\,ds\right |+\left|\int_{\dom}\tilde\vb_{0}(\xb) \,d\xb\right|
\nonumber\\
&\leq \int_0^t \left|\frac{d}{ds}\int_{\dom} \vb^{\eps,m}(s,\xb)\,d\xb\right |\,ds+\Vert\tilde{\vb}_{0}\Vert_{L^1}\nonumber\\
&\leq 2\int_{0}^{t}\alpha^\eps(s)\left ( \Vert \vb^{\eps,m}(s)\Vert_{\frac{1}{2}}^{2}+\Vert \wb_m^{\eps}(s)\Vert_{\frac{1}{2}}^{2} \right )\,ds+\Vert\tilde{\vb}_{0}\Vert_{L^1}.
\end{align}
The proof is complete.
\end{proof}

\begin{proof}[Proof of \autoref{cor-dominate L2}] Let us choose and fix  $ s>0$.\\
Define $\bar{\vb}(t)$ by
\begin{align}\label{eq-bar-vn}
 \bar{\vb}(t):= \frac{1}{(2\pi)^3}\int_{\dom}\vb^{\eps,m}(t,\xb)\,d\xb, t\in [0, \xi^{\eps,m}).
\end{align}
Then, using notation \eqref{eqn-u-real} and inequality \eqref{eq-ekj},  we infer that
 \begin{align*}
\bar{\vb}(t)=\frac{1}{(2\pi)^3}\int_{\dom}\sum_{(k,j)\in \mathcal{K},} (\vb^{\eps,m})_\kb^\jb\eb_k^j(t,\xb)\,d\xb=\frac{1}{(2\pi)^3}\int_{\dom}\sum_{j=1}^{3} (\vb^{\eps,m})_0^\jb\eb_0^j(t,\xb)\,d\xb=\sum_{j=1}^{3} (\vb^{\eps,m})_0^\jb\eb_0^j(t,\xb).
 \end{align*}
Thus, by inequality \eqref{eqn-1013} we infer that  (also for $s=0$)
\begin{align}\label{eqn-v-bar v}
\Vert \vb^{\eps,m}-\bar{\vb}\Vert_{L^2} \leq \Vert \vb^{\eps,m}\Vert_{s}.
\end{align}
Hence, by the triangle inequality, \eqref{eqn-v-bar v} and \autoref{lemma-2.1} we obtain
\begin{align}
\Vert \vb^{\eps,m}(t) \Vert_{L^2} &\leq \Vert \vb^{\eps,m}(t)-\bar{\vb}(t)\Vert_{L^2}+\Vert \bar{\vb}(t)\Vert_{L^2}= \Vert \vb^{\eps,m}(t)-\bar{\vb}(t)\Vert_{L^2}+
\frac{1}{(2\pi)^\frac{3}{2}}\left | \int_{\dom}\vb^{\eps,m}(t,\xb)\,d\xb \right |\nonumber\\
&\leq  \Vert \vb^{\eps,m}(t)\Vert_{s}+\frac{1}{(2\pi)^\frac{3}{2}}\left | \int_{\dom}\vb^{\eps,m}(t,\xb)\,d\xb \right |\nonumber\\
&\leq \Vert \vb^{\eps,m}(t)\Vert_{s}+\frac{2}{(2\pi)^\frac{3}{2}}\int_{0}^{t}\alpha^{\eps}(s) \left(\Vert \vb^{\eps,m}(s)\Vert_{\frac{1}{2}}^{2}+\Vert \wb_m^{\eps}(s)\Vert_{\frac{1}{2}}^{2}\right)\,ds+\frac{1}{(2\pi)^\frac{3}{2}}\Vert\tilde{\vb}_{0}\Vert_{L^1}.
\end{align}

Recall that
\begin{align}
\Vert \cdot\Vert_{H^s}^2=\Vert \cdot\Vert_{L^2}^2+\Vert \cdot\Vert_{s}^2.
\end{align}
It follows that
\begin{align}
\Vert \vb^{\eps,m}(t)\Vert_{s}&\leq \Vert \vb^{\eps,m}(t)\Vert_{H^{s}}\leq \Vert \vb^{\eps,m}(t)\Vert_{L^2}+\Vert \vb^{\eps,m}(t)\Vert_{s}\nonumber\\
&\leq 2\Big (\Vert \vb^{\eps,m}(t)\Vert_{s}+\frac{1}{(2\pi)^\frac{3}{2}}\int_{0}^{t}\alpha^{\eps}(s)\left (\Vert \vb^{\eps,m}(s)\Vert_{\frac{1}{2}}^{2}+\Vert \wb_m^{\eps}(s)\Vert_{\frac{1}{2}}^{2}\right )\,ds\Big )+\frac{1}{(2\pi)^\frac{3}{2}}\Vert\tilde{\vb}_{0}\Vert_{L^1},
\end{align}
which implies \eqref{eq-2.8} holds.

Fix $s\geq \frac12$ and $0<T<\xi^{\eps,m}$ such that $\sup_{t\in[0,T]}\Vert\vb^{\eps,m}(t)\Vert_s<\infty$. Since
\[
\Vert \vb^{\eps,m}(t)\Vert_{\frac12} \leq \Vert \vb^{\eps,m}(t)\Vert_{s},\;\;t\in[0,T]
\]
and  $\wb_m^{\eps}\in C^\infty([0,T]\times\dom;\mathbb R^3)$ and $\alpha^\eps$ is bounded on $[0,T]$, we infer that for every $t\in [0,T]$
\begin{align*}
&\int_0^t\alpha^\eps(r)\left(\Vert \vb^{\eps,m}(r)\Vert_{\frac12}^2+\Vert \wb_m^\eps(r)\Vert_{\frac12}^2\right)\,dr
\leq
T\sup_{r\in[0,T]}\alpha^\eps(r)\sup_{r\in[0,T]}\Vert \vb^{\eps,m}(r)\Vert_s^2+\int_0^T\alpha^\eps(r)\Vert \wb_m^\eps(r)\Vert_{\frac12}^2\,dr.
\end{align*}
Hence, in view of \eqref{eq-2.8},
\begin{align*}
\sup_{t\in[0,T]}
\Vert\vb^{\eps,m}(t)\Vert_{H^s}
&\leq
2\sup_{t\in[0,T]}
\Vert\vb^{\eps,m}(t)\Vert_s+\frac{2T}{(2\pi)^{\frac32}}
\sup_{r\in[0,T]}\alpha^\eps(r)\sup_{r\in[0,T]}\Vert\vb^{\eps,m}(r)\Vert_s^2
\\
&\quad+
\frac{2}{(2\pi)^{\frac32}}
\int_0^T
\alpha^\eps(r)
\Vert\wb_m^\eps(r)\Vert_{\frac12}^2\,dr
+
\frac{1}{(2\pi)^{\frac32}}
\Vert\tilde{\vb}_0\Vert_{L^1},
\end{align*}
which implies there exists a constant $C_T>0$ such that
\[
\sup_{t\in[0,T]}
\Vert\vb^{\eps,m}(t)\Vert_{H^s}
\leq
C_T
\left(
1+
\sup_{t\in[0,T]}
\Vert\vb^{\eps,m}(t)\Vert_s
+
\sup_{t\in[0,T]}
\Vert\vb^{\eps,m}(t)\Vert_s^2
\right).
\]

This proves \eqref{eqn-H^s bound}.
The proof of \autoref{cor-dominate L2} is complete.
\end{proof}

\subsection{Galerkin approximation and regularity}\label{app-galerkin-regularity}

In the following, we show that the strong solution $ \vb^{\eps,m}$ is classical, that is, $\vb^{\eps,m} \in C^{1,0}((0, \tau^{\eps,m}]\times \dom; \mathbb{R}^{3})\cap C^{0,2}((0, \tau^{\eps,m}]\times \dom; \mathbb{R}^{3})$. Here, for $n,r\in\mathbb{N}$, $C^{n,r} ((0, \tau^{\eps,m}]\times \dom; \mathbb{R}^{3})$ is defined by
\begin{align}\label{eq-Cnm}
C^{n,r} ((0, \tau^{\eps,m}]\times \dom; \mathbb{R}^{3}):=\Big{\{}\ub: \frac{\partial^{k}}{\partial t^{k}}\frac{\partial^{[\mathbf{j}]}}{\partial x_{1}^{j_{1}}...\partial x_{3}^{j_{3}}} \ub\in C((0, \tau^{\eps,m}]\times \dom; \mathbb{R}^{3}),\  k\leq n,\, [\mathbf{j}]\leq r \Big{\}},
\end{align}
where $[\mathbf{j}]=j_{1}+j_{2}+j_{3}, j_{i}\in \mathbb{N},i=1,2,3$.

For this purpose, we need to consider the Galerkin approximation of (\ref{ee-epsilon}).

For $n\in \mathbb{N}$, let $P_{n}:\mathbb{L}^{2}(\dom) \to \mathbb{L}^{2}(\dom) $ denote the projection onto the Fourier modes of order up to $n$, that is,
\begin{equation}\label{eqn-P_n}
P_n(\ub)= \sum_{(k,j)\in \mathcal{K},|\kb|\leq n} \hat{\ub}_{\kb}^{\jb}\eb_k^j, \mbox{ if }  \ub=\sum_{(k,j)\in \mathcal{K}} \hat{\ub}_{\kb}^{\jb}\eb_k^j \in \mathbb{L}^{2}(\dom) .
    \end{equation}

Let us denote by $H_n$  closed subspace of $\mathbb{L}^{2}(\dom)$ generated by Fourier modes of order up to $n$, i.e.
\begin{equation}\label{eqn-H_n}
H_n:=\left\{\ub \in  \mathbb{L}^{2}(\dom): \ub= \sum_{(k,j)\in \mathcal{K},|\kb|\leq n} \hat{\ub}_{\kb}^{\jb}\eb_k^j \right\}  .  \end{equation}
It follows that $P_n$ is the orthogonal projection from the space $\mathbb{L}^{2}(\dom)$ onto $H_n$.
Note that since $ H_n$ is finite dimensional, all norms on $H_n$ are equivalent. In what follows we will use the norm inherited from the Sobolev space $\mathbb{H}^{\frac{1}{2}}(\dom)$.

We also have
\begin{lemma} If $\ub \in H_n$ then $\Delta \ub \in H_n$ and
\[
P_n (\Delta \ub )=\Delta \ub.
\]
\end{lemma}

Note that $\wb_m^{\eps}\in C^{\infty}([0, T]\times \dom; \mathbb{R}^{3})$, we introduce the Faedo-Galerkin approximation of (\ref{ee-epsilon}) below.

\begin{align}\label{ee-epsilon approximation}
&\partial_{t} \vb_{n}^{\eps,m}-\Delta \vb_{n}^{\eps,m}+\alpha^{\eps}(t)P_{n}\Big{(}((\vb_{n}^{\eps,m}+\wb_m^{\eps})\cdot \nabla) (\vb_{n}^{\eps,m}+\wb_m^{\eps})\Big{)}=0,\ \mathrm{on}\ [0,\infty)\times \dom,\\
&\vb_{n}^{\eps,m}(0)=P_{n}\tilde\vb_{0}(\xb),\quad \xb=(x_{1},x_{2},x_{3})\in \dom.\nonumber
\end{align}

In view of the proof of Proposition 3.1 of \cite{Dong+Guo+Wu+Zhou_2023}, we have the following result.

\begin{lemma}\label{lem-vneps-regularization}
Assume that $\tilde\vb_{0}\in \mathbb{H}^{1}(\dom)$ and \autoref{A2} holds.

For every $n\in\mathbb N$, let $(\vb_n^{\eps,m},\tau_n)$ be the unique maximal solution to the Galerkin approximation \eqref{ee-epsilon approximation}, where $\tau_n\in(0,\infty]$ denotes its maximal existence time. Let $\xi_1^{\eps,m}$ denote the maximal existence time of \eqref{ee-epsilon},
as in \autoref{prop-maximal global-veps}.
Then for every $\tau^{\eps,m}\in (0,\xi_1^{\eps,m})$, there exists $n_0=n_0(\tau^{\eps,m},\eps,m,\omega)$ such that $\tau_n>\tau^{\eps,m}$ and $\{\vb_{n}^{\eps,m}\}_{n\ge n_0}$ is uniformly bounded in $L^\infty([0, \tau^{\eps,m}];\mathbb{H}^1(\dom))\cap L^2([0, \tau^{\eps,m}];\mathbb{H}^2(\dom))$, i.e., the following estimate holds almost surely,
\begin{align}\label{eq-uniform-H1-H2}
\sup\limits_{n\geq n_0} \sup\limits_{t\in [0, \tau^{\eps,m}]} \Vert \vb_{n}^{\eps,m} (t) \Vert_{H^1}
+\sup\limits_{n\geq n_0} \int_{0}^ {\tau^{\eps,m}} \Vert \vb_{n}^{\eps,m} (t) \Vert_{H^2}^{2}\,dt
<\infty .
\end{align}
Moreover, $\vb_n^{\eps,m}\to\vb^{\eps,m}$ in $C([0,\tau^{\eps,m}];\mathbb H^1(\dom))\cap L^2(0,\tau^{\eps,m};\mathbb H^2(\dom))$.
\end{lemma}
\begin{proof}
Let us choose and fix $\eps,m$ and $\omega$ outside a
$\mathbb P$-null set. Let $T=\tau^{\eps,m}<\xi_1^{\eps,m}$.
In view of \autoref{prop-maximal global-veps}, we infer that
\[
\vb^{\eps,m}\in C([0,T];\mathbb H^1(\dom))
\cap L^2([0,T];\mathbb H^2(\dom)).
\]

Let $I=[t_0,t_0+h]\subset[0,T]$, $h\in(0,1]$.
Recall that
\[
X_I^1:=C(I;\mathbb H^1(\dom))
\cap L^2(I;\mathbb H^2(\dom)).
\]
Following the proof of \autoref{prop-tildev-local strong},
by the H\"older inequality and the Sobolev embedding theorem,
we infer that for every $u,v\in X_I^1$,
\begin{align}\label{eq-Galerkin-local-Lipschitz}
&\Vert\alpha^\eps
\big[((u+\wb_m^\eps)\cdot\nabla)(u+\wb_m^\eps)
-((v+\wb_m^\eps)\cdot\nabla)(v+\wb_m^\eps)\big]
\Vert_{L^2(I;\mathbb L^2(\dom))}\nonumber\\
&\leq ch^{\frac14}
\Vert\alpha^\eps\Vert_{L^\infty([0,T])}
\Big(\Vert u\Vert_{X_I^1}+\Vert v\Vert_{X_I^1}
+2\Vert\wb_m^\eps\Vert_{C([0,T];\mathbb H^2(\dom))}\Big)
\Vert u-v\Vert_{X_I^1}.
\end{align}

Let us choose $M>0$ such that
\[
\Vert\vb^{\eps,m}\Vert_{X_{[0,T]}^1}
+\Vert\wb_m^\eps\Vert_{C([0,T];\mathbb H^2(\dom))}<M.
\]
Since $P_n$ commutes with $\Delta$ and
\[
\Vert P_nu\Vert_{\mathbb H^j(\dom)}
\leq\Vert u\Vert_{\mathbb H^j(\dom)},
\qquad j=0,1,2,
\]
by \autoref{thm-maximal regularity}, we infer that
for every $u_{n,0}\in H_n$ satisfying
$\Vert u_{n,0}\Vert_{\mathbb H^1(\dom)}\leq M+1$,
\begin{align*}
&\Vert e^{(\cdot-t_0)\Delta}u_{n,0}
\Vert_{X_I^1} \leq c\|u_{n,0}\|_{\mathbb H^1(\dom)}\leq c(M+1),\\
&\left\Vert\int_{t_0}^{\cdot}e^{(\cdot-s)\Delta}
P_n\left[\alpha^\eps
((u+\wb_m^\eps)\cdot\nabla)(u+\wb_m^\eps)\right]ds
\right\Vert_{X_I^1}\leq ch^{\frac14}
\Vert\alpha^\eps\Vert_{L^\infty([0,T])}
\Big(\Vert u\Vert_{X_I^1}
+\Vert\wb_m^\eps\Vert_{C([0,T];\mathbb H^2(\dom))}\Big)^2.
\end{align*}
Therefore, choosing $R>0$ sufficiently large
and $h>0$ sufficiently small, by
\eqref{eq-Galerkin-local-Lipschitz}
and the Banach fixed point theorem, we infer
that the Galerkin equation with initial value
$\vb_n^{\eps,m}(t_0)=u_{n,0}$ admits a unique
solution on $I$, satisfying
\[
\Vert\vb_n^{\eps,m}\Vert_{X_I^1}\leq R.
\]
Moreover, $R$ and $h$ are independent of
$n$ and $t_0$, and $h$ can be chosen such that
\[
ch^{\frac14}\Vert\alpha^\eps\Vert_{L^\infty([0,T])}
\Big(2R+2\Vert\wb_m^\eps
\Vert_{C([0,T];\mathbb H^2(\dom))}\Big)
\leq\frac12.
\]

In view of the mild formulations of
\eqref{ee-epsilon approximation} and
\eqref{ee-epsilon}, by \autoref{thm-maximal regularity}
and \eqref{eq-Galerkin-local-Lipschitz}, we infer that
\begin{align*}
\Vert\vb_n^{\eps,m}-\vb^{\eps,m}\Vert_{X_I^1}\leq &c\Vert\vb_n^{\eps,m}(t_0)
-\vb^{\eps,m}(t_0)\Vert_{\mathbb H^1(\dom)}
+\frac12\Vert\vb_n^{\eps,m}
-\vb^{\eps,m}\Vert_{X_I^1}\\
&+c\Vert(\mathrm{Id}-P_n)
\big[\alpha^\eps((\vb^{\eps,m}+\wb_m^\eps)\cdot\nabla)
(\vb^{\eps,m}+\wb_m^\eps)\big]
\Vert_{L^2(I;\mathbb L^2(\dom))}.
\end{align*}
Consequently,
\begin{align}\label{eq-Galerkin-local-convergence}
&\Vert\vb_n^{\eps,m}-\vb^{\eps,m}\Vert_{X_I^1}\nonumber\\
&\leq c\Big(
\Vert\vb_n^{\eps,m}(t_0)
-\vb^{\eps,m}(t_0)\Vert_{\mathbb H^1(\dom)}
+\Vert(\mathrm{Id}-P_n)
\big[\alpha^\eps((\vb^{\eps,m}+\wb_m^\eps)\cdot\nabla)
(\vb^{\eps,m}+\wb_m^\eps)\big]
\Vert_{L^2(I;\mathbb L^2(\dom))}
\Big).
\end{align}

Since $\vb_n^{\eps,m}(0)=P_n\tilde{\vb}_0$
and $\tilde{\vb}_0\in\mathbb H^1(\dom)$,
we infer that
\[
\Vert P_n\tilde{\vb}_0-\tilde{\vb}_0
\Vert_{\mathbb H^1(\dom)}
\longrightarrow0,\qquad n\rightarrow\infty.
\]
Moreover, by the H\"older inequality and
the Sobolev embedding theorem, we have
\[
\alpha^\eps((\vb^{\eps,m}+\wb_m^\eps)\cdot\nabla)
(\vb^{\eps,m}+\wb_m^\eps)
\in L^2([0,T];\mathbb L^2(\dom)).
\]
Hence,
\begin{align}\label{eq-Galerkin-projection-error}
&\Vert(\mathrm{Id}-P_n)
\big[\alpha^\eps((\vb^{\eps,m}+\wb_m^\eps)\cdot\nabla)
(\vb^{\eps,m}+\wb_m^\eps)\big]
\Vert_{L^2([0,T];\mathbb L^2(\dom))}
\longrightarrow0.
\end{align}

Let us divide $[0,T]$ into finitely many
subintervals
\[
0=t_0<t_1<\cdots<t_J=T,
\qquad t_j-t_{j-1}\leq h.
\]
Since
\[
\Vert\vb_n^{\eps,m}(0)\Vert_{\mathbb H^1(\dom)}
\leq\Vert\tilde{\vb}_0\Vert_{\mathbb H^1(\dom)}
<M,
\]
the uniform local existence result implies
that $\vb_n^{\eps,m}$ exists on $[0,t_1]$.
By \eqref{eq-Galerkin-local-convergence}
and \eqref{eq-Galerkin-projection-error},
we infer that
\[
\Vert\vb_n^{\eps,m}-\vb^{\eps,m}
\Vert_{X_{[0,t_1]}^1}\longrightarrow0.
\]
In particular,
\[
\Vert\vb_n^{\eps,m}(t_1)-\vb^{\eps,m}(t_1)
\Vert_{\mathbb H^1(\dom)}\longrightarrow0.
\]
Since
$\Vert\vb^{\eps,m}(t_1)\Vert_{\mathbb H^1(\dom)}<M$,
for sufficiently large $n$, we have
\[
\Vert\vb_n^{\eps,m}(t_1)\Vert_{\mathbb H^1(\dom)}
\leq M+1.
\]
Therefore, the uniform local existence result
also applies on $[t_1,t_2]$.

Repeating the above argument, for every
$j=2,\ldots,J$, by
\eqref{eq-Galerkin-local-convergence}, we infer that
\begin{align*}
&\Vert\vb_n^{\eps,m}(t_j)
-\vb^{\eps,m}(t_j)\Vert_{\mathbb H^1(\dom)}\\
&\leq c\Big(
\Vert\vb_n^{\eps,m}(t_{j-1})
-\vb^{\eps,m}(t_{j-1})\Vert_{\mathbb H^1(\dom)}
+\Vert(\mathrm{Id}-P_n)
\big[\alpha^\eps((\vb^{\eps,m}+\wb_m^\eps)\cdot\nabla)
(\vb^{\eps,m}+\wb_m^\eps)\big]
\Vert_{L^2([0,T];\mathbb L^2(\dom))}
\Big)
\longrightarrow0.
\end{align*}
Since $J$ is finite, we infer that
the Galerkin solutions exist on $[0,T]$
for sufficiently large $n$, and
\begin{align}\label{eq-Galerkin-convergence-H1-H2}
&\Vert\vb_n^{\eps,m}-\vb^{\eps,m}
\Vert_{C([0,T];\mathbb H^1(\dom))}
+\Vert\vb_n^{\eps,m}-\vb^{\eps,m}
\Vert_{L^2([0,T];\mathbb H^2(\dom))}
\longrightarrow0,\qquad n\rightarrow\infty.
\end{align}

Since the Galerkin equation is a finite-dimensional
ordinary differential equation, its solution
can be extended beyond $T$. Hence, there exists
$n_0=n_0(T,\eps,m,\omega)$ such that
\[
\tau_n>T,\qquad n\geq n_0.
\]
Finally, by \eqref{eq-Galerkin-convergence-H1-H2},
we infer that
\begin{align*}
&\sup_{n\geq n_0}\sup_{t\in[0,T]}
\Vert\vb_n^{\eps,m}(t)\Vert_{\mathbb H^1(\dom)}
+\sup_{n\geq n_0}\int_0^T
\Vert\vb_n^{\eps,m}(t)\Vert_{\mathbb H^2(\dom)}^2\,dt
<\infty.
\end{align*}
Therefore, \eqref{eq-uniform-H1-H2} holds.
The proof of \autoref{lem-vneps-regularization} is complete.
\end{proof}

Following the proofs of Lemmas \ref{lemma-2.1-aux} and \ref{lemma-2.1} and Corollary \ref{cor-dominate L2}, we obtain the following lemma.

\begin{cor}\label{cor-dominate L2-n}

Assume that $\tilde{\vb}_{0}\in \mathbb{H}^{1}(\dom)$ and \autoref{A2} holds. Let $\vb_{n}^{\eps,m}:[0,\tau_n) \to H_n$ be the local maximal solutions of (\ref{ee-epsilon approximation}). Then for every $s>0$ and $t\in [0, \tau_{n})$,
\begin{align}\label{eq-2.8-1}
 \Vert \vb_{n}^{\eps,m}(t)\Vert_{s}\leq \Vert \vb_{n}^{\eps,m}(t)\Vert_{H^{s}}\leq 2\Big (\Vert \vb_{n}^{\eps,m}(t)\Vert_{s}+\frac{1}{(2\pi)^\frac{3}{2}}\int_{0}^{t}\alpha^{\eps}(s)\left (\Vert \vb_{n}^{\eps,m}(s)\Vert_{\frac{1}{2}}^{2}+\Vert \wb_m^{\eps}(s)\Vert_{\frac{1}{2}}^{2}\right )\,ds\Big )+\frac{1}{(2\pi)^\frac{3}{2}}\Vert\tilde{\vb}_{0}\Vert_{L^1}.
\end{align}

Moreover, for every $s\geq \frac12$ and $\tau^{\eps,m}\in(0,\xi_1^{\eps,m})$, there exists a constant $C_{\tau^{\eps,m}}>0$, independent of $n$, such that
\begin{align}\label{eqn-H^s bound-n}
\Vert \vb_{n}^{\eps,m}(t)\Vert_{H^{s}}
\leq
C_{\tau^{\eps,m}}
\left(
\Vert \vb_{n}^{\eps,m}(t)\Vert_{s}+1
\right),
\qquad n\geq n_0,\;
t\in[0,\tau^{\eps,m}],
\end{align}
where $n_0$ is chosen as in \autoref{lem-vneps-regularization}.
\end{cor}
\begin{proof}
The estimate \eqref{eq-2.8-1} follows from the proof of
\autoref{cor-dominate L2}.

Fix $\tau^{\eps,m}\in(0,\xi_1^{\eps,m})$. By
\autoref{lem-vneps-regularization},
\[
\sup_{n\geq n_0}
\sup_{t\in[0,\tau^{\eps,m}]}
\Vert\vb_n^{\eps,m}(t)\Vert_{H^1}<\infty.
\]
Since
\[
\Vert\vb_n^{\eps,m}(t)\Vert_{\frac12}
\leq
\Vert\vb_n^{\eps,m}(t)\Vert_{H^1},
\]
and $\wb_m^\eps$ is smooth and $\alpha^\eps$ is bounded on
$[0,\tau^{\eps,m}]$, there exists a constant
$C_{\tau^{\eps,m}}^\prime>0$, independent of $n$, such that
\[
\sup_{n\geq n_0}
\int_0^{\tau^{\eps,m}}
\alpha^\eps(r)
\left(
\Vert\vb_n^{\eps,m}(r)\Vert_{\frac12}^2
+
\Vert\wb_m^\eps(r)\Vert_{\frac12}^2
\right)\,dr
\leq
C_{\tau^{\eps,m}}^\prime.
\]
Hence, \eqref{eqn-H^s bound-n} follows from
\eqref{eq-2.8-1}, with a constant
$C_{\tau^{\eps,m}}$ independent of $n$.
\end{proof}

In view of \autoref{lem-vneps-regularization} and \autoref{cor-dominate L2-n}, we have the following lemma.
\begin{lemma}\label{eq-strong convergence galerkin regularization}
Assume that $\tilde\vb_{0}\in \mathbb{H}^{1}(\dom)$ and \autoref{A2} holds. Let $\vb_{n}^{\eps,m}:[0,\tau_n) \to H_n$ be the local maximal solutions of (\ref{ee-epsilon approximation}) and $\xi_1^{\eps,m}$ denote the maximal existence time of \eqref{ee-epsilon}
as in \autoref{prop-maximal global-veps}.
Then for every $0<\tau_\ast<\tau^{\eps,m}<\xi_1^{\eps,m}$ and integer $k\geq 2$, it holds almost surely that

\begin{align}\label{eqn-vneps-HK}
\sup_{n\geq n_0} \sup\limits_{t\in [\tau_\ast, \tau^{\eps,m}]} \Vert \vb_{n}^{\eps,m} (t) \Vert_{H^k}^2+\sup_{n\geq n_0} \int_{\tau_\ast}^ {\tau^{\eps,m}} \Vert \vb_{n}^{\eps,m} (t) \Vert_{H^{k+1}}^{2}\,dt <\infty,
\end{align}
where $n_0$ is chosen as in \autoref{lem-vneps-regularization}. In particular,
\[
\sup_{n\geq n_0} \sup\limits_{t\in [\tau_\ast, \tau^{\eps,m}]} \Vert \vb_{n}^{\eps,m} (t) \Vert_{H^6}<\infty,\;\mathbb{P}\text{-a.s.}.
\]
\end{lemma}

\begin{proof}[Proof of \autoref{eq-strong convergence galerkin regularization}]
Fix $\omega\in \Omega$ outside a null set such that \eqref{eq-uniform-H1-H2} holds.
Let us choose and fix  $0<\tau_\ast<\tau^{\eps,m}<\xi_1^{\eps,m}$.

For $k\geq2$, set
\[
a_j:=2^{j-k}\tau_\ast,\qquad j=2,\ldots,k.
\]
Then
\[
0<a_2<\cdots<a_k=\tau_\ast.
\]

We prove by induction that, for $j=2,\ldots,k$,
\begin{equation}\label{eqn-vneps-HK-1}
\sup_{n\geq n_0}
\sup_{t\in[a_j,\tau^{\eps,m}]}
\Vert\vb_n^{\eps,m}(t)\Vert_{H^j}^2
+
\sup_{n\geq n_0}
\int_{a_j}^{\tau^{\eps,m}}
\Vert\vb_n^{\eps,m}(t)\Vert_{H^{j+1}}^2\,dt
<\infty.
\end{equation}

Firstly, we aim to show that \eqref{eqn-vneps-HK-1} holds for $j=2$, i.e.
\begin{equation}\label{eqn-429}
\begin{aligned}
\sup_{n\geq n_0} \sup\limits_{t\in [a_2, \tau^{\eps,m}]} \Vert \vb_{n}^{\eps,m} (t) \Vert_{H^2}^2+\sup_{n\geq n_0} \int_{a_2}^ {\tau^{\eps,m}} \Vert \vb_{n}^{\eps,m} (t) \Vert_{H^3}^{2}\,dt <\infty.
\end{aligned}
\end{equation}
Multiplying (\ref{ee-epsilon approximation}) with $(-\Delta)^{2} \vb_{n}^{\eps,m}$ and integrating on $\dom$, since $\Lambda$ is self-adjoint, by the H\"older inequality, the Young inequality, the Sobolev inequality we infer that there exist constants $c>0$ and $\delta_1\in(0,1)$ such that
\begin{align}\label{eqn-427}
\frac{1}{2}\partial_{t}\Vert  \vb_{n}^{\eps,m}\Vert_{2}^{2}+\Vert \vb_{n}^{\eps,m}\Vert_{3}^{2}&= -\alpha^{\eps}\langle ((\vb_{n}^{\eps,m}+\wb_m^{\eps})\cdot \nabla) (\vb_{n}^{\eps,m}+\wb_m^{\eps}) , (-\Delta)^{2} \vb_{n}^{\eps,m}\rangle\\
&\leq \alpha^\eps \Vert\nabla[((\vb_{n}^{\eps,m} + \wb_m^\eps)\cdot \nabla) (\vb_{n}^{\eps,m} + \wb_m^\eps)]\Vert_{L^2} \Vert\vb_{n}^{\eps,m}\Vert_{3}\nonumber\\
&\leq \alpha^{\eps}\Vert \nabla (\vb_{n}^{\eps,m}+\wb_m^{\eps})\Vert_{L^4}^{2}\Vert \vb_{n}^{\eps,m}\Vert_{3}+\alpha^{\eps}\Vert\vb_{n}^{\eps,m}+\wb_m^{\eps}\Vert_{L^\infty}\Vert  \vb_{n}^{\eps,m}+\wb_m^{\eps}\Vert_{2}\Vert \vb_{n}^{\eps,m}\Vert_{3}\nonumber\\
&\leq c(\alpha^{\eps})^{2}\Vert \nabla (\vb_{n}^{\eps,m}+\wb_m^{\eps})\Vert_{L^4}^{4}+c(\alpha^{\eps})^{2}\Vert\vb_{n}^{\eps,m}+\wb_m^{\eps}\Vert_{L^\infty}^{2}\Vert  \vb_{n}^{\eps,m}+\wb_m^{\eps}\Vert_{2}^{2}+\delta_1\Vert   \vb_{n}^{\eps,m}\Vert_{3}^{2}.\nonumber
\end{align}
Since $\mathbb{H}^1(\dom)\hookrightarrow \mathbb{L}^4(\dom)$ and $\Vert\cdot\Vert_{s_1}\leq \Vert\cdot\Vert_{s_2}$ for $0<s_1\leq s_2$, we infer there exists a constant $c>0$ such that
\[
\Vert \nabla (\vb_{n}^{\eps,m}+\wb_m^{\eps})\Vert_{L^4}^2\leq c\Vert \nabla (\vb_{n}^{\eps,m}+\wb_m^{\eps})\Vert_{H^1}^2=c\Vert \nabla (\vb_{n}^{\eps,m}+\wb_m^{\eps})\Vert_{L^2}^2+c\Vert \nabla (\vb_{n}^{\eps,m}+\wb_m^{\eps})\Vert_{1}^2\leq c\Vert \vb_{n}^{\eps,m}+\wb_m^{\eps}\Vert_{2}^2.
\]
Therefore, since $\mathbb{H}^2(\dom)\hookrightarrow \mathbb{L}^\infty(\dom)$, \eqref{eqn-427} implies
\begin{align*}
\frac{1}{2}\partial_{t}\Vert  \vb_{n}^{\eps,m}\Vert_{2}^{2}+\Vert \vb_{n}^{\eps,m}\Vert_{3}^{2}\leq c(\alpha^{\eps})^{2}\Vert \vb_{n}^{\eps,m}+\wb_m^{\eps}\Vert_{2}^{4} +c(\alpha^{\eps})^{2}\Vert \vb_{n}^{\eps,m}+\wb_m^{\eps}\Vert_{H^{2}}^{2}\Vert  \vb_{n}^{\eps,m}+\wb_m^{\eps}\Vert_{2}^{2}+\delta_1\Vert \vb_{n}^{\eps,m}\Vert_{3}^{2},
\end{align*}
then triangle inequality implies
\begin{align}\label{eqn-vneps-2}
\frac{1}{2}\partial_{t}\Vert  \vb_{n}^{\eps,m}\Vert_{2}^{2}+(1-\delta_1)\Vert \vb_{n}^{\eps,m}\Vert_{3}^{2}&\leq c(\alpha^{\eps})^{2}\Vert \vb_{n}^{\eps,m}+\wb_m^{\eps}\Vert_{H^{2}}^{2}\Vert  \vb_{n}^{\eps,m}+\wb_m^{\eps}\Vert_{2}^{2}\\
&\leq c(\alpha^{\eps})^{2}\Vert \vb_{n}^{\eps,m}+\wb_m^{\eps}\Vert_{H^{2}}^{2}\Vert  \vb_{n}^{\eps,m}\Vert_{2}^{2}+c(\alpha^{\eps})^{2}\Vert \vb_{n}^{\eps,m}+\wb_m^{\eps}\Vert_{H^{2}}^{2}\Vert  \wb_m^{\eps}\Vert_{2}^{2}.\nonumber
\end{align}

In view of \autoref{lem-vneps-regularization}
\[
\sup_{n\geq n_0} \sup\limits_{t\in [0, \tau^{\eps,m}]} \Vert \vb_{n}^{\eps,m} (t) \Vert_{H^1}<\infty,\  \ \sup_{n\geq n_0} \int_{0}^ {\tau^{\eps,m}} \Vert \vb_{n}^{\eps,m} (t) \Vert_{H^2}^{2}\,dt<\infty.
\]

For every $n\geq n_0$, there exists $t_n\in(0,a_2)$ such that
\[
\Vert\vb_n^{\eps,m}(t_n)\Vert_{H^2}^2
\leq \frac1{a_2}\int_0^{a_2}\Vert\vb_n^{\eps,m}(t)\Vert_{H^2}^2\,dt.
\]
Therefore,
\begin{align}\label{eq-vneps-H2bound}
\sup_{n\geq n_0} \Vert \vb_{n}^{\eps,m} (t_n) \Vert_{H^2}^{2}<\infty.
\end{align}
Since $\alpha^\eps$ is bounded on $[0,\tau^{\eps,m}]$ and $\wb_m^{\eps}\in  C^{\infty}([0, \tau^{\eps,m}]\times \dom;\mathbb{R}^3)$, we infer that

\begin{align*}
\sup_{n\geq n_0}\int_0^{\tau^{\eps,m}}(\alpha^\eps)^2\Vert\vb_n^{\eps,m}+\wb_m^\eps\Vert_{H^2}^2\left(1+\Vert\wb_m^\eps\Vert_2^2\right)\,dt<\infty.
\end{align*}
Hence, in view of \eqref{eqn-vneps-2} and \eqref{eq-vneps-H2bound}, using Gronwall's inequality on $[t_n,\tau^{\eps,m}]$, we deduce
\[
\sup_{n\geq n_0} \sup\limits_{t\in [a_2, \tau^{\eps,m}]} \Vert \vb_{n}^{\eps,m} (t) \Vert_{2}\leq \sup_{n\geq n_0} \sup\limits_{t\in [t_n, \tau^{\eps,m}]} \Vert \vb_{n}^{\eps,m} (t) \Vert_{2}<\infty.
\]
Moreover,
\[
 \sup_{n\geq n_0} \int_{a_2}^ {\tau^{\eps,m}} \Vert \vb_{n}^{\eps,m} (t) \Vert_{3}^{2}\,dt\leq  \sup_{n\geq n_0} \int_{t_n}^ {\tau^{\eps,m}} \Vert \vb_{n}^{\eps,m} (t) \Vert_{3}^{2}\,dt<\infty.
\]
In view of \autoref{cor-dominate L2-n}, for every $\tau_\ast\in (0, \tau^{\eps,m})$, we infer that there exists a constant $C>0$ such that
\begin{equation*}
\begin{aligned}
&\sup_{n\geq n_0} \sup\limits_{t\in [a_2, \tau^{\eps,m}]} \Vert \vb_{n}^{\eps,m} (t) \Vert_{H^2}\leq C(\sup_{n\geq n_0} \sup\limits_{t\in [a_2, \tau^{\eps,m}]}\Vert \vb_{n}^{\eps,m}\Vert_{2}+1) <\infty,\\
 &\sup_{n\geq n_0} \int_{a_2}^ {\tau^{\eps,m}} \Vert \vb_{n}^{\eps,m} (t) \Vert_{H^3}^{2}\,dt\leq C(\sup_{n\geq n_0} \int_{a_2}^ {\tau^{\eps,m}} \Vert \vb_{n}^{\eps,m} (t) \Vert_{3}^{2}\,dt+1)  <\infty,
\end{aligned}
\end{equation*}
which implies \eqref{eqn-429} holds.

Let $j\in\{3,\ldots,k\}$ and assume that \eqref{eqn-vneps-HK-1} holds with $j$ replaced by $j-1$, i.e.
\begin{equation}\label{eqn-vneps-HK-1-1}
\sup_{n\geq n_0}
\sup_{t\in[a_{j-1},\tau^{\eps,m}]}
\Vert\vb_n^{\eps,m}(t)\Vert_{H^{j-1}}^2
+
\sup_{n\geq n_0}
\int_{a_{j-1}}^{\tau^{\eps,m}}
\Vert\vb_n^{\eps,m}(t)\Vert_{H^{j}}^2\,dt
<\infty.
\end{equation}
In the following we aim to show that \eqref{eqn-vneps-HK-1} holds for $j$ by induction.

In view of \eqref{eqn-vneps-HK-1-1}, for every $n\geq n_0$, there exists
$t_{n,j}\in(a_{j-1},a_j)$ such that
\[
\Vert\vb_n^{\eps,m}(t_{n,j})\Vert_{H^j}^2\leq\frac1{a_j-a_{j-1}}\int_{a_{j-1}}^{a_j}\Vert\vb_n^{\eps,m}(t)\Vert_{H^j}^2\,dt.
\]
Therefore,
\begin{align}\label{eqn-vnem-hj}
\sup_{n\geq n_0}\Vert\vb_n^{\eps,m}(t_{n,j})\Vert_{H^j}^2<\infty.
\end{align}

We first give the product estimate used below.
For every integer $j\geq3$ and every smooth periodic vector field $h$, the Leibniz rule, the Sobolev embedding theorem and the Gagliardo--Nirenberg interpolation inequality give
\begin{align*}
\Vert\Lambda^{j-1}((h\cdot\nabla)h)\Vert_{L^2}
&\leq C_j\left(\Vert h\Vert_{L^\infty}\Vert h\Vert_j
+\sum_{\ell=1}^{j-1}\Vert D^\ell h\Vert_{L^4}\Vert D^{j-\ell}h\Vert_{L^4}\right)\\
&\leq C_j\Vert h\Vert_{H^2}\Vert h\Vert_j,
\end{align*}
where $D^\ell h$ denotes the collection of spatial derivatives of order $\ell$.
Indeed, writing $\hat h_0:=|\dom|^{-1}\int_\dom h\,dx$, Sobolev embedding bounds each factor $\Vert D^\ell h\Vert_{L^4}$ by $c\Vert h-\hat h_0\Vert_{H^{\ell+3/4}}$.
Interpolation between $\mathbb H^2$ and $\mathbb H^j$, together with the periodic Poincar\'e inequality, therefore gives
\[
\Vert D^\ell h\Vert_{L^4}\Vert D^{j-\ell}h\Vert_{L^4}
\leq C_j\Vert h-\hat h_0\Vert_{H^2}\Vert h-\hat h_0\Vert_{H^j}
\leq C_j\Vert h\Vert_{H^2}\Vert h\Vert_j,
\qquad 1\leq\ell\leq j-1.
\]
The endpoint term is bounded using $\mathbb H^2(\dom)\hookrightarrow\mathbb L^\infty(\dom)$.
The constant $C_j$ is independent of $n$.

Set $h=\vb_n^{\eps,m}+\wb_m^\eps$.
Taking the inner product of \eqref{ee-epsilon approximation} with $(-\Delta)^j\vb_n^{\eps,m}$, and using that $P_n$ is orthogonal and commutes with $\Lambda$, we infer that
\begin{align*}
\frac12\partial_t\Vert\vb_n^{\eps,m}\Vert_j^2+\Vert\vb_n^{\eps,m}\Vert_{j+1}^2
&=-\alpha^\eps\left\langle\Lambda^{j-1}((h\cdot\nabla)h),\Lambda^{j+1}\vb_n^{\eps,m}\right\rangle\\
&\leq C_j|\alpha^\eps|\Vert h\Vert_{H^2}\Vert h\Vert_j\Vert\vb_n^{\eps,m}\Vert_{j+1}.
\end{align*}

Therefore, by the Young inequality, for any $\delta_2\in(0,1)$,
\begin{align*}
\frac{1}{2}\partial_{t}\Vert  \vb_{n}^{\eps,m}\Vert_{j}^{2}+\Vert \vb_{n}^{\eps,m}\Vert_{j+1}^{2}
&\leq c(\alpha^{\eps})^{2}\Vert \vb_{n}^{\eps,m}+\wb_m^{\eps}\Vert_{H^{2}}^{2}\Vert \vb_{n}^{\eps,m}+\wb_m^{\eps}\Vert_{j}^{2}+\delta_2\Vert \vb_{n}^{\eps,m}\Vert_{j+1}^{2},
\end{align*}
which implies
\begin{align}\label{eqn-vnemj}
\frac{1}{2}\partial_{t}\Vert  \vb_{n}^{\eps,m}\Vert_{j}^{2}+(1-\delta_2)\Vert \vb_{n}^{\eps,m}\Vert_{j+1}^{2}
&\leq c(\alpha^{\eps})^{2}\Vert \vb_{n}^{\eps,m}+\wb_m^{\eps}\Vert_{H^{2}}^{2}\Vert \vb_{n}^{\eps,m}+\wb_m^{\eps}\Vert_{j}^{2}\\
&\leq c(\alpha^{\eps})^{2}\Vert \vb_{n}^{\eps,m}+\wb_m^{\eps}\Vert_{H^{2}}^{2}\Vert \vb_{n}^{\eps,m}\Vert_{j}^{2}+c(\alpha^{\eps})^{2}\Vert \vb_{n}^{\eps,m}+\wb_m^{\eps}\Vert_{H^{2}}^{2}\Vert \wb_m^{\eps}\Vert_{j}^{2}.\nonumber
\end{align}

Since $\alpha^\eps$ is bounded on $[0,\tau^{\eps,m}]$ and $\wb_m^\eps\in C^\infty([0,\tau^{\eps,m}]\times\dom;\mathbb R^3)$, by \eqref{eq-uniform-H1-H2} in \autoref{lem-vneps-regularization}, we infer that

\begin{align*}
\sup_{n\geq n_0}\int_0^{\tau^{\eps,m}}(\alpha^\eps)^2\Vert\vb_n^{\eps,m}+\wb_m^\eps\Vert_{H^2}^2\left(1+\Vert\wb_m^\eps\Vert_j^2\right)\,dt<\infty.
\end{align*}
Therefore, by \eqref{eqn-vnemj} and \eqref{eqn-vnem-hj}, applying the Gronwall inequality on $[t_{n,j},\tau^{\eps,m}]$ and using \autoref{cor-dominate L2-n}, we deduce
\[
\sup_{n\geq n_0} \sup\limits_{t\in [a_j, \tau^{\eps,m}]} \Vert \vb_{n}^{\eps,m} (t) \Vert_{H^j}^2+\sup_{n\geq n_0} \int_{a_j}^ {\tau^{\eps,m}} \Vert \vb_{n}^{\eps,m} (t) \Vert_{H^{j+1}}^{2}\,dt<\infty, \mathbb{P}\text{-a.s.},
\]
which implies \eqref{eqn-vneps-HK-1} holds for $j$.

By induction, \eqref{eqn-vneps-HK-1} holds for
$j=2,\ldots,k$. Since $a_k=\tau_\ast$, we infer that
\[
\sup_{n\geq n_0}\sup_{t\in[\tau_\ast,\tau^{\eps,m}]}\Vert\vb_n^{\eps,m}(t)\Vert_{H^k}^2+\sup_{n\geq n_0}\int_{\tau_\ast}^{\tau^{\eps,m}}\Vert\vb_n^{\eps,m}(t)\Vert_{H^{k+1}}^2\,dt<\infty.
\]
Since $k\geq2$ is arbitrary, \eqref{eqn-vneps-HK} holds.

In particular, taking $k=6$ in \eqref{eqn-vneps-HK} yields
\[
\sup_{n\geq n_0}\sup_{t\in[\tau_\ast,\tau^{\eps,m}]}\Vert\vb_n^{\eps,m}(t)\Vert_{H^6}<\infty.
\]
The proof of \autoref{eq-strong convergence galerkin regularization} is complete.
\end{proof}

\subsection{Classical regularity}\label{app-proof-classical-solution}

\begin{proof}[Proof of \autoref{eq-classical solution}]
By the Sobolev product estimate and the Sobolev embedding theorem, we infer that, for every $s>\frac32$, there exists a constant $C_s>0$ such that for every $a\in\mathbb{H}^s$ and $b\in\mathbb{H}^{s+1}$,
\begin{equation}\label{eqn-a nabla b}
\begin{aligned}
\Vert (a\cdot \nabla) b\Vert_{H^s}&\leq C_s\left(\Vert a\Vert_{H^s}\Vert \nabla b\Vert_{L^\infty}+\Vert a\Vert_{L^\infty}\Vert \nabla b\Vert_{H^s}\right)\\
&\leq C_s\Vert a\Vert_{H^s}\Vert \nabla b\Vert_{H^s}+C_s\Vert a\Vert_{H^s}\Vert \nabla b\Vert_{H^s}\leq C_s\Vert a\Vert_{H^s}\Vert b\Vert_{H^{s+1}}.
\end{aligned}
\end{equation}
Here, $C_s$ may change from line to line.
Let us choose and fix $0<\tau_\ast<\tau^{\eps,m}<\xi_1^{\eps,m}$ and $n\geq n_0$, where $n_0$ is chosen as in \autoref{lem-vneps-regularization}.
In view of \eqref{ee-epsilon approximation}, by triangle inequality and \eqref{eqn-a nabla b} we infer
\begin{align*}
\sup\limits_{t\in [\tau_\ast, \tau^{\eps,m}]} \Vert \partial_{t} \vb_{n}^{\eps,m}\Vert_{H^4} &\leq c \sup\limits_{t\in [\tau_\ast, \tau^{\eps,m}]} \Big{(}\Vert \Delta \vb_{n}^{\eps,m}\Vert_{H^4} + \alpha^{\eps}\Vert ((\vb_{n}^{\eps,m}+\wb_m^{\eps})\cdot \nabla)(\vb_{n}^{\eps,m}+\wb_m^{\eps})\Vert_{H^4}\Big{)}\\
&\leq c \sup\limits_{t\in [\tau_\ast, \tau^{\eps,m}]} \Big{(}\Vert \vb_{n}^{\eps,m}\Vert_{H^6} + \alpha^{\eps}(\Vert \vb_{n}^{\eps,m}\Vert_{H^4} +\Vert \wb_m^{\eps}\Vert_{H^4})(\Vert \vb_{n}^{\eps,m}\Vert_{H^5} +\Vert \wb_m^{\eps}\Vert_{H^5} )\Big{)},
\end{align*}
Moreover, taking the time derivative on both sides of \eqref{ee-epsilon approximation}, we have
\begin{equation}\label{eqn-partial 2 t}
\begin{aligned}
\partial_{t}^2 \vb_{n}^{\eps,m}&=\Delta \partial_{t}\vb_{n}^{\eps,m}-\left(\partial_{t}\alpha^{\eps}(t)\right)P_{n}\Big{(}(\vb_{n}^{\eps,m}+\wb_m^{\eps})\cdot \nabla (\vb_{n}^{\eps,m}+\wb_m^{\eps})\Big{)}\\
&-\alpha^{\eps}(t)P_{n}\Big{(}\partial_{t}(\vb_{n}^{\eps,m}+\wb_m^{\eps})\cdot \nabla (\vb_{n}^{\eps,m}+\wb_m^{\eps})\Big{)}-\alpha^{\eps}(t)P_{n}\Big{(}(\vb_{n}^{\eps,m}+\wb_m^{\eps})\cdot \nabla \partial_{t}(\vb_{n}^{\eps,m}+\wb_m^{\eps})\Big{)},\ \mathrm{on}\ [0,\tau_n)\times \dom.
\end{aligned}
\end{equation}
Recall that  $\alpha^{\eps}(t)=\exp{(bB^{\eps}(t))}$, we infer
\[
\partial_{t}\alpha^{\eps}(t)=\partial_{t}(\exp{(bB^{\eps}(t))})=b\alpha^{\eps}(t)(\partial_{t}B^{\eps}(t)),
\]
therefore, by \eqref{eqn-partial 2 t} and \eqref{eqn-a nabla b} we infer that
\begin{align*}
\sup\limits_{t\in [\tau_\ast, \tau^{\eps,m}]} \Vert \partial_{t}^{2} \vb_{n}^{\eps,m}\Vert_{H^2} &\leq c\sup\limits_{t\in [\tau_\ast, \tau^{\eps,m}]} \Big{(}\vert b\vert \alpha^{\eps}(t)\vert \partial_{t}B^{\eps}(t)\vert (\Vert \vb_{n}^{\eps,m}\Vert_{H^2} +\Vert \wb_m^{\eps}\Vert_{H^2} )(\Vert \vb_{n}^{\eps,m}\Vert_{H^3} +\Vert \wb_m^{\eps}\Vert_{H^3} )\Big{)}\\
&+c\sup\limits_{t\in [\tau_\ast, \tau^{\eps,m}]}\Big{(}\alpha^{\eps}(t)(\Vert \partial_{t}\vb_{n}^{\eps,m}\Vert_{H^2} +\Vert \partial_{t}\wb_m^{\eps}\Vert_{H^2} )(\Vert \vb_{n}^{\eps,m}\Vert_{H^3} +\Vert \wb_m^{\eps}\Vert_{H^3} )\Big{)}\\
&+c\sup\limits_{t\in [\tau_\ast, \tau^{\eps,m}]}\Big{(}\alpha^{\eps}(t)(\Vert \vb_{n}^{\eps,m}\Vert_{H^2} +\Vert \wb_m^{\eps}\Vert_{H^2} )(\Vert \partial_{t}\vb_{n}^{\eps,m}\Vert_{H^3} +\Vert \partial_{t}\wb_m^{\eps}\Vert_{H^3} )\Big{)}\\
&+\sup\limits_{t\in [\tau_\ast, \tau^{\eps,m}]} \Vert \partial_{t} \vb_{n}^{\eps,m}\Vert_{H^4}.
\end{align*}
The above a priori estimates in \autoref{eq-strong convergence galerkin regularization} imply that
\begin{align*}
\sup\limits_{n\geq n_0}\sup\limits_{t\in [\tau_\ast, \tau^{\eps,m}]} \Vert \partial_{t} \vb_{n}^{\eps,m}\Vert_{H^4} < \infty, \mathbb{P}\text{-a.s.}, \mathrm{and}\ \sup\limits_{n\geq n_0}\sup\limits_{t\in [\tau_\ast, \tau^{\eps,m}]} \Vert \partial_{t}^{2} \vb_{n}^{\eps,m}\Vert_{H^2}< \infty, \mathbb{P}\text{-a.s.}.
\end{align*}
It further means $\{ \vb_{n}^{\eps,m}\}_{n\geq n_0}$ is bounded in $\mathbb{H}^{1}( [\tau_\ast, \tau^{\eps,m}];\mathbb{H}^4(\dom) )\cap \mathbb{H}^{2}( [\tau_\ast, \tau^{\eps,m}];\mathbb{H}^2(\dom) )$.
Consequently,

by the strong convergence in \autoref{lem-vneps-regularization}, every weak limit equals $\vb^{\eps,m}$. Hence, $ \vb^{\eps,m} \in \mathbb{H}^{1}( [\tau_\ast, \tau^{\eps,m}];\mathbb{H}^4(\dom) )\cap \mathbb{H}^{2}( [\tau_\ast, \tau^{\eps,m}];\mathbb{H}^2(\dom) )$. By the Sobolev embedding theorem in both time and space, we infer that
\begin{align*}
\mathbb{H}^{1}( [\tau_\ast, \tau^{\eps,m}];\mathbb{H}^4(\dom) )\hookrightarrow C( [\tau_\ast, \tau^{\eps,m}];\mathbb{H}^4(\dom) ) \hookrightarrow C^{0,2}( [\tau_\ast, \tau^{\eps,m}]\times \dom; \mathbb{R}^{3}),
\\
 \mathbb{H}^{2}( [\tau_\ast, \tau^{\eps,m}];\mathbb{H}^2(\dom) )\hookrightarrow C^1( [\tau_\ast, \tau^{\eps,m}];\mathbb{H}^2(\dom) )
 \hookrightarrow C^{1,0}( [\tau_\ast, \tau^{\eps,m}]\times \dom; \mathbb{R}^{3}),
\end{align*}
where the space
$C^{n,r} ((0, \tau^{\eps,m}]\times \dom; \mathbb{R}^{3}),\;n,r\in\mathbb{N}$ is defined in \eqref{eq-Cnm}.
Thus, $\vb^{\eps,m}\in C^{1,0}( [\tau_\ast, \tau^{\eps,m}]\times \dom; \mathbb{R}^{3})\cap C^{0,2}( [\tau_\ast, \tau^{\eps,m}]\times \dom; \mathbb{R}^{3})$.
Since $\tau_\ast\in(0,\tau^{\eps,m}]$ is arbitrary, $ \vb^{\eps,m} \in C^{1,0}( (0, \tau^{\eps,m}]\times \dom; \mathbb{R}^{3})\cap C^{0,2}( (0, \tau^{\eps,m}]\times \dom; \mathbb{R}^{3})$. Therefore, $\vb^{\eps,m}$ is the classical solution to (\ref{ee-epsilon}).
\end{proof}

\end{document}